\documentclass[11pt]{article}

\usepackage{natbib}
\usepackage{geometry}
\usepackage{booktabs}
\usepackage{footmisc}
\usepackage{amsmath}
\usepackage{amssymb}
\usepackage{amsthm}
\usepackage{graphicx}
\usepackage{courier}
\usepackage{ulem}
\usepackage[ruled]{algorithm2e}
\usepackage{algorithmic}
\usepackage{multirow}
\usepackage{float}

\usepackage{tikz}
\usetikzlibrary{positioning}
\usetikzlibrary{arrows.meta}

\newtheorem{theorem}{Theorem}[section]
\newtheorem{lemma}{Lemma}[section]
\newtheorem{proposition}{Proposition}[section]
\newtheorem{corollary}{Corollary}[section]
\newtheorem{assumption}{Assumption}
\newtheorem{definition}{Definition}
\newtheorem{remark}{Remark}
\newtheorem{example}{Example}
\def\todo#1 {\textcolor{red}{Todo: #1}}

\allowdisplaybreaks[4]

\usepackage[utf8]{inputenc} 
\usepackage[T1]{fontenc}    
\usepackage{hyperref}       
\usepackage{url}            
\usepackage{booktabs}       
\usepackage{amsfonts}       
\usepackage{nicefrac}       
\usepackage{microtype}      
\usepackage{xcolor}         

\usepackage{amsmath,amsfonts,bm}
\usepackage{bbm}

\def\rmG{{\mathbf{G}}}

\def\rmI{{\mathbf{I}}}

\def\rmN{{\mathbf{N}}}

\def\rmV{{\mathbf{V}}}
\def\rmW{{\mathbf{W}}}
\def\rmX{{\mathbf{X}}}
\def\rmY{{\mathbf{Y}}}

\def\vzero{{\bm{0}}}

\def\vmu{{\bm{\mu}}}
\def\vnu{\bm{\nu}}

\def\va{{\bm{a}}}
\def\vb{{\bm{b}}}
\def\vc{{\bm{c}}}
\def\vd{{\bm{d}}}
\def\ve{{\bm{e}}}

\def\vh{{\bm{h}}}
\def\vi{{\bm{i}}}

\def\vp{{\bm{p}}}
\def\vq{{\bm{q}}}

\def\vu{{\bm{u}}}
\def\vv{{\bm{v}}}
\def\vw{{\bm{w}}}
\def\vx{{\bm{x}}}
\def\vy{{\bm{y}}}
\def\vz{{\bm{z}}}
\def\vA{{\bm{A}}}

\def\mA{{\bm{A}}}
\def\mB{{\bm{B}}}

\def\mD{{\bm{D}}}
\def\mE{{\bm{E}}}

\def\mI{{\bm{I}}}

\def\mM{{\bm{M}}}

\def\mP{{\bm{P}}}
\def\mQ{{\bm{Q}}}

\def\mS{{\bm{S}}}
\def\mT{{\bm{T}}}
\def\mU{{\bm{U}}}
\def\mV{{\bm{V}}}
\def\mW{{\bm{W}}}
\def\mX{{\bm{X}}}
\def\mY{{\bm{Y}}}

\DeclareMathAlphabet{\mathsfit}{\encodingdefault}{\sfdefault}{m}{sl}
\SetMathAlphabet{\mathsfit}{bold}{\encodingdefault}{\sfdefault}{bx}{n}

\def\gC{{\mathcal{C}}}
\def\gD{{\mathcal{D}}}
\def\gE{{\mathcal{E}}}
\def\gF{{\mathcal{F}}}
\def\gG{{\mathcal{G}}}

\def\gI{{\mathcal{I}}}
\def\gJ{{\mathcal{J}}}

\def\gL{{\mathcal{L}}}
\def\gM{{\mathcal{M}}}
\def\gN{{\mathcal{N}}}
\def\gO{{\mathcal{O}}}
\def\gP{{\mathcal{P}}}
\def\gQ{{\mathcal{Q}}}
\def\gR{{\mathcal{R}}}
\def\gS{{\mathcal{S}}}
\def\gT{{\mathcal{T}}}
\def\gU{{\mathcal{U}}}

\def\gW{{\mathcal{W}}}

\def\sE{{\mathbb{E}}}

\def\sN{{\mathbb{N}}}

\def\sP{{\mathbb{P}}}

\def\sR{{\mathbb{R}}}

\def\0{{\bf 0}}
\def\1{{\bf 1}}

\def\bfB{{\bf B}}

\def\bfN{{\bf N}}

\def\bfx{{\bf x}}

\def\EB{{\mathbb E}}

\def\PB{{\mathbb P}}

\def\RB{{\mathbb R}}

\def\lam{\mbox{\boldmath$\lambda$\unboldmath}}

\def\Xii{\mbox{\boldmath$\Xi$\unboldmath}}

\def\var{\mathrm{var}}

\def\tr{\mathrm{tr}}

\def\diag{\mathrm{diag}}

\newcommand{\ssum}[3]{\sum\limits_{{#1}={#2}}^{#3}}

\def\inner#1#2{\left\langle #1, #2 \right\rangle}
\newcommand{\norm}[1]{\left\|{#1}\right\|}
\newcommand{\abs}[1]{\left|{#1}\right|}

\begin{document}

\title{Asymptotic Behaviors and Phase Transitions in Projected Stochastic Approximation: A Jump Diffusion  Approach}

  



\maketitle
\begin{center}
\begin{tabular}{ccc}
 Jiadong Liang & Yuze Han \\
{\texttt{jdliang@pku.edu.cn}} &
{\texttt{hanyuze97@pku.edu.cn}} \\
Xiang Li & Zhihua Zhang \\
{\texttt{lx10077@pku.edu.cn}} &
{\texttt{zhzhang@math.pku.edu.cn}}
\end{tabular}

\begin{tabular}{ccc}
     School of Mathematical Sciences, Peking University, Beijing
\end{tabular}
\end{center}

\begin{abstract}
In this paper we consider linearly constrained optimization problems and
propose a loopless projection stochastic approximation (LPSA) algorithm.  It  performs the projection with probability $p_n$ at the $n$-th iteration to ensure feasibility.
Considering a specific family of the probability $p_n$ and step size $\eta_n$, we analyze our algorithm from an asymptotic and continuous perspective.
Using a novel jump diffusion approximation,
we show that the trajectories connecting those properly rescaled last iterates weakly converge to the solution of specific stochastic differential equations (SDEs).
By analyzing  SDEs, we identify the asymptotic behaviors  of LPSA for different choices of $(p_n, \eta_n)$.
We find that the algorithm presents an intriguing asymptotic bias-variance trade-off and yields phase transition phenomenons, according to the relative magnitude of $p_n$ w.r.t.\ $\eta_n$. 
This finding provides insights on selecting appropriate ${(p_n, \eta_n)}_{n \geq 1}$ to minimize the projection cost. Additionally, we propose the Debiased LPSA (DLPSA) as a practical application of our jump diffusion approximation result. DLPSA is shown to effectively reduce projection complexity compared to vanilla LPSA.
\end{abstract}

\section{Introduction}
Optimization problems involving uncertainty 
play a vital role
in many areas such as machine learning, artificial intelligence, and statistical inference. In this paper we are concerned with
linearly constrained optimization problems,  
which are  the most fundamental class of constrained optimization problems. That is, 
\begin{equation}\label{eq:proj_opt}
\begin{aligned}
\min\limits_{\vx}\sE_{\zeta\sim \gD} f(\vx,\zeta) ~~\text{subject to }\vA^\top \vx= \vzero,
\end{aligned}
\end{equation}
where $\zeta$ denotes the data point generated independently and identically from a distribution $\gD$ and $\vA$ is the constraint matrix.
Classic methods such as simplex methods, interior-point methods, and primal-dual methods have been widely used to provide simple and effective solutions for linearly constrained optimization problems~\citep{goldfarb1969extension, fletcher1972algorithm, murtagh1978large,boyd2011distributed}. 
These algorithms, supported by theoretical guarantees, have found wide applications in various fields such as resource allocation~\citep{arrow1954existence}, portfolio combination~\citep{best2010portfolio}, statistical classification~\citep{bartlett2006convexity}, and decentralized optimization~\citep{boyd2011distributed}. 
However, in the era of big data, the computational cost of these classic methods becomes prohibitively expensive due to the high dimensionality and large data volume.


In contrast, stochastic algorithms have become more attractive due to their lower per-iteration computational complexity and greater versatility. 
When using stochastic gradient descent (SGD) to solve linearly constrained optimization problems, projections are required to ensure feasibility. 
However, in large-scale cases, projections could become a bottleneck due to their high computational cost.
To improve projection efficiency, an idea is to decrease the frequency of projections by performing them after several steps of SGD periodically.
This idea has been explored in the context of distributed learning and has been shown to be effective in reducing communication complexity~\citep{mcmahan2017communication}.
In fact, in the context of linearly constrained optimization, communication is equivalent to projection. 
Among the algorithms that have been proposed, \textit{Local SGD} that alternates between multiple steps of parallel SGD and one communication is the most representative due to its simplicity and effectiveness~\citep{mcmahan2017communication}. 
Empirical investigations have shown that it outperforms other methods in terms of communication efficiency~\cite{lin2018don}.
Non-asymptotic convergence analyses are provided to illustrate its fast convergence which is robust to the underlying data generation mechanism~\citep{li2019communication, bayoumi2020tighter, koloskova2020unified, woodworth2020local, woodworth2020minibatch, koloskova2020unified}. 
However, for the general setting considered in this work, it is still unclear whether the lazy projection approach is still effective or not. 

Another limitation of existing works is the inaccessibility of asymptotic analysis.
The large-sample analysis that characterizes the asymptotic behavior of iterates from the given iterative algorithm provides many practical instructions.
It not only helps to conduct hypothesis testing for scientific decision-making but also constructs confidence intervals for statistical inference.
However, only a few works consider such an aspect.
\citet{duchi2016local,asi2019stochastic} established asymptotic normality for SGD variants but didn't take the cost of performing projections into account.
Similar asymptotic normality has been established for Local SGD by~\citet{li2021statistical}. They showed the averaged iterates from Local SGD achieve the optimal asymptotic variance when communications are performed at an appropriate frequency.
However, all of these results neither reveal the effect of projection frequency on the asymptotic normality nor uncover the mechanisms of lazy projections.

To address these gaps in the existing literature, we pose two questions in our work:
\begin{itemize}
\item[(Q1)] Is the lazy projection approach effective for solving general linearly constrained problems?
\item[(Q2)] If the answer to (Q1) is affirmative, what are its asymptotic behaviors, and how does it work?
\end{itemize}
These questions are interdependent, and answering (Q1) requires developing an iterative algorithm that adheres to the lazy projection approach. 
Subsequently, we aim to answer (Q2) by analyzing the algorithm's asymptotic convergence.


The difficulty of the theoretical analysis is decided by the specific form of the considered algorithm. 
Unfortunately, the lazy projection circumvents the convenience of theoretical analysis due to its double-loop structure.
The running iterates behave differently in the inner loop where no projection is performed, while in the outer loop where projection is utilized to maintain feasibility.
However, recent developments have simplified this structure for variants of SGD by introducing a  so-called "loopless" technique, which essentially replaces the hard loop with a soft scheme~\citep{kovalev2020don}. 
At iteration $n$, we toss a (possibly biased) coin $\omega_n$ with  head probability $p_n$. 
When we get a head $\omega_n=1$, we start a new loop and update the outer-loop intermediate variables. 
When we get a tail $\omega_n=0$, we stay in the same loop and keep the intermediate variables. 
We obtain a loopless counterpart algorithm in this way and no longer need to distinguish between inner and outer loops. 
This homogeneity in performing loops facilitates theoretical analysis and typically does not deteriorate the convergence rate~\citep{hanzely2020federated,li2021page,li2021anita,gargiani2022page}. 
It is worth mentioning that~\citet{hanzely2020federated} first introduced the loopless technique to FL and obtained many efficient FL algorithms.
In addition, \citet{li2021anita} used a dynamic $p_n$ that varies with $n$ to generalize the scope of the original methods. 
This motivates us to analyze a loopless version of SGD methods with lazy projection with decreasing $p_n$.

To capture the asymptotic behavior of the proposed algorithm, we take a continuous perspective and use diffusion approximation to characterize the convergence of the trajectory of the discrete stochastic algorithm. 
This continuous approach has several advantages, including a rich toolbox and the ability to provide intuitive explanations for complex phenomena that are difficultly analyzed in discrete cases.
Prior works applying diffusion approximation to stochastic optimization algorithms can be roughly divided into two classes.
The first class interprets the optimization algorithm as a numerical discretization of a specific stochastic differential equation (SDE) in a finite time interval $[0, T]$~\citep{kloeden1991numerical}.
When the step size $\eta$ is small and the length $T (= n\eta)$ of the interval is fixed (where $n$ is the total iterations), such an approximation is of high accuracy, and it is easy to analyze the geometric properties of our target algorithms~\citep{wibisono2016variational,li2017stochastic,he2018differential,feng2020uniform,raginsky2017non,orvieto2019continuous,fontaine2021convergence}.
However, this avenue is difficult to capture the asymptotic behaviors around the optimal point due to the fixed length $T$.\footnote{A finite value of $T$ implies not only the algorithm but also its corresponding SDE do not converge to the optimum.} 
The second class addresses the issue by rescaling the iterates by a proper power function of step sizes and showing that under certain conditions the rescaled iterates weakly converge to the stationary solution of specific SDEs~\citep{kushner1981asymptotic,pelletier1998weak,fan2018statistical,gadat2017stochastic,gadat2018stochastic}. 
However, to the best of our knowledge, in linearly constrained optimization schemes, no previous work has analyzed projected stochastic gradient methods from this perspective, which is the focus of our work.
Moreover, the analysis techniques used in the previous works can not be immediately borrowed to our concerned issue.

\subsection{Contribution}

Formally, we analyze a loopless version of SGD methods with lazy projection with decreasing $p_n$, but from an asymptotic and continuous perspective in this work.
Our first theoretical contribution is to provide a complete answer to both (Q1) and (Q2).

\begin{algorithm}[t!]
\caption{Loopless Projected Stochastic Approximation (LPSA)}
 	\label{alg:lp_proj_sa}
 	\begin{algorithmic}
 		\STATE {\bfseries Input:} function $f$, data distribution $\gD$, initial point $\vx_0$, step size $\eta_n$, projection probability $p_n$.
 		\FOR{$n=1$ {\bfseries to} $T-1$}
 		 \STATE{Sample $\zeta_n\sim \gD$ and $\omega_n\sim \mathrm{Bernoulli}(p_n)$}
 		 \STATE{$\vx_{n+\frac{1}{2}}=\vx_n-\eta_n\nabla f(\vx_n,\zeta_n)$}
 		 \IF{$\omega_n=1$}
 		  \STATE{$\vx_{n+1}=\gP_{\mA^\bot}\vx_{n+\frac{1}{2}}$ \quad \# $\gP_{\mA^\bot}$ denotes the projection onto the null space of $\mA^\top$.}
 		 \ELSE
 		  \STATE{$\vx_{n+1}=\vx_{n+\frac{1}{2}}$}
 		 \ENDIF
 		 \ENDFOR
 		\STATE {\bfseries Return: $\gP_{\mA^\bot}\vx_T$.}
 	\end{algorithmic}
\end{algorithm}

\paragraph{LPSA algorithm and convergence analysis}
To answer (Q1), we propose Loopless Projected Stochastic Approximation (LPSA) in Algorithm~\ref{alg:lp_proj_sa}, a loop-less algorithm for solving~(\ref{eq:proj_opt}).
LPSA decides whether to perform projection by observing the realization of a Bernoulli random variable $\omega_n$ with head probability $p_n$ after performing one SGD step. 
This algorithm has two hyperparameters, namely the step size ${\eta_n}$ and the projection probability ${p_n}$.

Let $\vx^\star$ be the solution of (\ref{eq:proj_opt}).  We decompose the running iterate $\vx_n$ into the sum of two orthogonal components: $\vx_n = \vu_n + \vv_n$ where $\vu_n = \gP_{\mA^\bot}\vx_n$ and $\vv_n = 
(\rmI - \gP_{\mA^\bot})
\vx_n $.
We show in Theorem~\ref{thm:converge} that for $\eta_n \propto n^{-\alpha}$ and $p_n \propto \min\{ \gamma \eta_n^\beta, 1\}$, the algorithm converges with a non-asymptotic convergence rate in the sense that $\EB\|\vu_n-\vx^{\star}\|^2 = \gO(n^{- \alpha \min\{1,2-2\beta\}})$ and $\EB\|\vv_n\|^2 = \gO(n^{-\alpha(1-\beta)})$. 
Our focus is on $\vu_n$ because it always satisfies the linear constraint, while $\vv_n$ can be viewed as the error incurred by lazy projections.
From this theorem, the convergence rate of $\EB\|\vu_n-\vx^{\star}\|^2$ is stable across different values of $\beta$ in the high-frequency regime (i.e., $\beta < 0.5$), but begins to deteriorate as $\beta$ increases beyond 0.5. On the other hand, the convergence rate of $\EB\|\vv_n\|^2$ monotonically increases with $\beta$, indicating that the error incurred by lazy projections gradually diminishes in the increasing-frequency regime.


\paragraph{Phase transitions in asymptotic behaviors}
To answer (Q2), we provide a complete characterization for the asymptotic behavior of the iterates $\{\vu_n\}_{n=1}^{\infty}$ and $\{\vv_n\}_{n=1}^{\infty}$.
According to the particular value of $\beta \in [0, 1)$, our analysis and results are divided into three cases.
\begin{enumerate}
    \item[(i)] At the frequent projection regime where $\beta \in [0, 1/2)$, we first scale and center $\vu_n$ to $\check{\vu}_{n}:=\frac{\vu_n - \vx^\star}{\sqrt{\eta_{n-1}}}$ and then construct  piecewise constant interpolated processes $\{\Bar{\vu}_t^{(n)}\colon t\ge 0\}_{n=1}^\infty$ that pass through these points $\{\check{\vu}_{k}\}_{k \ge n}$. 
    In Theorem~\ref{thm:diff_approx}, we show that under some regularity conditions, the random function $\Bar{\vu}_t^{(n)}$ weakly converges to the stationary distribution of an SDE driven by the Brownian motion.
    As a result of $\Bar{\vu}_0^{(n)} = \check{\vu}_{n}$, we know that the random vector $\check{\vu}_n$ weakly converges to a scaled Gaussian distribution $\gN (\vzero, \tilde{\Sigma})$ where the variance $\tilde{\Sigma}$ satisfies the specific Lyapunov equation~\eqref{eq:variance-lya}.
    \item[(ii)] 
    At the occasional projection regime where $\beta \in (1/2, 1)$, the asymptotic behavior of $\vu_n$ starts to change.
     We introduce a new variable $\hat{\vu}_n:= \frac{1}{\eta_{n-1}^{1-\beta}}(\vu_n - \vx^\star)$, obtained by centering and scaling $\vu_n$ appropriately, and observe that $\vv_n$ now affects the update of $\hat{\vu}_n$ in a non-trivial way, leading to a non-diminishing asymptotic bias.
     To characterize the behavior of $\vv_n$, we similarly construct a piecewise constant interpolated process $\Bar{\vv}_t^{(n)}$ that passes through all the points $\{\check{\vv}_k\}_{k \ge n}$, with $\check{\vv}_n = \eta_{n-1}^{\beta-1}\vv_{n}$. Theorem~\ref{thm:jump_approx} shows that the random function $\Bar{\vv}_t^{(n)}$ converges weakly to the stationary distribution of an SDE driven by a Poisson process. 
     Using this result, we derive the asymptotic bias of $\hat{\vu}_n$ in Theorem~\ref{thm:asy_bias_u} and show that it converges to a vector $\vmu$ in the $L_2$ sense. 
     The vector $\vu$ is defined in~\eqref{eq:bias-u} and assumed to be non-zero.
    \item[(iii)] At the moderate projection regime where $\beta = 1/2$, we have $\hat{\vu}_n = \check{\vu}_{n}=\frac{\vu_n - \vx^\star}{\sqrt{\eta_{n-1}}}$, which is expected to have an intermediate asymptotic behavior.
    In Theorem~\ref{thm:noncen_u}, we show that $\frac{\vu_n - \vx^\star}{\sqrt{\eta_{n-1}}}$ weakly converges to a non-centered Gaussian distribution $\gN (\vmu, \tilde{\Sigma})$ where $\vu$ is the same vector given in Case (ii) and $\tilde{\Sigma}$ is the covariance matrix given in Case (i).
\end{enumerate}

The above phase transition introduces a bias-variance trade-off on the mean squared error $\EB\|\vu_n-\vx^{\star}\|^2$. 
When $\beta \in [0, 1/2)$, the variance of $\vu_n$ dominates that error, while the bias of $\vu_n$ becomes dominant once $\beta 
\in [1/2, 1)$.
It provides a consistent and more precise characterization of the asymptotic behaviors than the scalar non-asymptotic bound on $\EB\|\vu_n-\vx^{\star}\|^2$ in Theorem~\ref{thm:converge}.

\paragraph{Degenerate bias algorithm and its theoretical analysis}
The above analysis assumes that the vector $\vmu$ in the asymptotic bias is non-zero, but to fully answer (Q2), we must also consider the case where $\vmu$ is a zero vector. 
To simplify matters, we focus on linearly constrained quadratic minimization and investigate the behavior of LPSA when $\vmu$ is zero. 
We introduce the degenerate order $k$ in Definition~\ref{def:degenerate}, which identifies the order at which a non-zero asymptotic bias can still be observed, and allows us to establish both non-asymptotic and asymptotic results. 
In particular, Theorem~\ref{thm:cov_dgnrt_bias} provides more detailed convergence rates for both $\EB\|\vu_n-\vx^\star\|^2$ and $\EB\|\gP_{\mA^\bot} \mS\vv_n\|^2$, while Theorem~\ref{thm:degenerate_bias_u} characterizes the new phase transition in asymptotic behaviors for the degenerate case. 
Note that all of these results depend on the degenerate order $k$. 
The previous paragraph analyzes the case where $k=1$. 
As $k$ increases, the problem becomes more degenerate, and both $\vu_n$ and $\gP_{\mA^\bot} \mS\vv_n$ converge at faster rates. 
Finally, at the end of Section~\ref{sec:dgnrt_bias&homo_fed}, we discuss how to extend the results from linear cases to nonlinear ones.

We have observed that a degenerated bias can expedite convergence, and we investigate how to leverage this finding to develop a faster algorithm that requires fewer projections to achieve a target accuracy. Specifically, we notice that the non-zero expectation of $\vv_n$ contributes to the asymptotic bias of $\vu_n$, which impedes the convergence rate of $\vu_n$ when $\beta > 1/2$. Based on this observation, we propose a novel debiased algorithm, named Debiased LPSA (DLPSA), that subtracts an independent expectation estimator of $\vv_n$ from the stochastic gradient estimator used to update $\vu_n$. We demonstrate that DLPSA achieves a strictly better convergence rate than vanilla LPSA. Moreover, we analyze the projection complexity of both algorithms and find that the projection complexity of DLPSA is $\gO(\varepsilon^{-1/3})$ to achieve an $\varepsilon$-accurate solution, which is better than the $\gO(\varepsilon^{-1/2})$ of LPSA.

\paragraph{Technical contribution}
Our theoretical results are different than previous attempts in the sense that we focus on the trajectory that connects the appropriately scaled last iterates rather than point estimates~\citep{duchi2016local,asi2019stochastic} or partial-sum processes~\citep{li2021statistical}.
We construct continuous random functions that pass through appropriately rescaled and centered $\vu_n$ and $\vv_n$.
The idea motivating the constructions is that when the iteration $n$ goes to infinity, the considered processes will approximate a mean SDE with increasing accuracy.
To that end, we make use of operator semigroup theories to support this approximation~\citep{trotter1958approximation,kurtz1969extensions,kurtz1970general,kurtz1975semigroups,kushner1980martingale}.
From a technical level, we propose several novel proof techniques to analyze the discontinuity brought by probabilistic projection.
We borrow more related operator semigroup theories for jump-diffusion and verify three sufficient conditions (e.g., stochastic tightness) to apply them.
Among them, we provide a novel result for the stochastic tightness of c\`adl\`ag random functions.
To disentangle the dependence of $\vv_n$ in the update of $\vu_n$, we use a coupling argument to eliminate the effect of $\vv_n$.
See Section~\ref{sec:proof-idea} for the proof sketches.
We believe our technique can extend to and help analyze other stochastic approximation algorithms which can be approximated by a jump-diffusion.

\section{Problem Formulation}

The
optimization problem with linear constraint
(\ref{eq:proj_opt})
has been extensively studied, notably in the seminal work of Boyd et al. \cite{boyd2011distributed}. Although it is undoubtedly significant in many settings, we will provide some examples to help readers better understand this problem.

\begin{example}[General LASSO Problem]

The framework of the General LASSO model was proposed in \cite{she2010sparse, tibshirani2011solution}. And it has the form as follows:
\[
\underset{\boldsymbol{\theta}}{\operatorname{argmin}} \frac{1}{2}\|\mathbf{y}-\mathbf{X} \boldsymbol{\theta}\|_2^2+\lambda\|\mathbf{D} \boldsymbol{\theta}\|_1,
\]
where $\mX$ is the matrix  of covariates, and $\mY$ is the response vector. This model covers many of the specific problems on sparse recovery such as the fused lasso, trend filtering, wavelet smoothing, and a method for outlier detection.  Lemma~1 of \cite{james2020penalized} shows that the general LASSO problem has an equivalent version by simple algebraic transformation, namely, 
\[
\underset{\boldsymbol{\beta}}{\operatorname{argmin}} \frac{1}{2}\|\mathbf{y}-\widetilde{\mathbf{X}} \boldsymbol{\beta}\|_2^2+\lambda\|\boldsymbol{\beta}\|_1 \text { s.t. } \mathbf{A}^{\top} \boldsymbol{\beta}=\mathbf{0}.
\]
This new form can be regarded as one particular setting of the linear constrained optimization problem we are interested in.
\end{example}

\begin{example}[Distributed Learning]
\label{eg:distributed}
Typical distributed optimization problems can be rewritten as a global consensus problem
\begin{equation}\label{eq:dis_opt}
 \min\limits_{\vx^{(1)},\vx^{(2)},\dots,\vx^{(N)}}\frac{1}{N}\ssum{k}{1}{N}\sE_{\zeta^{(k)}\sim\gD_k} g_k (\vx^{(k)},\zeta^{(k)}) ~~\text{ s.t. } \vx^{(1)}=\cdots=\vx^{(N)},
\end{equation}
where $N$ is the number of clients, 
each $g_k$ is a local objective that involves only the data $\zeta^{(k)} \sim \gD_k$ residing on the client $k$ and $\vx^{(k)} \in \sR^d$ is the local parameter at this client.
This type of problem is commonly encountered in various fields, including machine learning, signal processing, and mobile communication \cite{boyd2011distributed, bertsekas2015parallel, zhu2010distributed}. 
If we concatenate all the local parameters and data points as $\vx=\left[(\vx^{(1)})^\top,(\vx^{(2)})^\top,\dots,(\vx^{(N)})^\top\right]^\top\in\sR^{Nd}$ and 
$\zeta = ( \zeta^{(1)}, \zeta^{(2)}, \dots, \zeta^{(N)} )^\top$,
the problem (\ref{eq:dis_opt}) can be formulated as a linear equality-constrained optimization problem (\ref{eq:proj_opt}), where
$f(\vx, \zeta) = \frac{1}{N} \sum_{k=1}^N g( \vx^{(k)}, \zeta^{(k)} )$,
$\gD = \gD_1 \times \gD_2 \times \cdots \times \gD_N$
and $\mA$ is equipped with a particular structure
\begin{equation}\label{eq:mat_A_fl}
    \mA = \mB \otimes \rmI_d \in \sR^{Nd\times (N-1)d}
    \text{ and }
    \mB^{\top}=\begin{bmatrix}
    1 & -1 & 0 & \cdots & 0 & 0 \\
    0 & 1 & -1 & \cdots & 0 & 0 \\
    \vdots & \vdots & \vdots & \ddots & \vdots & \vdots\\
    0 & 0 & 0 & \cdots & 1 & -1
    \end{bmatrix}\in \sR^{(N-1)\times N}.
\end{equation}
In the expression of $\mA$, $\otimes$ denotes the Kronecker product and $\rmI_d \in \sR^{d \times d}$ is the identity matrix.
\end{example}


\subsection{Loopless Projected Stochastic Approximation}\label{sec:lpsa_alg}

In fact, the linear constrained optimization problem~(\ref{eq:proj_opt}) can be solved via a randomly (and infrequently) projected stochastic approximation algorithm.
In particular, at iteration $n$, we first perform one step  SGD via
\begin{equation}\label{eq:update_rule}
    \vx_{n + \frac{1}{2}} = \vx_n-\eta_n\nabla f(\vx_n) +\eta_n\xi_n,
\end{equation}
where $f(\vx)=\sE_{\zeta\sim \gD}f(\vx,\zeta)$ and $\xi_{n}=\nabla f(\vx_n)-\nabla f(\vx_n, \zeta_n)$.
Here $\{\xi_n\}$ is a martingale difference sequence (m.d.s.) under the natural filtration $\gF_{n+1}:=\sigma(\zeta_k,\omega_k; k\le n+1)$.
We then use the loopless trick introduced in the introduction, i.e.,
we independently cast a coin with the head probability $p_n$ and obtain the result $\omega_n\sim \mathrm{Bernoulli}(p_n)$.
If $\omega_n = 1$, we perform one step of projection to ensure $\vx_{n+1}$ fall into the feasible region: $\vx_{n+1} = \gP_{\mA^\bot}(\vx_{n + \frac{1}{2}})$ where $\gP_{\mA^\bot}$ denotes the projection onto the null space of $\mA^\top$.
If $\omega_n = 0$, we assign $\vx_{n+1}$ as the same value of $\vx_{n + \frac{1}{2}}$, i.e., $\vx_{n+1}=\vx_{n + \frac{1}{2}}$.
It is clear this algorithm~\eqref{eq:update_rule} mimics the behavior of Local SGD in distributed learning settings.

\subsection{Assumptions}\label{sec:formulation-assume}

For the linearly constrained convex optimization problem~\eqref{eq:proj_opt}, we make the following assumptions which are quite common in the literature.
Without special clarification, $\|\cdot\|$ denotes the Euclidean norm for vectors and the spectral norm for matrices.
\begin{assumption}[Smoothness]
\label{asp:smooth}
We assume that $f\colon \sR^d \to \sR$ is $L$-smooth, that is, 
\begin{equation*}
\|\nabla f(\vx)-\nabla f(\vy)\|\le L\|\vx -\vy \|, \quad \forall \, \vx, \vy \in \sR^d.
\end{equation*}
\end{assumption}

\begin{assumption}[Strong convexity]
\label{asp:str_cov}
We assume that $f\colon \sR^d \to \sR$ is $\mu$-strongly convex, that is,
\begin{equation*}
    f(\vx)- f(\vy) \geq\langle\nabla f(\vy), \vx - \vy \rangle+\frac{\mu}{2}\|\vx - \vy\|^{2}, \quad \forall\, \vx, \vy \in \mathbb{R}^{d}.
\end{equation*}
\end{assumption}

\begin{assumption}[Continuous Hessian matrix]
\label{asp:hess_lip}
We assume that $f\colon \sR^d \to \sR$ is Hessian Lipschitz, that is, there is a constant $\tilde{L}$ such that
\begin{equation*}
\begin{aligned}
\left\|\nabla^2 f(\vx)-\nabla^2 f(\vy )\right\| \le \tilde{L}\|\vx - \vy\|, \quad \forall ~ \vx, \vy \in \sR^d.
\end{aligned}
\end{equation*}
\end{assumption}

\begin{assumption}[Continuous covariance matrix]
\label{asp:noi_lip}
Given an m.d.s.\ $\{\xi_t\}$, we denote the conditional covariance as $\sE[\xi_t\xi_t^\top|\gF_t]=\Sigma(\vx_t)$ and assume it is $L$-Lipschitz continuous in the sense that
\[ 
    \|\Sigma(\vx)-\Sigma(\vy)\|_2\le L\|\vx - \vy\|, \quad \forall \, \vx, \vy \in \sR^d.
\]
\end{assumption}

\begin{assumption}
\label{asp:high_noi_mom}
For the m.d.s.\ $\{\xi_n\}$, we assume there exists a $p>2$ such that the $p$-th moment of every element in $\{\xi_n\}$ is uniformly bounded, that is,
\begin{equation*}
\sup\limits_{n\ge 0}\sE\|\xi_n\|^p < \infty.
\end{equation*}
\end{assumption}

The first three assumptions imply we consider the strongly convex case.
The last two assumptions help us identify the asymptotic variance.
Especially, the assumption of uniformly bounded $p~(p>2)$ moments is typically required to establish central limit theorems \citep{gadat2017stochastic,gadat2018stochastic,li2021statistical}.
Finally, we want to emphasize that the stationary condition for the problem~\eqref{eq:proj_opt} is different from unconstrained ones.
Specially, we assume $\nabla f(\vx^\star) \neq \vzero$
where $\vx^\star$ denotes the solution of the problem (\ref{eq:proj_opt}).
This implies the constrained solution does not coincide with the unconstrained solution.
Otherwise, there is no need to perform the projection step and it suffices to apply vanilla SGD.
However, the projection of $\nabla f(\vx^\star)$ onto the null space of $\mA^\top$ must be zero, as claimed by the following proposition.

\begin{proposition}[\cite{li2021delayed}, Corollary 2.1]\label{prop:proj_A_bot_nabla}
\label{prop:f(x)}
Let $\gP_{\mA}$ be the projection onto the column space of $\mA$ and $\gP_{\mA^\bot}$ the projection onto the null space of $\mA^\top$.
Under Assumption \ref{asp:str_cov}, the solution of \eqref{eq:proj_opt} is unique (denoted  $\vx^\star$). Moreover, we have $\gP_{\mA^\bot} (\nabla f(\vx^\star)) = \vzero$.
\end{proposition}

\subsection{Jump Diffusion}
Jump diffusion is a stochastic L\'evy process that involves jumps and diffusion.
Typically, the former is modeled by a Poisson process, while the latter is modeled as a Brownian motion.
It has wide and important applications in physics, finance\citep{platen2010numerical}, and computer vision.

We say a function $f$ defined on $\sR$ to be c\`adl\`ag when $f$ is right-continuous and has left limits everywhere. 
For a c\`adl\`ag process $(\rmV_s)_{s \ge 0}$, we denote $\rmV_{t-}$ as the left limit of $\rmV_{(\cdot)}$ at time $t$.
Let $\rmN_\gamma(t)$ denote the Poisson process with $\gamma$ the intensity, which quantifies the number of jumps up to the time $t$ and is clearly c\`adl\`ag.
We use $\rmN_{\gamma}(dt) = \rmN_\gamma(t)-\rmN_\gamma(t-) \in \{0, 1\}$ to indicate whether $\rmN_{\gamma}$ jumps at time $t$ and $\int_0^T g(t-)\rmN(dt) = \sum_{\{t:\rmN_\gamma(t)\neq \rmN_\gamma(t-)\}}g(t-)$ to denote the integral that drives for a measurable function $g(\cdot)$.
We will consider a special class of jump diffusion in the following form
\begin{equation}\label{eq:jump_diff}
    d\rmX_t = \alpha(t,\rmX_t)dt + \beta(t,\rmX_t)d\rmW_t + \varphi(t,\rmX_{t-})\rmN_\gamma(dt).
\end{equation}
When the coefficient functions $\alpha(t,\rmX_t)$ and $\beta(t,\rmX_t)$ satisfy conditions like linear growth and Lipschitz continuity, there exists a solution for the jump diffusion~\eqref{eq:jump_diff} (e.g., Theorem 1.19 in \citep{oksendal2005stochastic}).

In the last of this subsection, we introduce the definition of tightness, which is an important nature of a sequence of measures.
\begin{definition}
    A sequence of probability measures $\{\mu_n\}_{n=1}^\infty$ in the measurable space $(\Omega, \gF)$ equipped with metric $\mathrm{d}(\cdot, \cdot)$ is said to be tight, if for each $\varepsilon$, there exists a compact set $\mathbf{K}$ such that $\mu_n(\mathbf{K}) \ge 1-\varepsilon, \; \forall n \ge 1$. 
\end{definition}

\section{Main Results}\label{sec:main_rst}

In this section, we aim to examine the convergence behaviors of our projected stochastic approximation method~\eqref{eq:update_rule} from both non-asymptotic and asymptotic perspectives. 
In particular, we consider a specific family of step sizes $\eta_n$ and projection probabilities $p_n$, given by $\eta_n = \eta_0 n^{- \alpha}$ and $p_n = \min \{ \gamma \eta_n^\beta, 1 \}$, respectively, where $0 < \alpha \le 1$ and $0 \le \beta < 1$ are the indices. 
The step size $\eta_n$ has been used to establish central limit theorems (CLTs)~\citep{polyak1963gradient,li2021statistical}, whereas the choice of $p_n$ is relatively new. 
To provide a comprehensive overview of convergence, we will examine almost all combinations of $\alpha$ and $\beta$.

\subsection{Non-asymptotic Analysis}
\label{sec:converge}
To provide a convergence rate, it is natural to focus on the projection of $\vx_n$ into the column space of $\mA$ since it is the easiest feasible solution obtainable from $\vx_n$. 
Thus, we decompose the iterated $\vx_n$ into two orthogonal components, namely, $\vx_n = \vu_n + \vv_n$, where $\vu_n = \gP_{\mA^\bot}(\vx_n)$ and $\vv_n = \gP_{\mA}(\vx_n)$.\footnote{Proposition \ref{prop:proj} in Appendix \ref{sec:append_converge} can be checked to understand why $\vu_n$ is orthogonal to $\vv_n$.} 
The following theorem specifies the convergence rate of $\sE \norm{ \vu_n - \vx^\star }^2$ and $\sE \norm{ \vv_n }^2$ in terms of $\alpha$, $\beta$, and $n$.

\begin{theorem}
\label{thm:converge}
Suppose that Assumptions \ref{asp:smooth}, \ref{asp:str_cov} and \ref{asp:noi_lip} hold. Let $\eta_n = \eta_0 n^{-\alpha}$ and $p_n = \min\{ \gamma \eta_n^\beta, 1 \}$ with $0 < \alpha \le 1$,
$0 \le \beta < 1$
and $\eta_0 > 2 \mathbbm{1}_{ \{ \alpha=1 \} } / \mu$ with $\mu$ the strong convexity parameter of the objective function $f$.
Then we have
\begin{align*}
   \sE \norm{ \vu_n - \vx^\star }^2 = \gO ( n^{- \alpha \min \{ 1 , 2 - 2 \beta \}} )
   \quad \text{and} \quad \sE \norm{ \vv_n }^2 = \gO ( n^{ - 2 \alpha (1 - \beta) } ).
\end{align*}

\end{theorem}

According to Theorem~\ref{thm:converge}, as $\beta$ decreases which means more frequent occurrences of projection, the convergence of $\sE \norm{ \vu_n - \vx^\star }^2$ accelerates.
The convergence rate is $\gO(n^{-\alpha})$ when $\beta < 0.5$, while it is $\gO(n^{-2\alpha(1-\beta)})$ when $\beta > 0.5$. 
Therefore, a phase transition occurs when $\beta$ crosses $0.5$, and we should analyze the asymptotic performances for these two phases separately. 
An interesting observation is that when $\beta = 1$, the algorithm can disconverge in an artifact quadratic loss with a specific $\mA$, as illustrated in Theorem~\ref{thm:counter}. 


\begin{theorem}\label{thm:counter}
If $\eta_n = \eta_0 n^{- \alpha}$ and $p_n = \min\{ \gamma \eta_n, 1 \}$ with $0 < \alpha \le 1$, for a carefully designed $\mA$, there exists a quadratic function $f(\vx)$ so that $\nabla^2 f(\vx) \succeq \rmI_d$ and  $\sE \norm{ \vu_n {-} \vx^\star }^2$ does not converge to $0$. Here $\rmI_d \in \RB^{d \times d}$ is the identity matrix, and  $\nabla^2 f(\vx) \succeq \rmI_d$ means $\nabla^2 f(\vx) {-} \rmI_d$ is positive semidefinite.
\end{theorem}

\subsection{Asymptotic Behavior of the Rescaled Trajectory}\label{sec:dynamics_approx}

In this section, 
our goal is to obtain the asymptotic convergence of the sequence $\{\vu_n\}_{n=1}^{\infty}$.
As mentioned in the previous subsection, the asymptotic behavior of $\{\vu_n\}_{n=1}^{\infty}$ changes abruptly when $\beta$ crosses $0.5$.
Therefore, we divide our analysis and results into three cases based on the value of $\beta$: $\beta \in [0, 1/2)$, $\beta = 1/2$, and $\beta \in (1/2,1)$.


\subsubsection{Case 1: Frequent Projection where $\beta \in [0, 1/2)$}
\label{sec:freq_proj}
From an asymptotic perspective, the typical central limit theorem states that $\check{\vu}_{n}:=\frac{\vu_n - \vx^\star}{\sqrt{\eta_{n-1}}}$ weakly converges to a rescaled standard distribution \citep{kushner1981asymptotic}, which helps us capture large-sample convergence behaviors and provides avenues for future statistical inference. 
However, we can derive an even stronger result that captures the asymptotic behavior of the trajectory passing these points $\check{\vu}_{n}$'s. 
Specifically, we serialize the sequence $\{ \check{\vu}_{n}\}_{n=1}^{\infty}$ by constructing a continuous random function denoted as $\Bar{\vu}_t^{(n)}$. 
This function starts from $\Bar{\vu}_0^{(n)} = \check{\vu}_{n}$ and as $t$ increases, it passes through $\check{\vu}_{n+1}, \check{\vu}_{n+2}$ and so on. 
We will show that such a random function $\Bar{\vu}_t^{(n)}$ converges weakly to the solution of a specific SDE, from which we can derive the asymptotic variance of $\check{\vu}_{n}$ and the trajectory's evolution.


We now provide details on how to construct the continuous function $\Bar{\vu}_t^{(n)}$, which interpolates all the points $\{ \check{\vu}_{k}\}_{k \ge n}$ via piecewise linear functions. 
To begin, we derive the one-step relation between $\check{\vu}_{n}$ and $\check{\vu}_{n+1}$ for illustration purposes.
In particular, 
\begin{align}
    \check{\vu}_{n+1} &= \check{\vu}_n -\eta_n \vb_n + \sqrt{\eta_n}\xi^{(1)}_n, \label{eq:relation-u} \\
    \vb_n &:= \gP_{\mA^\bot} \left(\nabla^2 f(\vx^\star) - \frac{1}{2\eta_0}\mathbbm{1}_{\{\alpha = 1\}}\rmI_d \right)\check{\vu}_n + \frac{1}{\eta_n}\gR_n,  \label{eq:b}
\end{align}
where  $\gR_n$ stands for a high-order residual error, and $\xi^{(1)}_n$ denotes the  component of noise $\xi_n$ on the null space of $\mA$.
We defer the specific but tedious derivation of~\eqref{eq:relation-u} to Appendix \ref{sec:prf_of_c1}. 
In essence, equation~\eqref{eq:relation-u} is a first-order Euler-Maruyama discretization with timescale $\eta_n$ for an SDE that begins at $\check{\vu}_n$, with a local drift coefficient $\vb_n$ and a local diffusion coefficient of $\var(\xi^{(1)}_n)$.

\begin{definition}[\textbf{Time interpolation}]\label{def:time_itpl}
Given a positive sequence $\iota= \{\iota_n\}_n^\infty$ that decreases to zero, we define the time interpolation $N(n,t,\iota)$ as what follows.
 For $n\in \sN, ~ t\ge 0$,
 \[ 
 N(n,t,\iota) = \min\limits_{m\in\sN}\left\{ m\ge n: \ssum{k}{n}{m}\iota_k > t\right\},  \Gamma_n(\iota)=\ssum{k}{1}{n-1}\iota_k, \; \mbox{and} \; \underline{t}_n(\iota)= \Gamma_{N(n,t,\iota)}-\Gamma_n.
 \]
\end{definition}
We then introduce a time interpolation for the formal description of the continuous function and further analysis. 
Intuitively, $N(n,t,\eta)$ represents the number of iterations $m$ at which the sum of step sizes $\sum_{k=n+1}^{m+1} \eta_k$ is just larger than $t$, and $\underline{t}_n(\eta)$ is the approximation of $t$ when we only use step sizes $\{\eta_k\}_{k \ge n}$. 
Since $\eta_n \to 0$, $\underline{t}_n(\eta) \to t$ as $n$ approaches infinity. 
A key property of Definition \ref{def:time_itpl} is that $N(n, \Gamma_{m}(\eta) - \Gamma_{n}(\eta), \eta) = m$ for any $m \ge n$. 
With this, we can construct $\Bar{\vu}_t^{(n)}$. 
Given $n\in \sN$, let $\Bar{\vu}_0^{(n)} = \check{\vu}_n$ and define for $t \ge 0$,
\begin{equation}\label{eq:itpl_u}
\begin{aligned}
   \Bar{\vu}_t^{(n)} &= \check{\vu}_n + \left\{\ssum{k}{n}{N(n,t,\eta)-1}\eta_k\vb_k + (t-\underline{t}_n(\eta))\vb_{N(n,t,\eta)}\right\}\\
   & \quad \quad \; + \left\{\ssum{k}{n}{N(n,t,\eta)-1}\sqrt{\eta_k}\xi_k^{(1)}+\sqrt{t-\underline{t}_n(\eta)}\xi_{N(n,t,\eta)}^{(1)}\right\}.
\end{aligned}
\end{equation}


From the construction, we can see that $\Bar{\vu}_{\underline{t}_k(\eta)}^{(n)} = \check{\vu}_{n+k}$.

\begin{theorem}[\textbf{Diffusion Approximation}]
\label{thm:diff_approx}
Let Assumptions \ref{asp:smooth}-\ref{asp:high_noi_mom} hold and $\eta_0 > 2 \mathbbm{1}_{ \{ \alpha=1 \} 
} / \mu$.
The continuous stochastic processes $\{\Bar{\vu}_t^{(n)}: t\ge 0\}_{n=1}^\infty$ weakly converges to the stationary weak solution of the following SDE:
\begin{equation}\label{eq:diff_approx}
    d\rmX_t=-\gP_{\mA^\bot}\left(\nabla^2 f(\vx^\star) - \frac{1}{2\eta_0}\mathbbm{1}_{\{\alpha = 1\}}\rmI_d \right) \gP_{\mA^\bot} \rmX_t dt +\gP_{\mA^\bot}\Sigma(\vx^\star)^{\frac{1}{2}}d\rmW_t.
\end{equation}
Further, the rescaled sequence $\left\{\check{\vu}_n\right\}_{n=1}^\infty$ converges weakly to the invariant distribution of the dynamics \eqref{eq:diff_approx}, i.e., $\gN (\vzero, \tilde{\Sigma})$. Here the variance matrix $\tilde{\Sigma}$ satisfies the following Lyapunov equation
\begin{equation}
\label{eq:variance-lya}
    \gP_{\mA^\bot}\left(\nabla^2 f(\vx^\star) {-} \frac{1}{2\eta_0}\mathbbm{1}_{\{\alpha = 1\}}\rmI_d \right) \gP_{\mA^\bot} \tilde{\Sigma} + \tilde{\Sigma} \gP_{\mA^\bot} \left(\nabla^2 f(\vx^\star) {-} \frac{1}{2\eta_0}\mathbbm{1}_{\{\alpha = 1\}}\rmI_d \right)\gP_{\mA^\bot} = \gP_{\mA^\bot}\Sigma(\vx^\star)\gP_{\mA^\bot}.
\end{equation}
\end{theorem}

\begin{remark}
By applying the continuous-time version of the Lyapunov theorem (Lemma 1 in \citep{wei1987multivariate}), a unique positive semidefinite solution to the Lyapunov equation (denoted $\tilde{\Sigma}$) can be obtained. 
Using Theorem \ref{thm:diff_approx} and Theorem 4.1.1 in \citep{gadat2017stochastic}, it can be concluded that for $\beta \in [0, \frac{1}{2})$, our algorithm LPSA achieves the same asymptotic variance as SGD with the same step size. 
Notably, the typical projected SGD corresponds to the case where $\beta = 0$, whereas LPSA allows $\beta$ to vary within the interval $[0, \frac{1}{2})$. 
Additionally, the projection frequency can be reduced by increasing $\beta$, which is equivalent to decreasing the probability $p_n$. Therefore, LPSA is more efficient in performing projections than projected SGD.
The advantages brought by more flexible and moderate projection frequency are particularly apparent when projections are expensive.

\end{remark}

\subsubsection{Case 2: Occasional Projection where $\beta\in (1/2, 1)$}\label{sec:occa_proj}
Upon entering the low-frequency regime where $\beta \in \left(\frac{1}{2}, 1\right)$, the situation undergoes a complete transformation. Intuitively, as LPSA undergoes significantly fewer projections, the parameters are updated using infeasible $\vx_t$, leading to the accumulation of residual errors that mislead the update direction. These errors not only dominate and hinder the non-asymptotic convergence rate (refer to Theorem~\ref{thm:converge}) but also alter the asymptotic behavior. 
Consequently, in this scenario, determining the appropriate timescale is not enough; we must also ascertain how these errors accumulate. 
To resolve this issue, we develop a novel analysis routine. Specifically, we consider $p_t = \gamma \eta_t^\beta$ with $\gamma>0$.

Our approach involves monitoring another random process that is linked with $\{\vv_n\}_{n=1}^{\infty}$ and acts as a bridge to derive the asymptotic behavior of $\{\vv_n\}_{n=1}^{\infty}$. 
The right scale should ensure that the scaled sequence has a non-vanishing expected $L_2$ norm. 
To achieve this, we define $\check{\vv}_n = \eta_{n-1}^{\beta-1}\vv_{n}$, as per Theorem \ref{thm:converge}. 
Furthermore, given $\check{\vv}_n$, the candidate value of $\check{\vv}_{n+1}$, prior to tossing the coin $\omega_n$, can be determined from LPSA's Algorithm \ref{alg:lp_proj_sa};
we refer to this candidate as $\check{\vv}_{(n+1)-}$. 


\begin{equation}
\label{eq:relation-v}
\begin{aligned}
    \check{\vv}_{(n+1)-}
    & := \check{\vv}_n - \eta_n^\beta \vd_n + \eta^\beta_n \xi^{(2)}_n,
\end{aligned}
\end{equation}
where $\vd_n = \nabla f(\vx^\star) + \eta_n^{-\beta}\gS_n$ with $\gS_n$ a residual error satisfying $\eta_n^{-\beta}\gS_n= o_{\PB}(1)$ (see Appendix~\ref{sec:prf_of_c2} for more details) and $\xi^{(2)}_n$ stands for the component of the noise $\xi_n$ on the orthogonal complementary space $A^\bot$.
Owing to the probabilistic projection, $\check{\vv}_{n+1}$ takes the value of $\check{\vv}_{(n+1)-}$ with probability $1-\gamma\eta_n^\beta$ and takes value zero with probability $\gamma\eta_n^\beta$. 
Following the definition used in the previous section, we construct a c\`adl\`ag random process $\bar{\vv}_t^{(n)}$ that originates from $\check{\vv}_n$ and passes through $\{ \check{\vv}_{(k)-}\}_{k \ge n}$.
Connecting these discrete points with a step function leads to the construction described below: $\bar{\vv}_{\underline{t}_n(\eta^\beta)}^{(n)}= \check{\vv}_{N(n,t,\eta^\beta)}$ and
\begin{equation}\label{eq:itpl_v}
\begin{aligned}
    \bar{\vv}_t^{(n)} = \bar{\vv}_{\underline{t}_n(\eta^\beta)}^{(n)} {-} \left(t {-} \underline{t}_n(\eta^\beta)\right) (\vd_{N(n,t,\eta^\beta)} {-} \xi_{N(n,t,\eta^\beta)}^{(2)}) ~\text{for}~t\in \left(\underline{t}_n(\eta^\beta), \underline{t}_n(\eta^\beta) {+} \eta_{\underline{t}_n(\eta^\beta)}^\beta\right).
    \end{aligned}
\end{equation}

We can deduce from \eqref{eq:itpl_v} that $\bar{\vv}_{\underline{t}_n(\eta^\beta)-}^{(n)} = \check{\vv}_{N(n,t,\eta^\beta)-}$ holds true for any $t \ge 0$. 
At a probability of $p_{N(n,t,\eta^\beta)}$, the value of $\check{\vv}_{N(n,t,\eta^\beta)}$ is zero, which leads to abrupt changes in the process $\{\bar{\vv}^{(n)}_{t}:t \ge 0\}$ at the time $\underline{t}_n(\eta^\beta)$. 
The discontinuities in $\bar{\vv}^{(n)}_{\cdot}$ hinder the diffusion process from being applicable, as Theorem \ref{thm:diff_approx} is primarily concerned with continuous functions.
However, the following theorem demonstrates that we can still identify an appropriate process within the broader class of \textit{jump} diffusion processes to approximate $\{\bar{\vv}^{(n)}_{t}:t \ge 0\}$.

\begin{theorem}[\textbf{Jump Approximation}]\label{thm:jump_approx}
Let Assumptions \ref{asp:smooth}, \ref{asp:str_cov}, \ref{asp:noi_lip} and \ref{asp:high_noi_mom} hold and $\eta_0 > 2 \mathbbm{1}_{ \{ \alpha=1 \} 
} / \mu$.
The following c\`adl\`ag stochastic processes $\{\Bar{\vv}_t^{(n)}: t\ge 0\}_{n=1}^\infty$ weakly converges to the stationary weak solution of the following SDE
\begin{equation}\label{eq:jump_approx}
    d\rmY_t=-\nabla f(\vx^\star)dt - \rmY_{t-}\cdot\bfN_\gamma(dt).
\end{equation}
Here $\bfN_\gamma(t)$ represents Poisson process with intensity $\gamma$, and $\bfN_\gamma(dt) = \bfN_\gamma(t) - \bfN_\gamma(t-)$.
Further, the rescaled sequence $\left\{\check{\vv}_n\right\}_{n=1}^\infty$ weakly converges to the invariant distribution of the dynamics \eqref{eq:jump_approx}, i.e., $ -\frac{\nabla f(\vx^\star)}{\|\nabla f(\vx^\star)\|}\cdot\gE \left(\frac{\|\nabla f(\vx^\star)\|}{\gamma}\right)$
where $\gE(\theta)$ represents the exponential distribution with intensity $\frac{1}{\theta}$.
\end{theorem}

Theorem \ref{thm:jump_approx} demonstrates that the sequence $\{\bar{\vv}^{(n)}_t\}_{n=1}^{\infty}$, which is obtained by shifting initial points, will ultimately approximate a jump process with a constant drift as $n$ approaches infinity. 
The SDE~\eqref{eq:jump_approx} provides insight into the movement of $\bar{\vv}^{(n)}_t$ (or a rescaled version of $\vv_n$) as $t$ increases. 
Due to the error caused by infrequent projections, $\bar{\vv}^{(n)}_t$ will move towards the direction of $\nabla f(\vx^\star)$ (caused by the drift term $-\nabla f(\vx^\star)dt$) and will periodically be forced to reset as a zero vector (due to the correcting term $- \rmY_{t-}\cdot\bfN_\gamma(dt)$). 
From a qualitative perspective, the SDE~\eqref{eq:jump_approx} captures the periodic behavior of $\vv_n$, which shows that without projection, the residual error will accumulate along the direction of $\nabla f(\vx^\star)$. 
As discussed in Section~\ref{sec:formulation-assume},
$\nabla f(\vx^\star)$ is not
zero in our constrained problems.


The remaining issue is how to establish a link between $\{\Bar{\vv}_t^{(n)}: t\ge 0\}_{n=1}^{\infty}$ and our target sequence $\{\vu_n\}_{n=1}^{\infty}$.
Similar to the construction of $\check{\vv}_n$,  we need to consider a rescaled $\vu_n$, denoted as, $\hat{\vu}_n: = (\vu_n - \vx^\star) / \eta_{n-1}^{1-\beta}$.
Due to the technical complexity, we provide a high-level overview of our approach below. 
Firstly, we observe that the update rule of $\check{\vu}_{n+1}$ linearly depends on $\check{\vu}_{n}$ and $\check{\vv}_{n}$ up to a high order approximation error (see Appendix~\ref{sec:prf_of_cs3}). Secondly, we leverage Theorem~\ref{thm:jump_approx} to determine the asymptotic behavior of $\check{\vv}_{n}$, which is found to contribute to the asymptotic bias of $\hat{\vu}_n$. 
Finally, by iteratively applying the update rule for $\hat{\vu}_n$ and analyzing the particular form of the asymptotic bias, we complete the analysis for $\hat{\vu}_n$.
As a result, we have the following Theorem that is based on Theorem \ref{thm:jump_approx} and indicates that when $\beta \in \left(\frac{1}{2}, 1\right)$, $\hat{\vu}_n$ converges to a non-zero vector.


\begin{theorem}
\label{thm:asy_bias_u}
Let Assumptions \ref{asp:smooth}- \ref{asp:high_noi_mom} hold and set $\eta_0 > 2 \mathbbm{1}_{ \{ \alpha=1 \} } / \mu$. 
Assume the following vector $\vmu$ is non-zero where
\begin{equation}
\label{eq:bias-u}
\vmu = \frac{1}{\gamma}\left\{\gP_{\mA^\bot}\left(\nabla^2 f(\vx^\star) - \frac{1-\beta}{\eta_0}\mathbbm{1}_{\{\alpha = 1\}}\rmI\right)\gP_{\mA^\bot}\right\}^{\dag}\left(\gP_{\mA^\bot}\nabla^2 f(\vx^\star) \nabla f(\vx^\star)\right)
\end{equation}
 and $\rmG^\dag$ denotes the pseudoinverse of the symmetric matrix $\rmG$.
Then $\hat{\vu}_n := \frac{1}{\eta_{n-1}^{1-\beta}}(\vu_n - \vx^\star)$ converges to the non-zero vector $\vmu$ in the $L_2$ as $n\to \infty$.
\end{theorem}
Recall that $\check{\vu}_{n} =\frac{\vu_n - \vx^\star}{\sqrt{\eta_{n-1}}} = \eta_{n-1}^{0.5 - \beta} \hat{\vu}_n$.
By using this equation along with Theorem~\ref{thm:asy_bias_u}, it can be deduced that $\norm{\EB \check{\vu}_n } = \eta_{n-1}^{0.5 - \beta} \norm{ \EB \hat{\vu}_n } \rightarrow \infty$.
It implies the original $\check{\vu}_{n}$ diverges in the $L_2$ sense.
Furthermore, the bias in Theorem \ref{thm:asy_bias_u} rather than the Gaussian fluctuation in Theorem \ref{thm:diff_approx} becomes the leading term hindering the convergence.


\subsubsection{Case 3: Moderate Projection where $\beta = 1/2$}\label{sec:appr_proj}
The most subtle case is when $\beta$ falls exactly on the boundary of the two regimes discussed earlier that is $\beta = 1/2$.
In this moderate projection regime, we have $\check{\vu}_n =\hat{\vu}_n =\frac{\vu_n - \vx^\star}{\sqrt{\eta_{n-1}}}$. 
At the boundary where the phase transition occurs, we expect the behavior of $\frac{\vu_n - \vx^\star}{\sqrt{\eta_{n-1}}}$ to be similar to both cases. 

\begin{theorem}
    \label{thm:noncen_u}
Under the same conditions of Theorem~\ref{thm:asy_bias_u}, $\frac{\vu_n - \vx^\star}{\sqrt{\eta_{n-1}}}$ converges weakly to a non-central Gaussian distribution $\gN(\vmu, \tilde{\Sigma})$, where the mean vector $\vmu$ is defined in~\eqref{eq:bias-u} and the covraince matrix $\tilde{\Sigma}$ satisfies the Lyapunov equation~\eqref{eq:variance-lya}.
\end{theorem}

The asymptotic behavior of the sequence $\{\check{\vu}_n\}_{n=1}^{\infty}$ changes with the value of $\beta$. 
When $\beta \in [0, 1/2)$, $\check{\vu}_n=\frac{\vu_n - \vx^\star}{\sqrt{\eta_{n-1}}}$ converges weakly to the mean-zero Gaussian distribution $\gN(\vzero, \tilde{\Sigma})$, where the covariance matrix $\tilde{\Sigma}$ satisfies the Lyapunov equation~\eqref{eq:variance-lya}.
On the other hand, when $\beta \in (1/2, 1)$, $\hat{\vu}_n=\frac{\vu_n - \vx^\star}{\eta_{n-1}^{1-\beta}}$ converges to a non-zero vector $\vmu$ in the $L_2$ sense.
When $\beta = 1/2$, the limiting distribution of $\check{\vu}_n =\hat{\vu}_n =\frac{\vu_n - \vx^\star}{\sqrt{\eta_{n-1}}}$ weakly converges to a non-centered Gaussian distribution $\gN(\vmu, \tilde{\Sigma})$, as shown in Theorem~\ref{thm:noncen_u}. 
This can be viewed as an intermediate case of two asymptotic behaviors described above.

The phase transition in the asymptotic behavior of the sequence $\{\check{\vu}_n\}_{n=1}^{\infty}$ and $\{\hat{\vu}_n\}_{n=1}^{\infty}$ leads to an interesting bias-variance trade-off on the mean squared error $\EB\|\vu_n-\vx^{\star}\|^2$. 
When $\beta\in [0, 1/2)$, the variance of $\vu_n$ dominates the squared error because $\check{\vu}_n$ is asymptotically unbiased, resulting in an error of order $\gO( \tr(\tilde{\Sigma}) \eta_{n})$. 
However, when $\beta = 1/2$, $\check{\vu}_n$ starts to have a 
bias of the order $\gO( \norm{\vmu} \sqrt{\eta_n} )$,
whose square has the same order as the variance term.
Thus, the mean squared error $\EB\|\vu_n-\vx^{\star}\|^2$ remains of order $\gO( 
( \tr( \tilde{\Sigma} ) + \norm{\vmu}^2 ) \eta_{n})$. 
In contrast, when $\beta \in (1/2, 1)$, the bias of $\hat{\vu}_n$ increases to the order of $\gO( \norm{\vmu} \eta_n^{1-\beta})$ and dominates the error. 
Then the mean squared error is of the order $\gO( \norm{\vmu}^2 \eta_n^{2(1 - \beta)} )$.
These analyses are consistent with the non-asymptotic result Theorem~\ref{thm:converge}.

\subsection{High-Level Proof Sketch of Main Results}
\label{sec:proof-idea}
The proof of the theorems introduced in this section is provided in the appendix.
In this section, we outline the main ingredients of the proofs.

\begin{figure}[t!]
\begin{tikzpicture}[>=stealth,every node/.style={shape=rectangle,draw,rounded corners, minimum width=3.5cm,},]
    \node[very thick] (l2) {\begin{tabular}{c} Convergence of \\ operator semigroups \\
     \citep{kurtz1970general}\end{tabular}};
    \node[very thick] (l3) [right=of l2]{\begin{tabular}
    {c}Tightness of $\bar{\vu}_t^{(n)}$ \\ Theorem 7.3 in \\ \citep{billingsley2013convergence} \end{tabular}} ;
    \node[very thick, fill=orange!30] (t6) [right=of l3]{\begin{tabular}{c}
    Theorem~\ref{thm:noncen_u}
    \\ non-central normality \\ when $\beta = 0.5$ \end{tabular}};
    \node[very thick, fill=orange!30] (t1) [below=of l2]{\begin{tabular}{c} 
    Theorem~\ref{thm:converge}
    \\ convergence rate\end{tabular}};
    \node[very thick, fill=red!30] (t3) [right=of t1]{\begin{tabular}{c}
    Theorem~\ref{thm:diff_approx}\\ diffusion\\ approximation
    \end{tabular}};
    \node[very thick, fill=orange!30] (l4) [right=of t3]{\begin{tabular}{c}
    Lemma~\ref{lem:bd_innr_pdct}\\ asymptotic uncorrelation\\ between $\check{\vu}_n, \check{\vv}_n$ \end{tabular}};
    \node[very thick, fill=orange!30] (l5) [below=of t1]{\begin{tabular}{c}
    Lemma~\ref{lem:tight_v_smp}\\ tightness of \\ c\`adl\`ag process \end{tabular}};
    \node[very thick, fill=red!30](t4)[right=of l5]{\begin{tabular}{c}
         Theorem~\ref{thm:jump_approx}\\jump diffusion\\ approximation \end{tabular}};
    \node[very thick, fill=orange!30](t5)[right=of t4]{\begin{tabular}{c}
         Theorem~\ref{thm:asy_bias_u}\\ asymptotic bias of $\hat{\vu}_n$\\
         when $\beta > 0.5$ \end{tabular}};
    \node[very thick](l6)[below=of t4]{\begin{tabular}{c}
        Convergence of\\ operator semigroups\\ for jump diffusion\\ \citep{kushner1980martingale} \end{tabular}};
    \draw[->, line width=.4mm] (l2) to[out=345,in=155] (t3);
    \draw[->, line width=.4mm] (l3) to[out=270,in=90] (t3);
    \draw[->, line width=.4mm] (t1) to[out=0,in=180] (t3);
    \draw[->, dashed, line width=0.4mm] (t1) to[out=345,in=155] (t4);
    \draw[->, dashed, line width=0.4mm] (l5) to[out=0,in=180] (t4);
    \draw[->, dashed, line width=0.4mm] (l6) to[out=90,in=270] (t4);
    \draw[->, line width=.4mm] (t3) to[out=25,in=185] (t6);
    \draw[->, line width=.4mm] (t4) to[out=25,in=200] (t6);
    \draw[->, line width=.4mm] (t4) to[out=0,in=180] (t5);
    \draw[->, line width=.4mm] (l4) to[out=90,in=270] (t6);
    \draw[->, line width=.4mm] (l4) to[out=270,in=90] (t5);
\end{tikzpicture}
  \caption{Proof framework for results in Section~\ref{sec:main_rst}
  }
    \label{fig:proof_sketch}
\end{figure}
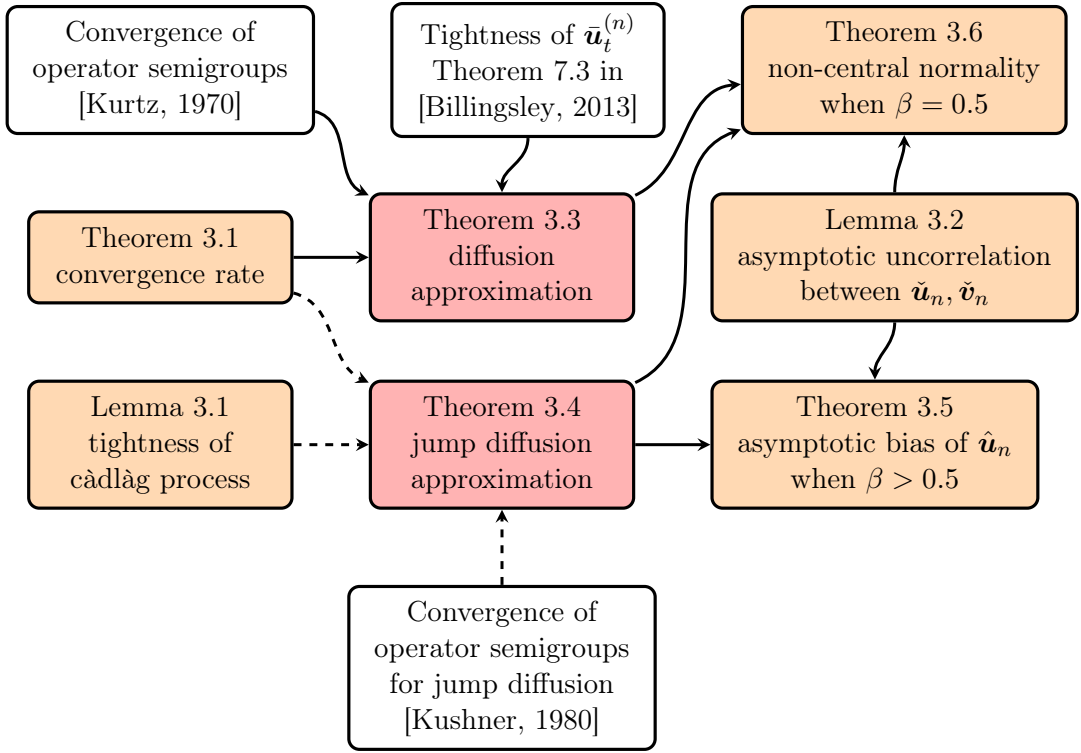

\subsubsection{Proof idea of Theorem~\ref{thm:diff_approx}}
By the definition of $\check{\vu}_n$ in~\eqref{eq:relation-u} and the construction of the process $\bar{\vu}_t^{(n)}$ in ~\eqref{eq:itpl_u}, for any test function $g: \RB^d \to \RB$ that is compactly supported and has Lipschitz continuous second derivatives, we can show that
\begin{equation}
\begin{aligned}
    \forall t \ge 0, \quad g\left(\bar{\vu}_{t}^{(n)}\right)-g\left(\bar{\vu}_{0}^{(n)}\right)-\int_{0}^{t} \mathcal{L} g\left(\bar{\vu}_{s}^{(n)}\right) d s=\mathcal{M}_{t}^{(n, g)}+\mathcal{R}_{t}^{(n, g)},
\end{aligned}
\end{equation}
where $\mathcal{L}$ is the infinitesimal generator of the stochastic process in~\eqref{eq:diff_approx}, $\gM_t^{(n,g)}$ is a $\gF_t^{(n)}$-martingale with $\gF^{(n)}_t$ the natural filtration of $\Bar{\vu}^{(n)}_t$, and $\gR_t^{(n,g)}$ denotes a reminder term.
Since we fix a test function $g$, both $\mathcal{M}_{t}^{(n, g)}$ and $\mathcal{R}_{t}^{(n, g)}$ depend on $g$ implicitly.

To prove Theorem~\ref{thm:diff_approx}, we mainly use the results of operator semigroups theories developed by Trotter and Kurtz~\citep{trotter1958approximation,kurtz1969extensions,kurtz1970general,kurtz1975semigroups}.
In brief, their results in our context imply that if the sequence of the stochastic processes $\{\bar{\vu}_{t}^{(n)}: t \ge 0\}_{n \ge 0}$ satisfies the following three conditions,
\begin{enumerate}
    \item The sequence of processes $\{\bar{\vu}_{t}^{(n)}: t \ge 0\}_{n \ge 0}$ itself is tight,
    \item The sequence of infinitesimal generators converges, that is, $\mathcal{R}_{t}^{(n, g)} = o_{\PB}(1)$ as $n \to \infty$,
    \item The sequence of the initial state converges weakly, that is $\bar{\vu}_{0}^{(n)}$ weakly converges to some distribution $\pi$,
\end{enumerate}
then the stochastic process $\bar{\vu}_{t}^{(n)}$ as a random function converges weakly to the stochastic process driven by the infinitesimal generator $\mathcal{L}$ with the initial distribution $\pi$.

Most efforts are then paid to validate the three conditions in our case.
First, the classic result from \citet{billingsley2013convergence} (see Theorem 7.3 therein) guarantees that $\{\bar{\vu}_{t}^{(n)}: t \ge 0\}_{n \ge 0}$ is tight once it has stochastic equicontinuity and its initial values are also tight, which we prove in Lemma~\ref{lem:tight_u}.
Second, by the smoothness of the test function $g$ and the non-asymptotic rate $\EB\|\vu_t-\vx^{\star}\|^2$ in Theorem~\ref{thm:converge}, we can show $\EB|\gR_t^{(n,g)}| \to 0$ as $n \to \infty$.
It implies $\mathcal{R}_{t}^{(n, g)} = o_{\PB}(1)$ indeed holds due to Markov's inequality.
Finally, we show $\bar{\vu}_{0}^{(n)}$ weakly converges the stationary distribution of the stochastic process in~\eqref{eq:diff_approx} by leveraging the Prokhorov’s theorem.
See Appendix~\ref{sec:prf_of_c1} for detailed proof.

Though some works also utilize the three conditions to analyze asymptotic behaviors of stochastic iterative methods such as stochastic gradient descent~\citep{gadat2017stochastic} or stochastic heavy-ball methods~\citep{gadat2018stochastic}, we are the first to analyze the loopless projected stochastic method which is designed for linearly constrained problems, while previous works~\citep{gadat2017stochastic,gadat2018stochastic} consider unconstrained problems.


\subsubsection{Proof idea of Theorem~\ref{thm:jump_approx}} 
Similar to the proof of Theorem~\ref{thm:diff_approx}, we again use the operator semigroups theories to establish the weak convergence in Theorem~\ref{thm:jump_approx}. 
By the definition of $\check{\vv}_n$ in~\eqref{eq:relation-v} and the construction of the process $\bar{\vv}_t^{(n)}$ in ~\eqref{eq:itpl_v}, for any fixed test function $g$, we also have the following martingale decomposition
\begin{equation}
    \forall t>0, \quad g(\bar{\vv}_t^{(n)})-g(\bar{\vv}_0^{(n)})-\int_0^t \gJ g\left(\bar{\vv}_s^{(n)}\right)ds = \gN_t^{(n,g)} + \gT_t^{(n,g)},
\end{equation}
where $\mathcal{J}$ is the infinitesimal generator of the stochastic process in~\eqref{eq:jump_approx},  $\gN_t^{(n,g)}$ is a $\gD_t^{(n)}$-martingale with $\gD_t^{(n)}$ the natural filtration of $\bar{\vv}_t^{(n)}$, and $\gT_t^{(n,g)}$ is a reminder term.
Since we fix a test function $g$, both $\mathcal{N}_{t}^{(n, g)}$ and $\mathcal{T}_{t}^{(n, g)}$ depend on $g$ implicitly.

We then check the three similar conditions introduced earlier.
However, due to the discontinuity of $\{\bar{\vv}_t^{(n)}: t \ge 0 \}_{n \ge 0}$, many previous arguments are not applicable.  
For example, Theorem 7.3 in \citep{billingsley2013convergence} fails to support the tightness of $\{\bar{\vv}_t^{(n)}: t \ge 0 \}_{n \ge 0}$ due to the discontinuous issue.
Instead, we provide original proof in Lemma~\ref{lem:tight_v_smp} for the tightness of any c\`adl\`ag processes with the help of Theorem 4.1 in \citep{kurtz1975semigroups}.

\begin{lemma}[Tightness of $\{\bar{\vv}_t^{(n)}\}$]\label{lem:tight_v_smp}
    Given Assumptions~\ref{asp:smooth}, \ref{asp:str_cov} and \ref{asp:noi_lip}, the sequence $\{\bar{\vv}_t^{(n)}\}$ is tight.
\end{lemma}

On the other hand, to show the convergence to the discontinuous infinitesimal generators $\mathcal{J}$, we make use of Theorem 1 in \citep{kushner1980martingale} that tailored the general results developed in~\citep{kurtz1975semigroups} to jump diffusions.
As a result, we only need to show that $\gT_t^{(n,g)}$ converges to zero in $L_1$, which is true by the smoothness of the test function $g$ and the non-asymptotic rate $\EB\|\vv_t\|^2$ in Theorem~\ref{thm:converge}.
For the last condition, we first show that there exists a unique invariant distribution for the Lévy process in~\eqref{eq:jump_approx} and then show that $\bar{\vv}_0^{(n)}$ weakly converges to this particular distribution by its tightness and the Prokhorov’s theorem.
See Appendix~\ref{sec:prf_of_c2} for detailed proof.

\subsubsection{Proof idea of Theorem~\ref{thm:asy_bias_u}}
The proof is composed of two steps; we first show that $\EB \hat{\vu}_{n} \to \vmu$ and then complete the proof with the statement that $\mathrm{Var}(\hat{\vu}_n) \to 0$.

In the first step, doing some algebra to the update rule of $\hat{\vu}_n$ and taking expectation, we have
\begin{equation*}
\begin{aligned}
 \sE \hat{\vu}_{n+1} &= \left(\rmI - \eta_n \left(\gP_{\mA^\bot}\nabla^2 f(\vx^\star) -\frac{1-\beta}{\eta_0}\mathbbm{1}_{\{\alpha=1\}}\rmI_d \right) \gP_{\mA^\bot} \right)\sE \hat{\vu}_n
 -\eta_n \gP_{\mA^\bot} \nabla^2 f(\vx^\star) \sE \check{\vv}_n + o(\eta_n)
\end{aligned}
\end{equation*}
where $o(1)$ denotes an infinitesimal term when $n \to \infty$.
Subtracting $\vmu$ from both sides, we arrive at
\begin{equation*}
\begin{aligned}
 \sE\hat{\vu}_{n+1} - \vmu &= \left(\rmI - \eta_n \gP_{\mA^\bot}\left(\nabla^2 f(\vx^\star)-\frac{1-\beta}{\eta_0}\mathbbm{1}_{\{\alpha=1\}}\rmI_d\right) \gP_{\mA^\bot}\right)(\sE \hat{\vu}_n -\vmu)\\ 
 & \quad - \eta_n \gP_{\mA^\bot} \nabla^2 f(\vx^\star) (\sE (\check{\vv}_n - \vnu) - o(\eta_n)
\end{aligned}
\end{equation*}
where $\vnu = \frac{1}{\gamma} \nabla f (\vx^{\star})$.
From Theorem~\ref{thm:jump_approx}, we have $\|\sE \check{\vv}_n - \vnu\| = o(1)$.
Disentangling this iteration, we have that $\RB \hat{\vu}_{n} \to \vmu$ (see Lemma~\ref{lem:converge_r_t_3} for the rationale). 

For the second step, we find an interesting intermediate result and present it in Lemma~\ref{lem:bd_innr_pdct}.
\begin{lemma}\label{lem:bd_innr_pdct}
When $\beta > 1/2$, under appropriate conditions, we have
\[
\lim\limits_{n \to \infty}\left| \EB \inner{\hat{\vu}_n - \EB \hat{\vu}_n}{\nabla^2 f(\vx^\star) (\check{\vv}_n - \EB \check{\vv}_n)}\right| = 0,
\]
where the definitions of $\hat{\vu}_n$ and $\check{\vv}_n$ are consistent with one in Theorem~\ref{thm:asy_bias_u}.
\end{lemma}
As Lemma~\ref{lem:bd_innr_pdct} shows, $\hat{\vu}_n$ and $\check{\vv}_n$ are asymptotically uncorrelated, i.e., $\mathrm{Cov}(\hat{\vu}_n, \check{\vv}_n)$ converges to zero.
This lemma helps us to single out the effect of $\check{\vv}_n$ on the update of $\hat{\vu}_n$, given the fact that $\check{\vu}_{n+1}$ linearly depends on $\check{\vu}_{n}$ and $\check{\vv}_{n}$ up to a high order approximation error.
As a result, we can show that
\[
\mathrm{Var}(\hat{\vu}_{n+1}) \le (1-\mu \eta_n/2) \mathrm{Var}(\hat{\vu}_n) + o(\eta_n)
\]
from which one can show that $\mathrm{Var}(\hat{\vu}_n) \to 0$.


\subsubsection{Proof idea of Theorem~\ref{thm:noncen_u}}
When $\beta= 0.5$, the orders of the bias term and fluctuation term are the same. 
Hence, $\check{\vu}_n$ will not converge to a Dirac measure at the point $\vmu = \lim\limits_{n\to \infty} \EB \check{\vu}_n$. 
To figure out its distribution, we apply a coupling technique by constructing a surrogate sequence $\{\vz_n\}$, which has almost the same update formula that shares the randomness of the gradient noise with $\{\check{\vu}_n\}$. 
The only difference is that the update of $\{\vz_n\}$ does not rely on $\{\check{\vv}_n\}$ anymore. 
Once eliminating the dependence on $\{\check{\vv}_n\}$, we know that $\vz_n$ converges weakly to $\mathcal{N}(\boldsymbol{0}, \Sigma)$ by a similar argument in proving Theorem~\ref{thm:diff_approx}.
Similar to what we did in analyzing $\mathrm{Var}(\hat{\vu}_n)$ in the proof of Theorem~\ref{thm:asy_bias_u}, we can show $\lim\limits_{n\to \infty}\EB \norm{(\check{\vu}_n - \vmu) - \vz_n}_2^2 = 0$ by leveraging Theorem~\ref{thm:jump_approx} and Lemma~\ref{lem:bd_innr_pdct}.
This result implies the Wasserstein-2 distance between $\check{\vu}_n - \vmu$ and $\vz_n$ is diminishing. 
Therefore, we conclude that $\check{\vu}_n$ has the distribution $\mathcal{N}(\vmu, \Sigma)$ and complete the proof.

\section{The Case of Bias Degeneration}
\label{sec:dgnrt_bias&homo_fed}
In the previous section, we have discussed the asymptotic behaviors of LPSA under three different cases according to the value of $\beta$.
An implicit assumption made in Theorems~\ref{thm:asy_bias_u} and \ref{thm:noncen_u} is that the vector $\vmu$ defined in~\eqref{eq:bias-u} that represents the asymptotic bias is non-zero.
However, once $\gP_{\mA^\bot}\nabla^2 f(\vx^\star)\nabla f(\vx^\star) = \0$, $\vmu$ would definitely become a zero vector.
In this degenerate case, the scaled iterate $\hat{\vu}_n = \frac{1}{ \eta_{n-1}^{1 - \beta} } (\vu_n - \vx^\star)$ defined in Theorem~\ref{thm:asy_bias_u} is no longer seriously biased as $\EB\|\hat{\vu}_n\|^2 \to 0$.
It implies that for $\beta \in \left( 1/2, 1 \right)$, the scale $\frac{1}{ \eta_{n-1}^{1 - \beta}}$ used in the definition of $\hat{\vu}_n$ is not the right scale to derive non-trivial asymptotic behaviors.
Furthermore, $1/2$ is no longer the change point where the phase transition occurs.
It would be totally unclear how the iterates $\{\vu_n\}_{n \ge 0}$ would behave when $\beta \in [1/2, 1)$.

This section is devoted to the degenerate case where $\gP_{\mA^\bot} \nabla^2 f(\vx^\star) \nabla f(\vx^\star) = \vzero$.
This case is intriguing even for the linear regression or equivalently quadratic loss functions.
In Section~\ref{sec:dgnrt_convergence}, we present a refined convergence analysis for the linear regression with linear constraints and analyze the optimal choice of projection frequency for the degenerate case.
In Section~\ref{sec:nonlinear}, we discuss how to extend our results for linear regressions to non-linear cases.

\subsection{Convergence Analysis}\label{sec:dgnrt_convergence}
In this subsection, we present a more complete theoretical result on this degenerate linear regression.
In this case, the objective function $f$ is quadratic, and our goal becomes to 
\begin{equation}\label{eq:quad_obj}
\min\limits_{\vx} \frac{1}{2} \vx^\top \mS \vx - \vb^\top \vx
\text{ s.t.} ~ \mA^\top \vx = \vzero,
\end{equation}
where $\mS$ is a positive definite matrix.
In this case, LPSA takes the following form
\begin{align}
\tag{\ref{eq:update_rule}}
\begin{split}
        \vx_{n + \frac{1}{2} } & = \vx_{n} - \eta_n (\mS \vx_n - \vb) + \eta_n \xi_n, \\
    \vx_{n + 1} & = \left\{
    \begin{matrix}
        \gP_{\mA^\bot} (\vx_{n + \frac{1}{2} } ), \hfill
        & \text{ with probability } p_n, \hfill \\
        \vx_{n + \frac{1}{2}}, \hfill
        & \text{ with probability } 1 - p_n \hfill.
    \end{matrix}
    \right.
\end{split}
\end{align}
The gradient function becomes $\nabla f(\vx) = \mS \vx - \vb$ and the Hessian matrix keeps constant as $\nabla^2 f(\vx) = \mS$.
From Proposition~\ref{prop:f(x)}, the optimal solution $\vx^{\star}$ of the problem \eqref{eq:quad_obj} satisfies $\gP_{\mA^\bot} (\mS \vx^\star - \vb) = \vzero$ and the asymptotic bias is reduced to $\vmu =  \frac{1}{\gamma}\left\{\gP_{\mA^\bot}\left(\mS - \frac{1-\beta}{\eta_0}\mathbbm{1}_{\{\alpha = 1\}}\rmI\right)\gP_{\mA^\bot}\right\}^{\dag}\left(\gP_{\mA^\bot} \mS (\mS \vx-\vb)\right)$.
Hence, the degenerate condition $\gP_{\mA^\bot} \nabla^2 f(\vx^\star) \nabla f(\vx^\star) = \vzero$ corresponds to $\gP_{\mA^\bot} \mS (\mS \vx-\vb) = \0$.
In order to establish a non-trivial asymptotic convergence, we introduce the \textit{degenerate order} $k$ which is the smallest integer making $\gP_{\mA^\bot} \mS^{k} (\mS\vx^\star - \vb)$ non-zero.

\begin{definition}[Degenerate order]
    \label{def:degenerate}
        We say $k$ is the degenerate order for the tuple $(\mA, \mS, \vb)$ if $\gP_{\mA^\bot} \mS^i (\mS \vx^\star - \vb)= \vzero, ~~\forall i = 0, 1,\dots, k-1$ and $\gP_{\mA^\bot} \mS^{k} (\mS\vx^\star - \vb)\neq \vzero$ where $\vx^{\star}$ is the minimizer of the constrained problem~\eqref{eq:quad_obj} instanced by the tuple $(\mA, \mS, \vb)$.
\end{definition}

\begin{assumption}[Degenerate condition]\label{asp:degenerate}
We assume the degenerate order $k$ for the tuple $(\mA, \mS, \vb)$ satisfies $1 \le k <  \infty$.
\end{assumption}

If $k=1$, it is reduced to the non-degenerate case which we have discussed in the previous section.
Assumption~\ref{asp:degenerate} remains the possibility of $k=1$ for generality.
If $k \ge 2$, here comes a truly non-trivial degenerate case. 
Lemma~\ref{lem:well_def_dgnt} below shows that Assumption~\ref{asp:degenerate} is not vacuous and indeed has instances satisfying it.
\begin{lemma}\label{lem:well_def_dgnt}
    For any positive integer $k,d$ satisfying $k + 1 < d$, we can construct two matrices $\mA_0, \mS_0$ and a vector $\vb_0$ such that the tuple $(\mA_0, \mS_0, \vb_0)$ satisfies Assumption~\ref{asp:degenerate}.
\end{lemma}




Recall that we define $\vu_n = \gP_{\mA^\bot}(\vx_n)$ and $\vv_n = \gP_{\mA} (\vx_n)$.
In the following, we provide a generalization of Theorem~\ref{thm:converge} under the degenerate condition.
In this particular problem \eqref{eq:quad_obj}, Assumptions~\ref{asp:smooth} and \ref{asp:str_cov} naturally hold.
Then under Assumptions~\ref{asp:noi_lip} and~\ref{asp:degenerate}, we establish the following convergence rate in Theorem~\ref{thm:cov_dgnrt_bias}.
By using a similar argument, we can further characterize the asymptotic behavior of the iterates $\{ \vu_n\}_{n \ge 0}$ in Theorem~\ref{thm:degenerate_bias_u}.

\begin{theorem}(Convergence rate for the degenerate case)\label{thm:cov_dgnrt_bias}
Suppose that Assumption~\ref{asp:noi_lip} is true and Assumption~\ref{asp:degenerate}
holds with parameter $k$.
Let $\eta_n = \eta_0 n^{-\alpha}$ and $p_n = \min\{ \gamma \eta_n^\beta, 1 \}$ with $0 < \alpha \le 1$, $0 \le \beta < 1$ and $\eta_0 > 2 k \mathbbm{1}_{ \{ \alpha=1 \} } / \lambda_{\mathrm{min}} (\mS)$ where $\lambda_{\mathrm{min}} (\mS)$ denotes the smallest eigenvalue of $\mS$.
Then,
\begin{align*}
    \EB\|\vu_n - \vx^\star \|^2 &= \gO\left( n^{-\alpha \min\{1, 2k(1-\beta)\}} \right), \\
    \EB\| \gP_{\mA^\bot} \mS \vv_n \|^2 &= \gO\left( n^{-2k\alpha(1-\beta)} + n^{- \alpha (2 - \beta) } \right).
\end{align*}
\end{theorem}

Theorem~\ref{thm:cov_dgnrt_bias} provides upper bounds for the expected squared errors $\EB\|\vu_n - \vx^\star \|^2$ and $\EB\| \gP_{\mA^\bot} \mS \vv_n \|^2$.
We emphasize that Theorem~\ref{thm:converge} still holds under the conditions of Theorem~\ref{thm:cov_dgnrt_bias}.
However, the latter differs from the former in two aspects.
First, Theorem~\ref{thm:cov_dgnrt_bias} studies $\EB\| \gP_{\mA^\bot} \mS \vv_n \|^2$ rather than $\EB\|\vv_n \|^2$.
This is because it is $\gP_{\mA^\bot} \mS \vv_n$ that directly affects the convergence of $\vu_n$ due to the following decomposition from~\eqref{eq:update_rule}
\begin{equation}
\label{eq:help0}
\vu_{n+1} - \vx^{\star} = \left[ \rmI - \eta_n \gP_{\mA^\bot} \mS \right](\vu_{n} - \vx^{\star}) - \eta_n \gP_{\mA^\bot} \mS \vv_n  + \eta_n  \gP_{\mA^\bot} \xi_n.
\end{equation}
Under the non-degenerate case where $k=1$, one can show that $\EB\|\vv_n \|^2$ is equivalent to $\EB\| \gP_{\mA^\bot} \mS \vv_n \|^2$ up to a constant factor.

Secondly, Theorem~\ref{thm:cov_dgnrt_bias} provides more detailed rates that depend on the degenerate order $k$. 
When $k=1$, we obtain the exact rates in Theorem~\ref{thm:converge}. 
As $k$ increases, these rates decrease, indicating that a degenerate bias with higher-order degeneration speeds up convergence. 
When $k$ approaches infinity, we have $\EB\|\vu_n - \vx^\star \|^2 = \gO(n^{-\alpha})$, which is the convergence rate of SGD on non-constrained problems. 
This is reasonable because, in the case of $k = \infty$, there is no positive integer $k$ such that $\gP_{\mA^\bot} \mS^{k} (\mS\vx^\star - \vb) \neq \vzero$. 
In other words, $\gP_{\mA^\bot} \mS^k (\mS \vx^\star - \vb) = \vzero$ for any positive integer $k$. 
Roughly speaking, this condition is almost equivalent to $\mS \vx^\star - \vb = \vzero$, which means the minimizer for the constrained optimization problem coincides with that of the non-constrained one.
Though this case is beyond the scope of our paper, one can solve it by treating the constrained problem as a non-constrained one and applying non-constrained optimization methods. 
As we shall see, when $k$ is sufficiently large, the convergence rate of LPSA on the degenerated problems indeed matches that of the non-constrained one.

\def\uhomo{{\hat{\vu}}}

\begin{theorem}
\label{thm:degenerate_bias_u}
Under the conditions of Theorem~\ref{thm:cov_dgnrt_bias}, it follows that as $n\to \infty$,
\begin{itemize}
    \item If $0 < \beta < 1 - \frac{1}{2 k}$, $\frac{\vu_n-\vx^{\star}}{\sqrt{\eta_{n-1}}}$ weakly converges to $\gN(\vzero, \tilde{\Sigma})$ where the variance matrix $\tilde{\Sigma}$ satisfies the following Lyapunov equation~\eqref{eq:variance-lya}.
    \item If $1 - \frac{1}{2 k} < \beta < 1$, $\frac{\vu_n - \vx^\star}{\eta_{n-1}^{ k (1-\beta)}}$ converges to the non-zero vector $\vmu_k$ in the $L_2$ sense with $\vmu_k$ given by
    \[
    \vmu_k = \frac{ (-1)^{k+1} }{ \gamma^{k} }
    \left[ \gP_{\mA^\bot} \left( \mS - \frac{ k (1{-}\beta) }{ \eta_0 } \mathbbm{1}_{ \{ \alpha=1 \} } \rmI_d  \right) \gP_{\mA^\bot} \right]^\dag \gP_{\mA^\bot} \mS^k 
    (\mS \vx^\star - \vb).
    \]
    \item If $ \beta = 1-\frac{1}{2 k}$, $\frac{\vu_n-\vx^{\star}}{\sqrt{\eta_{n-1}}}$ weakly converges to $\gN(\vmu_k, \tilde{\Sigma})$ with $\vmu_k$ and $\tilde{\Sigma}$ given in previous items.
\end{itemize}
\end{theorem}

Theorem~\ref{thm:degenerate_bias_u} characterizes the asymptotic behavior of the iterates $\{ \vu_n\}_{n \ge 0}$ for the degenerate case.
When $k=1$, the results of Theorems~\ref{thm:diff_approx},~\ref{thm:asy_bias_u} and~\ref{thm:noncen_u} on quadratic objectives are covered.
However, when $k \ge 2$, Theorem~\ref{thm:degenerate_bias_u} modifies several aspects to reveal a non-trivial behavior.
Specifically, the scale $\eta_{n-1}^{-k(1-\beta)}$, the non-degenerate bias $\vmu_k$, and the boundary point $1-\frac{1}{2k}$ start to depend on the degenerate order $k$.

\subsubsection{
\texorpdfstring{The structure of $\vmu_k$}{The structure of muk}}
In the following, we provide an intuitive explanation of how the non-degenerate bias $\vmu_k$ is related to the quantity $\gP_{\mA^\bot} \mS^k (\mS \vx^\star - \vb)$, which is non-zero due to Assumption~\ref{asp:degenerate}.
The existence of the asymptotic bias $\vmu_k$ can be attributed to the non-ignorable effect that $\vv_n$ has on $\vu_{n+1}$ via the relation~\eqref{eq:help0}. 
Specifically, from the update rule~\eqref{eq:update_rule}, we know that $\vv_{n + \frac{1}{2}} = \vv_n - \eta_n \gP_{\mA^\bot}\mS \vv_n$, and that $\vv_{n+1}$ takes the value $\vv_{n + \frac{1}{2}}$ with probability $1-p_n$ and becomes zero with the remaining probability.
To further explore this relation, we can multiply both sides of the update rule~\eqref{eq:update_rule} with $\gP_{\mA^\bot}\mS(\gP_{\mA}\mS)^{i}$, where $0 \le i\le k-1$, and take expectations. For any $n \ge i-1$, we can see that:
\begin{equation}
\label{eq:help1}
\EB \gP_{\mA^\bot}\mS(\gP_{\mA}\mS)^{i} \vv_{n-i+1} = (1-p_{n-i})(\EB  \gP_{\mA^\bot}\mS(\gP_{\mA}\mS)^{i} \vv_{n-i} - \eta_{n-i} \EB\gP_{\mA^\bot}\mS(\gP_{\mA}\mS)^{i+1} \vv_{n-i}).
\end{equation}
Theorem~\ref{thm:jump_approx} implies that $\vv_n$ weakly converges to a scaled exponential distribution whose expectation is approximately $\eta_n^{1-\beta} \nabla f(\vx^\star)/\gamma$, where $\nabla f(\vx^\star) = \mS \vx^\star - \vb$ in this case.
 As $\vv_n$ weakly converges, we can expect that $\EB \gP_{\mA^\bot}\mS(\gP_{\mA}\mS)^{i} \vv_{n-i+1} \approx \EB \gP_{\mA^\bot}\mS(\gP_{\mA}\mS)^{i} \vv_{n-i}$ when $n$ is sufficiently large. 
 By plugging this approximation into~\eqref{eq:help1}, we get that for sufficiently large $n$,
\begin{equation}
\label{eq:help2}
\EB  \gP_{\mA^\bot}\mS(\gP_{\mA}\mS)^{i} \vv_{n-i} \approx \frac{\eta_{n-i}}{p_{n-i}} \EB\gP_{\mA^\bot}\mS(\gP_{\mA}\mS)^{i+1} \vv_{n-i}.
\end{equation}
If $k=2$, we can set $i=0$ in~\eqref{eq:help2} and obtain that
\begin{align*}
\EB  \gP_{\mA^\bot}\mS \vv_{n} &\approx \frac{\eta_{n}}{p_{n}} \EB\gP_{\mA^\bot}\mS \gP_{\mA}\mS \vv_{n} \\
&\approx \frac{\eta_{n}}{p_{n}}  \gamma \EB\gP_{\mA^\bot}\mS \gP_{\mA}\mS \nabla f(\vx^\star)\\
& =\frac{\eta_{n}}{p_{n}}  \gamma \EB\gP_{\mA^\bot}\mS^2 \nabla f(\vx^\star) 
\neq \0,
\end{align*}
where we use the fact $\EB \vv_{n} \approx \gamma \nabla f(\vx^\star)$ and the following equation
\[
\gP_{\mA^\bot}\mS (\gP_{\mA}\mS)^{i-1} \nabla f(\vx^\star) = \gP_{\mA^\bot}\mS^i \nabla f(\vx^\star)~\forall i \ge 1,
\]
which is proved in Proposition~\ref{prop:degenerate} in the appendix.
However, this argument cannot be applied when $k>2$ because we would have $\EB \gP_{\mA^\bot}\mS \vv_{n-1} \approx \0$, which implies a failure in finding a non-degenerate asymptotic bias.
As a remedy, we can iterate \eqref{eq:help2} over $i=0,\cdots,k-1$, change the value of $n$ accordingly, and obtain that for sufficiently large $n$,
\begin{align*}
\EB  \gP_{\mA^\bot}\mS \vv_{n} 
&\approx \frac{\eta_{n}}{p_{n}} \EB\gP_{\mA^\bot}\mS \gP_{\mA}\mS \vv_{n} 
\approx \frac{\eta_{n}^2}{p_{n}^2} \EB\gP_{\mA^\bot}\mS (\gP_{\mA}\mS)^2 \vv_{n} \\
&\approx \cdots \approx  \frac{\eta_{n}^{k-1}}{p_{n}^{k-1}} \EB\gP_{\mA^\bot}\mS (\gP_{\mA}\mS)^{k-1} \vv_{n} \\
&\approx  \frac{\eta_{n}^{k-\beta}}{\gamma  p_{n}^{k-1}} \gP_{\mA^\bot}\mS (\gP_{\mA}\mS)^{k-1} \nabla f(\vx^\star) 
\approx  \frac{\eta_{n}^{k-\beta}}{\gamma  p_{n}^{k-1}} \gP_{\mA^\bot}\mS^k \nabla f(\vx^\star) 
\neq \0.
\end{align*}
The choice of $p_n = \gamma \eta_n^\beta$, together with the last inequality, implies that
\begin{equation}
\label{eq:help3}
\EB  \gP_{\mA^\bot}\mS \vv_{n} 
\approx  \frac{\eta_{n}^{k(1-\beta)}}{\gamma^k} \gP_{\mA^\bot} \mS^{k} \nabla f(\vx^\star).
\end{equation}
The approximation given by~\eqref{eq:help3} provides an explanation for the appearance of $\gP_{\mA^\bot} \mS^k (\mS \vx^\star - \vb)$ in the definition of $\vmu_k$. Furthermore, it clarifies the choice of scale $\eta_n^{-k(1-\beta)}$: in order to reveal a non-zero bias, it is necessary to normalize $\vv_{n}$ by $\eta_{n}^{k(1-\beta)}$, as suggested by~\eqref{eq:help3}, and similarly for $\vu_n$ due to the relation~\eqref{eq:help0}.

\subsubsection{When $k = \infty$}
Now, let us focus on the case where $k = \infty$. This situation arises when $\gP_{\mA^\bot} \mS^k (\mS \vx^\star - \vb) = \vzero$ for all positive integers $k$.

To support this claim, we can consider the distributed learning problem introduced in Example~\ref{eg:distributed}, which is a specific instance of the optimization problem we are analyzing. 
This problem involves minimizing a sum of convex functions, each of which is associated with a different worker node. 
More specifically, we consider the following quadratic objectives
\begin{align*}
    \min_{\vx} f(\vx) = \frac{1}{N} \sum_{i=1}^N \left[ \frac{1}{2} (\vx^{(i)})^\top \mS_i \vx^{(i)} - \vb_i^\top \vx^{(i)} \right]
    ~~\text{ s.t. } \vx^{(1)}=\cdots=\vx^{(N)},
\end{align*}
where $\vx$ denotes the concatenated variable $\left[(\vx^{(1)})^\top,(\vx^{(2)})^\top,\dots,(\vx^{(N)})^\top\right]^\top$ and
$\mS = \diag \{\mS_1, \mS_2, \dots, \mS_N \}$ is a positive definite block diagonal matrix composed of local Hessian matrices $\{S_i\}_{i=1}^N$.
Moreover, $\vb=\left[\vb_1^\top, \vb_2^\top, \dots, \vb_N^\top\right]^\top $ represents the concatenated vector of local responses.
The communication between the nodes takes place through a coordinator node, which aggregates the gradients computed by the workers. 
By construction, the communication matrix $\mA$ has a block structure and is given in~\eqref{eq:mat_A_fl}, where each block corresponds to the interactions between the coordinator and one of the workers. In this setting, it is possible to show that $\gP_{\mA^\bot} \mS^k (\mS \vx^\star - \vb) = \vzero$ for all $k$, indicating that $k = \infty$ is a valid scenario.

This is because the linear operator $\gP_{\mA^\bot}$ synchronizes all local parameters with their mean, that is, 
$\gP_{\mA^\bot} (\vx) = \left[ \bar{\vx}^\top, \bar{\vx}^\top, \dots, \bar{\vx}^\top \right]^\top $
where $\bar{\vx} = \frac{1}{N} \sum_{i=1}^N \vx^{(i)}$ is the average of all the local parameter vectors $\vx^{(i)}$. 
Consequently, the solution $\vx^\star$ takes the form
\[
\vx^\star = \left[ (\vx^{(\star)} )^\top, ( \vx^{(\star)} )^\top, \dots, ( \vx^{(\star)} )^\top   \right]^\top 
~\text{with}~\vx^{(\star)} = (\sum_{i=1}^N \mS_i)^{-1} \sum_{i=1}^N \vb_i .
\]
Then the condition $\gP_{\mA^\bot} \mS^k (\mS \vx^\star - \vb) = \vzero$ is equivalent to $\sum_{i=1}^N (\mS_i^{k+1} \vx^{(\star)} - \mS_i^k \vb_i) = \vzero$, which holds naturally if all the $\mS_i$'s have the same form, i.e., $\mS_i \equiv \mS_1$.
In this case, $\vx^{(\star)} = \frac{\mS_1^{-1}}{N} \sum_{i=1}^N \vb_i$ and $\gP_{\mA^\bot} \mS^k (\mS \vx^\star - \vb) = \vzero$ for any positive $k$.
Thus, non-constrained optimization methods that involve totally local training followed by a single round of communication at the end are sufficient, as argued before.

To analyze the case of $k=\infty$, we can replace Assumption~\ref{asp:degenerate} with the following: there is a positive integer $k$ such that $\gP_{\mA^\bot} \mS^i (\mS \vx^\star - \vb)= \vzero, ~~\forall i = 0, 1,\dots, k-1$.
Fortunately, this modification would not affect the validity of Theorem~\ref{thm:cov_dgnrt_bias} (See Appendix~\ref{sec:app:degenerate} for the rationale).

\subsection{Extension to Degenerate Non-linear Cases}
\label{sec:nonlinear}
In this subsection, our focus is on the degenerate non-linear cases where $\gP_{\mA^\bot}\nabla^2 f(\vx^\star) \nabla f(\vx^\star) = \vzero$ for a non-quadratic function $f$. To analyze this case, a natural method is approximating the non-linear loss function with a linear one and then studying the resulting approximation error. We accomplish this by making use of the update rule~\eqref{eq:update_rule} and the mean value theorem, which allows us to obtain a recursive relation between $\vu_{n+1}$ and $\vu_n$, as follows:
    \begin{align}
        \vu_{n+1} & = \vu_n - \eta_n \gP_{\mA^\bot} \nabla f(\vx_n) + \eta_n \xi_n^{(1)} \nonumber \\
        &= \vu_n - \eta_n \gP_{\mA^\bot}\nabla^2 f(\vx^\star) (\vx_n - \vx^\star) \nonumber \\
        & \quad + \eta_n \gP_{\mA^\bot} \left\{ \nabla^2 f(\vx^\star) - \int_0^1
        \nabla^2 f(t\vx_n + (1-t)\vx^\star)dt\right\}(\vx_n - \vx^\star) + \eta_n \xi_n^{(1)} \nonumber \\
        \begin{split}
        \label{eq:update_rule1}
     &= (\mI - \eta_n \gP_{\mA^\bot}\nabla^2 f(\vx^\star) \gP_{\mA^\bot})(\vu_n - \vx^\star) - \eta_n \underbrace{ \gP_{\mA^\bot} \nabla^2 f(\vx^\star) \vv_n }_{\gT_{1},~\text{the linear part}}\\
        & \quad + \eta_n \underbrace{\gP_{\mA^\bot}\int_0^1 \left\{ \nabla^2 f(\vx^\star) - \nabla^2 f(t\vx_n + (1-t)\vx^\star) \right\}(\vx_n - \vx^\star)dt}_{\gT_{2},~\text{the non-linear part}} + \eta_n \xi_n^{(1)}.
        \end{split}
    \end{align}
We refer to $\gT_{1}$ as the \textit{linear part} of the asymptotic bias of $\vu_n$ because it is present even in the linear problem (\ref{eq:quad_obj}), where the corresponding quantity is $\gP_{\mA^\bot} \mS \vv_n$. 
On the other hand, $\gT_2$ arises from the nonlinearity of the problem due to the fact that $\nabla^2 f$ is no longer constant. 
Consequently, $\gT_2$ is dependent on a higher-order expansion of $\nabla f(\vx_n) - \nabla f(\vx^\star)$ and we refer to it as the \textit{non-linear part} of the asymptotic bias of $\vu_n$.

Based on the intuitive bound~\eqref{eq:help3}, we can estimate that the magnitude of the linear part' of the bias is roughly $\eta_n^{k(1-\beta)}$, where $k \ge 2$ is defined in Assumption~\ref{asp:degenerate}. Regarding the non-linear part, assuming the loss function $f$ has adequate continuity and smoothness, we can infer that
    \begin{align}
    \label{eq:T2}
        \gT_2 &= \gP_{\mA^\bot}\int_0^1 \left\{ \int_0^1 \nabla^3 f(s\vx^\star + (1-s)[t \vx_n + (1-t)\vx^\star])ds \right\}t(\vx^\star - \vx_n)(\vx_n - \vx^\star)^\top dt \nonumber \\
        &= -\gP_{\mA^\bot}\int_0^1 t\inner{\nabla^3 f(\vx^\star)}{(\vx_n - \vx^\star)(\vx_n - \vx^\star)^\top} + \gR_n \nonumber \\
        &= - \gP_{\mA^\bot}\frac{1}{2}\inner{\nabla^3 f(\vx^\star)}{(\vx_n - \vx^\star)(\vx_n - \vx^\star)^\top} + \gR_n.
    \end{align}
    Here, $\gR_n$ represents a high-order infinitesimal of $\inner{\nabla^3 f(\vx^\star)}{(\vx_n - \vx^\star)(\vx_n - \vx^\star)^\top}$ under high-order smoothness conditions.
In this expression, the inner product is defined as 
\[
\inner{\mT}{\mM} := \left(\sum\limits_{j,k = 1}^d \mT[i,j,k]\mM(j,k)\right)_{i=1}^d \in \sR^d
\]
for any tensor $\mT \in \sR^{d\times d\times d}$ and any matrix $\mM \in \sR^{d\times d}$. 
    
    Using the equation
    \[
 \eta_n^{2(\beta - 1)}\EB (\vx_n - \vx^\star)(\vx_n - \vx^\star)^\top =
    \check{\vv}_n \check{\vv}_n^\top + \check{\vv}_n \hat{\vu}_n^\top + \hat{\vu}_n \check{\vv}_n^\top + \hat{\vu}_n \hat{\vu}_n^\top,
    \]
    we can observe that in the degenerated bias setting, the dominant term is $\check{\vv}_n \check{\vv}_n^\top$.\footnote{This is because $\check{\vv}_n = \gO_\PB(1)$ due to Theorem~\ref{thm:diff_approx} while $\hat{\vu}_n = o_{\PB}(1)$ due to Theorem~\ref{thm:degenerate_bias_u}.}
According to Theorem~\ref{thm:jump_approx}, we have $\EB\check{\vv}_n \check{\vv}_n^\top \approx \frac{2}{\gamma^2} \nabla f(\vx^\star)\nabla f(\vx^\star)^\top$ and thus
\[
\eta_n^{2(\beta - 1)}\EB \gT_2 \approx -\frac{1}{2}\gP_{\mA^\bot}\inner{\nabla^3 f(\vx^\star)}{\EB\check{\vv}_n \check{\vv}_n^\top} \approx -\frac{1}{\gamma^2}\gP_{\mA^\bot}\inner{\nabla^3 f(\vx^\star)}{\nabla f(\vx^\star)\nabla f(\vx^\star)^\top}.
\]
If $\gP_{\mA^\bot}\inner{\nabla^3 f(\vx^\star)}{\nabla f(\vx^\star)\nabla f(\vx^\star)} \neq 0$, the order of $\gT_2$ is approximately $\eta_n^{2(1-\beta)}$, which is greater than the order of $\gT_1$ as long as $k \ge 3$. 
In this scenario, $\gT_2$ dominates the update rule~\eqref{eq:update_rule1}, while the theoretical results presented in this subsection focus on the case $\gT_1$ is the dominant and thus are no longer valid.
As a result, we should normalize $\vu_n$ and $\vv_n$ by $\eta_n^{2(1-\beta)}$ to reveal a non-trivial bias and similarly derive its explicit form.

However, there is a possibility that $\gP_{\mA^\bot}\inner{\nabla^3 f(\vx^\star)}{\nabla f(\vx^\star) \nabla f(\vx^\star)} = \vzero$. 
If this is the case, the order of $\gT_2$ becomes unclear again.
By substituting the update rule~\eqref{eq:update_rule} into the decomposition~\eqref{eq:T2}, we can further decompose $\gT_2$ into a linear part and a non-linear part. 
Specifically, the linear part takes the form of $\gP_{\mA^\bot}\inner{\nabla^3 f(\vx^\star)}{\vv_n \vv_n^{\top}}$, while its non-linear part is the fourth-order term obtained from expanding $\nabla f(\vx_n) - \nabla f(\vx^\star)$.
Similarly, we need to carry out a similar discussion for the decomposition of $\gT_1$ as well. 
This process should be repeated until all the quantities converge to a non-trivial limit.
From the above analysis, we can infer that the influence on convergence caused by the nonlinear part and the linear part is independent, but they interfere with each other. 
As the degenerate order increases, the number of cases needing consideration grows exponentially and the resulting problem becomes much easier to solve.
It is tedious and time-consuming to identify all the cases, so we only discuss the case where only the linear part exists in detail because it is more fundamental and relatively harder to identify than the nonlinear part.


\section{Debiased Loopless Projection Stochastic Approximation}\label{sec:dlpsa}

In the previous section, we observed that a degenerated bias can speed up convergence. 
This raises the question of how we can use this observation to develop a fast algorithm that requires as few projections as possible while still meeting a given accuracy condition. 
It is also worth noting that increasing $\beta~(<1)$ leads to a decrease in the overall projection frequency. 
Thus, our goal is to find the largest possible value of $\beta$ while keeping the global convergence rate constant. 
However, as we have seen from the theorems developed in the previous sections, the convergence of LPSA, as measured by $\EB\|\vu_n-\vx^\star\|^2$, decreases if $\beta$ is greater than $1/2$. This is due to the shift of the dominant term from the fluctuation caused by gradient noises to the asymptotic bias incurred by low-frequency projection. 
Therefore, if we want to increase $\beta$ larger than $1/2$, we cannot continue to use LPSA. 
Nevertheless, we recall an interesting conclusion in the proof of Theorem~\ref{thm:asy_bias_u} (see Lemma~\ref{lem:asmp_mix_uv} for more details)
\[
\lim\limits_{n \to \infty}\left| \EB \inner{\hat{\vu}_n - \EB \hat{\vu}_n}{\nabla^2 f(\vx^\star) (\check{\vv}_n - \EB \check{\vv}_n)}\right| = 0,
\]
where $\hat{\vu}_n $ and $\check{\vv}_n$ are given in Theorem~\ref{thm:asy_bias_u}.
This fact tells us that the impact of residuals $\vv_n - \EB \vv_n$ on the convergence of $\EB \norm{\vu_n - \vx^\star}^2$ is negligible compared to the bias resulting from the expectation of $\vv_n$. Hence, intuitively, to improve the algorithm, we only need to eliminate the sequence of bias ${\EB \vv_n}$ from the update, without trying to remove the entire occurrence of $\vv_n$, which is much more difficult to implement than the former.

As a direct corollary of Theorem~\ref{thm:jump_approx}, we can approximate $\EB\vv_n$ as $\eta_n^{1-\beta}\EB\check{\vv}_n \approx \eta_n^{1-\beta}\vv$, where $\vv = - \gamma^{-1}\nabla f(\vx^\star)$. If we knew the exact value of $\vv$, we could use it to improve the stochastic gradient involved in the update formula (\ref{eq:update_rule}), replacing $\nabla f(\vx_n, \zeta_n)$ with $\nabla f(\vx_n - \eta_n^{1-\beta}\vv, \zeta_n)$. 
This substitution yields an improvement in the analysis of the mean square error $\EB \norm{\vu_n - \vx^\star}^2$.
We explain it in the following.
In the iterative relation between $\EB \norm{\vu_{n+1} - \vx^\star}^2$ and $\EB \norm{\vu_n - \vx^\star}^2$, we need to handle the inner product term between $\vu_n$ and $\vv_n$. 
However, due to the substitution of the stochastic gradient mentioned earlier, this inner product term changes from $\EB \inner{\vu_n - \vx^\star}{\nabla^2 f(\vx^\star)\vv_n}$ to 
\[
\EB \inner{\vu_n - \vx^\star}{\nabla^2 f(\vx^\star)(\vv_n - \eta_n^{1-\beta}\vv)}= \eta_n^{1-\beta}\EB \inner{\vu_n - \vx^\star}{\nabla^2 f(\vx^\star)(\check{\vv}_n - \vv)}.
\]
One can show that the former is $\gO(\eta_n^{2(1-\beta)})$, while the latter is only $o(\eta_n^{2(1-\beta)})$ due to a similar reason in Lemma~\ref{lem:bd_innr_pdct}.
Hence, the impact of the cross-term is negligible asymptotically. 

Unfortunately, it is nearly impossible to obtain the exact value of $\vv$ (or $\nabla f(\vx^\star)$). However, we can estimate it using the stochastic gradient at the $n$th iteration, i.e., $\nabla f(\vx_n, \zeta_n^\prime)$, where $\zeta_n^\prime$ is another independent random variable. Thanks to the independence, the noise generated by $\nabla f(\vx_n, \zeta_n^\prime) - \nabla f(\vx_n)$ has little impact on the magnitude of the inner product term between $\vu_n$ and $\vv_n$. Refer to Algorithm~\ref{alg:d_lp_proj_sa} for details.

\begin{algorithm}[h]
\caption{Debiased Loopless Projected Stochastic Approximation (DLPSA)}
 	\label{alg:d_lp_proj_sa}
 	\begin{algorithmic}
 		\STATE {\bfseries Input:} function $f$, data distribution $\gD$, initial point $\vx_0$, step size $\eta_n$, projection probability $p_n = \gamma \eta_n^\beta$.
 		\FOR{$n=0$ {\bfseries to} $T-1$}
 		 \STATE{Sample $\zeta_n^{(a)}, \zeta_n^{(b)}\sim \gD$ and $\omega_n\sim \mathrm{Bernoulli}(p_n)$ independently.}
 		 \STATE{$\vx_{n+\frac{1}{2}}=\vx_n-\eta_n\nabla f\left(\vx_n + \gamma^{-1}\eta_n^{1-\beta}\nabla f(\vx_n, \zeta_n^{(a)}),\zeta_n^{(b)}\right)$}
 		 \IF{$\omega_n=1$}
 		  \STATE{$\vx_{n+1}=\gP_{\mA^\bot}\vx_{n+\frac{1}{2}}$}
 		 \ELSE
 		  \STATE{$\vx_{n+1}=\vx_{n+\frac{1}{2}}$}
 		 \ENDIF
 		 \ENDFOR
 		\STATE {\bfseries Return: $\gP_{\mA^\bot}\vx_T$.}
 	\end{algorithmic}
\end{algorithm}

With some abuse of notation, we still use $\vu_n$ and $\vv_n$ to represent $\gP_{\mA^\bot}\vx_n$ and $\gP_\mA \vx_n$. And correspondingly we have a better convergence result about $\EB \norm{\vu_n - \vx^\star}^2$.

\begin{assumption}\label{asp:lip_4th_mmt}
For the random gradient noise $\nabla f(\vx, \zeta) - \nabla f(\vx)$, we assume that there exists a constant $\Sigma^\prime > 0$ such that
\[
\EB_\gD \norm{\nabla f(\vx, \zeta) - \nabla f(\vx)}^4 \le \Sigma^\prime ( 1 + \norm{\vx - \vx^\star}^4).
\]
\end{assumption}

\begin{theorem}\label{thm:dlpsa_cov_rt}
Assume that Assumptions~\ref{asp:smooth} - \ref{asp:noi_lip} and Assumption~\ref{asp:lip_4th_mmt} hold. Taking Algorithm~\ref{alg:d_lp_proj_sa} with step size $\eta_n = \eta_0 n ^ {-\alpha}$ and projection probability $p_n = \min\{ \gamma\eta_n^\beta, 1 \}$ which satisfy 
$0 < \alpha \le 1$,
$0\le \beta < 1$
and $\eta_0 > 2 \mathbbm{1}_{ \{ \alpha=1 \} } / \mu$
($\mu$ is the strong convexity parameter of $f$),
we have
\[
\EB \norm{ \vu_n - \vx^\star}^2 \precsim
\left(\frac{\Tilde{L}}{\mu^2 \gamma} + L(\tr(\Sigma) + \norm{\nabla f(\vx^\star)}^2)\right)\frac{L^2\eta_n^{3(1-\beta)}}{\mu\gamma^2} + \frac{\tr(\Sigma)}{\mu}\eta_n.
\]
\end{theorem}

Disregarding constant factor dependence, Theorem~\ref{thm:dlpsa_cov_rt} reveals a $\gO(n^{-\alpha \min\{1, 3(1-\beta)\}})$ convergence rate for DLPSA, demonstrating its superior convergence rate to LPSA. This highlights the power of debiasing, achieved by eliminating the expectation of $\vv_n$.

\subsection{Comparison on Projection Efficiency}
In practical settings, the projection operation can be computationally expensive, which emphasizes the importance of designing efficient algorithms and adjusting hyperparameters such as $\alpha$ and $\beta$ to minimize the projection complexity. To quantify the projection efficiency, we use the average projection complexity (APC), which measures the number of projections required to achieve a feasible solution with an accuracy of $\epsilon>0$. 

Table \ref{tab:table1} summarizes the derived convergence results and their corresponding APCs. When $\beta \in [0, 0.5)$, the APC is $\gO\left(\epsilon^{\beta - \frac{1}{\alpha}}\right)$, and when $\beta \in (0.5, 1)$, the APC is $\gO\left(\epsilon^{\frac{\alpha\beta - 1}{2\alpha(1-\beta)}}\right)$. It is observed that for LPSA, the APC is minimized as $\alpha \to 1$ and $\beta \to 0.5$, where APC approaches $\frac{1}{\sqrt{\epsilon}}$. In the scenario of degenerate bias discussed in Section~\ref{sec:dgnrt_bias&homo_fed}, the optimal APC is $\gO(\epsilon^{-\frac{1}{2k}})$ with $\alpha = 1$ and $\beta = \frac{2k-1}{2k}$. In contrast, the optimal APC for DLPSA is $\gO(\epsilon^{-1/3})$ with $\alpha = 1$ and $\beta = 2/3$, which is a factor of $\epsilon^{-1/6}$ improvement over vanilla LPSA. This improvement is due to the discovery that the asymptotic bias of $\check{\vv}_n$ plays a dominant role in the convergence process.


There is an interesting parallelism between LPSA and Local SGD in the context of FL, as projection complexity corresponds to communication complexity in FL. 
Synchronization in FL is essentially a projection in linearly constrained problems. 
In \citep{li2021statistical}, the authors analyzed the averaged communication complexity (ACC) for Local SGD with Polyak-Ruppert averaging, where ACC is defined as the number of communication required to obtain a $\epsilon$-accuracy global parameter. 
They considered a general case where the length of the $m$-th inner loop could be up to $E_m:=m^\nu$ with $\nu \in [0, 1)$. 
After $E_m$ steps of the inner loop, communication would perform to synchronize local models. 
Hence, $1/E_m$ plays a role similar to $p_m$ in our paper. 
\citet{li2021statistical} found that when $\nu\in [0,1)$, the averaged Local SGD iterates enjoy an optimal asymptotic normality up to a known constant scale and its ACC is $\left(\frac{1}{\epsilon}\right)^{\frac{1}{1+\nu}}$. 
When $\nu \to 1$, ACC approaches $\frac{1}{\sqrt{\epsilon}}$, similar to the LPSA case where APC converges to $\frac{1}{\sqrt{\epsilon}}$ when $\alpha \to 1$ and $\beta \to 0.5$. 
The equivalence suggests that the communication complexity of the federated counterpart of DLPSA can achieve $\epsilon^{-1/3}$, which is a significant improvement compared with classical Local SGD algorithms. 
We note that some work achieves better communication complexity using more advanced techniques like variance reduction~\citep{karimireddy2019scaffold} or primal-dual methods~\citep{zhang2021fedpd}.
However, the primary goal of our work in this section is to demonstrate that our new asymptotic analysis method can facilitate downstream works, and the design of new algorithms like DLPSA is a direct and simple application.

\begin{table}[t!]
\renewcommand\arraystretch{1.5}
\vspace{-0.1in}
\caption{(Non-)Asymptotic results and projection complexity under different choice of $\eta_n$ and $p_n$. The first two columns list the non-asymptotic and asymptotic results respectively, and the last column characterizes projection complexity.}
    \centering
    \begin{tabular}{p{2cm}|c|c|c|c}
         \toprule 
         Algorithm & $\beta$ & $\sE\|\vu_n - \vx^\star\|^2$ & Asymptotic behavior & APC \\
         \hline \multirow{2}{*}{LPSA}
         & $[0,1/2)$ & $\gO\left(n^{-\alpha}\right)$ & Normal & $\gO\left(\epsilon^{\beta - \frac{1}{\alpha}}\right)$ \\
         \cline{2-5} & $(1/2, 1)$ & $\gO\left(n^{-2\alpha(1-\beta)}\right)$ & Biased & $\gO\left(\epsilon^{\frac{\alpha\beta - 1}{2\alpha(1-\beta)}}\right)$\\
         \hline \multirow{2}{2cm}{LPSA (bias degenerated)}
         & $[0,1-\frac{1}{2k})$ & $\gO(n^{-\alpha})$ & Normal & $\gO\left(\epsilon^{\beta-\frac{1}{\alpha}}\right)$\\
         \cline{2-5} & $[1-\frac{1}{2k}, 1)$ & $\gO(n^{-2k(1-\beta)\alpha})$ & Biased & $\gO\left(\epsilon^{\frac{\alpha\beta - 1}{2k(1-\beta)\alpha}}\right)$ \\
         \hline \multirow{2}{*}{DLPSA} & $[0,2/3)$ & $\gO(n^{-\alpha})$ & NA & $\gO\left(\epsilon^{\beta - \frac{1}{\alpha}}\right)$ \\
         \cline{2-5} & $[2/3, 1)$ & $\gO(n^{-3(1-\beta)\alpha})$ & NA & $\gO\left(\epsilon^{\frac{\alpha\beta - 1}{3(1-\beta)\alpha}}\right)$ \\
         \bottomrule
    \end{tabular}
    \label{tab:table1}
\vspace{-0.1in}
\end{table}

\section{Experiments}\label{sec:exp}
In this section, we validate our theoretical results 
through 
experiments. 
We first introduce the experimental setup and then show the experimental results corresponding to Sections~\ref{sec:main_rst} to \ref{sec:dlpsa}.

\paragraph{Experimental setup}
We focus on the problem \eqref{eq:proj_opt} where the objective is a quadratic function of the form 
\begin{align}\label{eq:proj_opt_quadratic}
    f(\vx, \zeta) = \frac{1}{2} \vx^\top \mS \vx - \vb^\top \vx + \zeta^\top \vx \  \text{s.t.} \ \mA \vx = \vzero.
\end{align}
Here, the Hessian matrix is given by $\mS = \mU \mU^\top + \rmI_{d_1}$, where $\mU \in \sR^{d_1 \times d_1}$, $\vb \in \sR^{d_1}$ and $\zeta \in \sR^{d_1}$ are independently generated, and $\rmI_{d_1} \in \sR^{d_1 \times d_1}$ is the identity matrix. 
The entries of $\mU$ are independent $\gN(0, 1/d_1)$ random variables, and $\vb$ and $\zeta$ are modeled as $\gN(\vzero, \rmI_{d_1})$. Therefore, with high probability, the population objective $f(\vx) = \sE_\zeta f(\vx, \zeta)$ is $1$-strongly convex and $C$-smooth with $C$ being a positive constant number.\footnote{Note that the smoothness of $f$ is guaranteed by Theorem~4.4.5 in \citet{vershynin2018high}.} 
The constraint matrix $\mA \in \sR^{d_1 \times (d_1 - d_2)}$ also has i.i.d. standard normal entries. For our experiments, we  set $d_1 = 5$ and $d_2 = 2$. 

\begin{figure}[t!]
    \centering
    \includegraphics[width=\textwidth]{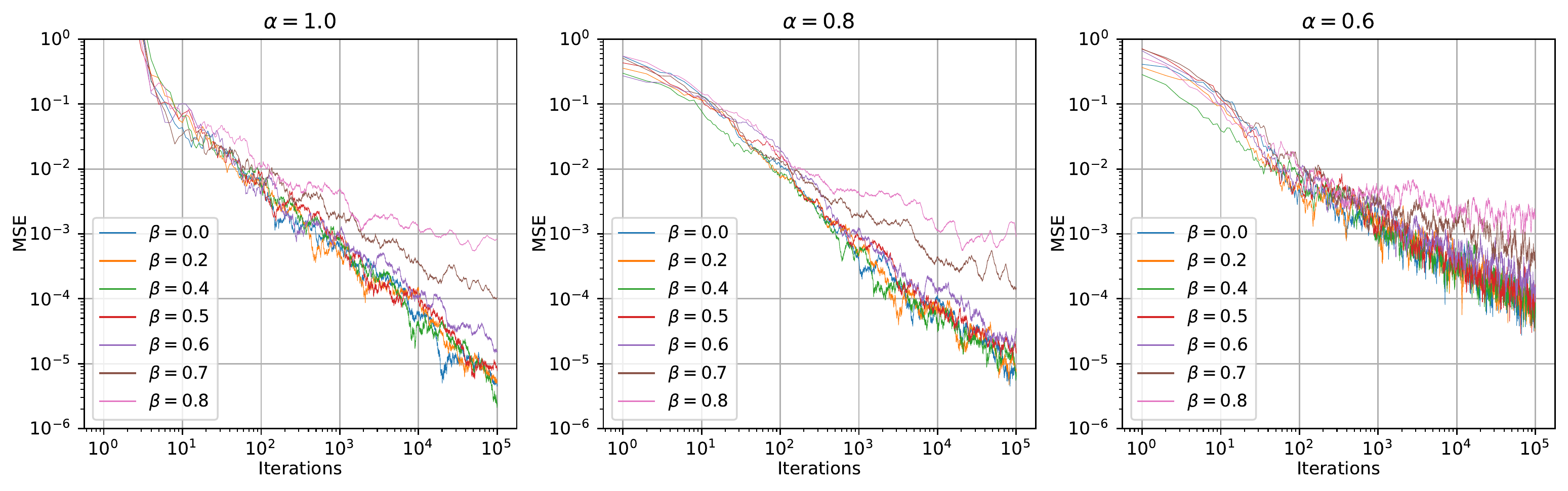}
    \caption{The log-log scale graphs of averaged MSEs versus iterations over $10$ repetitions.}    \label{fig:lc_converge_rate}
\end{figure}

\paragraph{Convergence rates}
Figure~\ref{fig:lc_converge_rate} shows log-log plots of the averaged mean squared errors (MSEs) over 10 repetitions versus iterations. 
We set $\alpha$ to $\{1, 0.8, 0.6 \}$ and $\beta$ to $\{0, 0.2, 0.4, 0.5, 0.6, 0.7, 0.8 \}$, with $d_1 = 5$ and $d_2 = 2$. Each repetition runs $10^5$ steps of LPSA.
According to Theorem~\ref{thm:converge}, the slope of the line in the log-log plot should be $-\alpha \min \{1, 2-2\beta \}$. We observe that the experimental results align with this prediction when the iteration count is larger than $1000$. Moreover, for $\beta \in [0, 1/2)$, the value of $\beta$ does not affect the slope. However, for $\beta \in (1/2, 1)$, larger values of $\beta$ and smaller values of $\alpha$ result in smoother lines. Notably, the point at which the phase transition occurs is when $\beta = 1/2$.


\begin{figure}[t!]
    \centering
    \includegraphics[width=\linewidth]{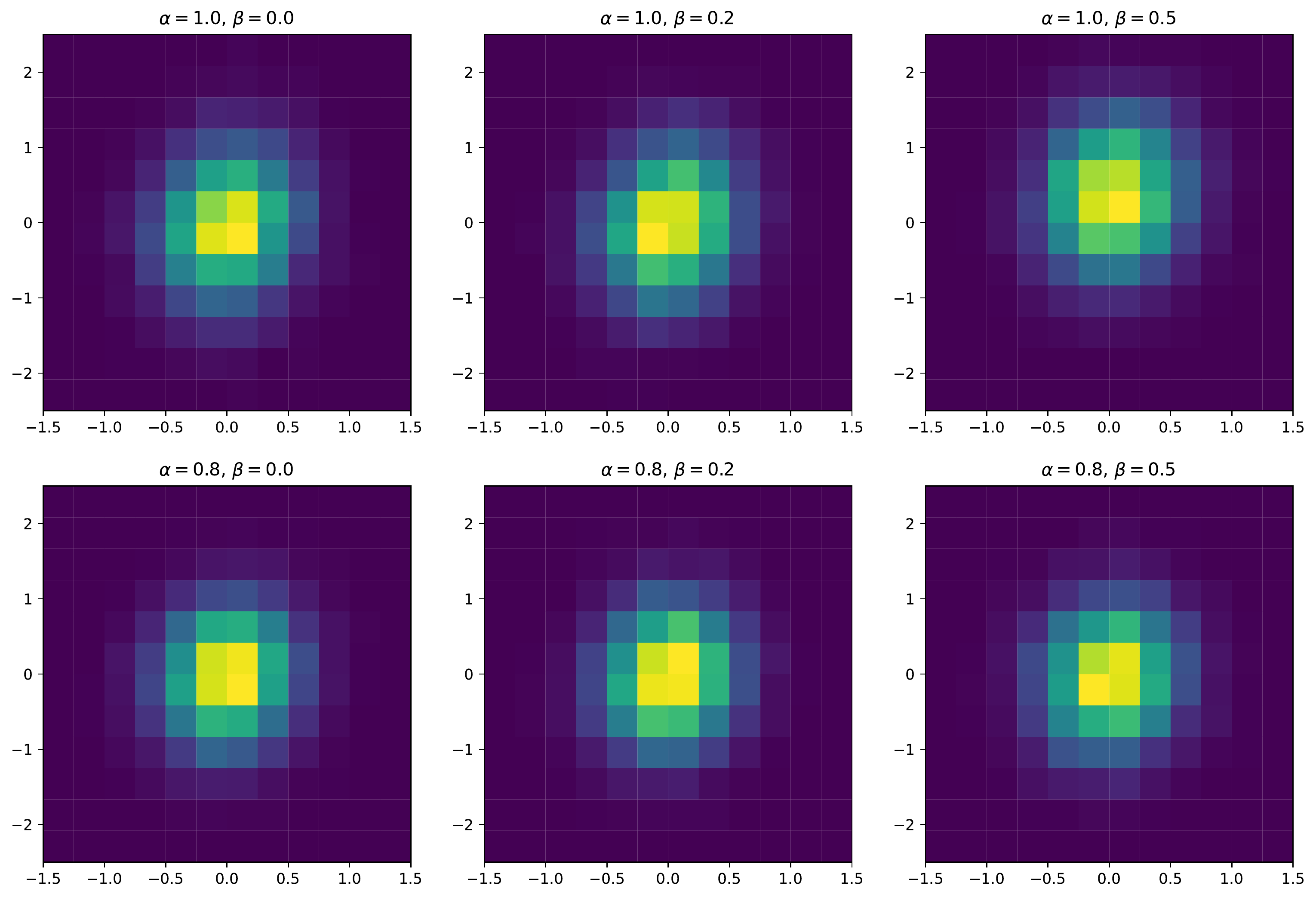}
        \caption{The heatmaps of (centralized) $\check{\vu}_{10001}$ along two orthogonal directions.
        }
        \label{fig:lc_heatmap}
\end{figure}

\paragraph{Frequent projection}

We conduct experiments with LPSA by running $10^5$ steps of the algorithm over $10^5$ repetitions for $\alpha=1$ and $\beta \in \{0, 0.2, 0.5\}$, and we pick the last iterates ${\vu}_{10001}$. 
We then compute the rescaled vectors $\check{\vu}_{10001}$ for these iterates, as defined in Section~\ref{sec:freq_proj}. 
Since $d_1=5$ and $d_2=2$, we find that the $\check{\vu}n$ sequence lies in a two-dimensional subspace of $\sR^5$. We plot the heatmaps of these $\check{\vu}_{10001}$'s across two orthogonal directions of the subspace in the left two columns of Fig.\ref{fig:lc_heatmap} for $\beta \in {0, 0.2}$. 
For $\beta=0.5$, we first centralize $\check{\vu}_{10001}$ by subtracting it by the bias vector $\vmu$ in Theorem\ref{thm:noncen_u} and then plot the heatmaps as before in the last column of Fig.~\ref{fig:lc_heatmap}. 
We observe that cells near the origin have lighter colors, and as we move away from the origin, the cell color becomes darker. The cells with lighter colors imply more frequencies, which agrees with Theorems \ref{thm:diff_approx} and \ref{thm:noncen_u}, where the limiting distribution of the (centralized) $\check{\vu}_n$ is Gaussian.

\begin{figure}[t!]
    \centering
    \begin{minipage}{0.49\textwidth}
        \centering
        \includegraphics[width=\linewidth]{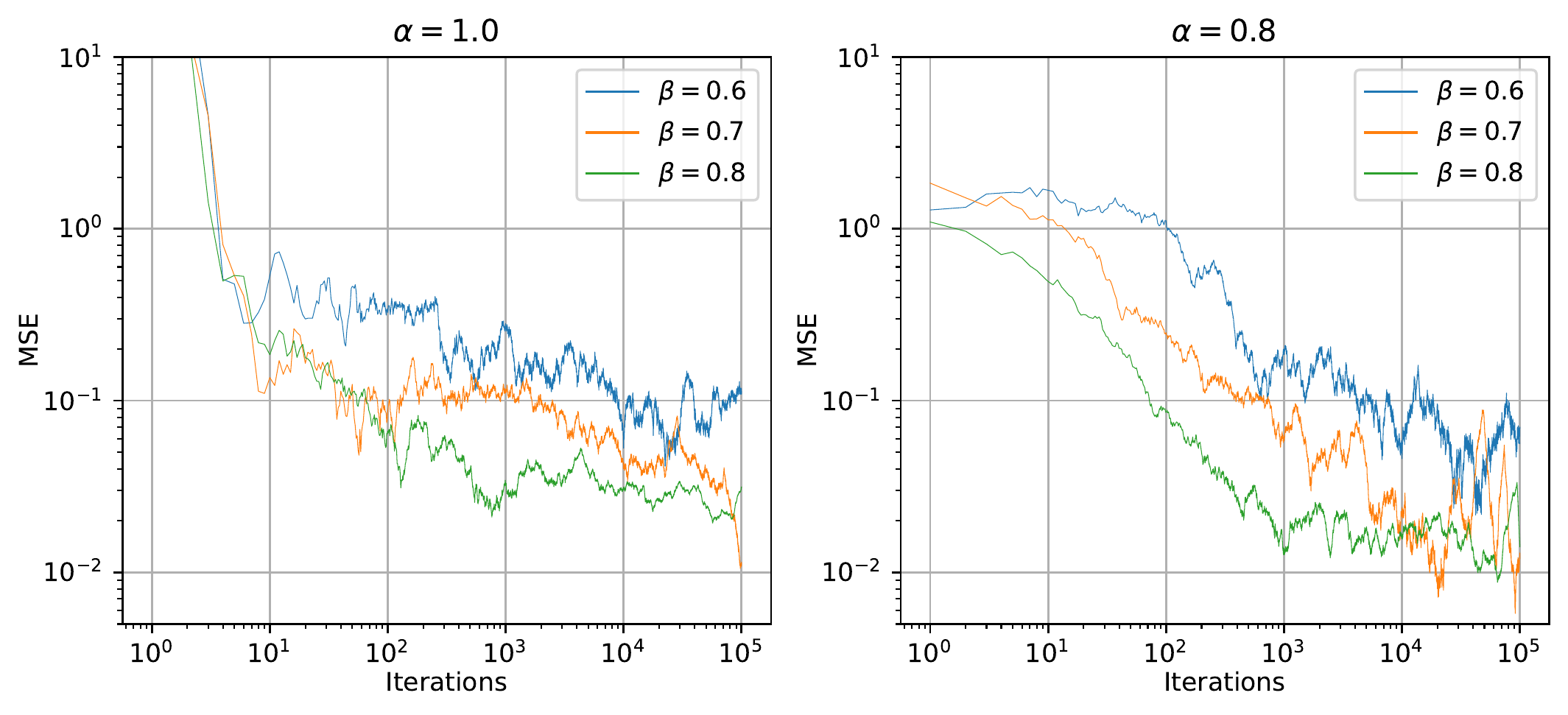}
        \caption{ The log-log scale graphs of averaged MSEs between $\check{\vu}_n$ and the bias vector $\vmu$ (defined in Theorem~\ref{thm:jump_approx}) versus iterations  over $10$ repetitions.}
        \label{fig:lc_bias_mse}
    \end{minipage}
    \hfill
    \begin{minipage}{0.49\textwidth}
        \centering
        \includegraphics[width=\linewidth]        {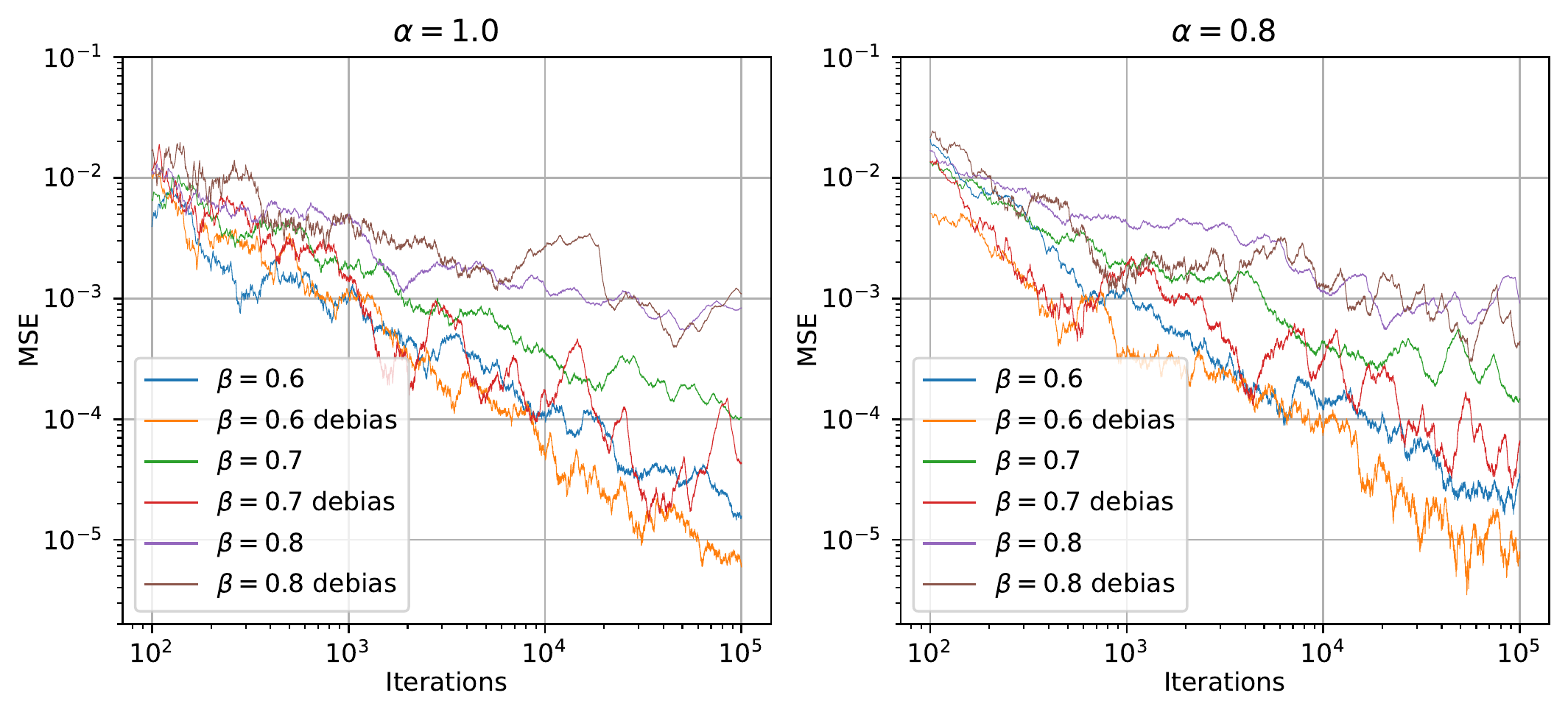}
        \caption{
        The log-log scale graphs of averaged MSEs of LPSA and DLPSA versus iterations over $10$ repetitions.
        The `debias' in the legend indicates the line for DLPSA.
        }
        \label{fig:lc_debias}
    \end{minipage}
\end{figure}

\paragraph{Occasional projection}

For $\alpha \in \{1.0, 0.8\}$ and $\beta \in \{0.6, 0.7, 0.8\}$, we conduct $10$ repetitions of LPSA with $10^5$ steps, and plot the average mean squared errors (MSEs) between the estimated vector $\hat{\vu}_n$ (defined in Section~\ref{sec:occa_proj}) and the bias vector $\vmu$ (defined in Theorem~\ref{thm:asy_bias_u}) in Fig.\ref{fig:lc_bias_mse}.
Although the MSEs have a non-negligible magnitude, they gradually decrease as the number of iterations increases.
Note that Theorem\ref{thm:asy_bias_u} guarantees asymptotic convergence without specifying the convergence rate, which can be extremely slow.
For the problem~\eqref{eq:proj_opt_quadratic}, it can be verified that the bias vector $\vmu$ is non-zero.
Therefore, we can conclude that $\hat{\vu}_n$ is indeed asymptotically biased for $\beta > 0.5$, which is consistent with the analysis in Section~\ref{sec:occa_proj}.\footnote{We also plot the trajectories of $\hat{\vu}_n$, which are deferred to Appendix \ref{sec:append_expe_complement}}


\paragraph{Debiased algorithm}
In Section~\ref{sec:dlpsa}, we introduced the debiased algorithm DLPSA, which enjoys superior convergence rates compared to LPSA, as shown in Theorem~\ref{thm:dlpsa_cov_rt}. 
Specifically, for $\beta > 0.5$, DLPSA converges at rates $\gO(n^{-\alpha\min\{1, 3(1-\beta)\}})$, which is faster than LPSA. To demonstrate this, we run $10^5$ steps of DLPSA for $\alpha \in \{1.0, 0.8\}$ and $\beta \in \{0.6, 0.7, 0.8\}$ over $10$ repetitions, and compare the averaged MSEs between LPSA and DLPSA in Fig.\ref{fig:lc_debias}. 
Note that the approximation of $\nabla f(\vx^\star)$ by $\nabla f(\vx_n, \zeta_n')$ is only valid for large $n$, so we warm up DLPSA by running the first $10^2$ steps with LPSA. 
We only plot the MSEs after this stage in Fig.\ref{fig:lc_debias}.
We observe that for $\beta \in \{0.6, 0.7\}$, DLPSA exhibits significantly faster convergence rates than LPSA, while for $\beta = 0.8$, the acceleration is not as obvious due to the less frequent projections.


\begin{figure}[t!]
    \centering
    \includegraphics[width=\textwidth]{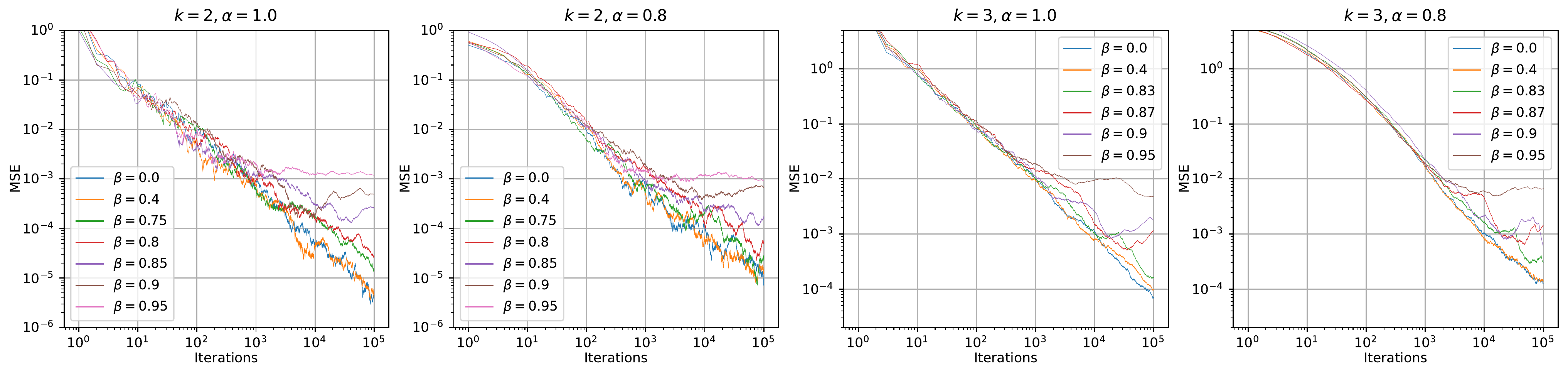}
    \caption{The log-log scale graphs of averaged MSEs on degenerate cases over $10$ repetition vs iterations.}    \label{fig:lc_degenerate}
\end{figure}

\paragraph{Degenerate cases}

In general, the problem~\eqref{eq:proj_opt_quadratic} in the \textit{Experimental setup} parapragh does not satisfy the degenerate condition $\gP_{\mA^\bot} \nabla^2 f(\vx^\star) \nabla f(\vx^\star) = \vzero$. 
Fortunately, Lemma~\ref{lem:well_def_dgnt} provides a way to find a degenerate example.\footnote{For more details, see the proof of Lemma~\ref{lem:well_def_dgnt} in Appendix~\ref{sec:app:degenerate}.}
To demonstrate the performance of LPSA on degenerate cases, we construct two problems with degenerate orders $k=2$ and $k=3$, respectively, as defined in Assumption~\ref{asp:degenerate}. 
For $k=2$, we set $d_1 = 5$ and $d_2 = 2$, and select $\alpha \in \{ 1.0, 0.8 \}$ and $\beta \in \{ 0, 0.4, 0.75, 0.8, 0.85, 0.9, 0.95 \}$. For each $(\alpha, \beta)$ pair, we perform $10^5$ steps of LPSA and plot the log-log scale graphs of the averaged MSEs over $10$ repetitions in the left half of Fig.\ref{fig:lc_degenerate}. 
Recall that Theorem~\ref{thm:degenerate_bias_u} implies the convergence rate for $k=2$ is of order $ \gO( n^{-\alpha \min\{ 1, 4(1-\beta) \}} )$. Therefore, $\beta = 0.75$ is the change point where the phase transition occurs, as shown in Fig.~\ref{fig:lc_degenerate}.

For $k=3$, we consider a problem with larger dimensions $d_1 = 20$ and $d_2 = 10$ and choose $\alpha \in \{ 1.0, 0.8 \}$ and $\beta \in \{ 0, 0.4, 0.83, 0.7, 0.9, 0.95\}$. For each $(\alpha, \beta)$ pair, we still perform $10^5$ steps of LPSA and plot the log-log scale graphs of the averaged MSEs over $10$ repetitions in the right half of Fig.\ref{fig:lc_degenerate}. Theorem~\ref{thm:degenerate_bias_u} implies the convergence rate for $k=3$ is of order $ \gO( n^{-\alpha \min\{ 1, 6(1-\beta) \}} )$. Fig.~\ref{fig:lc_degenerate} also shows that the phase transition occurs when $\beta$ crosses $5/6$.

\section{Concluding Remarks}\label{sec:ccld}
In this paper we study the linearly constrained optimization problem. 
We propose the LPSA algorithm that is inspired by Local SGD.
The probabilistic projection in LPSA follows the spirit of loopless methods \citep{kovalev2020don,hanzely2020federated,li2021anita} and simplifies the double-loop structure of original Local SGD, facilitating theoretical analysis.
We thoroughly analyze the (non-)asymptotic properties of properly scaled trajectories obtained from $\{\vu_n\}$ and discover an interesting phase transition where $\{\vu_n\}$ changes from asymptotically normal to asymptotically biased as the projection frequency decreases.
From a technical level, we generalize jump diffusion approximations to accommodate the particularity and discontinuity of LPSA.  



There are also some open problems.
It is unclear about the asymptotic behavior of $\vu_n$ when $\beta = 0.5$, i.e., $p_n = \Theta(\sqrt\eta_n)$.
The jump diffusion approach fails because we can't analyze $\{{\vu}_n\}$ via the length of $\{\vv_n\}$ anymore.
It accounts for failure that $\{\check{\vu}_n\}$ and $\{\check{\vv}_n\}$ are incompatible in the sense that they use different time scales and the time interpolation.
However, we speculate $\check{\vu}_n$ would finally converge weakly to a non-centred Gaussian distribution.
In addition, it is also interesting to analyze the performance of projection complexity of LPSA. 
From Theorem \ref{thm:asy_bias_u}, to achieve a better convergence rate at lower projection frequencies, we must overcome the asymptotically biased nature of $\vu_n$. 
One feasible approach is to build a `de-biasing' algorithm which attenuates the effect of $\vv_n$ during the update of $\vu_n$. 
We leave them as future work.


\newpage
\bibliographystyle{plainnat}
\bibliography{ref.bib,bib/optimization,bib/federated,bib/distributed}

\appendix

\newpage


\section{Proof of Section \ref{sec:converge}}\label{sec:append_converge}
In this section, we give the proof of Theorems \ref{thm:converge} and \ref{thm:counter}.

\subsection{Useful Propositions and Lemmas}
In this subsection, we present some existing results and auxiliary lemmas useful for our later analysis.
\begin{proposition}[\cite{nesterov2018lectures}, Theorem 2.1.9, property of strong convexity]\label{prop:sc}
If $f(\vx) $ is $\mu$-strongly convex, then we have 
\[
\inner{\nabla f(\vx) - \nabla f(\vy)}{ \vx - \vy} \ge \mu \norm{ \vx - \vy}^2 , \ \forall\, \vx, \vy \in \sR^d.
\]
\end{proposition}

\begin{proposition}[Cauchy–Schwarz Inequality]\label{prop:cauchy}
For any vectors $\va, \vb \in \sR^d$ and positive number $\gamma$, it holds that
\[
2 \inner{\va}{\vb} \le \gamma \norm{\va}^2 + \frac{1}{\gamma} \norm{\vb}^2.
\]
Moreover, for any positive integer $n$ and any vectors $\vx_1, \vx_2, \dots, \vx_n \in \sR^d$, it holds that
\[
\norm{ \sum_{i=1}^n \vx_i }^2 \le n \sum_{i=1}^n \norm{ \vx_i }^2.
\]

\end{proposition}

\begin{proposition}
[\cite{li2021delayed}, Proposition 2.1 and Lemma B.1, property of projection]
\label{prop:proj}
Suppose that $\mA$ is a $p \times q$ matrix.
Let $\gP_{\mA}$ be the projection onto the column space of $\mA$ and $\gP_{\mA^\bot}$ the projection onto the null space of $\mA^\top$. Then we have
\begin{enumerate}
    \item Linearity: $\gP_{\mA} (\alpha \vx + \beta \vy) = \alpha \gP_{\mA} (\vx) + \beta \gP_{\mA} (\vy)$ for any $\vx, \vy \in \sR^p$ and $\alpha, \beta \in \sR$.
    \item Non-expansiveness: $ \max\{ \| \gP_{\mA} (\vx) - \gP_{\mA} (\vy) \|, \| \gP_{\mA^\bot} (\vx) - \gP_{\mA^\bot} (\vy) \|  \} \le \| \vx - \vy \| $ for any $\vx, \vy \in \sR^p$.
    \item Orthogonality: any $\vx \in \sR^p$ can be decomposed uniquely into $\vx = \vu + \vv$ where $\vu = \gP_{\mA^\bot}(\vx)$ and $\vv = \mP_{\mA} (\vx)$ satisfying $\langle \vu, \vv \rangle = 0$.
\end{enumerate}
More specifically, we have $\gP_{\mA}(\vx) = \mA (\mA^\top \mA)^\dag \mA^\top \vx = ( \mA^\top )^\dag \mA^\top \vx$ and $\gP_{\mA^\bot}(\vx) = \mI_p - \gP_{\mA} (\vx) = \left( \mI_p - \mA (\mA^\top \mA)^\dag \mA^\top \right) \vx = \left( \mI_p - ( \mA^\top )^\dag \mA^\top \right) \vx $ with $\dag$ the pseudo inverse.
\end{proposition}

\begin{proposition}[Stolz–Cesàro theorem]\label{prop:stolz}
Let $\{a_n\}$ and $\{b_n\}$ be two sequences of real numbers such that
\begin{enumerate}
    \item $0 < b_1 < b_2 < \dots < b_n < \dots $ and $\lim_{t \rightarrow \infty} b_t = \infty$.
    \item $\lim_{n \rightarrow \infty } \frac{ a_{n+1} - a_n }{ b_{n+1} - b_n } = l \in \sR$.
\end{enumerate}
Then, $\lim_{n \rightarrow \infty} \frac{a_n}{b_n}$ exists and is equal to $l$.
\end{proposition}


\begin{lemma}\label{lem:converge_r_t}
Let $\{r_n\} \subset (0, 1)$ be a sequence of positive numbers that decays to zero monotonically.
If $\frac{r_n}{r_{n+1}} - 1 = o (r_{n})$, for $p \ge 1$, we have that
\[
\lim_{T \rightarrow \infty} \frac{ \sum_{n=1}^T r_n^p \prod_{s=n+1}^T (1 - r_s) }{ r_T^{p-1} } = 1.
\]
\end{lemma}

\begin{lemma}\label{lem:converge_r_t_2}
Let $\{r_n\} \subset (0, 1)$ be a sequence of positive numbers that decays to zero monotonically and $a$ is a positive number.  
If $ \frac{r_n}{r_{n+1}} - 1 = a r_{n} +  o(r_n)  $, for $p \ge 1$ and $1 / a > p - 1$, we have 
\[
\lim_{T \rightarrow \infty} \frac{ \sum_{n=1}^T r_t^p \prod_{s=n+1}^T (1 - r_s) }{ r_T^{p-1} } = \frac{1}{1 - a (p - 1)}.
\]
\end{lemma}

\begin{lemma}\label{lem:converge_r_t_3}
Let $\{ r_n \} \subset (0, 1)$ be a sequence of positive numbers that decays to zero monotonically and $\{ s_n \}$ is a sequence of positive numbers.
If $ \frac{r_n}{r_{n+1}} - 1 = a r_{n} +  o(r_n)  $ for $a \ge 0$ 
and
$s_{n+1} \le (1 - r_n) s_n + o(r_n)$.
Then we have $s_n = o(1)$.
\end{lemma}
The proof of the three lemmas are deferred to Appendix \ref{sec:append_lem_r_n}.

\subsection{Proof of Theorem \ref{thm:converge}}
In this subsection, we give the formal statement of Theorem \ref{thm:converge} and its proof.
Before that, we first present the one-step descent lemmas of $\sE \norm{ \vu_n - \vx^\star }^2$ and $\sE \norm{ \vv_n }^2$, whose proof is deferred to Appendix \ref{sec:append_one_step_des}.

\begin{lemma}[One-step descent of $\sE \norm{\vu_{n} - \vx^\star}^2$ ]\label{lem:des_u_n}
Suppose that Assumptions \ref{asp:smooth}, \ref{asp:str_cov} and \ref{asp:noi_lip} hold.
Then there exists a $n_0$ such that for any $n \ge n_0$,
\begin{equation}\label{eq:recur_u_n}
\begin{aligned}
   \sE\|\vu_{n+1}-\vx^\star\|^2\le (1-\mu\eta_n)\sE\|\vu_n-\vx^\star\|^2+\frac{3L^2}{\mu}\eta_n\sE\|\vv_n\|^2+ 2 \eta_n^2
    \Sigma_\star^{(1)},
\end{aligned}
\end{equation}
where $\Sigma_\star^{(1)}:= \sE \norm{ \gP_{\mA^\bot} \xi^\star }^2$ with $\xi^\star = \nabla f(\vx^\star) - \nabla f(\vx^\star, \zeta), \zeta \sim \gD$.

\end{lemma}

\begin{lemma}[One-step descent of $\sE \norm{\vv_{n} }^2$ ]\label{lem:des_v_n}
Suppose that Assumptions \ref{asp:smooth}, \ref{asp:str_cov} and \ref{asp:noi_lip} hold.
Then there exists a $n_0$ such that for any $n \ge n_0$
\begin{equation}\label{eq:recur_v_n}
\begin{aligned}
   \sE\|\vv_{n+1}\|^2\le \left( 1 - \frac{p_n}{2} \right)\sE\|\vv_n\|^2+\frac{7L^2\eta_n^2}{p_n}\sE\|\vu_n-\vx^\star\|^2+\frac{7L^2\eta_n^2}{p_n}\|\nabla f(\vx^\star)\|^2 +2 \eta_n^2\Sigma_\star^{(2)},
\end{aligned}
\end{equation}
where $\Sigma_\star^{(2)}:= \sE \norm{ \gP_{\mA} \xi^\star }^2$ with $\xi^\star = \nabla f(\vx^\star) - \nabla f(\vx^\star, \zeta), \zeta \sim \gD$.

\end{lemma}

Now we are prepared to give the formal statement of Theorem \ref{thm:converge}.


\begin{proof}[
{Proof of Theorem \ref{thm:converge}}]

Let $z_n = \| \vu_n - \vx^\star \|^2 + c_0 \sqrt{ \frac{p_n}{ \eta_n} } \| \vv_n \|^2 $ with $c_0 = \sqrt{3 / (7 \mu)}$.
By Lemmas \ref{lem:des_u_n} and \ref{lem:des_v_n}, there exists a $n_0$ such that for any $n \ge n_0$, we have
\begin{align*}
    \sE z_{n+1} 
    & \le \left( 1 - \min \left\{ \mu \eta_n, \frac{p_n}{2}  \right\} + 7 c_0 L^2 \frac{\eta_n^{3/2} }{p_n^{1/2} } 
    \right) \sE z_n + 2 \eta_n^2 \Sigma_\star^{(1)} \\
    & \quad \ \ + 7 c_0 L^2 \frac{\eta_n^{3/2} }{p_n^{1/2}} \| \nabla f(\vx^\star) \|^2 + 2 c_0 \eta_n^{3/2} p_n^{1/2} \Sigma_\star^{(2)}.
\end{align*}
With $p_n = \min\{ \gamma \eta_n^\beta, 1\}$ for some $0 \le \beta < 1$, there exists a $n_1 \ge n_0$ such that for any $n \ge n_1$, we have $p_n = \gamma \eta_n^\beta$ and
\begin{align}\label{eq:recur_z_n}
    \sE z_{n+1} 
    & \le \left( 1 - \frac{ \mu \eta_n }{2}
    \right) \sE z_n + 2 \eta_n^2 
    \Sigma^{(1)}_\star
    + \frac{7 c_0 L^2}{\sqrt{\gamma}} \eta_n^{3/2 - \beta / 2}  \| \nabla f(\vx^\star) \|^2 + 2 c_0 \sqrt{\gamma} \eta_n^{3/2 + \beta / 2} 
    \Sigma^{(2)}_\star.
\end{align}

For any $T \ge n_1$, applying the recursion \eqref{eq:recur_z_n} $(T - n_1)$ times yields 
\begin{equation}\label{eq:after_recur_z_n}
\begin{aligned}
    \sE z_T 
    & \le \sE z_{n_1} \prod_{n = n_1}^{T-1} \left( 1 - \frac{\mu \eta_n}{2} \right)  + 
    2 \Sigma^{(1)}_\star \sum_{n = n_1}^{T-1} \eta_n^2 \prod_{s=n+1}^{T-1} \left( 1 - \frac{\mu \eta_s}{2} \right)  \\
    & \quad \ \ + \frac{7 c_0 L^2}{\sqrt{\gamma}} \| \nabla f(\vx^\star) \|^2  \sum_{n = n_1}^{T-1} \eta_n^{3/2 - \beta / 2} \prod_{s=n+1}^{T-1} \left( 1 - \frac{\mu \eta_s}{2} \right) \\
    & \quad \ \ + 2 c_0 \sqrt{\gamma} \Sigma_\star^{(2)} \sum_{n = n_1}^{T-1} \eta_n^{3 / 2 + \beta / 2} \prod_{s=n+1}^{T-1} \left( 1 - \frac{\mu \eta_s}{2} \right).
\end{aligned}
\end{equation}
For case (i) where $0 < \alpha < 1$, we have
$\eta_n = \eta_0 n^{- \alpha} $.
Thus, for the first term, we have
\begin{align*}
    \sE z_{n_1} \prod_{n=n_1}^{T-1} \left( 1 - \frac{\mu \eta_n}{2} \right)
    & \le \sE z_{n_1} \exp \left( - \frac{\mu}{2} \sum_{n=n_1}^{T-1} \eta_n \right) \\
    & \le \sE z_{n_1} \exp \left( - \frac{\mu \eta_0 ( T^{1 - \alpha} - n_1^{1 - \alpha} ) }{2 (1 - \alpha) }   \right).
\end{align*}
For other terms,
one can check that $ \frac{\eta_n}{ \eta_{n+1} } - 1 = o (\eta_n)$.
Then by Lemma \ref{lem:converge_r_t}, we have
\begin{align}\label{eq:u_n_order_primal}
\sE \| \vu_n - \vx^\star \|^2 \le
\sE z_n 
& {=} 
\gO \left( \frac{c_0 L^2 \norm{\nabla f (\vx^\star)}^2}{ \sqrt{\gamma}\, \mu }  \eta_n^{1 / 2 - \beta / 2} \right) 
{=} \gO \left( \frac{L^2 \norm{\nabla f (\vx^\star)}^2}{ \sqrt{\gamma \mu} }  \eta_n^{1 / 2 - \beta / 2} \right) .
\end{align}
This implies that there exists a positive number $c_1$ such that 
$ \sE \| \vu_n - \vx^\star \|^2 \le c_1 \frac{L^2 \norm{\nabla f (\vx^\star)}^2}{ \sqrt{\gamma \mu}} \eta_n^{1 / 2 - \beta / 2} $ for any $n \ge n_1$.
Substituting this into \eqref{eq:recur_v_n} yields that
\begin{align*}
    \sE\|\vv_{n+1}\|^2
    & \le \left( 1 - \frac{ \gamma \eta_n^\beta }{2} \right)\sE\|\vv_n\|^2 +
    \frac{7 c_1 L^4 \norm{\nabla f (\vx^\star)}^2 }{\gamma^{3/2} \sqrt{\mu} }
    \eta_n^{5/2 - 3\beta / 2}
    \sE\|\vu_n-\vx^\star\|^2 \\
    & \quad \ +
    \frac{7 L^2}{\gamma} \eta_n^{2 - \beta}
    \|\nabla f(\vx^\star)\|^2 + 2 \eta_n^2 \Sigma^{(2)}_\star
\end{align*}
hold for any $n \ge n_1$.
Following the same argument as before, we can prove 
\begin{align}\label{eq:v_n_order_primal}
    \sE \| \vv_{n} \|^2 = \gO \left( \frac{L^2 \norm{\nabla f (\vx^\star)}^2 }{ \gamma^{2 } } \eta_n^{2 - 2\beta} \right).
\end{align}
Then there exists $c_2 > 0$ such that
 $\sE \| \vv_{n} \|^2 \le 
 \frac{c_2 L^2 \norm{ \nabla f(\vx^\star) }_2^2 }{ \gamma^2 }
 \eta_n^{2 - 2\beta}$ for any $n \ge n_1$.
Substituting this into \eqref{eq:recur_u_n} yields that
\begin{align*}
    \sE\|\vu_{n+1}-\vx^\star\|^2\le (1-\mu\eta_n)\sE\|\vu_n-\vx^\star\|^2+
    \frac{3 c_2 L^4 \norm{ \nabla f(\vx^\star) }_2^2 }{\mu \gamma^2} 
    \eta_n^{3 - 2 \beta}
    + 2 \eta_n^2 \Sigma^{(1)}_\star
\end{align*}
hold for any $n \ge n_1$.
Following the same procedure again, we can obtain 
\begin{align}\label{eq:u_n_order_fianl}
\sE \| \vu_n - \vx^\star \|^2 = \gO 
\left( \frac{\Sigma^{(1)}_\star}{\mu} \eta_n + \frac{L^4 \norm{\nabla f(\vx^\star)}^2 }{ \mu^2 \gamma^2 } \eta_n^{ 2 - 2\beta } \right)
.\end{align}
For case (ii) where $\alpha = 1$ with $\eta_0 > 2 / \mu$, we can still obtain \eqref{eq:after_recur_z_n}.
Since $\eta_n = \eta_0 t^{-1}$, for the first term on the right-hand side of \eqref{eq:after_recur_z_n}, we have
\begin{align*}
    \sE z_{n_1} \prod_{n=n_1}^{T-1} \left( 1 - \frac{\mu \eta_n}{2} \right)
    & \le \sE z_{n_1} \exp \left( - \frac{\mu}{2} \sum_{n=n_1}^{T-1} \eta_n \right) \\
    & \le \sE z_{n_1} \exp \left( - \frac{\mu \eta_0 ( \ln T - \ln n_1 ) }{2}   \right) \\
    & = \gO \left( T^{- \mu \eta_0 / 2} \right).
\end{align*}
For other terms,
one can check that $ \frac{\eta_n}{ \eta_{t+1} } - 1 = \frac{2}{\mu \eta_0} \cdot \frac{\mu \eta_n}{2} + o (\eta_n)$.
Then by Lemma \ref{lem:converge_r_t_2},
we have \eqref{eq:u_n_order_primal} holds
for $\eta_0 > 2 / \mu$. 
Following the same procedure as before, we can also obtain \eqref{eq:v_n_order_primal}. 
Substituting this into \eqref{eq:recur_u_n} yields that
\begin{align*}
    \sE\|\vu_{n+1}-\vx^\star\|^2\le (1-\mu\eta_n)\sE\|\vu_n-\vx^\star\|^2+\frac{3 c_2 L^2}{\mu} \eta_n^{3 - 2 \beta}
    + 2 \eta_n^2\Sigma_\star^{(1)}
\end{align*}
holds for any $n \ge n_1$.
Since $ \frac{\eta_n}{ \eta_{t+1} } - 1 = \frac{1}{\mu \eta_0} \cdot \mu \eta_n  + o (\eta_n)$, 
following the same procedure as before,
we can obtain \eqref{eq:u_n_order_fianl}
for $\eta_0 > 2 / \mu >  \max\{ 2 - 2 \beta, 1 \} / \mu$.
\end{proof}

%

\subsection{Proof of Theorem \ref{thm:counter}}\label{sec:append_counter}
We first give a more detailed statement of Theorem \ref{thm:counter}.
\begin{theorem}\label{thm:counter_formal}
If $\eta_n = \eta_0 n^{- \alpha}$ and $p_n = \min\{ \gamma \eta_n, 1 \}$ with $0 < \alpha \le 1$, for a specific $\mA \in \sR^{p \times r}$ with $r < p$, there exists a quadratic function $f(\vx)$ defined on $\sR^p$ so that $\nabla^2 f(\vx) \succeq \rmI_p$ and  $\sE \norm{ \vu_n {-} \vx^\star }^2$ does not converge to $0$. Here $\rmI_p \in \RB^{p \times p}$ is the identity matrix, and  $\nabla^2 f(\vx) \succeq  \rmI_p$ means $\nabla^2 f(\vx) {-} \rmI_p$ is positive semidefinite.
Moreover, if $\gP_{\mA}$ is not of the form $\gP_{\mA} = \sum_{i \in I} \ve_i \ve_i^\top$, where $I \subseteq \{1,2,\dots, p \}$ and $\ve_i$ is the unit vector in $\sR^p$ with the $i$-th element equal to $1$, $\nabla^2 f (\vx)$ can be chosen as a diagonal matrix such that $\nabla^2 f(\vx) \succeq \rmI_p $.
\end{theorem}


\begin{proof}[Proof of Theorem \ref{thm:counter_formal}]

Consider the quadratic function
$f(\vx) = \frac{1}{2} \vx^\top \mB \vx + \vc^\top \vx$ where the positive definite matrix
$\mB \in \sR^{p \times p}$ and the vector $\vc \in \sR^{p}$ are specified later.

\vspace{-0.2cm}
\paragraph{The exact solution to problem~\eqref{eq:proj_opt} }
We first compute the exact solution to problem \eqref{eq:proj_opt}, where $\mA \in \sR^{p \times r}$ for some positive integer $r < p$. Without loss of generality, we assume $\mathrm{rank} (\mA) = r$.

Suppose that the singular value decomposition (SVD) of $\mA$ is $\mA = \mU \mD_A \mV^\top$ where $\mU \in \sR^{p \times p} $ and $\mV \in \sR^{r \times r} $ are orthogonal matrices and $\mD_A \in \sR^{p \times r}$ is a rectangular diagonal matrix with diagonal entries in descending order. One can check that the solution to $\mA^\top \vx = \vzero$ has the form $\vx = \left(\rmI_p - (\mA^\top)^\dag \mA^\top \right) \vw = \gP_{\mA^\bot} (\vw)$ where $\vw$ is an arbitrary vector in 
$
\sR^p$ and
$(\mA^\top)^\dag$ is the pseudo inverse of $\mA^\top$.
From the SVD of $\mA$, we have 
\[
\rmI_p - (\mA^\top)^\dag \mA^\top
= \mU
\begin{bmatrix}
\vzero_{r} & \vzero_{r \times (p-r)} \\
\vzero_{(p-r) \times r} & 
\rmI_{p-r}
\end{bmatrix}
\mU^\top,
\]
where $\vzero_{m \times n} \in \sR^{m \times n}$ denote the zero matrix and reduces to $\vzero_{n} \in \sR^{n \times n}$ for $m=n$.
We denote the first $r$ columns of $\mU$ by $\mU_1$ and last $p-r$ columns of $\mU$ by $\mU_2$
for simplicity, 
Then the problem \eqref{eq:proj_opt} becomes the following unconstrained problem
\begin{align*}
    & \quad\; \min_{\vw \in \sR^p} \frac{1}{2} \vw^\top \left(\rmI_p - (\mA^\top)^\dag \mA^\top \right)^\top \mB \left(\rmI_p - (\mA^\top)^\dag \mA^\top \right) \vw + \vw^\top \left(\rmI_p - (\mA^\top)^\dag \mA^\top \right)^\top \vc.
     \\
     & = \min_{{\vw}_2 \in \sR^{p-r} } \frac{1}{2} {\vw}_2^\top \mB_2 {\vw}_2 + {\vw}_2^\top {\vc}_2,
\end{align*}
where $\vw_2 = \mU_2^\top \vw$, $\mB_2 = \mU_2^\top \mB \mU_2$ and $\vc_2 = \mU_2^\top \vc$.
The solution is $\vw_2^\star = -\mB_2^{-1} \vc_2$.  
From the expression of $\rmI_p - ( \mA^\top )^\dag \mA^\top $, we know that the first $r$ elements of $\mU^\top \vw$ will not affect the value of $\vx$. Thus, the solution to the original problem \eqref{eq:proj_opt} is $\vx^\star = - \mU_2 \mB_2^{-1} \vc_2$.

Moreover, one can check
\[
\gP_{\mA} =
\mU
\begin{bmatrix}
\rmI_r & \vzero_{r \times (p-r)} \\
\vzero_{(p-r) \times r} & \vzero_{(p-r) \times (p-r)]}
\end{bmatrix}
 \mU^\top
 = \mU_1 \mU_1^\top
\]
and 
\[
\gP_{\mA^\bot} = 
\mU
\begin{bmatrix}
\vzero_{ r \times r}
& \vzero_{r \times (p-r)} \\
\vzero_{(p-r) \times r} & 
\rmI_{p-r}
\end{bmatrix}
 \mU^\top
 = \mU_2 \mU_2^\top.
\]


\vspace{-0.2cm}
\paragraph{Recursions of $\sE\vu_n$ and $\sE\vv_n$}
From the definition of $\vu_n$ and the linearity of $\gP_{\mA^\bot}$, we have
\begin{align*}
    \vu_{n+1} - \vx^\star
    & = \gP_{\mA^\bot} ( \vx_n - \eta_n \mB \vx_n - \eta_n \vc + \eta_n \xi_n ) - \vx^\star \\
    & = \vu_n - \vx^\star - \eta_n \gP_{\mA^\bot} (\mB \vx_n + \vc) + \eta_n \gP_{\mA^\bot} \xi_n \\
    & = \vu_n - \vx^\star - \eta_n \gP_{\mA^\bot} \mB ( \vu_n - \vx^\star) - \eta_n \gP_{\mA^\bot} \mB \vv_n - \eta_n \gP_{\mA^\bot} (\mB \vx^\star + \vc) + \eta_n \gP_{\mA^\bot} \xi_n \\
    & = \vu_n - \vx^\star - \eta_n \gP_{\mA^\bot} \mB \gP_{\mA^\bot} ( \vu_n - \vx^\star) - \eta_n \gP_{\mA^\bot} \mB \vv_n - \eta_n \gP_{\mA^\bot} (\mB \vx^\star + \vc) + \eta_n \gP_{\mA^\bot} \xi_n.
\end{align*}
The optimality of $\vx^\star$ implies that $ \gP_{\mA^\bot} ( \mB \vx^\star + \vc ) = \vzero $.
Taking expectation yields
\begin{align}\label{eq:counter_u_ite}
    \sE \vu_{n+1} - \vx^\star
    = (\rmI_{p} - \eta_n \gP_{\mA^\bot} \mB \gP_{\mA^\bot}) 
    (\sE \vu_n - \vx^\star)
    - \eta_n \gP_{\mA^\bot} \mB \, \sE \vv_n.
\end{align}
As for the iteration of $\sE \vv_{t}$.
From the definition of $\vv_n$, with probability $1 - p_n $ we have
\begin{align*}
    \vv_{t+1}
    & = \gP_{\mA} (\vx_n - \eta_n \mB \vx_n - \eta_n \vc + \eta_n \xi_n) \\
    & = \vv_n - \eta_n \gP_{\mA} (\mB \vx_n + \vc) + \eta_n \gP_{\mA} \xi_n \\
    & = \vv_n - \eta_n \gP_{\mA} \mB (\vu_n - \vx^\star) - \eta_n \gP_{\mA} \mB \vv_n - \eta_n \gP_{\mA} (\mB \vx^\star + \vc) + \eta_n \gP_{\mA} \xi_n \\
    & = (\rmI_{p} - \eta_n \gP_{\mA} \mB \gP_{\mA} ) \vv_n - \eta_n \gP_{\mA} \mB (\vu_n - \vx^\star) - \eta_n \gP_{\mA} (\mB \vx^\star + \vc) + \eta_n \gP_{\mA} \xi_n,
\end{align*}
and with probability $p_n$ we have $\vv_{n+1} = \vzero$.
Taking expectation yields 
\begin{align}\label{eq:counter_v_ite}
    \sE \vv_{n+1} {=} (1 {-} p_n) ( \rmI_{p} {-} \eta_n \gP_{\mA} \mB \gP_{\mA} ) \sE \vv_{n} - (1 {-} p_n) \eta_n \left[ \gP_{\mA} \mB (\sE \vu_n {-} \vx^\star ) + \gP_{\mA} (\mB \vx^\star {+} \vc) \right] .
\end{align}

\vspace{-0.2cm}
\paragraph{Simultaneous diagonalization of $\gP_{\mA} \mB \gP_{\mA}$ and $\gP_{\mA^\bot} \mB \gP_{\mA^\bot}$  }
We first express the two matrices as follows:
\begin{align*}
    \gP_{\mA} \mB \gP_{\mA}
    & = \mU 
    \begin{bmatrix}
    \rmI_{r} & \vzero_{ r \times (n-r) } \\
    \vzero_{(n-r)\times r}  & 
    \vzero_{n-r}
    \end{bmatrix}
    \mU^\top \mB \mU
    \begin{bmatrix}
    \rmI_{r} & \vzero_{ r \times (n-r) } \\
    \vzero_{(n-r)\times r}  & 
    \vzero_{n-r}
    \end{bmatrix}
    \mU^\top \\
    & = \mU 
    \begin{bmatrix}
    \mB_1 & \vzero_{ r \times (n-r) } \\
    \vzero_{(n-r)\times r}  & 
    \vzero_{n-r}
    \end{bmatrix}
    \mU^\top, \\
    \gP_{\mA^\bot} \mB \gP_{\mA^\bot}
    & = \mU 
    \begin{bmatrix}
    \vzero_{r } & \vzero_{ r \times (n-r) } \\
    \vzero_{(n-r)\times r}  & 
    \rmI_{n-r}
    \end{bmatrix}
    \mU^\top \mB \mU
    \begin{bmatrix}
    \vzero_{r} & \vzero_{ r \times (n-r) } \\
    \vzero_{(n-r)\times r}  & 
    \rmI_{n-r}
    \end{bmatrix}
    \mU^\top \\
    & = \mU 
    \begin{bmatrix}
    \vzero_{r} & \vzero_{ r \times (n-r) } \\
    \vzero_{(n-r)\times r}  & 
    \mB_2
    \end{bmatrix}
    \mU^\top,
\end{align*}
where $\mB_1 = \mU_1^\top \mB \mU_1$ and $\mB_2 = \mU_2^\top \mB \mU_2$ are positive definite.
We 
suppose the eigenvalue decomposition of $\mB_1$ and $\mB_2$ is $\mB_1 = \mQ_1 \mD_{B_1} \mQ_1^\top$ and $\mB_2 = \mQ_2 \mD_{B_2} \mQ_2^\top$. With 
$\mQ := 
\begin{bmatrix}
\mQ_1 & \vzero_{ r \times (n-r) } \\
\vzero_{ (n-r) \times r} & \mQ_2
\end{bmatrix}
$ and $\mP := \mU \mQ$,
we obtain the eigenvalue decomposition of $\gP_{\mA^\bot} \mB \gP_{\mA^\bot}$ and $\gP_{\mA} \mB \gP_{\mA}$ as follows
\begin{align*}
    \gP_{\mA} \mB \gP_{\mA}
    & = 
    \mP
    \begin{bmatrix}
    \mD_{B_1} & \vzero_{r \times (n-r)} \\
    \vzero_{ (n-r) \times r } & \vzero_{ n-r }
    \end{bmatrix}
    \mP^\top
    =: \mP \widetilde{\mD}_{B_1} \mP^\top
\end{align*}
and
\begin{align*}
    \gP_{\mA^\bot} \mB \gP_{\mA^\bot}
    & = 
    \mP
    \begin{bmatrix}
    \vzero_{ r  } & \vzero_{ r \times (n-r) } \\
    \vzero_{ (n-r) \times r } & \mD_{B_2}
    \end{bmatrix}
    \mP^\top 
    =: \mP \widetilde{\mD}_{B_2} \mP^\top.
\end{align*}

\vspace{-0.2cm}
\paragraph{Proof by contradiction}
Left multiplication of 
\eqref{eq:counter_v_ite} by
$\mP^\top$ 
yields
\begin{align}\label{eq:counter_v_ite_2}
    \sE \tilde{\vv}_{n+1} 
    & = (\rmI_{p} - \widetilde{\mD}_n)
    \sE \tilde{\vv}_{n} - \eta_n (1 - p_n) \mB_0 (\sE {\vu}_n - \vx^\star) - (1 - p_n) \eta_n \vc_0.
\end{align}
where $\tilde{\vv}_{n} := \mP^\top \vv_{n} $, 
$\mB_0 := \mP^\top \gP_{\mA} \mB$,
$\widetilde{\mD}_n := \eta_n \widetilde{\mD}_{B_1} + p_n \rmI_{p} - \eta_n p_n \widetilde{\mD}_{B_1}$
and $\vc_0 := \mP^\top \gP_{\mA} ( \mB \vx^\star + \vc )$.
Adding $(\gamma \rmI_{p} + \widetilde{\mD}_{B_1})^{-1} \vc_0 $ to both sides of \eqref{eq:counter_v_ite_2}, we obtain
\begin{align}
    & \quad \ \sE \tilde{\vv}_{n+1} + (\gamma \rmI_{p} + \widetilde{\mD}_{B_1})^{-1} \vc_0 \nonumber \\
    & =  (\rmI_{p} - \widetilde{\mD}_n)  \sE \tilde{\vv}_{n} - \eta_n (1 - p_n) \mB_0 (\sE {\vu}_n - \vx^\star) - (1 - p_n) \eta_n \vc_0 + (\gamma \rmI_{p} + \widetilde{\mD}_{B_1})^{-1} \vc_0 \nonumber \\
    & = (\rmI_{p} - \widetilde{\mD}_n) [ \sE \tilde{\vv}_{n} + (\gamma \rmI_{p} + \widetilde{\mD}_{B_1})^{-1} \vc_0 ] - \eta_n (1 - p_n) \mB_0 (\sE {\vu}_n - \vx^\star) \nonumber \\
    & \quad \ + [ p_n - \gamma \eta_n (1 - p_n) ] (\gamma \rmI_{p} + \widetilde{\mD}_{B_1})^{-1} \vc_0. \label{eq:counter_v_ite_3}
\end{align}
Suppose $\sE \| \vu_n - \vx^\star \|^2 = o(1)$, which implies $\sE \vu_n - \vx^\star = o(1)$.
Let $\widetilde{\mD}_n = \mathrm{diag} \left\{ \tilde{d}_{n,1}, \tilde{d}_{n,2}, \dots, \tilde{d}_{n,p} \right\} $ amd $\widetilde{\mD}_{B_1} = \mathrm{diag} \left\{ \tilde{d}_{B_1, 1}, \tilde{d}_{B_1, 2}, \dots, \tilde{d}_{B_1, p}  \right\} $. Left multiplication of 
\eqref{eq:counter_v_ite_3} by $\ve_i^\top$ gives
\begin{align*}
    \left| \sE \ve_i^\top \tilde{\vv}_{n+1} + \frac{1}{\gamma + \tilde{d}_{B_1, i} } \ve_i^\top \vc_0  \right| 
    \le (1 -  \tilde{d}_{n,i}) \left|  \sE  \ve_i^\top \tilde{\vv}_{n} + \frac{1}{\gamma + \tilde{d}_{B_1, i}  } \ve_i^\top \vc_0  \right| + o(\eta_n),
\end{align*}
where $\ve_i$ is the unit vector with the $i$-th element equal to $1$.
Since $\tilde{d}_{n,i} = \eta_n d_{B_1, i} (1 - p_n) + p_n$ and $p_n = \min \{ \gamma \eta_n, 1 \}$, $o(\eta_n) = o(\tilde{d}_{n,i})$.
Lemma \ref{lem:converge_r_t_3} implies $\sE \ve_i^\top \tilde{\vv}_{n+1} = - \frac{1}{\gamma + d_{B_1, i} } \ve_i^\top \vc_0 + o(1) $.
It follows that $ \sE \tilde{\vv}_{n} = - (\gamma \rmI_p + \widetilde{\mD}_{B_1})^{-1} \vc_0 + o(1) $.
Thus we have
\begin{align}\label{eq:counter_v_limit}
    \sE \vv_{n} = - \mP (\gamma \rmI_p + \widetilde{\mD}_{B_1})^{-1} \mP^\top \gP_{\mA} ( \mB \vx^\star + \vc ) + o(1).
\end{align}
Denote the limit of $\sE \vv_n$ by $\vv_\infty$ and we come back to the iteration \eqref{eq:counter_u_ite}.
Left multiplication
of \eqref{eq:counter_u_ite} by $\mP^\top$ yields
\begin{align*}
    \sE \tilde{\vu}_{n+1}
    = (\rmI_p - \eta_n \widetilde{\mD}_{B_2}) \sE \tilde{\vu}_n - \eta_n \mP^\top \gP_{\mA^
    \bot} \mB \vv_\infty + o(\eta_n),
\end{align*}
where $\tilde{\vu}_n = \mP^\top (\vu_n - \vx^\star)$.
Similar to the above argument, adding $\mP^\top \gP_{\mA^\bot} \mB \vv_{\infty} $ to both sides and using Lemma \ref{lem:converge_r_t_3}, we can obtain 
\begin{align*}
    \sE {\tilde{\vu}_n} 
    & = - \mP^\top \gP_{\mA^\bot} \mB \vv_{\infty} + o(1) \\
    & = \mP^\top \gP_{\mA^\bot} \mB \mP (\gamma \rmI_p + \widetilde{\mD}_{B_1})^{-1} \mP^\top \gP_{\mA} ( \mB \vx^\star + \vc ) + o(1).
\end{align*}
It remains to
prove that there exists a positive definite matrix $\mB \in \sR^{p \times p}$ and a vector $\vc \in \sR^p$ such that the limit is nonzero.

\vspace{-0.2cm}
\paragraph{Specification of $\mB$ and $\vc$}
From the expression of $\vx^\star$, we have
\begin{align*}
    \mB \vx^\star + \vc
    = \vc - \mB \mU_2 ( \mU_2^\top \mB \mU_2 )^{-1} \mU_2^\top \vc
    = (\rmI_p - \mB \mU_2 ( \mU_2^\top \mB \mU_2 )^{-1} \mU_2^\top) \vc.
\end{align*}
Define $\widetilde{\mB} := \rmI_p - \mB \mU_2 ( \mU_2^\top \mB \mU_2 )^{-1} \mU_2^\top$ for short. We examine the column space of $\widetilde{\mB}$, which is denoted by $\gR (\widetilde{\mB})$.
We can easily find $\mU_2^\top \widetilde{\mB} = \vzero_{ (p-r) \times r }$. Thus $\gR (\widetilde{\mB}) \subseteq \gR (\mU_1)$.
On the other hand, we have $\widetilde{\mB} \mU_1 = \mU_1 $, which implies $\mathrm{rank} (\widetilde{\mB} ) \ge \mathrm{rank} (\mU_1)$.
As a result, $\gR( \widetilde{\mB}) = \gR ( \mU_1 ) $.
Then for any $\vz \in \sR^r$, there exists a $\vc \in \sR^p$ such that $\widetilde{\mB} \vc = \mU_1 \vz$. It suffices to prove that there exists a positive definite matrix $\mB \in \sR^{p \times p}$ and a vector $\vz \in \sR^r$ such that $\mP^\top \gP_{\mA^\bot} \mB \mP (\gamma \rmI_p + \widetilde{\mD}_{B_1})^{-1} \mP^\top \gP_{\mA} \mU_1 \vz$ is nonzero.

Since $\mP = \mU \mQ$, $\gP_{\mA} = \mU_1 \mU_1^\top$, $\gP_{\mA^\bot} = \mU_2 \mU_2^\top$, 
$\mQ = 
\begin{bmatrix}
\mQ_1 & \vzero_{ r \times (n-r) } \\
\vzero_{ (n-r) \times r} & \mQ_2
\end{bmatrix}$,
$\widetilde{\mD}_{B_1} =
\begin{bmatrix}
\mD_{B_1} & \vzero_{r \times (n-r)} \\
\vzero_{ (n-r) \times r } & \vzero_{ n-r }
\end{bmatrix}$
and
$\mB_1 = \mQ_1 \mD_{B_1} \mQ_1^\top$,
we have
\begin{align*}
    & \quad \ \mP^\top \gP_{\mA^\bot} \mB \mP (\gamma \rmI_p + \widetilde{\mD}_{B_1})^{-1} \mP^\top \gP_{\mA} \mU_1 \vz \\
    & = \mQ^\top \mU^\top \mU_2 \mU_2^\top \mB \mU \mQ ( \gamma \rmI_p + \widetilde{\mD}_{B_1} )^{-1} \mQ^\top \mU^\top \mU_1 \vz \\
    & = \mQ^\top 
    \begin{bmatrix}
    \vzero \\
    \mU_2^\top
    \end{bmatrix}
    \mB 
    [\mU_1 \ \mU_2 ]
    \begin{bmatrix}
    ( \gamma \rmI_r + \mB_1 )^{-1} & \vzero_{r \times (p-r)} \\
    \vzero_{ (p-r) \times r } & \frac{1}{\gamma} \rmI_{p-r}
    \end{bmatrix}
    \begin{bmatrix}
    \vz \\
    \vzero
    \end{bmatrix} \\
    & = \mQ^\top 
    \begin{bmatrix}
    \vzero \\
    \mU_2^\top \mB \mU_1 (\gamma \rmI_p + \widetilde{\mD}_{B_1})^{-1} \vz
    \end{bmatrix}.
\end{align*}
Then it suffices to prove that there exist a positive matrix $\mB \in \sR^{p \times p}$ such that $\mU_2^\top \mB \mU_1$ is nonzero.
Suppose that 
$\mU_1 = (\vp_1, \vp_2, \dots, \vp_r) = (p_{ki})_{p \times r}$ and $\mU_2 = (\vq_1, \vq_2, \dots, \vq_{p-r}) = (q_{kj})_{ p \times (p-r) }$. Then the column vectors of $\mU_1$ and $\mU_2$ form an orthonormal basis of $\sR^p$.

If there exist $i$,$j$ and $k_0$ such that  $p_{k_0i} q_{k_0j} \neq 0$,
Then we can take $\mB$ as a diagonal matrix $\rmI_p + \mE_{k_0 k_0}$ where $\mE_{ij}$ is the $p \times p$ matrix with $(i,j)$ entry equal to $1$ and others equal to $0$.
The $(j,i)$ entry of $\mU_2^\top \mB \mU_1$ is $\sum_{k=1}^{p} p_{ki} q_{kj} + p_{k_0 i} q_{k_0 j}  = p_{k_0 i} q_{k_0 j}  \neq 0 $.
And one can check $\mB \succeq \rmI_p$.

Otherwise, there must exist $i$, $j$, $k_0$ and $l_0$ such that $p_{k_0 i} q_{l_0 j} \neq 0$ and $k_0 \neq l_0$.
Since in this case $p_{ki} q_{kj} = 0$ for any $k$, then we have $q_{k_0 j} = p_{l_0 i} = 0$. We take $\mB = 2 \rmI_p + \mE_{k_0 l_0} + \mE_{l_0 k_0}$.
Then the $(j,i)$ entry of $\mU_2^\top \mB \mU_1$ is $2 \sum_{k=1}^p p_{ki} q_{kj} + p_{k_0 i} q_{l_0 j} + p_{l_0 i} q_{k_0 j} = p_{k_0 i} q_{l_0 j} \neq 0 $.
And one can check $\mB \succeq \rmI_p$.

As a result, there always exists $\mB$ and $\vc$ such that the limit of $ \norm{\sE\tilde{\vu}_n}$ is nonzero. This implies $\sE \norm{\vu_n - \vx^\star}^2 \neq o(1)$, which induces a contradiction.

In the latter case, for any  $\ve_i \in \sR^p$, either $\ve_i^\top \mU_1$ or $\ve_i^\top \mU_2$ is zero. Note that $\mU_1$ and $\mU_2$ are of full column rank. Then we have $\mU_1 = \sum_{i \in I_1} \ve_i \tilde{\vp}_i^\top $ and $\mU_2 = \sum_{j \in I_2} \ve_j \tilde{\vq}_j^\top$ where $|I_1| = r$, $|I_2| = p-r$, $I_1 \cup I_2 = \{1, 2, \dots, p \}$, $I_1 \cap I_2 = \varnothing$, $\tilde{\vp}_{i}~ (i \in I_1)$ are orthonormal basis of $\sR^r$ and $\tilde{\vq}_{j}~ (j \in I_2)$ are orthonormal basis of $\sR^{p-r}$.
As a consequence,
$\gP_{\mA} = \mU_1 \mU_1^\top = \sum_{i \in I_1} \ve_i \ve_i^\top$ and $\gP_{\mA^\bot} = \mU_2 \mU_2^\top = \sum_{j \in I_2} \ve_j \ve_j^\top$.
This implies that if $\gP_{\mA}$ is not of this form 
there must exist $i,j,k_0$ such that $p_{k_0 i} q_{k_0 j} \neq 0$. Then we can choose $\mB$ as a diagonal matrix such that $\mB \succeq \rmI_p$.
\end{proof}

Now we give a corresponding result of Theorem~\ref{thm:counter} in the context of distributed learning.
\begin{corollary}
\label{cor:counter}
Consider the problem \eqref{eq:dis_opt}.
If $\eta_n = \eta_0 n^{- \alpha}$ and $p_n = \min\{ p_0 \eta_n, 1 \}$ with $0 < \alpha \le 1$, then there exists a quadratic function $f(\vx)$ such that $ \nabla^2 f(\vx) $ is a diagonal matrix, $\nabla^2 f(\vx) \succeq \rmI_d$ and $\sE \norm{ \vu_n - \vx^\star }^2$ does not converge to $0$. 
\end{corollary}


\begin{proof}[Proof of Corollary \ref{cor:counter}]

With $\mA$ defined in \eqref{eq:mat_A_fl}, we have $p = Nd$ and $r = (N-1) d$.
Recall that for
\[\vx=\left[(\vx^{(1)})^\top,(\vx^{(2)})^\top,\cdots,(\vx^{(N)})^\top\right]^\top\in\sR^{Nd},\]
 we have
 \begin{align*}
\gP_{\mA^\bot} (\vx) = \left[ \bar{\vx}^\top, \bar{\vx}^\top, \cdots, \bar{\vx}^\top \right]^\top 
\end{align*}
where $\bar{\vx} = \frac{1}{N} \sum_{k=1}^N \vx^{(k)}$.
As a result, we have
$ \gP_{\mA^\bot} \ve_1 = \frac{1}{N} \sum_{k=0}^{N-1} \ve_{1 + kd } $.
This implies that $\gP_{\mA^\bot}$ can not be of the form $\gP_{\mA^\bot} = \sum_{i \in I} \ve_i \ve_i^\top$.
Thus, $\nabla^2 f(\vx)$ can be chosen as a diagonal matrix.
\end{proof}

\subsection{Proof of Lemmas \ref{lem:converge_r_t}, \ref{lem:converge_r_t_2} and \ref{lem:converge_r_t_3} }\label{sec:append_lem_r_n}

\begin{proof}[Proof of Lemma \ref{lem:converge_r_t}]

Define $a_T = \sum_{n=1}^T r_n^p \prod_{s=1}^n \frac{1}{1 - r_s}$ and $b_T = r_T^{p-1} \prod_{s=1}^T \frac{1}{1 - r_s}$.
We first prove that $b_{T+1} > b_T$ for sufficiently large $T$ and $\lim_{T \rightarrow \infty} b_T = \infty$.
Since
\begin{align*}
    \frac{ b_{T+1} }{ b_T }
    & = \left( \frac{ r_{T+1} }{ r_T } \right)^{p-1} \cdot \frac{1}{1 - r_{T+1}} \\
    & = \frac{1}{  (1 + o(r_{T}) )^{p-1} } \cdot (1 + r_{T+1} + o (r_{T+1}) )  \\
    & = (1 + o(r_{T+1})) (1 + r_{T+1} + o (r_{T+1}) ) \\
    & = 1 + r_{T+1} + o(r_{T+1}),
\end{align*}
then we have $b_{T+1} > b_T$ for sufficiently large $T$.
Besides, $\frac{r_n}{r_{n+1}} - 1 = o (r_{n}) $ implies $\frac{1}{r_n} - \frac{1}{r_{n+1}} = o(1)$. By Stolz–Cesàro theorem, we have $\lim_{n \rightarrow \infty} \frac{1}{n r_n} = 0$ and $\lim_{n \rightarrow \infty} \frac{ \sum_{s=1}^n 1 / s  }{ \sum_{s=1}^n r_s } = 0 $.
As a consequence, 
\begin{align*}
    b_T 
    & \ge r_T^{p-1} \exp\left( \sum_{s=1}^T r_s \right) \\
    & = r_T^{p-1} \exp \left( \frac{ \sum_{s=1}^T r_s }{ \sum_{s=1}^T 1 / s } \sum_{s=1}^T 1/s  \right) \\
    & \ge r_T^{p-1} \exp \left( \frac{ \sum_{s=1}^T r_s }{ \sum_{s=1}^T 1 / s } \log T \right) \\
    & = (T r_T)^{p-1} \exp \left[ \left( \frac{ \sum_{s=1}^T r_s }{ \sum_{s=1}^T 1 / s } - p + 1   \right) \log T \right].
\end{align*}
Thus $\lim_{T \rightarrow \infty} b_T = \infty$.
Now we use Stolz–Cesàro theorem to prove
$\lim_{T \rightarrow \infty} \frac{a_T}{b_T} = 1$.
With the definition of $a_T$ and $b_T$, we have
\[ a_{T+1} - a_T = r_{T+1}^p \prod_{s=1}^{T+1} \frac{1}{1 - r_s} \]
and 
\begin{align*}
b_{T+1} - b_T 
& = r_{T+1}^{p-1} \prod_{s=1}^{T + 1} \frac{1}{1 - r_s} - r_T^{p-1} \prod_{s=1}^T \frac{1}{1 - r_s} \\
& = (r_{T+1}^{p-1} - r_T^{p-1} ) \prod_{s=1}^{T+1} \frac{1}{1 - r_s} + r_T^{p-1} r_{T+1} \prod_{s=1}^{T+1} \frac{1}{1 - r_s}.  
\end{align*}
It follows that
\begin{align*}
    \frac{a_{T+1} - a_T}{b_{T+1} - b_T}
    & = 
    \frac{r_{T+1}^p}{r_{T+1}^{p-1} - r_T^{p-1} + r_T^{p-1} r_{T+1} } \\
    & = \frac{ r_{T+1} }{ 1 - ( r_T / r_{T+1} )^{p-1} + r_{T+1} ( r_T / r_{T+1} )^{p-1}  } \\
    & = \frac{ r_{T+1} }{ 1 - (1 + o(r_T))^{p-1} + r_{T+1} (1 + o(1)) } \\
    & = \frac{ r_{T+1} }{ o(r_T) + r_{T+1} (1 + o(1))  } \\
    & = \frac{1}{1 + o(1)},
\end{align*}
which implies $\lim_{T \rightarrow \infty} \frac{ a_{T+1} - a_T }{ b_{T+1} - b_T } = 1 $. By Stolz–Cesàro theorem, we obtain what we want.
\end{proof}

\begin{proof}[Proof of Lemma \ref{lem:converge_r_t_2}]

Define $a_T = \sum_{n=1}^T r_t^p \prod_{s=1}^n \frac{1}{1 - r_s}$ and $b_T = r_T^{p-1} \prod_{s=1}^T \frac{1}{1 - r_s}$.
We first prove that $b_{T+1} > b_T$ for sufficiently large $T$ and $\lim_{T \rightarrow \infty} b_T = \infty$.
Since
\begin{align*}
    \frac{ b_{T+1} }{ b_T }
    & = \left( \frac{ r_{T+1} }{ r_T } \right)^{p-1} \cdot \frac{1}{1 - r_{T+1}} \\
    & = \frac{1}{  (1 + a r_{T} + o(r_T) )^{p-1} } \left( 1 + r_{T+1} + o( r_{T+1} ) \right) \\
    & = \left( 1 - a(p-1) r_T + o (r_T) \right) \left( 1 + \frac{ r_{T} }{ 1 + a r_T + o(r_T) } + o( r_{T} ) \right) \\
    & = \left( 1 - a(p-1) r_T + o (r_T) \right) \left( 1 + r_T + o(r_T) \right) \\
    & = 1 + [1 - a(p-1)] r_T + o(r_T),
\end{align*}
then we have $b_{T+1} > b_T$ for sufficiently large $T$.
Besides, $\frac{r_n}{r_{n+1}} - 1 = a r_n +  o (r_{n}) $ implies $\frac{1}{r_n} - \frac{1}{r_{n+1}} = a + o(1)$. By Stolz–Cesàro theorem, we have $\lim_{t \rightarrow \infty} \frac{1}{n r_n} = a$ and $\lim_{n \rightarrow \infty} \frac{ \sum_{s=1}^n 1 / s  }{ \sum_{s=1}^n r_s } = a $.
As a consequence, 
\begin{align*}
    b_T 
    & \ge r_T^{p-1} \exp\left( \sum_{s=1}^T r_s \right) \\
    & = r_T^{p-1} \exp \left( \frac{ \sum_{s=1}^T r_s }{ \sum_{s=1}^T 1 / s } \sum_{s=1}^T 1/s  \right) \\
    & \ge r_T^{p-1} \exp \left( \frac{ \sum_{s=1}^T r_s }{ \sum_{s=1}^T 1 / s } \log T \right) \\
    & = (T r_T)^{p-1} \exp \left[ \left(
    \frac{ \sum_{s=1}^T r_s }{ \sum_{s=1}^T 1 / s }
    - p + 1   \right) \log T \right] \\
    & = (1 / a + o(1) )^{p-1} \exp \left[ \left( 1 / a + o(1)- p + 1   \right) \log T \right].
\end{align*}
Thus $\lim_{T \rightarrow \infty} b_T = \infty$.

Now we use Stolz–Cesàro theorem to prove $\lim_{T \rightarrow \infty} \frac{a_T}{b_T} = \frac{1}{1 - a(p-1)}$.
With the definition of $a_T$ and $b_T$, we have
\[ a_{T+1} - a_T = r_{T+1}^p \prod_{s=1}^{T+1} \frac{1}{1 - r_s} \]
and 
\begin{align*}
b_{T+1} - b_T 
& = r_{T+1}^{p-1} \prod_{s=1}^{T + 1} \frac{1}{1 - r_s} - r_T^{p-1} \prod_{s=1}^T \frac{1}{1 - r_s} \\
& = (r_{T+1}^{p-1} - r_T^{p-1} ) \prod_{s=1}^{T+1} \frac{1}{1 - r_s} + r_T^{p-1} r_{T+1} \prod_{s=1}^{T+1} \frac{1}{1 - r_s}.  
\end{align*}
It follows that
\begin{align*}
    \frac{a_{T+1} - a_T}{b_{T+1} - b_T}
    & = 
    \frac{r_{T+1}^p}{r_{T+1}^{p-1} - r_T^{p-1} + r_T^{p-1} r_{T+1} } \\
    & = \frac{ r_{T+1} }{ 1 - ( r_T / r_{T+1} )^{p-1} + r_{T+1} ( r_T / r_{T+1} )^{p-1}  } \\
    & = \frac{ r_{T+1} }{ 1 - ( 1 + a r_T + o(r_T) )^{p-1} + r_{T+1} ( 1 + a r_T + o(r_T) )^{p-1}  } \\
    & = \frac{ r_{T+1} }{ 1 - 1 - a (p-1) r_T + o(r_T) + r_{T+1} ( 1 +  o(1) )  } \\
    & = \frac{ 1 }{ - a (p-1) r_T / r_{T+1}  + o(1) +  1 + o(1)  } \\
    & = \frac{1}{1 - a(p-1) + o(1)},
\end{align*}
which implies $\lim_{T \rightarrow \infty} \frac{ a_{T+1} - a_T }{ b_{T+1} - b_T } = \frac{1}{1 - a(p-1)} $. By Stolz–Cesàro theorem, we obtain what we want.
\end{proof}

\begin{proof}[Proof of Lemma \ref{lem:converge_r_t_3}]

Suppose that $s_n = o(1)$ does not hold. Then for any positive number $\varepsilon > 0$, there exists a sequence of positive integers $\{ n_i \}$ that increases to $\infty$ such that $s_{n_i} \ge \varepsilon$. 
From the recursion of $s_n$, there exists a positive integer $T$ such that \begin{align}\label{eq:proof_recur_s}
    s_{n+1} \le (1 - r_n) s_n + \frac{\varepsilon}{2} r_n
\end{align}
for any $n \ge T$.
For $n_i > T$, we have
\[
\varepsilon \le s_{n_i} \le (1 - r_{n_i - 1}) s_{n_i-1} + \frac{\varepsilon}{2} r_{n_i - 1} 
\le (1 - r_{n_i - 1}) s_{n_i-1} + \varepsilon r_{n_i - 1}   .
\]
It follows that $s_{n_i-1} \ge \varepsilon$.
Since $n_i$ increases to $\infty$, we have $s_{n} \ge \varepsilon$ for any $n \ge T$ by induction.
For any $T_1 > T$, summing \eqref{eq:proof_recur_s} from $T$  to $T_1 - 1$ , we have
\begin{align*}
    \sum_{n=T+1}^{T_1} s_n \le \sum_{n=T}^{T_1-1} s_n - \sum_{n=T}^{T_1-1} s_n r_n + \frac{\varepsilon}{2} \sum_{n = n_T}^{T_1-1} r_n.
\end{align*}
Rearranging the terms yields
\begin{align*}
    s_T
    \ge s_{T_1} + \sum_{n=T}^{T_1-1} s_n r_n - \frac{\varepsilon}{2} \sum_{n = T}^{T_1-1} r_n 
    \ge s_{T_1} + \frac{\varepsilon}{2} \sum_{n = T}^{T_1-1} r_n.
\end{align*}
From the proofs of Lemmas \ref{lem:converge_r_t} and \ref{lem:converge_r_t_2}, we have $\lim_{n \rightarrow \infty} \frac{ \sum_{s=1}^n 1/s }{ \sum_{s=1}^n r_s } = a$. Thus $\lim_{T_1 \rightarrow \infty} \sum_{n=T}^{T_1} r_n = \infty$, which induces a contradiction.
As a consequence, we have $s_n = o(1)$.
\end{proof}

\subsection{Proof of Lemmas \ref{lem:des_u_n} and \ref{lem:des_v_n}}\label{sec:append_one_step_des}

\begin{proof}[Proof of Lemma \ref{lem:des_u_n}]

From the update rule and the linearity of $\gP_{\mA^\bot}$, we have
\begin{align}
    \sE \| \vu_{n+1} - \vx^\star \|^2
    & = \sE \| \gP_{\mA^\bot} ( \vx_n - \eta_n \nabla f(\vx_n) + \eta_n \xi_n ) - \vx^\star ) \|^2 \nonumber \\
    & = \sE \| \vu_n - \vx^\star - \eta_n \gP_{\mA^\bot} (\nabla f(\vx_n) - \nabla f(\vx^\star)) + \eta_n \gP_{\mA^\bot} \xi_n \|^2 \nonumber \\
    & = \sE \| \vu_n - \vx^\star \|^2 + \eta_n^2 \sE \| \gP_{\mA^\bot} (\nabla f(\vx_n) - \nabla f(\vx^\star)) \|^2 \nonumber \\
    & \quad \ - 2 \eta_n \sE \left\langle \vu_n - \vx^\star, \gP_{\mA^\bot} (\nabla f(\vx_n) - \nabla f(\vx^\star))  \right \rangle + \eta_n^2 \Sigma_n^{(1)}, \label{eq:proof_u_n_recur}
\end{align}
where the last equality is due to that $\{\xi_n\}$ is a m.d.s.
and $\Sigma_n^{(1)}:= \sE \| \gP_{\mA^\bot} \xi_n \|^2$.

For the second term of \eqref{eq:proof_u_n_recur}, we have
\begin{align*}
    \| \gP_{\mA^\bot} ( \nabla f(\vx_n) - \nabla f(\vx^\star)) \|^2
    & = \| \gP_{\mA^\bot} ( \nabla f(\vx_n) - \nabla f(\vu_n)) + \gP_{\mA^\bot} ( \nabla f(\vu_n) - \nabla f(\vx^\star))  \|^2 \\
    & \overset{(a)}{\le} 2 \| \gP_{\mA^\bot} ( \nabla f(\vx_n) - \nabla f(\vu_n)) \|^2 + 2 \| \gP_{\mA^\bot} ( \nabla f(\vu_n) - \nabla f(\vx^\star)) \|^2 \\
    & \overset{(b)}{\le} 2 L^2 \norm{\vv_n}^2 + 2 L^2 \norm{\vu_n - \vx^\star}^2,
\end{align*}
where (a) is by Proposition \ref{prop:cauchy} and (b) is due to non-expansiveness of $\gP_{\mA^\bot}$ and smoothness of $f$.
For the third term of \eqref{eq:proof_u_n_recur}, we have
\begin{align*}
    & \quad \ - \inner{\vu_n - \vx^\star}{ \gP_{\mA^\bot} (\nabla f(\vx_n) - \nabla f(\vx^\star) ) } \\
    & \overset{(a)}{=} - \inner{\vu_n - \vx^\star}{ \nabla f(\vx_n) - \nabla f(\vx^\star) } \\
    & = - \inner{\vu_n - \vx^\star}{ \nabla f (\vx_n) - \nabla f(\vu_n) } - \inner{ \vu_n - \vx^\star }{ \nabla f(\vu_n) - \nabla f(\vx^\star) } \\
    & \overset{(b)}{\le} \frac{\mu}{4} \norm{  \vu_n - \vx^\star }^2 + \frac{1}{\mu} \norm{ \nabla f(\vx_n) - \nabla f(\vu_n) }^2 - \mu \norm{ \vu_n - \vx^\star }^2 \\
    & \overset{(c)}{\le} - \frac{3 \mu}{4} \norm{ \vu_n - \vx^\star }^2 + \frac{L^2}{\mu} \norm{\vv_n}^2,
\end{align*}
where (a) follows from the orthogonality between $\gP_{\mA}$ and $\gP_{\mA^\bot}$, (b) is by Propositions \ref{prop:sc} and \ref{prop:cauchy} and (c) is due to the smoothness of $f$.
For the last term of \eqref{eq:proof_u_n_recur},
we first show that
$ | \Sigma_n^{(1)} - \Sigma_\star^{(1)} | \le dL \sE \norm{ \vx_n - \vx^\star } $.
From the definition of $ \Sigma_n^{(1)}$ and $\Sigma_\star^{(1)}$, we have
\begin{align*}
    | \Sigma_n^{(1)} - \Sigma_\star^{(1)} |
    & = \left| \sE\, \mathrm{trace} \left( \gP_{\mA^\bot} (\xi_n \xi_n^\top - \xi^\star (\xi^\star)^\top) \gP_{\mA^\bot}   \right) \right| \\
    & = \left| \mathrm{trace} \left( \gP_{\mA^\bot} ( \sE \xi_n \xi_n^\top - \sE \xi^\star (\xi^\star)^\top ) \gP_{\mA^\bot} \right) \right| \\
    & \le d \norm{ \gP_{\mA^\bot} \sE (\Sigma (\vx_n) - \Sigma(\vx^\star) ) \gP_{\mA^\bot} } \\
    & \le d L \sE \norm{ \vx_n - \vx^\star },
\end{align*}
where the last inequality is due to the non-expansiveness of $\gP_{\mA^\bot}$ and Assumption \ref{asp:hess_lip}.
It follows that
\begin{align*}
    \Sigma_n^{(1)}
    & \le \Sigma_\star^{(1)} +| \Sigma_n^{(1)} - \Sigma_\star^{(1)} | \\
    & \le \Sigma_\star^{(1)} + dL \sE \norm{ \vx_n - \vx^\star } \\
    & \le \Sigma_\star^{(1)} + dL \sE ( \norm{ \vu_n - \vx^\star } + \norm{ \vv_n } ) \\
    & \overset{(a)}{\le} \Sigma_\star^{(1)} + dL \left[ \frac{\Sigma_\star^{(1)}}{dL} + \frac{dL}{ 2\Sigma_\star^{(1)}} \sE ( \norm{ \vu_n - \vx^\star }^2 + \norm{ \vv_n }^2 )   \right] \\
    & = 2 \Sigma_\star^{(1)} + \frac{d^2 L^2}{2 \Sigma_\star^{(1)}} \sE \norm{ \vu_n - \vx^\star }^2 + \frac{d^2 L^2}{2 \Sigma_\star^{(1)}} \norm{ \vv_n }^2, 
\end{align*}
where (a) follows from Proposition \ref{prop:cauchy}.

By substituting these inequalities, we obtain
\begin{align*}
    \sE \norm{ \vu_{n+1} - \vx^\star }^2
    & \le \left( 1 - \frac{3\mu}{2} \eta_n + 2 L^2 \eta_n^2 + \frac{d^2 L^2}{2 \Sigma_\star^{(1)}} \eta_n^2 \right) \sE \norm{ \vu_n - \vx^\star }^2 \\
    & \quad \ + \left( \frac{2 L^2}{\mu} \eta_n + 2 L^2 \eta_n^2 + \frac{d^2 L^2}{2 \Sigma_\star^{(1)}} \eta_n^2 \right) \sE \norm{\vv_n}^2 + 2 \eta_n^2 \Sigma_n^{(1)} \\
    & \overset{(a)}{\le} (1-\mu\eta_n)\sE\|\vu_n-\vx^\star\|^2+\frac{4L^2}{\mu}\eta_n\sE\|\vv_n\|^2+  2\eta_n^2\Sigma_n^{(1)},
\end{align*}
where (a) holds if $n$ is large enough.
\end{proof}

\begin{proof}[Proof of Lemma \ref{lem:des_v_n}]

From the update rule and the linearity of $\gP_{\mA}$, we have
\begin{align}
    \sE \norm{ \vv_{n+1} }^2
    & = (1 - p_n) \sE \norm{ \gP_{\mA} ( \vx_n - \eta_n \nabla f (\vx_n) \eta_n \xi_n ) }^2 \nonumber \\
    & = (1 - p_n) \sE \norm{ \vv_n - \eta_n \gP_{\mA} \nabla f(\vx_n) + \eta_n \gP_{\mA} \xi_n }^2 \nonumber \\
    & = (1 - p_n) \sE \norm{ \vv_n }^2 + (1 - p_n) \eta_n^2 \sE \norm{ \gP_{\mA} \nabla f(\vx_n) }^2 \label{eq:proof_v_n_recur} \\
    & \quad \ - 2 (1 - p_n) \eta_n \sE \inner{\vv_n}{\gP_{\mA} \nabla f(\vx_n) } \nonumber 
     + (1 - p_n) \eta_n^2 \Sigma_{n}^{(2)}, 
\end{align}
where the last equality is due to that $\{ \xi_n \}$ is a m.d.s. and  $\Sigma_n^{(2)} := \sE \norm{ \gP_{\mA} \xi_n }^2$.

For the second term of \eqref{eq:proof_v_n_recur}, we have
\begin{align*}
    \norm{ \gP_{\mA} \nabla f (\vx_n) }^2
    & = \norm{ \gP_{\mA} ( \nabla f (\vx_n) - \nabla f(\vu_n) + \nabla f (\vu_n) - \nabla f (\vx^\star) + \nabla f (\vx^\star)  ) }^2 \\
    & \overset{(a)}{\le} 3 \norm{ \gP_{\mA} ( \nabla f (\vx_n) - \nabla f(\vu_n) ) }^2 + 3 \norm{ \gP_{\mA} ( \nabla f (\vu_n) - \nabla f(\vx^\star) ) }^2 + 3 \norm{ \gP_{\mA} \nabla f (\vx^\star) }^2 \\
    & \overset{(b)}{\le} 3 L^2 \norm{ \vv_n }^2 + 3 L^2 \norm{ \vu_n - \vx^\star }^2 + 3 \norm{ \nabla f (\vx^\star) }^2,
\end{align*}
where (a) is by Proposition \ref{prop:cauchy} and (b) follows from non-expansiveness of $\gP_{\mA}$ and smoothness of $f$.

For the third term of \eqref{eq:proof_v_n_recur}, we have
\begin{align*}
    & \quad \ - \inner{ \vv_n }{ \gP_{\mA} \nabla f(\vx_n) } \\
    & \overset{(a)}{=} - \inner{ \vv_n }{ \nabla f(\vx_n) } \\
    & = - \inner{ \vv_n }{ \nabla f (\vx_n) - \nabla f (\vu_n) } - \inner{ \vv_n }{ \nabla f (\vu_n) - \nabla f (\vx^\star) } - \inner{ \vv_n }{ \nabla f (\vx^\star) } \\
    & \overset{(b)}{\le} - \mu \norm{ \vv_n }^2 + \frac{p_n}{8 \eta_n} \norm{ \vv_n }^2 + \frac{2 \eta_n }{ p_n} \norm{ \nabla f (\vu_n) - \nabla f (\vx^\star) }^2 + \frac{p_n}{8 \eta_n} \norm{ \vv_n }^2 + \frac{2 \eta_n}{p_n} \norm{ \nabla f (\vx^\star) }^2 \\
    & \overset{(c)}{\le} - \mu \norm{ \vv_n }^2 + \frac{p_n}{4 \eta_n} \norm{ \vv_n }^2 + \frac{2 L^2 \eta_n}{p_n} \norm{ \vu_n - \vx^\star }^2 + \frac{2 \eta_n}{p_n} \norm{ \nabla f (\vx^\star) }^2,
\end{align*}
where (a) follows from the orthogonality between $\gP_{\mA}$ and $\gP_{\mA^\bot}$, (b) is by Propositions \ref{prop:sc} and \ref{prop:cauchy} and (c) is due to the smoothness of $f$.
For the last term of \eqref{eq:proof_v_n_recur},
we can obtain
\begin{align*}
    \Sigma_n^{(2)}
    \le  2 \Sigma_\star^{(2)} + \frac{d^2 L^2}{2 \Sigma_\star^{(2)}} \sE \norm{ \vu_n - \vx^\star }^2 + \frac{d^2 L^2}{2 \Sigma_\star^{(2)}} \norm{ \vv_n }^2
\end{align*}
by following similar procedure in the proof of Lemma \ref{lem:des_u_n}.
By substituting these inequalities, we obtain
\begin{align*}
    \sE \norm{ \vv_{n+1} }^2
    & \le (1 - p_n) \left( 1 - 2 \mu \eta_n + 3 L^2 \eta_n^2 + \frac{p_n}{2} + \frac{d^2 L^2}{ 2 \Sigma_\star^{(2)}  } \eta_n^2 \right) \sE \norm{ \vv_n }^2 \\
    & \quad \ + \left( \frac{4 L^2 \eta_n^2}{p_n} + 3 L^2 \eta_n^2 + \frac{d^2 L^2}{2 \Sigma_\star^{(2)}} \eta_n^2  \right) \sE \norm{ \vu_n - \vx^\star }^2 \\
    & \quad \ + \left( \frac{4 \eta_n^2}{p_n} + 3 \eta_n^2 \right) \norm{ \nabla f (\vx^\star) }^2 + 2 \eta_n^2 \Sigma_n^{(2)} \\
    & \overset{(a)}{\le} \left( 1 - \frac{p_n}{2} \right) \sE \norm{\vv_n}^2 + \frac{7 L^2 \eta_n^2}{p_n} \sE \norm{\vu_n - \vx^\star}^2 + \frac{7 \eta_n^2}{p_n} \norm{ \nabla f (\vx^\star) }^2 + 2 \eta_n^2 \Sigma_n^{(2)},
\end{align*}
where (a) holds if $n$ is large enough.
\end{proof}

\section{Proof of Section \ref{sec:dynamics_approx}}
For convenience, we assume that $\eta_0 > 2 \mathbbm{1}_{ \{ \alpha=1 \} } / \mu$ holds throughout this section and will not emphasize this condition in the auxiliary lemmas.

\subsection{Proof of Case 1} \label{sec:prf_of_c1}

We can deduce the recursive relationship of $\check{\vu}_n$ by the definition of $\vu_n$ and the update rule \eqref{eq:update_rule}.
\begin{equation}\label{eq:up_check_u}
\begin{aligned}
  \check{\vu}_{n+1}&=\gP_{\mA^\bot}\frac{\vx_{n+\frac{1}{2}}-\vx^\star}{\sqrt{\eta_n}}\\
  &= \gP_{\mA^\bot}\frac{1}{\sqrt{\eta_n}}(\vx_n-\eta_n\nabla f(\vx_n)+\eta_n\xi_n - \vx^\star)\\
  &= \frac{\sqrt{\eta_{n-1}}}{\sqrt{\eta_n}}\check{\vu}_n - \sqrt{\eta_n}\gP_{\mA^\bot}\left\{(\nabla f(\vx_n)-\nabla f(\vu_n))+(\nabla f(\vu_n)-\nabla f(\vx^\star))\right\} +\sqrt{\eta_n}\gP_{\mA^{\bot}}\xi_n\\
  &=  \frac{\sqrt{\eta_{n-1}}}{\sqrt{\eta_n}}\check{\vu}_n -\sqrt{\eta_n\eta_{n-1}}\gP_{\mA^\bot}\left\{\int^1_0\nabla^2 f\left(t\vx^\star + (1-t)\vu_n\right)dt \right\}\check{\vu}_n +\gR_n^{(1)} + \sqrt{\eta_n}\gP_{\mA^{\bot}}\xi_n\\
  &= \check{\vu}_n - \eta_n \gP_{\mA^\bot}\left(\nabla^2 f(\vx^\star) - \frac{1}{2\eta_0}\mathbbm{1}_{\{\alpha = 1\}}\rmI\right)\check{\vu}_n + \gR_n^{(1)} + \gR_n^{(2)} + \gR_n^{(3)} + \sqrt{\eta_n}\xi_n^{(1)},
\end{aligned}
\end{equation}
where $\gR^{(i)}_n,~ i=1,2,3$ are higher order terms with respect to $\eta_n$ of the form:
\begin{equation}\label{eq:diff_res1-3}
\begin{aligned}
  &\gR_n^{(1)} = -\sqrt{\eta_n}\gP_{\mA^\bot}(\nabla f(\vx_n)-\nabla f(\vu_n))\\
  &\gR_n^{(2)} = -\left(1-\sqrt{\frac{\eta_{n-1}}{\eta_n}} + \frac{\eta_n}{2\eta_0}\mathbbm{1}_{\{\alpha=1\}}\right)\check{\vu}_n + (\eta_n-\sqrt{\eta_n\eta_{n-1}})\gP_{\mA^\bot}\nabla^2 f(\vx^\star)\check{\vu}_n\\
  &\gR_n^{(3)} = \sqrt{\eta_n\eta_{n-1}}\gP_{\mA^\bot}\left(\nabla^2 f(\vx^\star) - \int^1_0\nabla^2 f\left(t\vx^\star + (1-t)\vu_n\right)dt\right)\check{\vu}_n.
\end{aligned}
\end{equation}
Lemma \ref{lem:bound_res} shows that the $\frac{1}{\eta_n}\gR^{(i)}_n$ are $o(1)$ in some sense.
\begin{lemma}\label{lem:bound_res}
When Assumptions \ref{asp:smooth}, \ref{asp:str_cov} and \ref{asp:noi_lip} hold, and let $p_t=\eta_t^\beta \text{ where } \beta\in [0,1/2)$, then for any $i\in \{1,2\}$, $ \sE \|\gR_n^{(i)}\|^2= o(\eta_n^2)$. For $\gR^{(3)}_n$, we have $\sE\|\gR^{(3)}_n\|^2 = \gO(\eta_n^2)$ and $\sE \|\gR^{(3)}_n\| = o(\eta_n)$.
\end{lemma}

We then show the tightness of the rescaling sequence we built. Actually, we make use of a classical criterion (Theorem 7.3 in \cite{billingsley2008probability}) to prove this property of $\Bar{\vu}_t^{(n)}$.

\begin{proposition}\label{prop:tight_cri}
The sequence $\Bar{\vu}^{(n)}$ is tight if these two conditions hold.
\begin{enumerate}
    \item For each positive $\zeta$, there exists an $a$ and an $n_0$ such that
    \begin{equation}\label{cond:init}
        \sP(\|\check{\vu}_n\|\ge a)\le \zeta \quad \forall n\ge n_0.
    \end{equation}
    \item For any $T>0$ and any positive $\epsilon, \zeta$,
    there exists
    a $\delta$ 
    and an integer $n_0$ 
    such that:
    \begin{equation}\label{cond:eq_cont}
        \mathbb{P}\left(\sup _{s \in[t, t+\delta]}\left\|\bar{\vu}_{s}^{(n)}-\bar{\vu}_{t}^{(n)}\right\| \geq \varepsilon\right) \leq \zeta \delta ;\quad \forall t \in[0, T] \quad \forall n \geq n_{0}.
    \end{equation}
\end{enumerate}
\end{proposition}

\begin{lemma}[\textbf{Tightness of $\bar{\vu}^{(n)}$}]\label{lem:tight_u}
Suppose that Assumptions \ref{asp:smooth}, \ref{asp:str_cov} and \ref{asp:noi_lip} holds, and assume that there exists a positive number $p>2$ such that $\sup\limits_{n\ge 0}\sE\|\xi_n\|^p<\infty$. Then the sequence of random processes $\{\Bar{\vu}^{(n)}\}$ is tight under the Skorokhod topology in finite interval.
\end{lemma}

\begin{lemma}\label{lem:generator_u}
Suppose Assumptions \ref{asp:smooth}, \ref{asp:str_cov} and \ref{asp:noi_lip} holds, and assume that there exists a positive number $p>2$ such that $\sup\limits_{n\ge 0}\sE\|\xi_n\|^p<\infty$. And suppose

\begin{equation}
    \sE[\xi_t\xi_t^\top|\gF_{t}]\overset{n\to\infty}{\longrightarrow} \Sigma ~~~\text{in probability},
\end{equation}
where $\Sigma$ is a positive definite $d\times d$-matrix. Then for any $C^2$ function $g:\sR^d\to \sR$,  compactly supported with Lipschitz continuous second derivatives, we have
\begin{equation}
    \sE [g(\check{\vu}_{n+1})-g(\check{\vu}_n)|\gF_n]= \eta_n\gL g(\check{\vu}_n) + \gR^g_n,
\end{equation}
where $\frac{1}{\eta_n}\gR^g_n \to 0$ in $L_1$ and $\gL$ is the infinitesimal generator defined by
\begin{equation}\label{eq:def_L}
    \forall \phi \in \mathcal{C}^{2}\left(\mathbb{R}^{p}\right) \quad \mathcal{L} \phi(\vx)=\left\langle-\gP_{\mA^\bot}\left(\nabla^2 f(\vx^\star) - \frac{1}{2\eta_0}\mathbbm{1}_{\{\alpha= 1\}}\rmI_d\right)\gP_{\mA^\bot}\vx, \nabla \phi\right\rangle+\frac{1}{2}\tr\left( \nabla^{2} \phi(\vx) \Sigma\right).
\end{equation}
\end{lemma}

All the proofs of above lemmas can be found in Appendix~\ref{sec:aux_prf_1}

\begin{proof}[
{Proof of Theorem \ref{thm:diff_approx}}]
The proof of this main results is divided into two steps. At first, we prove that every weak limit of sequence of random process $\{\Bar{\vu}^{(n)}\}$ is a solution of the martingale problem $(\gL,\gC)$, where $\gC$ denotes the class of $\gC^2$-functions with compact support and Lipschitz continuous second derivatives. $\gL$ is defined by \eqref{eq:def_L}. Then, from the property of Langevin dynamics, we know that \eqref{eq:diff_approx} converges to a unique invariant distribution $\pi^\star$. Further, by proving that the limit of every weakly converged subsequence equals to $\pi^\star$ and combining it with the Prokhorov's theorem, we conclude $\check{\vu}_n$ converges to $\pi^\star$ weakly. Finally,  repeat the first step of this proof, and we have $\{\Bar{\vu}^{(n)}_t\}$ converges to the solution of equation \eqref{eq:diff_approx} with initial distribution $\pi^\star$.

\paragraph{Step 1}  Let $g$ belong to $\gC$ and let $\gF^{(n)}_t$ denote the natural filtration of $\Bar{\vu}^{(n)}_t$. We aim to derive the following equation, which can guarantee that every sub-limit of $\{\Bar{\vu}^{(n)}_t\}$ is a weak solution of the martingale problem $(\gL,\gC)$.

\begin{equation}
\begin{aligned}
    \forall t \ge 0, \quad g\left(\bar{\vu}_{t}^{(n)}\right)-g\left(\bar{\vu}_{0}^{(n)}\right)-\int_{0}^{t} \mathcal{L} g\left(\bar{\vu}_{s}^{(n)}\right) d s=\mathcal{M}_{t}^{(n, g)}+\mathcal{R}_{t}^{(n, g)},
\end{aligned}
\end{equation}
where $\gM_t^{(n,g)}$ is a $\gF_t^{(n)}$-martingale and $\gR_t^{(n,g)}$ converges to zero in $L_1$.

In fact, we set
\begin{equation}
\begin{aligned}
     &\mathcal{M}_{t}^{(n, g)}&&=\sum_{k=n+1}^{N(n, t, \eta)-1} \left\{g\left(\check{\vu}_{k+1}\right)-g\left(\check{\vu}_{k}\right)-\mathbb{E}\left[g\left(\check{\vu}_{k+1}\right)-g\left(\check{\vu}_{k}\right) \mid \mathcal{F}_{k}\right]\right\},\\
     &\mathcal{R}_{t}^{(n, g)}&&=g\left(\bar{\vu}_{t}^{(n)}\right)-g\left(\bar{\vu}_{\underline{t}_{n}}^{(n)}\right)-\int_{\underline{t}_{n}}^{t} \mathcal{L} g\left(\bar{\vu}_{s}^{(n)}\right) d s\\ & && \quad +\int_{0}^{\underline{t}_{n}}\left(\mathcal{L} g\left(\bar{\vu}_{\underline{s}_{n}}^{(n)}\right)-\mathcal{L} g\left(\bar{\vu}_{s}^{(n)}\right)\right) d s+\sum_{k=n}^{N(n, t, \eta)-1} \gR_{k}^{g}.
\end{aligned}
\end{equation}

From the definition of $\Bar{\vu}_t^{(n)}$ \eqref{eq:itpl_u}, we can get
\begin{equation}
\begin{aligned}
     \bar{\vu}_t^{(n)}-\bar{\vu}_{\underline{t}_n(\eta)}^{(n)}= (t-\underline{t}_n(\eta))\vb_{N(n,t,\eta)} + \sqrt{t-\underline{t}_n(\eta)}\xi_{N(n,t,\eta)},
\end{aligned}
\end{equation}
which satisfies
\begin{equation}
\begin{aligned}
    \sE\left\|\Bar{\vu}_t^{(n)}-\bar{\vu}_{\underline{t}_n(\eta)}^{(n)}\right\|&\le (t-\underline{t}_n(\eta))\sE\|\vb_{N(n,t,\eta)}\| + \sqrt{t-\underline{t}_n(\eta)}\sE\|\xi_{N(n,t,\eta)}\|\\
    & \precsim \sqrt{\eta_{N(n,t,\eta)}}.
\end{aligned}
\end{equation}

Plug the above bound into the residual $\gR_t^{(n,g)}$, and note the Lipschitz continuity and boundedness of $g,~ \nabla g\text{ and } \nabla^2 g$,
\begin{equation}
\begin{aligned}
    &\sE\left|g(\bar{\vu}_t^{(n)})-g(\bar{\vu}_{\underline{t}_n(\eta)}^{(n)})\right| \precsim \sE\left\|\Bar{\vu}_t^{(n)}-\bar{\vu}_{\underline{t}_n(\eta)}^{(n)}\right\|
    =o(1),\\
    &\sE\left|\int_{\underline{t}_n(\eta)}^t\gL g(\bar{\vu}_s^{(n)})ds\right|\precsim \int_{\underline{t}_n(\eta)}^t \gC \precsim \eta_{N(n,t,\eta)}=o(1),\\
    &\sE\left|\int_0^{\underline{t}_n(\eta)}\gL g(\bar{\vu}_{\underline{s}_n(\eta)}^{(n)})-\gL g(\bar{\vu}_s^{(n)})ds\right|\precsim \sE\int_0^{\underline{t}_n(\eta)}\left\|\bar{\vu}_{\underline{s}_n(\eta)}^{(n)}-\bar{\vu}_s^{(n)}\right\|ds\\ 
    & \precsim \int_0^{\underline{t}_n(\eta)}\sqrt{\eta_{N(n,s,\eta)}}ds \precsim \sqrt{\eta_n} =o(1).
\end{aligned}
\end{equation}

Further, attributed to Lemma \ref{lem:generator_u},
\begin{equation}
\begin{aligned}
    \sE\left|\sum_{k=n}^{N(n, t, \eta)-1} \gR_{k}^{g}\right|&\le \ssum{k}{n}{N(n,t,\eta)}\eta_k\sE\left|\frac{\gR_k^g}{\eta_k}\right|\\
    &\le\sup\limits_{k\ge n}\sE\left|\frac{\gR_k^g}{\eta_k}\right|\ssum{k}{n}{N(n,t,\eta)}\eta_k \precsim o(1)t = o(1).
\end{aligned}
\end{equation}

So far we can say that $\sE|\gR_t^{(n,g)}|\to 0 ,~~ n\to \infty$. Finally, by using the Proposition~3.1 in \cite{kurtz1975semigroups}, we can get that, if a subsequence of index $\{n_k\}$ satisfies that $\{\bar{\vu}_0^{(n_k)}\}_{k=0}^\infty = \{\check{\vu}_{n_k}\}_{k=0}^\infty$ has a weak limitation $\vu^\prime_0$, then for any $t>0$, $\{\bar{\vu}_t^{n_k}\}$ converges weakly to $\vu^\prime_t$, where $\vu^\prime_t$ is the solution of the diffusion process with infinitesimal generator $\gL$ and initial distribution $\vu^\prime_0$.

\paragraph{Step 2} Now we suppose that there exists a weakly convergent subsequence $\{\check{\vu}_{n_k}\}_{k=1}^\infty$ with limit distribution $\Tilde{\pi}$. We should introduce some new notations. For $n\in \sN$ and $t\ge0$, we define $M(n,t,\eta) = \min\left\{m\ge 0; \ssum{i}{m}{n-1}\eta_i\le t\right\}$ and $\tilde{t}_n(\eta) = \Gamma_n - \Gamma_{M(n,t,\eta)}$. For the properties of step size sequence $\eta_n$, we can affirm $t-\tilde{t}_n(\eta)\to 0$ when $n\to \infty$. 

By leveraging the Prokhorov's theorem, for any $T>0$, we know that $\left\{\bar{\vu}^{(M(n_k,T,\eta))}_t\right\}$ has a weakly convergent subsequence. Without loss of generality, we can assume that the subsequence $\left\{\bar{\vu}^{(M(n_k,T,\eta))}_t\right\}$ itself converges weakly to a solution $\bar{\vu}^{\tilde{\nu}^{(T)}}_t$ of the SDE \eqref{eq:diff_approx} with initial distribution $\tilde{\nu}^{(T)}$. Owing to the tightness of the whole sequence $\{\bar{\vu}^{n}\}$, for any given $\epsilon>0$, there is a compact set $K_\epsilon\subset \sR^d$ only depending on $\epsilon$ such that $\sup\limits_{n}\sP(\check{\vu}_n \in K_\epsilon^c) \le \epsilon$. This makes us find the following holds: $\tilde{\nu}^{(T)}(K_\epsilon)\ge 1-\epsilon$ for any $T>0$.

By the geometrical ergodicity of the dynamics \eqref{eq:diff_approx}, we can choose $T_\epsilon$ such that 
\begin{equation}\label{eq:erg_diff}
    \sup\limits_{\vx\in K_\epsilon}\sup\limits_{g\in\gC}\left|\gP^{T_\epsilon}g(\vx)-\langle\pi^\star, g\rangle\right| \le \epsilon,
\end{equation}
where $\gP$ represents the Markov semigroup induced by the SDE \eqref{eq:diff_approx}. In virtue of the approximation of $\widetilde{(T_\epsilon)}_n(\eta) ( = \Gamma_n - \Gamma_{M(n,T_\epsilon, \eta)})$ to $T_\epsilon$ and the tightness of the sequence $\bar{\vu}^{(n)}$, we are able to deduce that $\check{\vu}_{n_k} (= \bar{\vu}^{M(n_k,T_\epsilon,\eta)}_{\widetilde{(T_\epsilon)}_n(\eta)})$ converges weakly to the limit random variable of the sequence $\bar{\vu}^{M(n_k,T_\epsilon,\eta)}_{T_\epsilon}$, i.e., $\bar{\vu}^{\tilde{\nu}^{(T_\epsilon)}}_{T_\epsilon}$. On the other hand, by assumption, $\check{\vu}_{n_k}$ converges weakly to $\tilde{\pi}$. Thus $\bar{\vu}^{\tilde{\nu}^{(T_\epsilon)}}_{T_\epsilon}\sim \tilde{\pi}$.

Given any $g\in\gC$, it is not difficult to derive the following bounds
\begin{equation}\label{eq:w_conv_bound}
\begin{aligned}
    \left|\langle\tilde{\pi}, g\rangle - \langle \pi^\star, g\rangle\right|&=\left|\sE g\left(\bar{\vu}^{(T_\epsilon)}_{T_\epsilon}\right)-\sE_{\pi^\star}g\right|
    = \left|\int \left(\gP^{T_\epsilon} g(\vx) - \sE_{\pi^\star}g\right)d\tilde{\nu}^{(T_\epsilon)}(\vx)\right|\\
    &\le \int \left|\gP^{T_\epsilon} g(\vx) - \sE_{\pi^\star}g\right|d\tilde{\nu}^{(T_\epsilon)}(\vx)\\ 
    &= \int_{K_\epsilon} \left|\gP^{T_\epsilon} g(\vx) - \sE_{\pi^\star}g\right|d\tilde{\nu}^{(T_\epsilon)}(\vx) + \int_{K_\epsilon^c} \left|\gP^{T_\epsilon} g(\vx) - \sE_{\pi^\star}g\right|d\tilde{\nu}(\vx)\\
    &\le \int_{K_\epsilon} \left|\gP^{T_\epsilon} g(\vx) - \sE_{\pi^\star}g\right|d\tilde{\nu}^{(T_\epsilon)}(\vx) + 2\|g\|_\infty \tilde{\nu}^{(T_\epsilon)}(K_\epsilon^c)\\ 
    &\overset{(a)}{\le} \epsilon + 2\|g\|_\infty \epsilon,
\end{aligned}
\end{equation}
where (a) holds for sake of $\tilde{\nu}^{(T_\epsilon)}(K_\epsilon)\ge 1-\epsilon$ and \eqref{eq:erg_diff}. We obtain $\tilde{\pi} = \pi^\star$ by taking $\epsilon \to 0$. Finally, owing to the Prokhorov's theorem, we have proved that $\check{\vu}_{n}$ converges weakly to $\pi^\star$. Further, the sequence of random process $\bar{\vu}^{(n)}_t$ converges weakly to the dynamics \eqref{eq:diff_approx} with stationary distribution $\pi^\star$ as initialization.

\end{proof}

\subsubsection{Proof of Lemma~\ref{lem:bound_res}, \ref{lem:tight_u}
and \ref{lem:generator_u}}\label{sec:aux_prf_1}

\begin{proof}[
{Proof of Lemma \ref{lem:bound_res}}]
By Theorem~\ref{thm:converge}, we have
\begin{equation}
\begin{aligned}
   &\sE\|\gR_n^{(1)}\|^2 &&\le\eta_n\sE\|\nabla f(\vx_n)-\nabla f(\vu_n)\|^2\le \eta_nL^2\sE\|\vv_n\|^2\precsim L^2\eta_n\times \eta_n^{2-2\beta}=o(\eta_n^2),\\
   &\sE\norm{\gR_n^{(3)}}^2 &&\precsim \eta_n^2\sE\left\|\left\{\int^1_0\nabla^2 f\left(t\vx^\star + (1-t)\vu_n\right)dt - \nabla^2 f(\vx^\star)\right\}\check{\vu}_n\right\|^2\\
   & && \precsim \eta_n^2\sE\|\check{\vu}_n\|^2 = \gO (\eta_n^2),\\
   & \sE\norm{\gR_n^{(3)}} && \precsim \eta_n \sE\norm{\int^1_0\nabla^2 f\left(t\vx^\star + (1-t)\vu_n\right)dt - \nabla^2 f(\vx^\star)}\cdot \norm{\check{\vu}_n}\\
   & && \precsim \eta_n\sE \norm{\check{\vu}_n}\int^1_0 \norm{\nabla^2 f(t\vx^\star - (1-t)\vu_n) - \nabla^2 f(\vx^\star)}dt\\
   & && \precsim \eta_n \sE \norm{\check{\vu}_n}\int^1_0 (1-t)\norm{\vu_n -\vx^\star}dt \precsim \eta_n^{3/2}\sE \norm{\check{\vu}_n}^2 = o(\eta_n).
\end{aligned}
\end{equation}
As for $\sE \norm{ \gR_n^{(2)} }^2$, when $\alpha < 1$, we have
\begin{align*}
    &\sE\|\gR_n^{(2)}\|^2 &&
   \precsim
   \left(1-\sqrt{\frac{\eta_{n-1}}{\eta_n}}\right)^2 + (\eta_n-\sqrt{\eta_n\eta_{n-1}})^2
   =\left(\frac{1}{\eta_n}+\eta_n\right)(\sqrt{\eta_n}-\sqrt{\eta_{n-1}})^2\\
   & &&\precsim \frac{(\eta_n-\eta_{n-1})^2}{\eta_n^2}=[1-(1+o(\eta_n))]^2=o(\eta_n^2);
\end{align*}
when $\alpha = 1$, we have
\begin{align*}
  \sE\norm{\gR^{(2)}_n}^2 \precsim \left( 1 - \sqrt{1 + \frac{1}{n-1}} + \frac{1}{2n}\right)^2 + \frac{1}{n} \left( \frac{1}{\sqrt{n}} - \frac{1}{\sqrt{n-1}} \right)^2 = o \left( \frac{1}{n^2} \right).
\end{align*}
\end{proof}

\begin{proof}[
{Proof of Lemma \ref{lem:tight_u}}]

From the construction of $\Bar{\vu}_t^{(n)}$ we know it is a continuous process. What remains to do is to verify two conditions \eqref{cond:init} and \eqref{cond:eq_cont}. 

For the first condition about initialization of the process, it is easy to check by the convergence rate result for $\vu_n-\vx^\star$.

For the condition \eqref{cond:eq_cont}, note that we have
\begin{equation}
\begin{aligned}
   \Bar{\vu}_s^{(n)}-\Bar{\vu}_t^{(n)}&= \left\{\ssum{k}{N(n,t,\eta)}{N(n,s,\eta)-1}\eta_k\vb_k - \left[(t-\underline{t}_n(\eta))\vb_{N(n,t,\eta)}-(s-\underline{s}_n(\eta)\vb_{N(n,s,\eta)})\right]\right\}\\
   & \quad + \left\{\ssum{k}{N(n,t,\eta)}{N(n,s,\eta)-1}\sqrt{\eta_k}\xi_k - \sqrt{t-\underline{t}_n(\eta)}\xi_{N(n,t,\eta)}+\sqrt{s-\underline{s}_n(\eta)}\xi_{N(n,s,\eta)}\right\}\\
   &=: \bfB + \Xii.
\end{aligned}
\end{equation}

From the discussion following Lemma \ref{lem:bound_res}, we can see that $\sE\|\vb_n\|^2$ is uniformly bounded. So
\begin{equation}\label{eq:tight_B}
\begin{aligned}
   \sP\left( \sup\limits_{s\in [t, t+\delta]} \|\bfB\|\ge \frac{\epsilon}{2}\right)&\le \sP\left( \sup\limits_{s\in [t, t+\delta]} \ssum{k}{N(n,t,\eta)}{N(n,s,\eta)}\eta_k\|\vb_k\|\ge \frac{\epsilon}{2}\right)\\
   &\le \frac{4}{\epsilon^2} \sE \sup\limits_{s\in [t, t+\delta]} \left(  \ssum{k}{N(n,t,\eta)}{N(n,s,\eta)}\eta_k\|\vb_k\|\right)^2\\
   &\le \frac{4}{\epsilon^2}\left( \sup\limits_{s\in [t, t+\delta]} \ssum{k}{N(n,t,\eta)}{N(n,s,\eta)}\eta_k\right) \sE\left( \sup\limits_{s\in [t, t+\delta]} \ssum{k}{N(n,t,\eta)}{N(n,s,\eta)}\eta_k\|\vb_k\|^2\right)\\
   &\le \frac{4}{\epsilon^2} \sup\limits_k \sE\|\vb_k\|^2\left( \sup\limits_{s\in [t, t+\delta]} \ssum{k}{N(n,t,\eta)}{N(n,s,\eta)}\eta_k\right)^2\\
   &\le \gC\frac{(\delta + 2\eta_n)^2}{\epsilon^2}\overset{(a)}{\le} \frac{\zeta(\delta+ 2\eta_n)}{4} \overset{(b)}{\le} \frac{\zeta\delta}{2},
\end{aligned}
\end{equation}
where (a) and (b) hold when we take $\delta + 2 \eta_{n_0} < \frac{\zeta\epsilon^2}{4\gC}$, $\eta_{n_0} < \delta / 2$ and $n \ge n_0$.

On the other hand, thanks to the property of monotone interpolation, we have
\begin{equation}
    \|\Xii\|\le \max\limits_{j\in \{0,1\}}\left\|\ssum{k}{N(n,t,\eta)}{N(n,s,\eta)-j}\sqrt{\eta_k}\xi_k-\sqrt{t-\underline{t}_n(\eta)}\xi_{N(n,t,\eta)}\right\|.
\end{equation}

By leveraging the Doob's inequality and the assumption of bounded p-th moment of $\xi_k$, we can get
\begin{equation}\label{eq:tight_Xi}
\begin{aligned}
   \sP\left(\sup\limits_{s\in[t,t+\delta]}\|\Xii\|\ge \frac{\epsilon}{2}\right)&\le \sP\left(\max\limits_{j\le N(n,t+\delta, \eta)}\left\|\ssum{k}{N(n,t,\eta)}{j}\sqrt{\eta_k}\xi_k-\sqrt{t-\underline{t}_n(\eta)}\xi_{N(n,t,\eta)}\right\|\ge \frac{\epsilon}{2}\right)\\
   &\le \frac{2^p}{\epsilon^p}\sE\left\|\ssum{k}{N(n,t,\eta)}{N(n,t+\delta, \eta)}\sqrt{\eta_k}\xi_k-\sqrt{t-\underline{t}_n(\eta)}\xi_{N(n,t,\eta)}\right\|^p\\
   &\overset{(a)}{\le}\frac{\gC_p2^p}{\epsilon^p}\ssum{k}{N(n,t,\eta)}{N(n,t+\delta,\eta)}\eta_k^{p/2}\sE\|\xi_k\|^p\\
   & \le\frac{\gC}{\epsilon^p}\eta_{n}^{\frac{p}{2}-1}\ssum{k}{N(n,t,\eta)}{N(n,s,\eta)}\eta_k\\ 
   & \overset{(b)}{\le} \frac{\gC}{\epsilon^p}\eta_{n}^{\frac{p}{2}-1} (\delta+ 2\eta_n) \overset{(c)}{\le} \frac{\zeta (\delta+ 2\eta_n)}{4}\le \frac{\zeta\delta}{2},
\end{aligned}
\end{equation}
where (a) holds by the Burkholder's inequality and (b), (c) hold when we choose $\eta_{n_0}^{\frac{p}{2}-1}\le \frac{\epsilon^p\zeta}{4\gC}$, $\eta_{n_0}\le \delta / 2$ and $n \ge n_0$.

Combining \eqref{eq:tight_B} and \eqref{eq:tight_Xi}, finally we can derive that

\begin{equation}
\begin{aligned}
   \mathbb{P}\left(\sup _{s \in[t, t+\delta]}\left\|\bar{\vu}_{s}^{(n)}-\bar{\vu}_{t}^{(n)}\right\| \geq \varepsilon\right)&\le \sP\left(\sup\limits_{s\in [t, t+\delta]}\|\bfB\|\ge \frac{\epsilon}{2}\right) + \sP\left(\sup\limits_{s\in[t,t+\delta]}\|\Xii\|\ge \frac{\epsilon}{2}\right)\\
   &\le \frac{\zeta\delta}{2}+\frac{\zeta\delta}{2}=\zeta\delta.
\end{aligned}
\end{equation}
So far, we conclude the proof of Lemma \ref{lem:tight_u}.
\end{proof}

\begin{proof}[
{Proof of Lemma \ref{lem:generator_u}}]

$\gC$ will represent a universal constant whose value may change from line to line, for the
sake of convenience. We use a Taylor expansion between $\vu_n$ and $\vu_{n+1}$

\begin{equation}
\begin{aligned}
   g(\check{\vu}_{n+1})-g(\check{\vu}_n) &= \langle\nabla g(\check{\vu}_n), \check{\vu}_{n+1}-\check{\vu}_n\rangle + \frac{1}{2}(\check{\vu}_{n+1}-\check{\vu}_n)^\top \nabla^2 g(\check{\vu}_{n}) (\check{\vu}_{n+1}-\check{\vu}_n)\\ 
   & \quad + \underbrace{\frac{1}{2}(\check{\vu}_{n+1}-\check{\vu}_n)^\top \left(\nabla^2 g(\lam_n) - \nabla^2 g(\check{\vu}_{n})\right) (\check{\vu}_{n+1}-\check{\vu}_n)}_{\gR_n^{(4)}}.
\end{aligned}
\end{equation}

Since $\nabla^2 g$ is Lipschitz continuous and compactly supported, $\nabla^2 g$ is also $\epsilon$-H\"older continuous for all $\epsilon\in (0,1]$. Then combining the equation \eqref{eq:up_check_u} and Lemma \ref{lem:bound_res}, we can control the order of $\gR_n^{(4)}$.

\begin{align*}
    \sE\|\gR_n^{(4)}\|&\precsim \sE\|\check{\vu}_{n+1}-\check{\vu}_n\|^{2+\epsilon}\\
    &\le \sE\left\|\eta_n \gP_{\mA^\bot}\left(\nabla^2 f(\vx^\star) - \frac{1}{2\eta_0}\mathbbm{1}_{\{\alpha= 1\}}\rmI_d\right)\gP_{\mA^\bot}\check{\vu}_n + \gR_n^{(1)} + \gR_n^{(2)} + \gR_n^{(3)} + \sqrt{\eta_n}\xi_n\right\|^{2+\epsilon}\\
    &\precsim \eta_n^{1+\frac{\epsilon}{2}}.
\end{align*}

So we deduce $\frac{1}{\eta_n}\gR_n^{(4)}\to 0$ in $L_1$. Further, we make use of the update formula \eqref{eq:up_check_u} again

\begin{equation}
\begin{aligned}
   &\sE[\langle\nabla g(\check{\vu}_n), \check{\vu}_{n+1}-\check{\vu}_n\rangle|\gF_n]\\ 
   =&\sE\left[\langle\nabla g(\check{\vu}_n), -\eta_n \gP_{\mA^\bot}\left(\nabla^2 f(\vx^\star) - \frac{1}{2\eta_0}\mathbbm{1}_{\{\alpha= 1\}}\rmI_d\right)\gP_{\mA^\bot}\check{\vu}_n + \gR_n^{(1)} + \gR_n^{(2)} + \gR_n^{(3)} + \sqrt{\eta_n}\xi_n\rangle|\gF_n\right]\\
   =& -\eta_n\sE\left[\langle \nabla g(\check{\vu}_n), \gP_{\mA^\bot}\left(\nabla^2 f(\vx^\star) - \frac{1}{2\eta_0}\mathbbm{1}_{\{\alpha= 1\}}\rmI_d\right)\gP_{\mA^\bot}\check{\vu}_n \rangle|\gF_n\right] + \ssum{i}{1}{3}\sE[\langle \nabla g(\check{\vu}_n), \gR_n^{(i)}\rangle|\gF_n].
\end{aligned}    
\end{equation}

Note by  Lemma \ref{lem:bound_res}, we have

\begin{equation}
\begin{aligned}
    &\sE\left|\sE\left[\left.\left\langle \nabla g(\check{\vu}_n), \frac{\gR_n^{(i)}}{\eta_n}\right\rangle\right|\gF_n\right]\right|\le \sE\left|\left\langle \nabla g(\check{\vu}_n), \frac{\gR_n^{(i)}}{\eta_n}\right\rangle\right|\\ 
    \precsim &\sE\left\|\frac{\gR_n^{(i)}}{\eta_n}\right\|\le \left(\sE\left\|\frac{\gR_n^{(i)}}{\eta_n}\right\|^2\right)^{\frac{1}{2}}=o(1).
\end{aligned}
\end{equation}

And at last, 

\begin{equation}\label{eq:quad_term}
\begin{aligned}
    &\frac{1}{2}\sE\left[(\check{\vu}_{n+1}-\check{\vu}_n)^\top \nabla^2 g(\check{\vu}_{n}) (\check{\vu}_{n+1}-\check{\vu}_n)|\gF_n\right]\\
    =& \frac{\eta_n}{2}\sE\left[\xi_n^\top\nabla^2 g(\check{\vu}_n) \xi_n|\gF_n\right]
    + \frac{\eta_n^{\frac{3}{2}}}{2}\sE\left\langle\vb_n, \nabla^2 g(\check{\vu}_n)\xi_n\right\rangle\\
    & + \frac{\eta_n^2}{2}\sE\langle\vb_n,\nabla^2 g(\check{\vu}_n)\vb_n\rangle.
\end{aligned}
\end{equation}

Cause $g$ is compactly supported, the norm of $\nabla^2 g$ is bounded. And by Lemma \ref{lem:bound_res}, we can deduce $\sE\|\vb_n\|^2$ is uniformly bounded. Therefore, the last two terms of \eqref{eq:quad_term} are $o(\eta_n)$. Combining the above analysis, we have

\begin{equation}
\begin{aligned}
    &\sE[g(\check{\vu}_{n+1})-g(\check{\vu}_n)|\gF_n]\\
    =& -\eta_n \left\langle \nabla g(\check{\vu}_n), \gP_{\mA^\bot}\left(\nabla^2 f(\vx^\star) - \frac{1}{2\eta_0}\mathbbm{1}_{\{\alpha= 1\}}\rmI_d\right)\gP_{\mA^\bot}\check{\vu}_n \right\rangle + \frac{\eta_n}{2}\left\langle\nabla^2 g(\check{\vu}_n), \Sigma\right\rangle + \gR_n^g
\end{aligned}
\end{equation}

with $\sE\|\gR_n^g\|=o(\eta_n)$.

\end{proof}

\subsection{Proof of Case 2}\label{sec:prf_of_c2}
 We first complete the formulation of the recursive relation for $\vv_n$ that was omitted from the main text
\begin{equation}\label{eq:updt_ck_v}
\begin{aligned}
    \check{\vv}_{(n+1)-} &= \eta^{\beta-1}_n\gP_{\mA}(\vx_n-\eta_n\nabla f(\vx_n)+\eta_n\xi_n)\\ 
    &= \left(\frac{\eta_n}{\eta_{n-1}}\right)^{\beta-1}\check{\vv}_n-\eta_n^\beta\gP_{\mA}\nabla f(\vx_n)+\eta_n^\beta \gP_\mA\xi_n\\ 
    &= \left(\frac{\eta_n}{\eta_{n-1}}\right)^{\beta-1}\check{\vv}_n- \eta_n^\beta \nabla f(\vx^\star) - \eta_n^\beta\gP_{\mA}\left(\nabla f(\vx_n)-\nabla f(\vu_n)\right)\\
    &~~ - \eta_n^\beta \gP_\mA (\nabla f(\vu_n)-\nabla f(\vx^\star))+\eta_n^\beta \gP_\mA\xi_n\\ 
    & = \check{\vv}_n - \eta_n^\beta\nabla f(\vx^\star) - \gS_n^{(1)} - \gS_n^{(2)} - \gS_n^{(3)} + \eta_n^\beta \xi_n^{(2)}\\ 
    & = : \check{\vv}_n - \eta_n^\beta\vd_n + \eta_n^\beta \xi_n^{(2)},
\end{aligned}
\end{equation}
where $\vd_n = \nabla f(\vx^\star) + \frac{1}{\eta_n^\beta}\gS_n^{(1)} + \frac{1}{\eta_n^\beta}\gS_n^{(2)} + \frac{1}{\eta_n^\beta}\gS_n^{(3)}$ with higher order terms
\begin{equation}\label{eq:updt_ck_v_res}
\begin{aligned}
    &\gS_n^{(1)} = \left(1-\frac{\eta_n^{\beta - 1} }{\eta_{n-1}^{\beta - 1} }\right)\check{\vv}_n,\\
    &\gS_n^{(2)} = \eta_n^\beta\gP_{\mA}(\nabla f(\vx_n)-\nabla f(\vu_n)),\\
    &\gS_n^{(3)} = \eta_n^\beta\gP_{\mA}(\nabla f(\vu_n) - \nabla f(\vx^\star)).
\end{aligned}
\end{equation}
We can see from Theorem \ref{thm:converge} that both $\frac{\gS_n^{(2)}}{\eta_n^\beta}$ and $\frac{\gS_n^{(3)}}{\eta_n^\beta}$ are of order $\gO(\eta_n^{1-\beta})$ in $L_1$. Moreover, owing to the slow diminishing property of step size $\{\eta_n\}$, the following bound holds
\begin{equation}
\begin{aligned}
    1-\frac{\eta_n^{\beta - 1} }{\eta_{n-1}^{\beta - 1} }& =1-\left(1+\frac{\eta_n-\eta_{n-1}}{\eta_{n-1}}\right)^{\beta - 1} \\ 
    &= 1 - \left(1 + \gO(\eta_n)\right)^{\beta - 1} = 1- (1+ (1 - \beta) \gO(\eta_n)) = \gO(\eta_n).
\end{aligned}
\end{equation}
So 
\begin{equation}\label{eq:T1_bound}
\begin{aligned}
    \frac{1}{\eta_n^\beta}\sE\left|\gS_n^{(1)}\right| \precsim \frac{\gO(\eta_n)}{\eta_n^\beta} = \gO(\eta_n^{1-\beta}).
\end{aligned}
\end{equation}
As in the derivation process of the first case, we first need to focus our attention on the discussion of the tightness of the sequence of random processes $ \{\bar{\vv}^{(n)}_t\}_{n=1}^\infty$. While the discontinuity of these processes constructed by \eqref{eq:itpl_v} prevents Proposition \ref{prop:tight_cri} from being used to verify the tightness of $\{\bar{\vv}^{(n)}_t\}$. Hence, we will leverage the following more general criterion for tightness proposed in \cite{kurtz1975semigroups}.

\begin{proposition}\label{prop:tight_cri_gen}
Let $\{\bfx_n(t)\}$ be a sequence of $\sR^d$-valued processes whose sample paths are c\`adl\`ag. Let 
\begin{equation}\label{eq:part_var}
    \omega(\bfx_n,\delta, T) = \inf\limits_{\{t_i\}}\max\limits_i \sup\limits_{t_{i-1}\le s<t<t_i}\|\bfx_t-\bfx_s\|.
\end{equation}
Where $\{t_i\}$ ranges over all finite partitions of the form $0=t_0 \le t_1< t_2< \cdots < t_{r-1}< T\le t_{r}$ with $\min\limits_{1\le i\le r}(t_i-t_{i-1})\ge \delta$. Then the sequence of processes $\{\bfx_n(t)\}$ is tight if and only if,
\begin{enumerate}
    \item for every $T>0$ and $\zeta>0$, there is a compact set $K$ such that
    \begin{equation}
        \liminf\limits_{n\to \infty}\sP(\bfx_n(t)\in K; ~~\forall t\in [0,T])> 1-\zeta ;
    \end{equation}
    \item for every $\epsilon,\zeta >0$, and $T>0$, there is a $\delta > 0$ such that
    \begin{equation}
        \limsup\limits_{n\to \infty}\sP(\omega(\bfx_n, \delta, T)\ge \epsilon) < \zeta .
    \end{equation}
\end{enumerate}
\end{proposition}

Denote $\gI_n(T) = \left\{N(n,t,\eta^\beta): t\in [0,T]\right\}\subset \sN$ and $\gL_n(T)=\left\{\Gamma_{n+k}-\Gamma_n: ~k\in\gI_n(T)\right\}$. From the update rule of parametric sequence $\{\vx_n\}$ and the construction of the rescaling process $\{\bar{\vv}^{(n)}_t\}$, an intuitive fact is that the discontinuous points of $\bar{\vv}^{(n)}_t$ in the interval $[0,T]$ belong to $\gL_n(T)$. The following lemma displays the lower bound of the distance between two adjoint discontinuous time point of $\bar{\vv}^{(n)}_{t}$,  which supports the proof of tightness.
\begin{lemma}\label{lem:sparse_jump}
For the sequence of c\`adl\`ag processes $\{\bar{\vv}^{(n)}\}$, consider the time point set $\gJ_n(T)$ such that
\begin{equation}
    \gJ_n(T) = \left\{N(n,t,\eta^\beta): t\in [0,T] \text{ and }\bar{\vv}^{(n)}_{\underline{t}_n(\eta^\beta)}\neq \bar{\vv}^{(n)}_{\underline{t}_n(\eta^\beta)-}\right\}
\end{equation}
and let 
\begin{equation}
    \Delta(\gJ_n(T)) = \min\left\{|\Gamma_{n+k}(\eta^\beta)-\Gamma_{n+l}(\eta^\beta)|:k,l\in \gJ_n(T)\text{ and }k\neq l\right\}.
\end{equation}
Then there is a universal constant $\gC$ and an $n_0\in \sN$ such that for any $\delta > 0$
\begin{equation}
    \sP(\Delta(\gJ_n(T)) < \delta) \le \gC \delta \quad \forall n\ge n_0.
\end{equation}
\end{lemma}

\begin{proof}[
{Proof of Lemma \ref{lem:sparse_jump}}]

By the sub-additivity of probability and the jump scheme proposed in the Algorithm \ref{alg:lp_proj_sa}, we have
\begin{equation}
\begin{aligned}
    &\sP(\Delta(\gJ_n(T)) < \delta)\\
    \le& \sum\limits_{k\in \gI_n(T)}\sP\left\{\exists ~ (k\le)~ l \in \gJ_n(T) \text{ s.t. } |\Gamma_{n+l}(\eta^\beta)-\Gamma_{n+k}(\eta^\beta)|<\delta; ~ k \in \gJ_n(T)\right\}\\
    \le& \sum\limits_{k\in \gI_n(T)}\sum\limits_{0< \Gamma_{n+l}(\eta^\beta)-\Gamma_{n+k}(\eta^\beta)< \delta} \sP\left\{l\in\gJ_n(T);k\in\gJ_n(T)\right\}\\ 
    \le& \sum\limits_{k\in\gI_n(T)}p_{n+k}\left(\sum\limits_{0< \Gamma_{n+l}(\eta^\beta)-\Gamma_{n+k}(\eta^\beta)< \delta}p_{n+l}\right)\\ 
    \le& \gamma^2 \sum\limits_{k\in\gI_n(T)}\eta^\beta_{n+k}\left(\sum\limits_{0< \Gamma_{n+l}(\eta^\beta)-\Gamma_{n+k}(\eta^\beta)< \delta}\eta^\beta_{n+l}\right)\\ 
    \overset{(a)}{\le}& \gamma^2(\delta+ \eta_n^\beta) \sum\limits_{k\in\gI_n(T)}\eta^\beta_{n+k}\le 2\gamma^2T\delta,
\end{aligned}
\end{equation}
where (a) holds when we let $\eta_{n_0}^\beta \le \delta$ and $n \ge n_0$. We conclude the proof by letting $\gC=2\gamma^2 T$.
\end{proof}

\begin{lemma}\label{lem:tight_v}
Suppose that Assumptions \ref{asp:smooth}, \ref{asp:str_cov} and \ref{asp:noi_lip} holds. Then the sequence of random processes $\{\Bar{\vv}^{(n)}\}$ is tight under the Skorokhod topology in a finite interval.
\end{lemma}

\begin{lemma}\label{lem:generator_v}
Suppose Assumptions \ref{asp:smooth}, \ref{asp:str_cov} and \ref{asp:noi_lip} holds, and assume that there exists a positive number $p>2$ such that $\sup\limits_{n\ge 0}\sE\|\xi_n\|^p<\infty$. When $p_n = \gamma \eta_n^\beta$ with $\gamma > 0$,
for any $C^2$ function $g:\sR^d\to \sR$,  compactly supported with Lipschitz continuous second derivatives, we have
\begin{equation}
    \sE [g(\check{\vv}_{n+1})-g(\check{\vv}_n)|\gF_n]= \eta_n^\beta\gJ g(\check{\vv}_n) + \gT^g_n,
\end{equation}
where $\frac{1}{\eta_n^\beta}\gT^g_n \to 0$ in $L_1$ and $\gJ$ is the infinitesimal generator defined by
\begin{equation}\label{eq:def_J}
    \forall \phi \in \mathcal{C}^{2}\left(\mathbb{R}^{p}\right), \quad \mathcal{J} \phi(\vx)=\left\langle-\nabla f(\vx^\star), \nabla \phi(\vx)\right\rangle+\gamma (\phi(\mathbf{0})-\phi(\vx)).
\end{equation}
\end{lemma}

From the It\^o's formula for the semimartingales, we know that the infinitesimal generator $\gJ$ defined in Lemma \ref{lem:generator_v} corresponds to the following stochastic differential equation driven by the Poisson process with intensity $\gamma$
\begin{equation}\label{eq:jump_approx_append}
\begin{aligned}
d\rmY_t = - \nabla f(\vx^\star) dt - \rmY_t\cdot \rmN_\gamma(dt).
\end{aligned}
\end{equation}

\begin{lemma}\label{lem:jump_ergo_inv}
There exists a unique invariant measure $\vu^\star$ for the L\'evy process \eqref{eq:jump_approx_append}. Further, for any initial distribution $\nu_0$, we have $\gW_2(\gG^t \nu_0, \nu^\star)\to 0$ as $t\to \infty$, where $\gW_2$ represents the Wasserstain-2 distance and $\{\gG^t\}$ is the Markovian semigroup generated by the infinitesimal generator $\gJ$.
\end{lemma}

The proofs of Lemma~\ref{lem:tight_v}, \ref{lem:generator_v} and \ref{lem:jump_ergo_inv} can be found in Appendix~\ref{sec:aux_prf_2}

\begin{proof}[
{Proof of Theorem \ref{thm:jump_approx}}]

This proof is basically modeled after the proof of Theorem \ref{thm:diff_approx}. Therefore, for the sake of narrative simplicity, we will omit some details that overlap with the previous proofs. All symbols follow the meaning in the proof of Theorem \ref{thm:diff_approx} without special specification. And we denote $\gD_k$ as the natural filtration generated by $\{\check{\vv}_i\}_{i=1}^{k}$. Analogously, two steps are split to complete to proof.
\paragraph{Step 1} Let $g$ belongs to $\gC$ and let $\gD_t^{(n)}$ denote the natural filtration of $\bar{\vv}_t^{(n)}$. We aim to find the following martingale decomposition,
\begin{equation}
    \forall t>0, \quad g(\bar{\vv}_t^{(n)})-g(\bar{\vv}_0^{(n)})-\int_0^t \gJ g\left(\bar{\vv}_s^{(n)}\right)ds = \gN_t^{(n,g)} + \gT_t^{(n,g)},
\end{equation}
where $\gN_t^{(n,g)}$ is a $\gD_t^{(n)}$-martingale and $\gT_t^{(n,g)}$ converges to zero in $L_1$.
In fact, let
\begin{equation}
\begin{aligned}
& \gN_t^{(n,g)} &&= \ssum{k}{n+1}{N(n,t,\eta^\beta)}\{g(\check{\vv}_{k+1})-g(\check{\vv}_k)- \sE[g(\check{\vv}_{k+1})-g(\check{\vv}_k)|\gD_k]\}, \\ 
& \gT_t^{(n,g)} &&= g(\bar{\vv}_t^{(n)})-g(\bar{\vv}_{\underline{t}_n(\eta^\beta)}^{(n)}) - \int_{\underline{t}_n(\eta^\beta)}^t \gJ g(\bar{\vv}_s^{(n)})ds\\ 
& && \quad + \int_0^{\underline{t}_n(\eta^\beta)}\left(\gJ g\left(\bar{\vv}_{\underline{s}_n(\eta^\beta)}^{(n)}\right) - \gJ g\left( \bar{\vv}_s^{(n)}\right)\right)ds + \ssum{k}{n}{N(n,t,\eta^\beta)-1} \gT_k^g.
\end{aligned}
\end{equation}
Using the definition formula of $\bar{\vv}_t^{(n)}$ \eqref{eq:itpl_v}, when $t \notin \gL_n(T)$, we have
\begin{equation}
\begin{aligned}
\sE \left\| \bar{\vv}_t^{(n)} - \bar{\vv}_{\underline{t}_n(\eta^\beta)}\right\|^2 &= 
(t-\underline{t}_n(\eta^\beta))^2\sE\|\vd_{N(n,t,\eta^\beta)}-\xi^{(2)}_{N(n,t,\eta^\beta)}\|^2 \\ 
&\le \gC \eta_n^{2\beta}.
\end{aligned}
\end{equation}
This inequality combined with the Lipschitz continuity of $g$ and its derivatives implies that the
first three terms in the definition of $\gT_t^{(n,g)}$ tend to $0$ when $n\to \infty$. Further, by Lemma \ref{lem:generator_v},
\begin{equation}
\begin{aligned}
\sE\left|\ssum{k}{n}{N(n,t,\eta^\beta)-1}\gT_k^g\right|&\le \ssum{k}{n}{N(n,t,\eta^\beta)-1}\eta_k^{\beta}\sE\left|\frac{\gT_k^g}{\eta_k^\beta}\right|\\ 
&\le \sup\limits_{k\ge n}\sE\left|\frac{\gT_k^g}{\eta_k^\beta}\right|(t+\eta_n^\beta)\xrightarrow[]{n\to \infty} 0.
\end{aligned}
\end{equation}

\paragraph{Step 2} Suppose that there is a weakly convergent subsequence $\{\check{\vv}_{n_k}\}_{k=1}^\infty$ with limit distribution $\tilde{\nu}$. The definition of $M(n,t,\eta^\beta)$ and $\tilde{t}_n(\eta^\beta)$ can be extended intuitively from the first paragraph in the second step of Theorem \ref{thm:diff_approx}'s proof.

Owing to Prokhorov's theorem and Lemma \ref{lem:tight_v}, for any $T>0$, there is a weakly convergent subsequence in $\left\{\bar{\vv}_t^{M(n_k,T,\eta^\beta)}\right\}$. By Theorem 1 in \citep{kushner1980martingale}, we know that the weak limit of this sequence is a solution of the stochastic differential equation \eqref{eq:jump_approx_append}. And WLOG, we assume the sequence itself converges weakly to a solution $\bar{\vv}_t^{\tilde{\pi}^{(T)}}$ of \eqref{eq:jump_approx_append} with initial distribution $\tilde{\pi}^{(T)}$ (the notation of $\tilde{\nu}$ and $\tilde{\pi}$ is independent with that in proof of Theorem \ref{thm:diff_approx}). By Lemma \ref{lem:tight_v}, for any given $\epsilon > 0$, a compact set $K_\epsilon$ can be found such that $\sup\limits_{n}\sP(\check{\vv}_n \in K^c_\epsilon)\le \epsilon$. Therefore, for all $T>0$, $\tilde{\pi}^{(T)}(K_\epsilon)\ge 1-\epsilon$.

Due to Lemma \ref{lem:jump_ergo_inv}, a $T_\epsilon$ can be found such that
\begin{equation}
    \sup\limits_{\vi\in K_\epsilon}\sup\limits_{g\in \gC}\left|\gG^{T_\epsilon}g(\vx)-\langle \nu^\star, g\rangle\right| \le \epsilon,
\end{equation}
where $\gG$ is the Markov semigroup induced by the SDE \eqref{eq:jump_approx_append}. Because $\widetilde{(T_\epsilon)}_n(\eta^\beta)$ converges to $T_\epsilon$ when $n\to \infty$, we have $\check{\vv}_{n_k}\left(=\bar{\vv}_{\widetilde{(T_\epsilon)}_n(\eta^\beta)}^{M(n_k,T_\epsilon,\eta^\beta)}\right)$ converges weakly to the limit random variable of the sequence $\bar{\vv}_{T_\epsilon}^{M(n_k,T_\epsilon, \eta^\beta)}$,i.e., $\bar{\vv}_{T_\epsilon}^{\tilde{\pi}^{(T_\epsilon)}}$. On the other hand, by assumption, $\check{\vv}_{n_k}$ converges weakly to $\tilde{\nu}$. Thus $\bar{\vv}_{T_\epsilon}^{\tilde{\pi}^{(T_\epsilon)}}\sim \tilde{\nu}$. Combining all the results we have obtained, an inequality corresponding to \eqref{eq:w_conv_bound} can be derived. Consequently, we obtain $\tilde{\nu} = \nu^\star$. Finally, by Prokhorov's theorem, $\check{\vv}_n$ convergence weakly to $\nu^\star$. Further, $\{\bar{\vv}_t^{(n)}\}$ converges weakly to the dynamics \eqref{eq:jump_approx_append} with stationary distribution $\nu^\star$ as initialization.
\end{proof}



\begin{proof}[
{Proof of Theorem \ref{thm:asy_bias_u}}]

We prove the target conclusion in two steps. First, we show that the mean of $\hat{\vu}_n$ converges to a constant non-zero vector. Next, we will see that the asymptotic variance of $\{\hat{\vu}_n\}$ is zero for $\beta \in \left(\frac{1}{2},1\right)$.

\paragraph{Step 1} The first thing we need to do is to derive the recursive relation for $\hat{\vu}_n$,
\begin{equation}\label{eq:hat_u_rec}
\begin{aligned}
\hat{\vu}_{n+1}&=\gP_{\mA^\bot}\frac{\vx_{(n+1)-}-\vx^\star}{\eta_n^{1-\beta}} = \frac{\gP_{\mA^\bot}(\vx_n-\vx^\star -\eta_n \nabla f(\vx_n) + \eta_n \xi_n)}{\eta_n^{1-\beta}}\\ 
& = \left(\frac{\eta_{n-1}}{\eta_n}\right)^{1-\beta}\hat{\vu}_n -\eta_n^\beta \gP_{\mA^\bot} \nabla f(\vx_n) + \eta_n^\beta \xi^{(1)}_n\\ 
& = \hat{\vu}_n -\eta_n^\beta \gP_{\mA^\bot}\left\{\nabla ^2 f(\vx^\star)(\vx_n-\vu_n) + [\nabla^2 f(\vartheta^{\vv}_n)-\nabla ^2 f(\vx^\star)](\vx_n-\vu_n)\right\}\\ 
& \quad + \left(\left(\frac{\eta_{n-1}}{\eta_n}\right)^{1-\beta}-1\right)\hat{\vu}_n + \eta_n^\beta \xi_n^{(1)}\\ 
& \quad - \eta_n^\beta \gP_{\mA^\bot}\left\{\nabla ^2 f(\vx^\star) (\vu_n - \vx^\star) + [\nabla ^2 f(\vartheta^{\vu}_n)-\nabla ^2 f(\vx^\star)](\vu_n - \vx^\star)\right\}\\ 
& = \left(\rmI - \eta_n \gP_{\mA^\bot}\nabla ^2 f(\vx^\star) \gP_{\mA^\bot}\right)\hat{\vu}_n -\eta_n \gP_{\mA^\bot}\nabla ^2 f(\vx^\star) \check{\vv}_n  + \eta_n^\beta \xi_n^{(1)} \\ 
& \quad + (\eta_n - \eta_n^\beta \eta_{n-1}^{1-\beta}) \gP_{\mA^\bot}\nabla ^2 f(\vx^\star) \hat{\vu}_n  + (\eta_n - \eta_n^\beta\eta_{n-1}^{1-\beta})\gP_{\mA^\bot}\nabla ^2 f(\vx^\star)\check{\vv}_n \\ 
& \quad - \left\{\eta_n^\beta  \gP_{\mA^\bot} \left(\nabla^2 f(\vartheta^{\vv}_n) - \nabla^2 f(\vx^\star)\right)(\vx_n - \vu_n)\right\}\\ 
& \quad - \left\{\eta_n^\beta \gP_{\mA^\bot}[\nabla ^2 f(\vartheta^{\vu}_n)- \nabla ^2 f(\vx^\star)](\vu_n -\vx^\star)\right\}\\
& \quad + \left(\left(\frac{\eta_{n-1}}{\eta_n}\right)^{1-\beta}-1\right)\hat{\vu}_n \\ 
& = : \left(\rmI - \eta_n \gP_{\mA^\bot}\left(\nabla ^2 f(\vx^\star)            -\frac{1-\beta}{\eta_0}\mathbbm{1}_{\{\alpha=1\}}\rmI_d\right)\gP_{\mA^\bot}\right)\hat{\vu}_n\\
& \quad -\eta_n \gP_{\mA^\bot}\nabla ^2 f(\vx^\star) \check{\vv}_n + \eta_n^\beta \xi_n^{(1)} + \eta_n \gR_n^{\vu},
\end{aligned}
\end{equation}
where $\vartheta^{\vu}_n$ and $\vartheta^{\vv}_n$ are two entrywise interpolation point between $\vu_n$ and $\vx_n$, $\vu_n$ and $\vx^\star$ respectively. And,
\begin{equation}\label{eq:Ru_def}
\begin{aligned}
\gR_n^{\vu} &= \left(1-\left(\frac{\eta_{n-1}}{\eta_n}\right)^{1-\beta}\right)\gP_{\mA^\bot}\nabla ^2 f(\vx^\star) (\hat{\vu}_{n}+ \check{\vv}_n) + \frac{1}{\eta_n}\left(\left(\frac{\eta_{n-1}}{\eta_n}\right)^{1-\beta} - 1 -\frac{1-\beta}{n}\mathbbm{1}_{\{\alpha=1\}}\right)\hat{\vu}_n\\ 
& \quad - \left(\frac{\eta_{n-1}}{\eta_n}\right)^{1-\beta} \gP_{\mA^\bot} \left\{(\nabla^2 f(\vartheta^{\vv}_n)-\nabla ^2 f(\vx^\star))\check{\vv}_n + (\nabla^2 f(\vartheta^{\vu}_n) - \nabla^2 f(\vx^\star))\hat{\vu}_n\right\}.
\end{aligned}
\end{equation}
The properties of the step size sequence $\{\eta_n\}$ tell us that, when $\alpha<1$,
\begin{equation}
\begin{aligned}
 1 - \left(\frac{\eta_{n-1}}{\eta_{n}}\right)^{1-\beta} &= 1 - \left(1 + \frac{\eta_{n-1}-\eta_n}{\eta_n}\right)^{1-\beta}\\ 
 & = 1 - (1+ o(\eta_n))^{1-\beta} = (1-\beta)o(\eta_n) = o(1)\eta_n.
\end{aligned}
\end{equation}
And when $\alpha = 1$,
\begin{align*}
    1 + \frac{1-\beta}{n} - \left(\frac{n}{n-1}\right)^{1-\beta} & = 1 + \frac{1-\beta}{n} - \left(1 + \frac{1-\beta}{n} + \gO(n^{-2})\right) = o(\eta_n), \\
    1 - \left( \frac{n}{n-1} \right)^{1-\beta} & = - \frac{1-\beta}{n} + o(\eta_n) = \gO(\eta_n).
\end{align*}
The result together with Theorem~\ref{thm:converge} can be used to guarantee the first line of \eqref{eq:Ru_def} being $o(1)$ in $L_2$.
By Assumptions \ref{asp:hess_lip} and \ref{asp:smooth}, $\nabla ^2 f(\cdot)$ is Lipschitz continuous and uniformly bounded. Then for any $\delta \in (0,1)$, $\nabla ^2 f(\cdot)$ is $\delta$- H\"oder continuous. By Proposition 3.1 in \citep{chen2021online}, we have
\begin{equation}\label{eq:p_mmt_conv}
\begin{aligned}
 &\sE\|\vu_n -\vx^\star\|^p \precsim \eta_n ^{p(1-\beta)},\\ 
 &\sE\norm{\vv_n}^ p \precsim \eta_n ^{p(1 -\beta)}
\end{aligned}
\end{equation}
when $\beta \in [\frac{1}{2}, 1)$. Based on these preparations, take $\delta = p/2 -1$ and use the H\"oder continuity of $\nabla^2 f(\cdot)$ and Young's inequality,
\begin{equation}
\begin{aligned}
 \sE\norm{(\nabla ^2 f(\vartheta^{\vv}_n) - \nabla^2 f(\vx^\star)) \check{\vv}_n}^2 &\precsim \sE\norm{\vartheta^{\vv}_n - \vx^\star}^{2\delta}\norm{\check{\vv}_n}^2\\ 
 &\precsim \sE\left(\norm{\vv_n}^{2\delta} + \norm{ \vu_n - \vx^\star}^{2\delta}\right)\norm{ \check{\vv}_n}^2 \\ 
 & = \frac{1}{\eta_{n-1}^{2(1-\beta)}}\sE\norm{\vv_n}^p + \frac{1}{\eta_{n-1}^{2(1-\beta)}}\sE \norm{\vu_n -\vx^\star}^{2\delta}\norm{\vv_n}^2\\ 
 & \precsim \frac{1}{\eta_{n-1}^{2(1-\beta)}}\left(\sE\norm{\vv_n}^p + \sE \norm{\vu_n -\vx^\star}^ p \right) \precsim \eta_{n}^{(p-2)(1-\beta)}.
\end{aligned}
\end{equation}
The same bound can be derived for $(\nabla^2 f(\vartheta^\vu_n)-\nabla ^2 f(\vx^\star))\check{\vu}_n$. These two results enable the second line of \eqref{eq:Ru_def} to be $o(1)$ in $L_2$. 
To simplify our writing, we denote $\vnu = -\frac{1}{\gamma}\nabla f(\vx^\star)$ and $\vmu = \frac{1}{\gamma}\left(\gP_{\mA^\bot}\left(\nabla^2 f(\vx^\star) -\frac{1-\beta}{\eta_0}\mathbbm{1}_{\{\alpha=1\}}\rmI_d\right) \gP_{\mA^\bot}\right)^\dag \gP_{\mA^\bot} \nabla ^2 f(\vx^\star) \nabla f(\vx^\star)$. Taking the expectation on both sides of \eqref{eq:hat_u_rec} yields
\begin{equation}\label{eq:mean_hat_u_1}
\begin{aligned}
 \sE \hat{\vu}_{n+1} &= \left(\rmI - \eta_n \left(\gP_{\mA^\bot}\nabla^2 f(\vx^\star) -\frac{1-\beta}{\eta_0}\mathbbm{1}_{\{\alpha=1\}}\rmI_d \right) \gP_{\mA^\bot} \right)\sE \hat{\vu}_n\\ 
 & \quad -\eta_n \gP_{\mA^\bot} \nabla^2 f(\vx^\star) \sE \check{\vv}_n + \eta_n \sE \gR^{\vu}_n.
\end{aligned}
\end{equation}
Subtract $\vmu$ 
from both sides of \eqref{eq:mean_hat_u_1} and we have
\begin{equation}\label{eq:mean_hat_u_2}
\begin{aligned}
 \sE\hat{\vu}_{n+1} - \vmu &= \left(\rmI - \eta_n \gP_{\mA^\bot}\left(\nabla^2 f(\vx^\star)-\frac{1-\beta}{\eta_0}\mathbbm{1}_{\{\alpha=1\}}\rmI_d\right) \gP_{\mA^\bot}\right)(\sE \hat{\vu}_n -\vmu)\\ 
 & \quad - \eta_n \left(\gP_{\mA^\bot} \nabla^2 f(\vx^\star) (\sE (\check{\vv}_n - \vnu) + \sE \gR^{\vu}_n\right).
\end{aligned}
\end{equation}
According to Theorem \ref{thm:jump_approx} we know $\norm{\sE\check{\vv}_n - \vnu} = o(1)$. As a result of this and $\norm{\sE \gR^{\vu}_n}\le \sE\norm{\gR^{\vu}_n} = o(1)$, the last term in the right hand of \eqref{eq:mean_hat_u_2} is $o(\eta_n)$. Therefore, using Lemma \ref{lem:converge_r_t_3}, it holds that $\norm{\sE \hat{\vu}_n - \vmu} = o(1)$, which means $ \sE \hat{\vu}_n \xrightarrow[]{n\to\infty}\vmu$.

\paragraph{Step 2} Before we calculate the asymptotic variance of $\hat{\vu}_n$, we introduce the following lemma which reveals the 'mixing' nature between two sequence $\{\hat{\vu}_n\}_{n=1}^\infty$ and $\{\check{\vv}_n\}_{n=1}^\infty$.
\begin{lemma}\label{lem:asmp_mix_uv}
Under the Assumptions of Theorem~\ref{thm:asy_bias_u}
\[
\lim\limits_{n\to\infty}\abs{\sE\inner{\hat{\vu}_n - \sE\hat{\vu}_n}{\nabla^2 f(\vx^\star)(\check{\vv}_n -\sE\check{\vv}_n)}} = 0.
\]
\end{lemma}

\begin{proof}[Proof of Lemma~\ref{lem:asmp_mix_uv}]

Actually, we can derive the following recursive relationship
\begin{equation}
\begin{aligned}
 &\abs{\sE \inner{\hat{\vu}_{n+1}- \sE \hat{\vu}_{n+1}}{\nabla^2 f(\vx^\star)(\check{\vv}_{n+1}-\sE\check{\vv}_{n+1})}}\\ 
 =& (1 - \gamma \eta_n^\beta)\abs{\sE \inner{\hat{\vu}_{n+1}- \sE \hat{\vu}_{n+1}}{\nabla^2 f(\vx^\star)(\check{\vv}_{(n+1)-}-\sE\check{\vv}_{n+1})}} + \gamma\eta_n^\beta \abs{\sE \inner{\hat{\vu}_{n+1} - \sE \hat{\vu}_{n+1}}{\nabla^2 f(\vx^\star)\sE\check{\vv}_{n+1}}}\\ 
 =& (1 - \gamma \eta_n^\beta)\abs{\sE \inner{\hat{\vu}_{n+1}- \sE \hat{\vu}_{n+1}}{\nabla^2 f(\vx^\star)(\check{\vv}_{(n+1)-}-\sE\check{\vv}_{n+1})}}\\ 
  =& (1 - \gamma \eta_n^\beta)\abs{\sE \inner{\hat{\vu}_{n+1}- \sE \hat{\vu}_{n+1}}{\nabla^2 f(\vx^\star)(\check{\vv}_{(n+1)-}-\sE\check{\vv}_{(n+1)-})}}\\ 
 =& (1-\gamma \eta_n^\beta) \left| \sE\left\langle \left(\rmI - \eta_n \gP_{\mA^\bot} \left(\nabla ^2 f(\vx^\star) -\frac{1-\beta}{\eta_0}\mathbbm{1}_{\{\alpha=1\}}\rmI_d\right)\gP_{\mA^\bot}\right)(\hat{\vu}_n -\sE \hat{\vu}_n) - \eta_n \gP_{\mA^\bot} \nabla^2 f(\vx^\star) (\check{\vv}_n -\sE \check{\vv}_n)\right.\right.\\
 &\left.\left.+ \eta_n^\beta \xi^{(1)}_n + \eta_n(\gR^{\vu}_n-\sE\gR^{\vu}_n), \nabla^2 f(\vx^\star)\left\{(\check{\vv}_n -\sE\check{\vv}_n) + \eta_n^\beta \xi^{(2)}_n + \eta_n^\beta( \vd_n
 - \sE \vd_n
 )\right\}\right\rangle\right|\\ 
 \le & (1 - \gamma \eta^\beta_n)\abs{\sE \inner{\hat{\vu}_n - \sE \hat{\vu}_n}{\nabla^2 f(\vx^\star)(\check{\vv}_n - \sE\check{\vv}_n)}} + \eta^\beta_n \abs{\sE\inner{\hat{\vu}_n - \sE\hat{\vu}_n}{\nabla^2 f(\vx^\star)( \vd_n
 - \sE \vd_n
 )}} + \gO(\eta_n),
\end{aligned}
\end{equation}
where $\vd_n$ is defined as (\ref{eq:updt_ck_v}).
From the fact $\sE\| \vd_n
\|^2 = o(1)$, we have
\begin{equation}\label{eq:asymp_uncor_uv}
\begin{aligned}
 &\abs{\sE \inner{\hat{\vu}_{n+1}- \sE \hat{\vu}_{n+1}}{\nabla^2 f(\vx^\star)(\check{\vv}_{n+1}-\sE\check{\vv}_{n+1}})}
 \le (1 - \gamma \eta^\beta_n)\abs{\sE \inner{\hat{\vu}_n - \sE \hat{\vu}_n}{\nabla^2 f(\vx^\star)(\check{\vv}_n - \sE\check{\vv}_n})} + o(\eta^\beta_n).
\end{aligned}
\end{equation}
We can obtain from Lemma \ref{lem:converge_r_t_3} that $\abs{\sE \inner{\hat{\vu}_n - \sE \hat{\vu}_n}{\nabla^2 f(\vx^\star)(\check{\vv}_n - \sE\check{\vv}_n)}} \xrightarrow{n\to\infty} 0$.
\end{proof}

Back to the main result's proof, we can write down the recursive rule for the variance of $\hat{\vu}_n$, 
\begin{equation}
\begin{aligned}
 &\sE \norm{\hat{\vu}_{n+1}-\sE\hat{\vu}_{n+1}}^2\\ 
 =& \sE \left\|\left(\rmI - \eta_n \gP_{\mA^\bot}\left(\nabla^2 f(\vx^\star)-\frac{1-\beta}{\eta_0}\mathbbm{1}_{\{\alpha=1\}}\rmI_d\right)\gP_{\mA^\bot}\right)(\hat{\vu}_n -\sE \hat{\vu}_n)\right.\\ 
 &\left.- \eta_n \gP_{\mA^\bot} \nabla ^2 f(\vx^\star)(\check{\vv}_n -\sE \check{\vv}_n) + \eta_n^\beta \xi^{(1)}_n + \eta_n (\gR^\vu_n - \sE\gR^\vu_n)\right\|^2\\ 
 \le & (1- \mu \eta_n / 2) \sE\norm{\hat{\vu}_n - \sE\hat{\vu}_n}^2 - \eta_n \sE\inner{\hat{\vu}_n - \sE\hat{\vu}_n}{ \nabla^2 f(\vx^\star) (\check{\vv}_n - \sE\check{\vv}_n)} + \eta_n^{2\beta}\sE\inner{\xi^{(1)}_n}{\xi^{(1)}_n} + o(\eta_n)\\
 \le& (1- \mu \eta_n / 2) \sE\norm{\hat{\vu}_n - \sE\hat{\vu}_n}^2 + o(\eta_n),
\end{aligned}
\end{equation}
where the last equation follows from the diminished correlation we derived just now and the precondition $\beta > \frac{1}{2}$. Finally, the Lemma is completed from Lemma \ref{lem:converge_r_t_3} and the fact $\sE\norm{\hat{\vu}_n - \vmu}^2 = \sE\norm{\hat{\vu}_n -\sE \hat{\vu}_n}^2 + \norm{\sE\hat{\vu}_n - \vmu}^2 \to 0$.
\end{proof}

\subsubsection{Proof of Lemma~\ref{lem:tight_v}, \ref{lem:generator_v} and \ref{lem:jump_ergo_inv}}\label{sec:aux_prf_2}

\begin{proof}[
{Proof of Lemma \ref{lem:tight_v}}]
What we need to do now is to verify the conditions in Proposition \ref{prop:tight_cri_gen} one by one.

For a given path $\bar{\vv}^{(n)}$, denote $t^\prime_n= \max\{ s\in\gL_n(T)\cap [0, t]\}\cup \{0\}$. Then $\bar{\vv}^{(n)}_{t^\prime_n}=0$ whenever $t^\prime_n > 0$. Let $R> 0$, we have the following inequalities
\begin{equation}
\begin{aligned}
    &\sP\left(\sup\limits_{t\in [0,T]}\left\|\bar{\vv}^{(n)}_t\right\|\ge R\right)\\
    \le& \sP\left(\sup\limits_{t\in[0,T]}\left\{\left\|\bar{\vv}^{(n)}_t-\bar{\vv}^{(n)}_{\underline{t}_n(\eta^\beta)}\right\| + \left\|\bar{\vv}^{(n)}_{\underline{t}_n(\eta^\beta)}-\bar{\vv}^{(n)}_{t^\prime_n}\right\| + \left\|\bar{\vv}^{(n)}_{t^\prime_n}\right\|\right\}\ge R\right)\\
    \le& \frac{2}{R}\sE\sup\limits_{t\in [0,T]}\left\{\left(t-\underline{t}_n(\eta^\beta)\right) \left\|\vd_{\underline{t}_n(\eta^\beta)}- \xi_{\underline{t}_n(\eta^\beta)}^{(2)}\right\| + \left\|\bar{\vv}^{(n)}_{\underline{t}_n(\eta^\beta)}-\bar{\vv}^{(n)}_{t^\prime_n}\right\|\right\}
    + \frac{2}{R}\sE\left\|\bar{\vv}^{(n)}_0\right\|\\ 
    \le & \frac{2}{R}\sE\sup\limits_{t\in \gL_n(T)}\left\|\bar{\vv}^{(n)}_{t-}-\bar{\vv}^{(n)}_{(t-)^\prime_n}\right\| + \frac{2}{R}\sE\left\|\check{\vv}_n\right\|\\ 
    \le & \frac{2}{R}\sE\sup\limits_{k\in\gI_n(T)}\ssum{i}{n}{n+k-1}\eta_i^\beta \left\|\vd_i + \xi^{(2)}_i\right\| + \frac{2}{R}\sE\|\check{\vv}_n\|\\ 
    =& \frac{2}{R}\ssum{i}{n}{n+\sup\gI_n(T)-1}\eta_i^\beta \sE\left\|\vd_i + \xi^{(2)}_i\right\| + \frac{2}{R}\sE\|\check{\vv}_n\|\\ 
    \le & \frac{2}{R}(T+\eta_n)\sup\limits_{n\in\sN}\sE\left\|\vd_i+\xi^{(2)}_i\right\| + \frac{2}{R}\sE\|\check{\vv}_n\|\\ 
    \overset{(a)}{\le} & \frac{\gC (1+T)}{R} \le \zeta,
\end{aligned}
\end{equation}
where (a) holds for the uniform boundedness of $\sE\|\vd_n+\xi_n^{(2)}\|+\sE\|\check{\vv}_n\|$. And the final inequality holds when we take $R > \frac{\gC(1+T)}{\zeta}$. Thus, the first condition of Proposition \ref{prop:tight_cri_gen} holds for $\bar{\vv}^{(n)}$.

As for the second condition, what we should do is to construct an appropriate partition that makes $\omega(\bar{\vv}^{(n)},\delta, T)$ defined as \eqref{eq:part_var} as small as possible.

For a given $\epsilon, \eta$ pair, let $\delta < \frac{\zeta}{2\gC}$, then from Lemma \ref{lem:sparse_jump} it can be seen $\sP(\Delta(\gJ_n(T))< \delta) < \frac{\zeta}{2} $. Now given the event $\gE = \{\Delta(\gJ_n(T))\ge \delta\} $, we choose the partition points $\{\tau_k\}\in [0,T]$ recursively from the set $\gL_n(T)$ such that the partition satisfies the following properties:
\begin{enumerate}
    \item $\min\limits_{k}\{\tau_k-\tau_{k-1}\} \in [\delta,3\delta)$,
    \item $\gJ_n(T) \subset \{\tau_k\}$.
\end{enumerate}
Let $\tau_0 = 0$ and suppose we have constructed the partition points $\tau_0,\dots,\tau_k\in [0,T]$ with inductive assumptions:
\begin{enumerate}
    \item $\min\limits_{i \le k-1}\{\tau_{i+1}-\tau_i\} \in [\delta,3\delta)$,
    \item $\gJ_n(\tau_k) \subset \{\tau_0,\tau_1,\dots, \tau_k\}$,
    \item there is no discontinuous point in $(\tau_k, \tau_k + \delta)$, i.e., $\gJ_n(T)\cap (\tau_k, \tau_k+ \delta)= \varnothing$.
\end{enumerate}
We will use these results to find the next partition point $\tau_{k+1}$. Define $\tilde{\tau}_{k+1} = \min\{t: t-\tau_k\ge \delta, t\in \gL_n(T)\}$, we use the following scheme:
\begin{equation}\label{eq:recur_tau}
    \tau_{k+1} = \left\{
    \begin{array}{cc}
         s,  & \exists ~s\in \left(\tilde{\tau}_{k+1}, \tilde{\tau}_{k+1}+\delta\right)\cap \gJ_n(T), \\
         \tilde{\tau}_{k+1},  & \text{Otherwise}.
    \end{array}\right.
\end{equation}
From the property of the event $\gE$ we know there is at most one discontinuous point in $(\tilde{\tau}_{k+1}, \tilde{\tau}_{k+1}+ \delta)$, which means the $\tau_{k+1}$ is always well-defined. Then we have $\delta \le \tau_{k+1}-\tau_k \le \tau_{k+1}-\tilde{\tau}_{k+1} + \tilde{\tau}_{k+1}-\tau_k\le \delta + \eta_n + \delta \le 3\delta$, where the last inequality holds when we choose $n_0$ such that $ \eta_{n_0}< \delta$. Thus $\tau_{k+1}$ satisfies the first inductive assumption.

By the third inductive assumption of $\tau_k$, we know there is no discontinuous point in $(\tau_k, \tilde{\tau}_{k+1})$. On the other hand, if $\tau_{k+1}\in \gJ_n(T)$, then $(\tau_{k+1}-\delta,\tau_{k+1})\cap \gJ_n(T)=\varnothing$, and especially we have $[\tilde{\tau}_{k+1},\tau_{k+1})\cap \gJ_n(T) = \varnothing$. Hence, the second inductive assumption of $\tau_{k+1}$ has been proved.

Finally, if $\tau_{k+1}=\tilde{\tau}_{k+1}$, then from the recursive construction scheme \eqref{eq:recur_tau} we know $(\tau_{k+1},\tau_{k+1}+\delta)\cap \gJ_n(T)=(\tilde{\tau}_{k+1},\tilde{\tau}_{k+1}+\delta)\cap \gJ_n(T)=\varnothing$. Otherwise, $\tau_{k+1}$ must belong to $\gJ_n(T)$. Combining this with the definition of $\gE$ we can make sure $(\tau_{k+1}, \tau_{k+1}+\delta)\cap \gJ_n(T)=\varnothing$. At this point, the proofs of the three inductive assumptions on $\tau_{k+1}$ are all complete.

With the partition $\{ \tau_k \}$, we can bound the tail probability of $\omega(\bar{\vv}^{(n)},\delta, T)$ as follows.
\begin{equation}\label{cond_2_back}
\begin{aligned}
    \sP(\omega(\bar{\vv}^{(n)}, \delta, T) \ge \epsilon) &\le \sP(\omega(\bar{\vv}^{(n)},\delta, T)\ge \epsilon ;~ \gE) + \sP(\gE^c)\\
    &\le \sP\left(\max\limits_{k}\sup\limits_{\tau_k\le t< s< \tau_{k+1}}\left\|\bar{\vv}^{(n)}_t-\bar{\vv}^{(n)}_s\right\|\ge \epsilon; ~ \gE\right) + \frac{\zeta}{2}\\ 
    &\le 2\sP\left(\sup\limits_{t\in [0,T]}\left\|\bar{\vv}^{(n)}_t-\bar{\vv}^{(n)}_{\underline{t}_n(\eta^\beta)}\right\|\ge \frac{\epsilon}{2}\right) + \frac{\zeta}{2} \\ 
    & \quad + \sP\left(\max\limits_{k}\sup\limits_{\tau_k\le t< s< \tau_{k+1}}\left\|\bar{\vv}^{(n)}_{\underline{t}_n(\eta^\beta)}-\bar{\vv}^{(n)}_{\underline{s}_n(\eta^\beta)}\right\|\ge \frac{\epsilon}{2}; ~ \gE\right).
\end{aligned}
\end{equation}
We will give the bound of two probabilities respectively. First,
\begin{equation}\label{eq:tight_margin_term}
\begin{aligned}
    \sP&\left( \sup\limits_{t\in [0,T]}\left\|\bar{\vv}^{(n
    )}_t-\bar{\vv}^{(n)}_{\underline{t}_n(\eta^\beta)}\right\|\ge \frac{\epsilon}{2}\right) = \sP\left(\sup\limits_{k\in \gI_n(T)}\left\|\bar{\vv}^{(n)}_{(\Gamma_{k+1}-\Gamma_n)-}-\bar{\vv}^{(n)}_{\Gamma_k-\Gamma_n}\right\|\ge \frac{\epsilon}{2}\right)\\
    &\le \sum\limits_{k\in \gI_n(T)}\sP\left(\left\|\check{\vv}_{(n+k+1)-}-\check{\vv}_{n+k}\right\|\ge\frac{\epsilon}{2}\right)
    \le \sum\limits_{k\in\gI_n(T)}\frac{4}{\epsilon^2}\sE\left\|\check{\vv}_{(n+k+1)-}-\check{\vv}_{n+k}\right\|^2\\
    &\le \frac{4}{\epsilon^2}\sum\limits_{k\in\gI_n(T)}\eta_{n+k}^{2\beta}\sE\|\vd_{n+k}+\xi^{(2)}_{n+k}\|^2 \le \frac{4\eta_n^\beta\sup\limits_i\sE\|\vd_i+\xi^{(2)}_i\|^2}{\epsilon^2}\sum\limits_{k\in\gI_n(T)}\eta_{n+k}^\beta\\ 
    &\le \frac{\gC(T+\eta_n^\beta)}{\epsilon^2}\eta_n^\beta < \frac{2\gC T}{\epsilon}\eta_n^\beta < \frac{\zeta}{8},
\end{aligned}
\end{equation}
where the last inequality holds when we take $\eta_{n_0}^\beta < \frac{\epsilon\zeta}{16\gC T}$ and $n \ge n_0$.

It is easy to see that we can use a bijection to link the elements in $\gI_n(T)$ and those in $\gL_n(T)$. Because the partition points $\{\tau_k\}$ are in $\gL_n(T)$, we assume that every $\tau_k$ corresponds to an index $\varsigma_k \in \gI_n(T)$. Then we have $\varsigma_{k+1} > \varsigma_k$. Denote $\mathfrak{S}_k = \gI_n(T)\cap [\varsigma_k,\varsigma_{k+1})$. So far we are ready to bound the last term in \eqref{cond_2_back}.
\begin{equation}\label{eq:tight_cent_term}
\begin{aligned}
&\sP \left(\max\limits_{k}\sup\limits_{\tau_k\le t< s< \tau_{k+1}}\left\|\bar{\vv}^{(n)}_{\underline{t}_n(\eta^\beta)}-\bar{\vv}^{(n)}_{\underline{s}_n(\eta^\beta)}\right\|\ge \frac{\epsilon}{2}; ~ \gE\right)\\ 
\le& \sum\limits_{k}\sP\left(\sup\limits_{l,h\in \mathfrak{S}_k}\left\|\check{\vv}_{n+l}-\check{\vv}_{n+h}\right\|\ge \frac{\epsilon}{2}; ~ \gE\right)\\
\le& \sum\limits_{k}\sP\left(\sup\limits_{l,h\in \mathfrak{S}_k}\ssum{i}{l}{h-1}\left\|\check{\vv}_{n+i+1}-\check{\vv}_{n+i}\right\|\ge \frac{\epsilon}{2};~ \gE\right) =\sum\limits_{k}\sP\left(\ssum{i}{\varsigma_k}{\varsigma_{k+1}-1}\left\|\check{\vv}_{n+i+1}-\check{\vv}_{n+i}\right\|\ge\frac{\epsilon}{2}; ~\gE\right)\\ 
\overset{(a)}{\le}& \sum\limits_{k}\sP\left(\ssum{i}{\varsigma_k}{\varsigma_{k+1}-1}\eta^\beta_{n+i}\left\|\vd_{n+i}+\xi^{(2)}_{n+i}\right\|\ge \frac{\epsilon}{2}\right) 
\le \sum\limits_{k}\frac{4}{\epsilon^2}\sE\left(\ssum{i}{\varsigma_k}{\varsigma_{k+1}-1}\eta_{n+i}^\beta\left\|\vd_{n+i}+\xi^{(2)}_{n+i}\right\|\right)^2\\
\le& \frac{4}{\epsilon^2}\sum\limits_{k}\left\{\left(\ssum{i}{\varsigma_k}{\varsigma_{k+1}-1}\eta_{n+i}^\beta\right)\left(\ssum{i}{\varsigma_k}{\varsigma_{k+1}-1}\eta^\beta_{n+i}\sE \left\|\vd_{n+i}+\xi^{(2)}_{n+i}\right\|^2\right)\right\}\\ 
\le& \frac{4\sup\limits_{i}\sE\|\vd_i+\xi^{(2)}_i\|^2}{\epsilon^2}\sum\limits_k\left(\ssum{i}{\varsigma_k}{\varsigma_{k+1}-1}\eta_{n+i}^\beta\right)^2
\le \frac{\gC}{\epsilon^2}\sum\limits_{k}(\tau_{k+1}-\tau_k)^2\\ 
\overset{(b)}{\le}& \frac{3\gC \delta}{\epsilon^2}\sum\limits_{k}(\tau_{k+1}-\tau_{k}) 
\le \frac{3\gC T\delta}{\epsilon^2} < \frac{\zeta}{4},
\end{aligned}
\end{equation}
where (a) follows from the combination of the fact that the path $\bar{\vv}^{(n)}$ is continuous in any interval $[\tau_{k+1},\tau_k)$ when $\gE$ holds and the update formula \eqref{eq:itpl_v}. And (b) is true by the property of the partition $\{\tau_k\}$ listed above. The last inequality holds when we take $\delta < \frac{\zeta\epsilon^2}{12\gC T}$.

Bring \eqref{eq:tight_margin_term} and \eqref{eq:tight_cent_term} into \eqref{cond_2_back}, we have
\begin{equation}
\begin{aligned}
\sP (\omega(\bar{\vv}^{(n)},\delta, T)\ge \epsilon)\ge 2\cdot\frac{\zeta}{8} + \frac{\zeta}{2} + \frac{\zeta}{4}=\zeta.
\end{aligned}
\end{equation}
At this point, we have checked the two sufficient conditions in Proposition \ref{prop:tight_cri_gen}. Hence, the tightness of $\{\bar{\vv}^{(n)}\}$ has been proved.
\end{proof}

\begin{proof}[
{Proof of Lemma \ref{lem:generator_v}}]
We would like to say that the overall proof framework is similar to the proof of Lemma \ref{lem:generator_u}. However, since $\check{\vv}_{n+1}$ may suddenly jump to $0$, we cannot directly use Taylor expansion to get the desired result. First, by the scheme on $\check{\vv}_{n+1}$ jumping to zero, we have 
\begin{equation}\label{eq:gen_v_back}
\begin{aligned}
\sE[g(\check{\vv}_{n+1})-g(\check{\vv}_n)|\gF_n] = p_{n}(g(\mathbf{0})-g(\check{\vv}_n)) + (1-p_n)\sE[g(\check{\vv}_{(n+1)-})-g(\check{\vv}_n)|\gF_n].
\end{aligned}
\end{equation}
Then we make use of the Taylor expansion between $\check{\vv}_{(n+1)-}$ and $\check{\vv}_n$.
\begin{equation}
\begin{aligned}
&g(\check{\vv}_{(n+1)-})-g(\check{\vv}_n) = \langle \nabla g(\check{\vv}_{n}), \check{\vv}_{(n+1)-}-\check{\vv}_n \rangle \\ 
& \quad + \frac{1}{2}(\check{\vv}_{(n+1)-}-\check{\vv}_n)^\top \nabla^2 g(\varrho_n) (\check{\vv}_{(n+1)-}-\check{\vv}_n)\\
&= \eta_n^\beta\langle \nabla g(\check{\vv}_{n}), \nabla f(\vx^\star) + \xi^{(2)}_n\rangle + \eta_n^\beta\left\langle \nabla g(\check{\vv}_n) , \frac{1}{\eta_n^\beta}\ssum{i}{1}{3}\gS_n^{(i)} \right\rangle \\ 
& \quad + \frac{\eta_n^{2\beta}}{2}(\vd_n+\xi^{(2)}_n)^\top\nabla^2 g(\varrho_n)(\vd_n+ \xi^{(2)}_n).
\end{aligned}
\end{equation}
Substitute this equation into the second term of the right hand of the equation \eqref{eq:gen_v_back}. It follows that
\begin{equation}
\begin{aligned}
&\sE[g(\check{\vv}_{n+1})-g(\check{\vv}_n)|\gF_n] = \eta_n^\beta \left(\gamma(g(\mathbf{0})-g(\check{\vv}_n)) + \langle \nabla g(\check{\vv}_n), -\nabla f(\vx^\star) \rangle\right) \\
& \quad +\eta_n^\beta \sE\left[\left.\left\langle \nabla g(\check{\vv}_n), \frac{1}{\eta_n^\beta}\ssum{i}{1}{3}\gS_n^{(i)}  \right\rangle\right| \gF_n\right] +  \frac{\eta_n^{2\beta}}{2}\sE\left[\left. (\vd_n+\xi^{(2)}_n)^\top \nabla^2 g(\varrho_n) (\vd_n+\xi^{(2)}_n) \right|\gF_n\right]\\ 
& \quad - \gamma\eta_n^{2\beta} \sE\left[\left. \langle \nabla g(\check{\vv}_n), \vd_n + \xi^{(2)}_n\rangle + \frac{\eta_n^\beta}{2}(\vd_n + \xi^{(2)}_n)^\top \nabla^2 g(\varrho_n) (\vd_n + \xi^{(2)}_n)\right|\gF_n\right] \\
&=: \eta_n^\beta \{\gamma (g(\mathbf{0})-g(\check{\vv}_n))- \langle \nabla g(\check{\vv}_n), \nabla f(\vx^\star)\rangle\} + \underbrace{ \gT_n^{(1)} + \gT_n^{(2)} + \gT_n^{(3)}}_{\gT_n^g}.
\end{aligned}
\end{equation}
From Theorem \ref{thm:converge} and the equation \eqref{eq:T1_bound}, we have $\sE |\gT_n^{(1)}|=o(\eta_n^\beta)$. And $\sE|\gT_n^{(2)}| = \gO(\eta_n^{2\beta})$ by leveraging that $\|\nabla^2 g (\vx)\|$ is bounded for all $\vx\in \sR^d$ and that $\sE \|\vd_n + \xi^{(2)}_n\|^2$ is bounded. Similar approaches can be used to show that $\sE|\gT_n^{(3)}|= \gO (\eta_n^{2\beta})$. At this point, the result has been proved.
\end{proof}

\begin{proof}[
{Proof of Lemma \ref{lem:jump_ergo_inv}}]

Consider the set of probability density functions 
\[\left\{h(\vx)=p(t)\mathbbm{1}_{\left\{\vx=\frac{\nabla f(\vx^\star)}{\| 
\nabla f(\vx^\star)\|}t\right\}}: p(t) \text{ is a p.d.f on }\sR\right\}
\]
and denote it as $\gM$. Then the distribution of any $\rmY_t$ only has mass on the line $\left\{\frac{\nabla f(\vx^\star)}{\| \nabla f(\vx^\star)\|}t: t\in \sR\right\}$ if we choose the initial distribution in $\gM$. In this case, we can suppose $\rmY_t = -\frac{\nabla f(\vx^\star)}{\|\nabla f(\vx^\star)\|}\upsilon_t$. Consequently, $\upsilon_t$ satisfies the following one-dimensional stochastic differential equation
\begin{equation}\label{eq:norm_sde}
\begin{aligned}
d\upsilon_t = \|\nabla f(\vx^\star)\|dt - \upsilon_t \rmN_\gamma(dt).
\end{aligned}
\end{equation}
Let $\varphi_t(\lambda) = \sE_{p^\star} e^{i\lambda \upsilon_t}$ be the characteristic function of $\upsilon_t$ with stationary initialization $p^\star$. Then we have $\varphi_t(\lambda) = \varphi_s(\lambda),~ \forall t\neq s$. On the other hand, consider the martingale problem corresponding to \eqref{eq:norm_sde}. It says that $e^{i\lambda \upsilon_t}-e^{i\lambda \upsilon_0} - \int_0^t \left(i\lambda\|\nabla f(\vx^\star)\|e^{i\lambda \upsilon_s} + \gamma (1-e^{i\lambda \upsilon_s})\right)ds$ is a martingale with respect to the natural filtration generated by $\upsilon_t$. Taking expectations we have
\begin{equation}
\begin{aligned}
0&= \varphi_t(\lambda) -\varphi_0(\lambda) - \int_0^t \{i\lambda\|\nabla f(\vx^\star)\|\varphi_s(\lambda)+ \gamma (1-\varphi_s(\lambda))\}ds\\ 
&=  - \int_0^t \{i\lambda\|\nabla f(\vx^\star)\|\varphi_s(\lambda)+ \gamma (1-\varphi_s(\lambda))\}ds,
\end{aligned}
\end{equation}
which means that
\begin{equation}
\begin{aligned}
i\lambda\|\nabla f(\vx^\star)\|\varphi_s(\lambda) + \gamma (1-\varphi_s(\lambda)) = 0 ,\quad \forall s>0,
\end{aligned}
\end{equation}
i.e. $\varphi_s(\lambda) = \frac{1}{1-\frac{i\|\nabla f(\vx^\star)\|}{\gamma}\lambda}$. So the invariant distribution of $\upsilon_t$ is $\gE\left(\frac{\|\nabla f(\vx^\star)\|}{\gamma}\right)$. As a result, the invariant distribution of $\rmY_t$ is $- \frac{\nabla f(\vx^\star)}{\|\nabla f(\vx^\star)\|}\cdot \gE\left(\frac{\|\nabla f(\vx^\star)\|}{\gamma}\right)$.

To show the mixing result, it is enough to prove the following fact,
\begin{equation}\label{eq:dirac_init_mix}
\begin{aligned}
\forall \vy_0\neq \vy_1,\ 
\frac{1}{\|\vy_0-\vy_1\|}\gW_2(\gG^t \delta_{\vy_0}, \gG^t \delta_{\vy_1})\longrightarrow 0 \text{ as }  t\to \infty, \quad
\end{aligned}
\end{equation}
where $\delta_\vy$ represents the Dirac measure at the point $\vy$. Let $\rmY_t^0$ and $\rmY_t^1$ be the stochastic process generated by \eqref{eq:jump_approx_append} with initial distribution $\delta_{\vy_0}$ and $\delta_{\vy_1}$ respectively. To give a bound for the Wasserstain-2 distance between $\rmY_t^0$ and $\rmY_t^1$, we compute the $L_2$ norm under the identical coupling, i.e., the two dynamics share all randomness in the Poisson process $\rmN_\gamma(s), ~ s\in [0,t]$. Owing to the property of the corresponding martingale problem of \eqref{eq:jump_approx_append}, we have
\begin{equation}
\begin{aligned}
0 &= \sE\|\rmY_t^0-\rmY_t^1\|^2 - \|\vy_0-\vy_1\|^2 \\ 
& \quad - \int_0^t \sE\left\{-\left(\nabla f(\vx^\star)^\top, \nabla f(\vx^\star)^\top\right)\begin{bmatrix} \rmI & -\rmI \\ -\rmI & \rmI\end{bmatrix}\begin{bmatrix} \rmY_s^0 \\ \rmY_s^1\end{bmatrix} -\gamma \|\rmY_s^0-\rmY_s^\1\|^2\right\}ds \\ 
& = \sE\|\rmY_t^0 - \rmY_t^0\|^2 - \|\vy_0-\vy_1\|^2 + \gamma \int_0^1 \sE\|\rmY_s^0 - \rmY_s^1\|^2 ds.
\end{aligned}
\end{equation}
Solving above integral equation we finally get $\sE \|\rmY_t^0-\rmY_t^1\|^2 = \|\vy_0-\vy_1\|^2 e^{-\gamma t}$. Hence, the equation \eqref{eq:dirac_init_mix} has been proved.

\end{proof}

\subsection{Proof of Case 3}\label{sec:prf_of_cs3}
\begin{proof}[
{Proof of Theorem~\ref{thm:noncen_u}}]
Before the major proof, let's recall the update formula of $\{\check{\vu}_n\}$ and $\{\check{\vv}_n\}$.
\begin{equation}\label{eq:u_updt_fix}
\begin{aligned}
    \check{\vu}_{n+1} &= \check{\vu}_n - \eta_n\gP_{\mA^{\bot}}\left(\nabla^2f(\vx^{\star}) - \frac{1}{2\eta_0}\mathbbm{1}_{\alpha = 1}\rmI_d\right)\check{\vu}_n + \gR_n^{(1)} + \gR_n^{(2)} + \gR_n^{(3)} + \sqrt{\eta_n}\xi^{(1)}_n\\
    & =\check{\vu}_n - \eta_n\gP_{\mA^{\bot}}\left(\nabla^2f(\vx^{\star}) - \frac{1}{2\eta_0}\mathbbm{1}_{\alpha = 1}\rmI_d\right)\check{\vu}_n - \sqrt{\eta_n}\gP_{\mA^{\bot}}(\nabla f(\vx_n) - \nabla f(\vu_n))\\
    & \quad + \gR_n^{(2)} + \gR_n^{(3)} + \sqrt{\eta_n}\xi^{(1)}_n\\
    &= \check{\vu}_n - \eta_n\gP_{\mA^{\bot}}\left(\nabla^2f(\vx^{\star}) - \frac{1}{2\eta_0}\mathbbm{1}_{\alpha = 1}\rmI_d\right)\check{\vu}_n - {\eta_n}\gP_{\mA^{\bot}}\nabla^2 f(\vx^\star)\check{\vv}_n\\
    & \quad + \gR_n^{(4)} + \gR_n^{(2)} + \gR_n^{(3)} + \sqrt{\eta_n}\xi^{(1)}_n,
\end{aligned}
\end{equation}
where the definition of $\gR_n^{(i)}, i = 1,2,3$ can be referred to (\ref{eq:diff_res1-3}). And
\[
\gR^{(4)}_n = \sqrt{\eta_n}\gP_{\mA^{\bot}}\left(\nabla^2 f(\vx^\star) - \int_0^1 \nabla^2 f((1-\lambda)\vu_n - \lambda \vx_n)d\lambda\right)\vv_n.
\]

By making use of Lemma~\ref{lem:bound_res} and (\ref{eq:p_mmt_conv}), we can assert that all of $\gR_n^{(i)}, i = 2,3,4$ are $o(\eta_n)$ in $L_1$ sense. Denote $ - \frac{1}{\gamma}
f(\vx^\star)$ as $\vv$. Then from Theorem~\ref{thm:jump_approx} we know that $\EB \check{\vv}_n \to \vv$ as $n \to \infty$. Let $\tilde{\vu}_n = \check{\vu}_n - \vu$. Then (\ref{eq:u_updt_fix}) can be rewritten as
\begin{align*}
    \tilde{\vu}_{n+1}  &= \tilde{\vu}_n - \eta_n\gP_{\mA^{\bot}}\left(\nabla^2 f(\vx^\star) - \frac{1}{2\eta_0}\mathbbm{1}_{\alpha = 1}\rmI_d\right)\tilde{\vu}_n - \eta_n \gP_{\mA^\bot}\nabla^2 f(\vx^\star)(\check{\vv}_n - \vv)\\
    & \quad + \gR_n^{(2)} + \gR_n^{(3)} + \gR_n^{(4)} + \sqrt{\eta_n}\xi_n^{(1)}.
\end{align*}
The main idea of our proof is to reduce this special case to one of the first two cases. We construct a surrogate random sequence $\{\vz_n\}$ which satisfies the following relation
\begin{align*}
    \vz_{n+1} = \vz_n - \eta_n\gP_{\mA^{\bot}}\left(\nabla^2 f(\vx^\star) - \frac{1}{2\eta_0}\mathbbm{1}_{\alpha = 1}\rmI_d\right)\vz_n
    + \gR_n^{(2)} + \gR_n^{(3)} + \gR_n^{(4)} + \sqrt{\eta_n}\xi_n^{(1)}.
\end{align*}
Then,
\[
\tilde{\vu}_{n+1} - \vz_{n+1} = (\tilde{\vu}_{n} - \vz_{n}) - \eta_n \gP_{\mA^\bot}\left(\nabla^2 f(\vx^\star) - \frac{1}{2\eta_0}\mathbbm{1}_{\alpha = 1}\rmI_d\right)(\tilde{\vu}_n - \vz_n) - \eta_n\gP_{\mA^\bot}\nabla^2 f(\vx^\star)(\check{\vv}_n - \vv).
\]
As we can see, the update rule of $\vz_n$ is the same as that of $\check{\vu}_n$ in the first case. And it is not hard to verify that $\{\vz_n\}$ are located in the null space of matrix $\mA$. So we can assert that $\{\vz_n\}$ converges weakly to the Gaussian distribution $\gN(0,\tilde{\Sigma})$. What we aim to do next is to control the magnitude of the Wasserstein-2 distance between $\tilde{\vu}_n$ and $\vz_n$. And we will accomplish this goal in three steps.

\paragraph{Step 1}
First we will prove the fact $\lim_{n\to \infty}\left\|\EB\tilde{\vu}_n - \EB \vz_n\right\| = 0$. Actually, we have
\begin{align*}
&\|\EB\tilde{\vu}_{n+1} - \EB\vz_{n+1}\| \\
=& \left\|\EB(\tilde{\vu}_{n} - \vz_{n}) - \eta_n\gP_{\mA^\bot}\left(\nabla^2 f(\vx^\star) - \frac{1}{2\eta_0}\mathbbm{1}_{\alpha=1}\rmI_d\right)\EB(\tilde{\vu}_n - \vz_n) - \eta_n\gP_{\mA^\bot}\nabla^2 f(\vx^\star)(\EB \check{\vv}_n -\vv)\right\|\\
\overset{(a)}{\le}& \left\|\left\{\rmI_d - \eta_n \gP_{\mA^\bot}\left(\nabla^2 f(\vx^\star) - \frac{1}{2\eta_0}\mathbbm{1}_{\alpha = 1}\rmI_d\right)\gP_{\mA^\bot}\right\}\EB(\tilde{\vu}_{n} - \vz_n)\right\|+ \eta_n\left\|\nabla^2 f(\vx^\star)\right\|\|\EB\check{\vv}_n - \vv\|\\
\overset{(b)}{\le}&\left(1 - \left(\mu - \frac{1}{2\eta_0}\right)\eta_n\right)\|\EB\tilde{\vu}_n - \EB\vz_n\| + \eta_n \left\|\nabla^2 f(\vx^\star)\right\|\|\EB \check{\vv}_n -\vv\|\\
\overset{(c)}{\le}& \left(1 - \left(\mu - \frac{1}{2\eta_0}\right)\eta_n\right)\|\EB\tilde{\vu}_n - \EB\vz_n\| + o(\eta_n),
\end{align*}
where (a) holds for the triangular inequality, (b) holds since $\gP_{\mA^\bot}\left(\nabla^2 f(\vx^\star) - \frac{1}{2\eta_0}\mathbbm{1}_{\alpha = 1}\rmI_d\right)\gP_{\mA^\bot}$ is positive definite on the null space of $\mA^\top$, and (c) holds for the fact $\lim_{n\to \infty}\EB\check{\vv}_n = \vv$. Combine this and Lemma~\ref{lem:converge_r_t_3} and we verify the assertion.

\paragraph{Step 2}
Next we will show that $\tilde{\vu}_n - \vz_n$ and $\check{\vv}_n - \vv$ are asymptotically uncorrelated, i.e., 
\[
\left|\EB\langle \tilde{\vu}_n - \vz_n , \nabla^2 f(\vx^\star)(\check{\vv}_n - \vv) \rangle\right|\to 0,
\]
as $n\to \infty$. Here we need to use the update formula of $\check{\vv}_n$ ~(\ref{eq:updt_ck_v}) and~(\ref{eq:updt_ck_v_res}).
Let $\gS_n = (\gS_n^{(1)} + \gS_n^{(2)} + \gS_n^{(3)} ) / \eta_n^\beta $ with $\gS_n^{(i)}$ defined in (\ref{eq:updt_ck_v_res}). One can check $\gS_n = o(1)$. Then we have
\begin{equation}
\begin{aligned}
    &\left|\EB \left\langle \tilde{\vu}_{n+1} - \vz_{n+1}, \nabla^2 f(\vx^\star)(\check{\vv}_{n+1} - \EB\check{\vv}_{n+1})\right\rangle\right|\\
    =& \left|(1-\gamma\sqrt{\eta_n})\EB \left\langle \tilde{\vu}_{n+1} - \vz_{n+1}, \nabla^2 f(\vx^\star)(\check{\vv}_n - \EB\check{\vv}_{n+1})\right\rangle\right.\\
    &\left. -\sqrt{\eta_n}(1-\gamma\sqrt{\eta_n})\EB\left\langle \tilde{\vu}_{n+1} - \vz_{n+1}, \nabla^2 f(\vx^\star) ( \nabla f(\vx^\star) - \xi_n^{(2)} + \gS_n ) \right\rangle\right.\\
    &\left. - \gamma\sqrt{\eta_n}\EB\left\langle \tilde{\vu}_{n+1} - \vz_{n+1}, \nabla^2 f(\vx^\star) \EB\check{\vv}_{n+1} \right\rangle \right|\\
    \overset{(a)}{\le}& (1 - \gamma\sqrt{\eta_n})\left|\EB\left\langle \tilde{\vu}_{n+1} - \vz_{n+1}, \nabla^2 f(\vx^\star)(\check{\vv}_{n} - \EB\check{\vv}_{n}) \right\rangle\right| + \sqrt{\eta_n}\left|\left\langle \EB\tilde{\vu}_{n+1} - \EB\vz_{n+1}, \gC_n \right\rangle\right|\\
    &+ \sqrt{\eta_n}|\EB\langle \tilde{\vu}_{n+1} - \vz_{n+1}, \xi_n^{(2)} \rangle| + \sqrt{\eta_n}|\EB\langle \tilde{\vu}_{n+1} - \vz_{n+1}, \gS_n \rangle| + o( \sqrt{\eta_n} ) \\
    \overset{(b)}{\le}& (1 - \gamma\sqrt{\eta_n})\left|\EB\left\langle \tilde{\vu}_{n+1} - \vz_{n+1}, \nabla^2 f(\vx^\star)(\check{\vv}_{n} - \EB\check{\vv}_{n}) \right\rangle\right| + \sqrt{\eta_n}\gC |\EB\tilde{\vu}_{n+1} - \EB\vz_{n+1}|\\
    &+ \sqrt{\eta_n}\left( \EB\|\tilde{\vu}_{n+1} - \vz_{n+1}\|^2 \EB \|\gS_n\|^2 \right)^{1/2} + o ( \sqrt{\eta_n} ) \\
    \overset{(c)}{\le}& (1 - \gamma\sqrt{\eta_n})\left|\EB\left\langle \tilde{\vu}_{n+1} - \vz_{n+1}, \nabla^2 f(\vx^\star)(\check{\vv}_{n} - \EB\check{\vv}_{n}) \right\rangle\right| + o(\sqrt{\eta_n}) \\
    \overset{(d)}{\le}& (1 - \gamma\sqrt{\eta_n})\left|\EB\left\langle \tilde{\vu}_{n} - \vz_{n}, \nabla^2 f(\vx^\star)(\check{\vv}_{n} - \EB\check{\vv}_{n}) \right\rangle\right| + o(\sqrt{\eta_n}) + \gO(\eta_n),
\end{aligned}
\end{equation}
where $\gC_n \triangleq \nabla^2 f (\vx^\star) [ (1 - \gamma \sqrt{\eta_n} ) \nabla f (\vx^\star) + \gamma \sE \check{\vu}_{n+1} ] $ is a sequence of constant vector whose norm is bounded by $\gC$. And (a) holds by simply using the triangular inequality and the recursive formula of $\sE \check{\vv}_n$, (b) holds by Cauchy-Schwartz's inequality, (c) holds for the results discussed in the first step and the fact that $\gS_n$ is $o(1)$, and (d) holds for the update rule of $\tilde{\vu}_n - \vz_n$. Up to now, the main result of this step is showed by using Lemma~\ref{lem:converge_r_t_3} again.

\paragraph{Step 3}
Back to the corollary we aim to prove. We tackle it by the following inequality
\begin{equation}
\begin{aligned}
    &\EB\|\tilde{\vu}_{n+1} - \vz_{n+1}\|^2\\
    =& \EB\left\|\left\{\rmI_d - \eta_n \gP_{\mA^\bot}\left(\nabla^2 f(\vx^\star) - \frac{\mathbbm{1}_{\alpha = 1}}{2\eta_0}\rmI_d\right)\gP_{\mA^\bot}\right\}(\tilde{\vu}_n - \vz_n) - \eta_n\gP_{\mA^\bot}\nabla^2 f(\vx^\star)(\check{\vv}_n - \vv)\right\|^2\\
    \overset{(a)}{\le}& \left(1 - \left(\mu - \frac{\mathbbm{1}_{\alpha = 1}}{2\eta_0}\right)\eta_n\right)\EB\|\tilde{\vu}_n - \vz_n\|^2 + 2\eta_n \left|\EB\langle\tilde{\vu}_n - \vz_n, \nabla^2 f(\vx^\star)(\check{\vv}_n - \EB\check{\vv}_n)\rangle\right|\\
    &+ \eta_n^2 L^2 \EB\|\check{\vv}_n - \vv\|^2 + 2\eta_n |\langle \EB\tilde{\vu}_n - \EB\vz_n, \nabla^2 f(\vx^\star) ( \EB\check{\vv}_n - \vv ) \rangle| + o(\eta_n)\\
    \overset{(b)}{\le}& \left(1 - \left(\mu - \frac{\mathbbm{1}_{\alpha = 1}}{2\eta_0}\right)\eta_n\right)\EB\|\tilde{\vu}_n - \vz_n\|^2 + o(\eta_n),
\end{aligned}
\end{equation}
where (a) holds for the positive definiteness of the Hessian matrix $\nabla^2 f(\vx^\star)$ and the choice of $\eta_0$, and (b) follows from the fact that $|\EB\langle \tilde{\vu}_n - \vz_n , \nabla^2 f(\vx^\star) (\check{\vv}_n - \EB\check{\vv}_n) \rangle|\to 0$ as $n\to \infty$ and $\EB\check{\vv}_n \to \vv$ as $n\to \infty$.
Owing to Lemma~\ref{lem:converge_r_t_3}, we can claim that the Wasserstein distance between $\tilde{\vu}_n$ and $\vz_n$ decays to zero as the iteration goes to infinity. What's more, we have known that $\vz_n$ converges weakly to $\gN(0,\tilde{\Sigma})$, so we can say $\tilde{\vu}_n$ converges weakly to $\gN(0,\tilde{\Sigma})$. Recall the definition of $\tilde{\vu}_n (= \check{\vu}_n - \vu)$. Therefore, we finally obtain $\check{\vu}_n$ weakly converges to $\gN(\vu,\tilde{\Sigma})$.
\end{proof}

\section{Proof of Section~\ref{sec:dgnrt_bias&homo_fed}}\label{sec:app:degenerate}
In this section, we give the proof for the results in Section~\ref{sec:dgnrt_bias&homo_fed}.

\subsection{Proof of Lemma~\ref{lem:well_def_dgnt}}
\begin{proof}[Proof of Lemma~\ref{lem:well_def_dgnt}]
    We pick a positive definite matrix $\mS_0 \in \sR^{d\times d}$ and a vector $\vc_0 \in \sR^d$ such that the vectors $\left\{ \vc_0, \mS_0\vc_0, \cdots, \mS_0^{k-1}\vc_0, \mS_0^k\vc_0 \right\}$ are linearly uncorrelated.
    Then we choose a full rank matrix $\mA_0$ such that $\mathrm{span}(\mA_0) \supseteq \mathrm{span}\left\{ \mS_0^i\vc_0 | i = 0, \cdots, k-1\right\}$ and $\mS_0^k\vc_0 \notin \mathrm{span}(\mA_0)$. Consequently, we have $k \le m \le d-1$ with $\mA \in \sR^{d\times m}$.

    Then we have $\gP_{\mA_0^\bot}\mS_0^i \vc_0 = 0, ~~\forall i = 0, \cdots, k-1$ and $\gP_{\mA_0^\bot}\mS_0^k \vc_0 \neq 0$. After getting $\mA_0, \mS_0$ and $\vc_0$, we consider the following equations of $\vx$ and $\vb$,
    \[
    \vc_0 = \mS_0\vx - \vb; ~~ \mA_0\vx = 0.
    \]
    The solution is easy to find. Indeed, by taking $\vx_0\in \mA^\bot$ and letting $\vb_0 = \mS_0\vx_0 - \vc_0$, then we get a solution pair $(\vx_0, \vb_0)$.

    Now, we focus on the linearly constrained optimization problem corresponding to the triplet $(\mA_0; \mS_0; \vb_0)$. Suppose the unique solution of this problem is $\vx_1 \in \sR^d$. The by the KKT condition, $\vx_1$ is supposed to satisfy $\mS_0\vx_1 - \vb_0 + \mA_0 \boldsymbol{\lambda} = 0$ with some $\boldsymbol{\lambda} \in $ and $\mA_0^\top \vx_1 = 0$. The first equality is equivalent to $\gP_{\mA_0^\bot}(\mS_0 \vx_1 - \vb_0) = 0$. While it is true that $\gP_{\mA_0^\bot}(\mS_0 \vx_0 - \vb_0) = \gP_{\mA_0^\bot} \vc_0 = 0$. Hence, $\vx_0 = \vx_1$ by the uniqueness of the solution. Further, the solution of this problem satisfies that $\gP_{\mA_0^\bot}\mS_0^i (\mS_0 \vx_0 - \vb_0) = 0, ~\forall i = 1,\cdots, k-1$ and $\gP_{\mA_0^\bot}\mS_0^k (\mS_0 \vx_0 - \vb_0)\neq 0$. So we have constructed a triplet whose corresponding linearly constrained optimization problem satisfies Assumption~\ref{asp:degenerate}.
\end{proof}

\subsection{Proof of Theorem~\ref{thm:cov_dgnrt_bias}}
For ease of notation, we denote $\gP_{\mA^\bot}\mS$ and $\gP_{\mA}\mS$ as $\gS_\bot$ and $\gS$ respectively.
Then we have the following property, which is useful for the proof of Theorem~\ref{thm:cov_dgnrt_bias}.
\begin{proposition}
\label{prop:degenerate}
    The following conditions are equivalent to Assumption~\ref{asp:degenerate}:
    $\gS_\bot\gS^i (\mS \vx^\star - \vb)= \vzero, ~~\forall i = 0, 1,\dots, k-2$ and $\gS_\bot\gS^{k-1} (\mS\vx^\star - \vb)\neq \vzero$. Moreover, $\gS_\bot\gS^{k-1} (\mS\vx^\star - \vb) = \gP_{\mA^\bot} \mS^k (\mS \vx^\star - \vb) $.
\end{proposition}

\begin{proof}

\textbf{Necessity} 
If $k=1$, the conclusion naturally holds. So we assume $k \ge 2$.
For $i=0$, 
$\gS_\bot (\mS \vx^\star - \vb) = \gP_{\mA^\bot} \mS (\mS \vx^\star - \vb) = \vzero$.
Suppose that for $0 \le i \le j-1$ with $1 \le j \le k-1$, we have verified $\gS_\bot \gS^i (\mS \vx^\star - \vb) - \vzero $.
Since $\gP_\mA = \rmI_d - \gP_{\mA^\bot}$,
we have 
\begin{equation}\label{eq:pr_degenerate_example}
\begin{aligned}
    \gS_\bot \gS^j (\mS \vx^\star - \vb)
    & = \gP_{\mA^\bot} \mS \gP_\mA \mS \gS^{j-1} (\mS \vx^\star - \vb) \\
    & = \gP_{\mA^\bot} \mS (\rmI_d - \gP_{\mA^\bot} ) \mS \gS^{j-1} (\mS \vx^\star - \vb) \\
    & = \gP_{\mA^\bot} \mS^2 \gS^{j-1} (\mS \vx^\star - \vb) - \gS_\bot \gS_\bot \gS^{j-1} (\mS \vx^\star - \vb) \\
    & = \gP_{\mA^\bot} \mS^2 (\rmI_d - \gP_{\mA^\bot}) \mS \gS^{j-2} (\mS \vx^\star - \vb) \\
    & = \dots = \gP_{\mA^\bot} \mS^{j+1} (\mS \vx^\star - \vb).
\end{aligned}
\end{equation}
By Assumption~\ref{asp:degenerate}, if $j < k-1$, we have $\gS_\bot \gS^j (\mS \vx^\star - \vb) = \vzero$; if $j = k-1$, we have $\gS_\bot \gS^j (\mS \vx^\star - \vb) = \vzero$.
The proof is completed by induction.

\textbf{Sufficiency} The sufficiency is a natural consequence of (\ref{eq:pr_degenerate_example}). 
\end{proof}

Now we can give the proof of Theorem~\ref{thm:cov_dgnrt_bias}.
\begin{proof}[Proof of Theorem~\ref{thm:cov_dgnrt_bias}]
Denote $\gS_\bot\gS^{k-i-1} \vv_{n}$ as $\vv_n^{(k-i-1)}$.
Without loss of generality, we assume $p_n = \min \{ \gamma \eta_n^\beta, 1 \} = \gamma \eta_n^\beta$.
For any $i = 1,\dots, k-1$, we can derive the following recursive relationship,
\begin{align*}
    &\EB\left\|\vv^{(k-i-1)}_{n+1}\right\|^2 = (1 - \gamma\eta_n^\beta)\EB\left\| \gS_\bot\gS^{k-i-1}\gP_{\mA}(\vx_n - \eta_n(\mS \vx_n - \vb) + \eta_n\xi_n) \right\|^2 + \gamma\eta_n^\beta\EB \|\mathbf{0}\|^2\\
    &\overset{(a)}{=} (1-\gamma \eta_n^\beta)\EB\left\| \vv_n^{(k-i-1)} - \eta_n\gS_\bot\gS^{k-i-1}\gP_\mA\mS(\vx_n -\vx^\star) + \eta_n \gS_\bot \gS^{k-i-1} \gP_\mA \xi_n \right\|^2\\
    &= (1- \gamma \eta_n^\beta)\EB\left\| \vv_n^{(k-i-1)} - \eta_n\gS_\bot\gS^{k-i-1}\gP_\mA \mS(\vx_n - \vx^\star) + \eta_n\gS_\bot\gS^{k-i-1}\xi^{(2)}_n \right\|^2\\
    &= (1-\gamma\eta_n^\beta)\EB\left\| \vv_n^{(k-i-1)}\right\|^2 - 2\eta_n(1- \gamma\eta_n^\beta)\EB \left\langle \vv_n^{(k-i-1)}, \gS_\bot\gS^{k-i}(\vu_n - \vx^\star) + \vv_n^{(k-i)} \right\rangle\\
    & \quad + \eta_n^2 (1-\gamma \eta_n^\beta)\EB\norm{\gS_\bot\gS^{k-i}(\vx_n - \vx^\star)}^2 + \eta_n^2(1-\gamma\eta_n)\EB\norm{\gS_\bot\gS^{k-i-1}\xi^{(2)}_n}^2\\
    &\le (1-\gamma\eta_n^\beta)\EB\norm{\vv_n^{(k-i-1)}}^2 + \left( \frac{\gamma}{2}\eta_n^\beta\EB\norm{\vv_n^{(k-i-1)}}^2 + \frac{2}{\gamma}\eta_n^{2-\beta}\EB\norm{\gS_\bot\gS^{k-i}(\vu_n - \vx^\star) + \vv_n^{(k-i)}}^2 \right)+  \gO(\eta_n^2)\\
    &\precsim \left(1-\frac{\gamma}{2}\eta_n^\beta\right)\EB\norm{\vv_n^{(k-i-1)}}^2 + \eta_n^{2-\beta}\left[\EB\norm{\vu_n - \vx^\star}^2 + \EB\norm{\vv_n^{(k-i)}}^2\right] + \gO(\eta_n^2).
\end{align*}
Here (a) follows from Assumption~\ref{asp:degenerate} and Proposition~\ref{prop:proj_A_bot_nabla}. On the other hand, we have
\begin{align*}
    &\EB\norm{\vu_{n+1} - \vx^\star}^2 = \EB\norm{(\vu_n - \vx^\star) - \eta_n \gP_{\mA^\bot}\mS( \vx_n - \vx^\star) + \eta_n \xi_n^{(1)}}^2\\
    &= \EB\norm{\vu_n - \vx^\star}^2 - 2\eta_n \left\langle \vu_n -\vx^\star,
    \gP_{\mA^\bot} \mS \gP_{\mA^\bot}
    (\vu_n - \vx^\star) \right\rangle - 2\eta_n\left\langle \vu_n - \vx^\star, \gS_\bot \vv_n \right\rangle + \gO( \eta_n^2 )\\
    &\le (1 - 2\mu \eta_n) \EB \norm{\vu_n - \vx^\star}^2 + \eta_n \left( \mu \EB \norm{\vu_n - \vx^\star}^2 + \mu^{-1} \EB \norm{\vv_n^{(0)}}^2\right) + \gO ( \eta_n^2 )\\
    &\precsim (1 - \mu \eta_n)\EB \norm{\vu_n - \vx^\star}^2 + \eta_n \EB \norm{\vv_n^{(0)}}^2 + \gO( \eta_n^2 ).
\end{align*}

And from Theorem~\ref{thm:converge}, we know that $\EB\norm{\vv^{(k-i-1)}_n}^2 = \gO(\eta_n^{2-2\beta})$ for all $i = 0,1,\dots, k-1$. What's next is to show the following convergence relationships inductively,
\begin{equation}\label{eq:dgnrt_v_conv_rt}
\EB\norm{\vv^{(k-i-1)}_n}^2 \precsim \eta_n^{2(i+1)(1-\beta)} + \eta_n^{2-\beta} ~~\forall i= 0,1,\dots, k-1.
\end{equation}

When $i = 0$, the conclusion is mundane. Suppose for $i =0,\dots, j-1$, the equation~(\ref{eq:dgnrt_v_conv_rt}) holds. Then we sum up over the update inequality of $\EB \norm{\vv_{n}^{(k-i-1)}}^2$ from $i= j$ to $k-1$ and then plus the recursive inequality of $\EB\norm{\vu_n - \vx^\star}^2$ and get
\begin{align*}
    &\left(\ssum{i}{j}{k-1}\EB \norm{\vv_{n+1}^{(k-i-1)}}^2 + \EB \norm{\vu_{n+1} - \vx^\star}^2\right) \precsim
    \left(1 - \frac{\gamma}{2}\eta_n^\beta + \eta_n^{2-\beta}\right)\ssum{i}{j}{k - 2}\EB \norm{\vv_n^{(k-i-1)}}^2 \\
    & \quad + \left(1 - \frac{\gamma}{2}\eta_n^\beta + \eta_n\right)\EB \norm{\vv_n^{(0)}}^2 + \left(1 - \mu \eta_n + (k-j) \eta_n^{2-\beta}\right) \EB \norm{\vu_n - \vx^\star}^2\\
    & \quad + \eta_n^{2-\beta} \EB \norm{\vv_n^{(k - j)}}^2 + \gO (\eta_n^2)\\
    &\overset{(a)}{\precsim} \left(1 - \frac{\gamma}{4} \eta_n^\beta \right) \ssum{i}{j}{k}\EB \norm{\vv_n^{(k-i-1)}}^2 + \left( 1 - \frac{\mu}{2}\eta_n\right)\EB \norm{\vu_n - \vx^\star}^2\\
    & \quad + \eta_n^{2 - \beta}\left( \eta_n^{2j(1-\beta)} + \eta_n^{2-\beta}\right) + \gO (\eta_n^2)\\
    &\precsim \left( 1 - \frac{\mu}{2}\eta_n \right)\left( \ssum{i}{j}{k}\EB \norm{\vv_n^{(k-i-1)}}^2 + \EB \norm{\vu_n - \vx^\star}^2\right) + \eta_n\left(\eta_n^{(2j + 1)(1 - \beta)} + \eta_n\right),
\end{align*}
where $(a)$ holds for the inductive assumption. And the last inequality holds for $\eta_n^{2(2-\beta)} = \eta_n^{2(1-\beta) } \cdot\eta_n^2 \precsim \eta_n^2$. By making use of Lemma~\ref{lem:converge_r_t}, we obtain the temporary bound
\[
\ssum{i}{j}{k}\EB \norm{\vv_n^{(k-i-1)}}^2 + \EB \norm{ \vu_n - \vx^\star}^2 \precsim \eta_n^{(2j+1)(1-\beta)} + \eta_n.
\]
Specifically, we have $\EB\norm{\vu_n - \vx^\star}^2 \precsim \eta_n^{(2j+1)(1 - \beta)} + \eta_n$. Bring this result back to the recursive bound of $\EB \norm{\vv_n^{(k - j-1)}}^2$, and get
\begin{align*}
\EB \norm{\vv_n^{(k - j - 1)}}^2 &\precsim \left( 1 - \frac{\gamma}{2}\eta_n^\beta \right)\EB \norm{\vv_n^{(k-j-1)}}^2 + \eta_n^{2-\beta}\left(\eta_n^{(2j+1)(1 - \beta)} + \eta_n + \eta_n^{2j(1-\beta)}\right) + \gO(\eta_n^2)\\
&\precsim \left(1 - \frac{\gamma}{2}\eta_n^\beta \right)\EB \norm{\vv_n^{(k - j - 1)}}^2 + \eta_n^{ 2 - \beta + 2j(1-\beta)} + \eta_n^2.
\end{align*}
By leveraging Lemma~\ref{lem:converge_r_t} again, we finally have $\EB \norm{\vv_n^{(k - j - 1)}}^2 \precsim \eta_n^{2(j+1)(1 - \beta)} + \eta_n^{2-\beta}$. Up to now, the inductive process is completed. So equation~(\ref{eq:dgnrt_v_conv_rt}) holds for all $i = 0, \dots, k-1$. Therefore, we can show that
\[
\EB \norm{\vu_{n+1} - \vx^\star}^2 \precsim (1 - \mu\eta_n)\EB \norm{\vu_n - \vx^\star}^2 + \eta_n(\eta_n^{2k(1-\beta)} + \eta_n).
\]
And the proof is concluded by using Lemma~\ref{lem:converge_r_t} once more.
\end{proof}

\subsection{Proof of Theorem~\ref{thm:degenerate_bias_u}}
In this subsection, we give the proof of Theorem~\ref{thm:degenerate_bias_u}. 
\def\vhomo{{\check{\vv}}}

\begin{proof}[
{Proof of 
Theorem~\ref{thm:degenerate_bias_u}}]
As claimed in the main text, the proof is divided into three parts.

\textbf{Case 1: $0 < \beta < 1 - \frac{1}{2k}$.} The proof is the same as that of Theorem~\ref{thm:diff_approx} except the range of $\beta$.

\textbf{Case 2: $1 - \frac{1}{2k} < \beta < 1$.} The proof is similar to that
of Theorem~\ref{thm:asy_bias_u}
and is divided into two parts.

With some abuse of notation, define $\uhomo_n = \frac{\vu_n - \vx^\star}{\eta_{n-1}^{ k (1-\beta) } }$, and
$\vhomo_{n}^{ (k-i-1) } = \frac{ 
\gS_\bot \gS^{k-i-1} \vv_n
}{ \eta_{n-1}^{ (i+1) (1-\beta) } }$ for $i = 0, \dots, k-1$. 
First, we show that $\sE \uhomo_n$ convergence to a constant non-zero vector.
Then we prove that the asymptotic variance of $\uhomo_n$ is zero.
Since we focus on the asymptotic convergence of $\uhomo_n$, we assume $p_n = \gamma \eta_n^\beta$ without loss of generality. 


\paragraph{Step 1}
We first derive the recursive relation of $\uhomo_n$.
From the update rule and the linearity of $\gP_{\mA^\bot}$, we have
\begin{align*}
    \uhomo_{n+1}
    & = \gP_{\mA^\bot} \frac{ \vx_{n+1} - \vx^\star }{ \eta_n^{ k (1-\beta) } } \\
    & = \gP_{\mA^\bot} \frac{ \vx_n - \vx^\star - \eta_n (\mS \vx_n - \vb) + \eta_n \xi_n }{ \eta_n^{ k (1 - \beta) } } \\
    & = \left( \frac{ \eta_{n-1} }{ \eta_n } \right)^{ k (1-\beta) } \uhomo_n - \eta_n^{ 1 - k (1-\beta) } \gP_{\mA^\bot} (\mS \vx_n - \vb) + \eta_n^{ 1 - k (1-\beta) } \xi_n^{(1)},
\end{align*}
where $\xi_n^{(1)} = \gP_{\mA^\bot} \xi_n$.
For the second term, we have
\begin{align*}
    \gP_{\mA^\bot} (\mS \vx_n - \vb)
    & = \gP_{\mA^\bot} ( \mS (\vu_n - \vx^\star) + \mS \vv_n + \mS \vx^\star - \vb ) \\
    & \overset{(a)}{=} \gP_{\mA^\bot} \mS \gP_{\mA^\bot} (\vu_n - \vx^\star) + \gS_\bot \vv_n + \gS_\bot (\mS \vx^\star - \vb ) \\
    & \overset{(b)}{=} \eta_{n-1}^{ k (1-\beta) } \gP_{\mA^\bot} \mS \gP_{\mA^\bot}  \uhomo_{n} + \eta_{n-1}^{ k (1-\beta) } \vhomo_{n}^{(0)},
\end{align*}
where (a) is because $\vu_n$ and $\vx_\star$ lie in the null space of $\mA^\top$ and (b) is by Assumption~\ref{asp:degenerate}.
This implies 
\begin{align*}
    \uhomo_{n+1}
    & = \left( \frac{ \eta_{n-1} }{ \eta_n } \right)^{ k (1-\beta) } \uhomo_n - \eta_n \left( \frac{ \eta_{n-1} }{ \eta_n } \right)^{ k (1-\beta) } \gP_{\mA^\bot} \mS \gP_{\mA^\bot}  \uhomo_{n} \\
    & \quad \ -  \eta_n \left( \frac{ \eta_{n-1} }{ \eta_n } \right)^{ k (1-\beta) } \vhomo_n^{(0)} + \eta_n^{ 1 - k (1-\beta) } \xi_n^{(1)}.
\end{align*}
Similar to the analysis below Equation~(\ref{eq:Ru_def}), when $\alpha < 1$, we have $ 1 - \left( \frac{ \eta_{n-1} }{ \eta_n }\right) ^{ k (1 - \beta) }  = o (\eta_n) $ and when $\alpha = 1$, we have $ 1 - \left( \frac{ \eta_{n-1} }{ \eta_n } 
\right)^{ k (1 - \beta) } = \gO (\eta_n) $ and $ 1 + \frac{ k (1-\beta) }{n} - \left( \frac{ \eta_{n-1} }{ \eta_n } \right)^{ k (1 - \beta) } = o(\eta_n) $.
As a result,
\begin{equation}\label{eq:pf_degenerate_u_n}
\begin{aligned}
    \uhomo_{n+1}
    & = \left[ \rmI_d - \eta_n \gP_{\mA^\bot} \left( \mS - \frac{ k (1-\beta) }{ \eta_0 } \mathbbm{1}_{ \{ \alpha=1 \} } \rmI_d  \right) \gP_{\mA^\bot}  \right] \uhomo_n - \eta_n \vhomo_n^{(0)} 
    \\
    & \quad \ 
    + \eta_n^{ 1 - k (1-\beta) } \xi_n^{(1)} + o(\eta_n) \uhomo_n + o(\eta_n) \vhomo_n^{(0)}. 
\end{aligned}
\end{equation}
Here we emphasize that the $o(\cdot), \gO(\cdot), \Theta(\cdot)$ notation can be matrix-valued in terms of the spectral norm  throughout this proof.
Theorem~\ref{thm:cov_dgnrt_bias} guarantees $\norm{\sE \uhomo_{n} } + \norm{ \sE \vhomo_{n}^{(0)} } \le \sqrt{ \sE \norm{ \uhomo_n }^2 } + \sqrt{ \sE \norm{ \vhomo_n^{(0)} }^2 } =\gO(1)$ for $\beta > 1 - \frac{1}{2k}$. Then taking expectation yields
\begin{align}\label{eq:pf_degenerate_u_n_expect}
    \sE \uhomo_{n+1} = \left[ \rmI_d - \eta_n \gP_{\mA^\bot} \left( \mS - \frac{ k (1-\beta) }{ \eta_0 } \mathbbm{1}_{ \{ \alpha=1 \} } \rmI_d  \right) \gP_{\mA^\bot}  \right] \sE \uhomo_n - \eta_n \sE \vhomo_n^{(0)} + o(\eta_n),
\end{align}
which implies the limit of $\sE \uhomo_n$ depends on that of $\sE \vhomo_{n}^{(0)}$.
Thus, we also need to derive the recursive relation of $\vhomo_{n}^{(0)}$, or more generally, $\vhomo_{n}^{(k-i-1)}$.
For $i = 1, \dots, k-1$, we have
\begin{align*}
    \vhomo_{n+1}^{(k-i-1)} 
    & = (1 - \gamma \eta_n^\beta) \gS_\bot \gS^{k-i-1} \mP_\mA \frac{ \vx_n - \eta_n (\mS \vx_n - \vb) + \eta_n \xi_n }{ \eta_n^{ (i+1) (1-\beta) } } \\
    & = (1 - \gamma \eta_n^\beta ) \left( \frac{ \eta_{n-1} }{ \eta_n } \right)^{ (i+1) (1 - \beta) } \vhomo_{n}^{(k-i-1)} - (1 - \gamma \eta_n^\beta) \eta_n^{1 - (i+1) (1-\beta) } \gS_\bot \gS^{k-i-1} \mP_\mA (\mS \vx_n - \vb) \\
    & \quad \ + (1 - \gamma \eta_n^\beta ) \eta_n^{ 1 - (i+1) (1-\beta) } \gS_\bot \gS^{k-i-1} \xi_n^{(2)},
\end{align*}
where $\xi_n^{(2)} = \gP_\mA \xi_n$.
For the second term, we have
\begin{align*}
    \gS_\bot \gS^{k-i-1} \gP_\mA (\mS \vx_n - \vb) 
    & = \gS_\bot \gS^{k-i-1} \gP_\mA ( \mS (\vu_n - \vx^\star) + \mS \vv_n + \mS \vx^\star - \vb ) \\
    & = \gS_\bot \gS^{k-i} (\vu_n - \vx^\star) + \gS_\bot \gS^{k-i} \vv_n+ \gS_\bot \gS^{k-i-1} \gP_{\mA} (\mS \vx^\star - \vb) \\
    & = \eta_{n-1}^{ k(1 - \beta) } \gS_\bot \gS^{k-i} \uhomo_n + \eta_{n-1}^{ i (1 - \beta) } \vhomo_n^{(k-i)},
\end{align*}
where the last equality is by Proposition~\ref{prop:proj_A_bot_nabla} and Assumption~\ref{asp:degenerate}.
This implies
\begin{equation*}
\begin{aligned}
     \vhomo_{n+1}^{(k-i-1)} 
     & = (1 - \gamma \eta_n^\beta ) \left( \frac{ \eta_{n-1} }{ \eta_n } \right)^{ (i+1) (1 - \beta) } \vhomo_{n}^{(k-i-1)} - (1 - \gamma \eta_n^\beta ) \eta_n^{\beta} \left( \frac{ \eta_{n-1} }{ \eta_n } \right)^{ i (1-\beta) } \vhomo_{n}^{ (k-i) } \\
     & \quad \ - (1 - \gamma \eta_n^\beta ) \eta_n^{ 1 + (k-i-1) (1 - \beta) } \left( \frac{ \eta_{n-1} }{ \eta_n } \right)^{ k (1-\beta) } \gS_\bot \gS^{(k-i)}  \uhomo_{n}  \\
     & \quad \ + (1 - \gamma \eta_n^\beta ) \eta_n^{ 1 - (i+1) (1-\beta) } \gS_\bot \gS^{k-i-1} \xi_n^{(2)}.
\end{aligned}
\end{equation*}
Similar to the analysis below Equation~(\ref{eq:Ru_def}), we have $ 1 - \left( \frac{ \eta_{n-1} }{ \eta_n }\right) ^{ j (1 - \beta) }  = \gO (\eta_n) = o(\eta_n^\beta) $ for $j \ge 0$. As a result,
\begin{equation}\label{eq:pf_degenerate_v_k-i}
\begin{aligned}
     \vhomo_{n+1}^{(k-i-1)} 
     & = (1 - \gamma \eta_n^\beta )  \vhomo_{n}^{(k-i-1)} - \eta_n^{\beta}  \vhomo_{n}^{ (k-i) } + (1 - \gamma \eta_n^\beta 
     ) \eta_n^{ 1 - (i+1) (1-\beta) } \gS_\bot \gS^{k-i-1} \xi_n^{(2)} \\
     & \quad \ + o(\eta_n^\beta) \uhomo_n + o(\eta_n^\beta) \vhomo_n^{(k-i-1)} + o(\eta_n^\beta) \vhomo_n^{(k-i)} .
\end{aligned}
\end{equation}
Theorem~\ref{thm:cov_dgnrt_bias} guarantees $\norm{\sE \uhomo_n } +\norm{ \sE \vhomo_{n}^{ (k-j) } } = \gO(1) $ for $0 \le j \le k$.
Then taking expectation yields
\begin{align}\label{eq:pf_degenerate_v_k-i_expect}
     \sE \vhomo_{n+1}^{(k-i-1)} 
     & = (1 - \gamma \eta_n^\beta )  \sE \vhomo_{n}^{(k-i-1)} - \eta_n^{\beta}  \sE \vhomo_{n}^{ (k-i) } + o(\eta_n^\beta),
\end{align}
which implies when $1 \le i \le k-1 $, the limit of $\sE \vhomo_n^{(k-i-1)}$ is related to that of $\sE \vhomo_n^{(k-i)}$.
For the sequence $\{ \vhomo_n^{(k-1)} \}$, we have
\begin{align*}
    \vhomo_{n+1}^{(k-1)} 
    & = (1 - \gamma \eta_n^\beta) \gS_\bot \gS^{k-1} \mP_\mA \frac{ \vx_n - \eta_n (\mS \vx_n - \vb) + \eta_n \xi_n }{ \eta_n^{  1-\beta } } \\
    & = (1 - \gamma \eta_n^\beta ) \left( \frac{ \eta_{n-1} }{ \eta_n } \right)^{ 1 - \beta } \vhomo_{n}^{(k-i-1)} - (1 - \gamma \eta_n^\beta) \eta_n^{ \beta } \gS_\bot \gS^{k-1} \mP_\mA (\mS \vx_n - \vb) 
    \\
    & \quad \ 
    + (1 - \gamma \eta_n^\beta ) \eta_n^{ \beta } \gS_\bot \gS^{k-1}  \xi_n^{(2)}.
\end{align*}
For the second term, we have
\begin{align*}
    \gS_\bot \gS^{k-1} \gP_\mA (\mS \vx_n - \vb) 
    & = \gS_\bot \gS^{k-1} \gP_\mA ( \mS (\vu_n - \vx^\star) + \mS \vv_n + \mS \vx^\star - \vb ) \\
    & = \gS_\bot \gS^{k} (\vu_n - \vx^\star) + \gS_\bot \gS^{k} \vv_n + \gS_\bot \gS^{k-1} (\mS \vx^\star - \vb),
\end{align*}
where Proposition~\ref{prop:proj_A_bot_nabla} leads to the fact $\gP_\mA (\mS \vx^\star - \vb) = \mS \vx^\star - \vb $.
Then we have
\begin{align*}
    \vhomo_{n+1}^{(k-1)} 
    & = (1 - \gamma \eta_n^\beta ) \left( \frac{ \eta_{n-1} }{ \eta_n } \right)^{ 1 - \beta } \vhomo_{n}^{(k-1)} 
    - (1 - \gamma \eta_n^\beta) \eta_n^{ \beta 
    } \gS_\bot \gS^{k} (\vu_n - \vx^\star)
    \\
    & \quad \ 
    - (1 - \gamma \eta_n^\beta) \eta_n^{ \beta } 
    \gS_\bot \gS^{k} \vv_n
    - (1 - \gamma \eta_n^\beta) \eta_n^{ \beta } 
     \gS_\bot \gS^{k-1} (\mS \vx^\star - \vb) 
    + (1 - \gamma \eta_n^\beta ) \eta_n^{ \beta } \gS_\bot \gS^{k-1} \xi_n^{(2)}.
\end{align*}
Recall that $ 1 - \left( \frac{ \eta_{n-1} }{ \eta_n }\right) ^{ j (1 - \beta) } = o(\eta_n^\beta) $.
As a result,
\begin{align*}
    \vhomo_{n+1}^{(k-1)} 
    & = (1 - \gamma \eta_n^\beta )  \vhomo_{n}^{(k-1)} 
    - \eta_n^{ \beta } 
     \gS_\bot \gS^{k-1} (\mS \vx^\star - \vb) 
    + (1 - \gamma \eta_n^\beta 
    ) \eta_n^{ \beta } \gS_\bot \gS^{k-1} \xi_n^{(2)} \\
    & \quad \ + o(\eta_n^\beta) \vhomo_n^{(k-1)} + \gO(\eta_n^\beta) (\vu_n - \vx^\star) +  \gO(\eta_n^\beta) \vv_n + o(\eta_n^\beta).
\end{align*}
Theorems~\ref{thm:cov_dgnrt_bias} and \ref{thm:converge} guarantees $\norm{\sE \vhomo_n^{(k-1)} } = \gO(1)$ and $\sE \norm{ \vu_n - \vx^\star } + \sE \norm{ \vv_n } = o(1)$.
Then taking expectation yields
\begin{align*}
    \sE \vhomo_{n+1}^{(k-1)} 
    & = (1 - \gamma \eta_n^\beta )  \sE \vhomo_{n}^{(k-1)} 
    - \eta_n^{ \beta } 
     \gS_\bot \gS^{k-1} (\mS \vx^\star - \vb) + o(\eta_n^\beta).
\end{align*}
Adding both sides by $ \gamma^{-1} \gS_\bot \gS^{k-1} (\mS \vx^\star - \vb)  $ and taking the norm, we obtain
\begin{align*}
    \norm{ \sE \vhomo_{n+1}^{(k-1)} + \gamma^{-1} \gS_\bot \gS^{k-1} (\mS \vx^\star - \vb) } 
    \le (1 - \gamma \eta_n^\beta) \norm{ \sE \vhomo_{n}^{(k-1)} + \gamma^{-1} \gS_\bot \gS^{k-1} (\mS \vx^\star - \vb) } + o(\eta_n^\beta).
\end{align*}
By Lemma~\ref{lem:converge_r_t_3}, we have $\lim_{n \rightarrow \infty} \sE \vhomo_{n}^{(k-1)} = - \gamma^{-1} \gS_\bot \gS^{k-1} (\mS \vx^\star - \vb)$. 
Similarly, given that $\sE \vhomo_n^{(k-i)} $ converges,  adding both sides of (\ref{eq:pf_degenerate_v_k-i_expect}) by $\gamma^{-1} \lim_{n \rightarrow \infty} \sE \vhomo_n^{ (k-i) } $ and taking the norm yield
\begin{align*}
    \norm{ \sE \vhomo_{n+1}^{ (k-i-1) } + \gamma^{-1} \lim_{n \rightarrow \infty} \sE \vhomo_n^{ (k-i) } } 
    \le (1 - \gamma \eta_n^\beta) \norm{ \sE \vhomo_{n}^{ (k-i-1) } + \gamma^{-1} \lim_{n \rightarrow \infty} \sE \vhomo_n^{ (k-i) } } + o(\eta_n^\beta).
\end{align*}
Then by Lemma~\ref{lem:converge_r_t_3}, we have $\lim_{n \rightarrow \infty } \sE \vhomo_n^{(k-i-1)} = - \gamma^{-1} \lim_{n \rightarrow \infty } \sE \vhomo_n^{(k-i)} $ for $1 \le i \le k-1$, which implies $ \lim_{n \rightarrow \infty} \sE \vhomo_n^{(0)} = (-1)^{k} \gamma^{-k} \gS_\bot \gS^{k-1} (\mS \vx^\star - \vb) $.

Now we return to the sequence $\{ \sE \uhomo_n \}$ and derive its limit.
Since $\eta_0 > 2 k \mathbbm{1}_{ \{ \alpha = 1 \} } / \lambda_{\mathrm{min}} (\mS) $,
the matrix $\gP_{\mA^\bot} \left( \mS - \frac{ k (1-\beta) }{ \eta_0 } \mathbbm{1}_{ \{ \alpha=1 \} } \rmI_d  \right) \gP_{\mA^\bot} $ is positive definite on the null space of $\mA^\top$ where $\uhomo$ lies in. 
Define $\vmu := - \left[ \gP_{\mA^\bot} \left( \mS - \frac{ k (1-\beta) }{ \eta_0 } \mathbbm{1}_{ \{ \alpha=1 \} } \rmI_d  \right) \gP_{\mA^\bot} \right]^\dag \lim_{n \rightarrow \infty} \sE \vhomo_n^{(0)}$.
Subtracting both sides of (\ref{eq:pf_degenerate_u_n_expect}) by 
$\vmu$
and taking the norm yield
\begin{align*}
    \norm{ \sE \uhomo_{n+1} - \vmu } 
    & \le \norm{ \left[ \rmI_d - \eta_n \gP_{\mA^\bot} \left( \mS - \frac{ k (1-\beta) }{ \eta_0 } \mathbbm{1}_{ \{ \alpha=1 \} } \rmI_d  \right) \gP_{\mA^\bot} \right] ( \sE \uhomo_{n} - \vmu ) } + o (\eta_n) \\
    & \le (1 - \Theta(\eta_n) ) \norm{ \sE \uhomo_n - \vmu } + o(\eta_n).
\end{align*}
Combining this with Lemma~\ref{lem:converge_r_t_3}, we obtain that 
\begin{equation}\label{eq:pf_degenerate_mu}
\begin{aligned}
\lim_{n \rightarrow \infty} \sE \uhomo_{n} = \vmu & = 
\frac{ (-1)^{k+1} }{ \gamma^{k} }
\left[ \gP_{\mA^\bot} \left( \mS - \frac{ k (1{-}\beta) }{ \eta_0 } \mathbbm{1}_{ \{ \alpha=1 \} } \rmI_d  \right) \gP_{\mA^\bot} \right]^\dag \gS_\bot \gS^{k-1} (\mS \vx^\star - \vb) \\
& = \frac{ (-1)^{k+1} }{ \gamma^{k} }
\left[ \gP_{\mA^\bot} \left( \mS - \frac{ k (1{-}\beta) }{ \eta_0 } \mathbbm{1}_{ \{ \alpha=1 \} } \rmI_d  \right) \gP_{\mA^\bot} \right]^\dag
\gP_{\mA^\bot} \mS^k 
(\mS \vx^\star - \vb),
\end{aligned}
\end{equation}
where the last equality is from the proof of Proposition~\ref{prop:degenerate}.


\paragraph{Step 2}

Before calculating the asymptotic variance of $\uhomo_n$, we first consider the inner product $ \left| \sE \inner{ \uhomo_n - \sE \uhomo_n }{ \vhomo_n^{(0)} - \sE \vhomo_n^{(0)} } \right| $.
For convenience, we define the right-hand side of (\ref{eq:pf_degenerate_v_k-i}) divided by $1 - \gamma \eta_n^\beta$ as $\vhomo_{( n+1) -}^{(k-i-1)} $. That is, we have $\vhomo_{n+1}^{(k-i-1)} = \vhomo_{(n+1)-}^{(k-i-1)}$ with probability $1 - \gamma \eta_n^\beta$ and $\vhomo_{n+1}^{(k-i-1)} = \vzero$ with probability $\gamma \eta_n^\beta$. This implies
\begin{align*}
    & \quad \ \left| \sE \inner{ \uhomo_{n+1} - \sE \uhomo_{n+1} }{ \vhomo_{n+1}^{(0)} - \sE \vhomo_{n+1}^{(0)} } \right| \\
    & = (1 - \gamma \eta_n^\beta) \left| \sE \inner{ \uhomo_{n+1} - \sE \uhomo_{n+1} }{ \vhomo_{ (n+1)- }^{(0)} - \sE \vhomo_{n+1}^{(0)} } \right| + \gamma \eta_n^\beta \left| \sE \inner{ \uhomo_{n+1} - \sE \uhomo_{n+1} }{\sE \vhomo_{n+1}^{(0)} } \right| \\
    & =  (1 - \gamma \eta_n^\beta) \left| \sE \inner{ \uhomo_{n+1} - \sE \uhomo_{n+1} }{ \vhomo_{ (n+1)- }^{(0)} - \sE \vhomo_{(n+1)-}^{(0)} + \sE \vhomo_{(n+1)-}^{(0)}  - \sE \vhomo_{n+1}^{(0)} } \right| \\
    & = (1 - \gamma \eta_n^\beta) \left| \sE \inner{ \uhomo_{n+1} - \sE \uhomo_{n+1} }{ \vhomo_{ (n+1)- }^{(0)} - \sE \vhomo_{(n+1)-}^{(0)} } \right|.
\end{align*}
Thus, it suffices to consider the inner product $ \left| \sE \inner{ \uhomo_{n+1} - \sE \uhomo_{n+1} }{ \vhomo_{ (n+1)- }^{(0)} - \sE \vhomo_{(n+1)-}^{(0)} } \right| $.
Theorem~\ref{thm:cov_dgnrt_bias} and Cauchy's inequality imply that 
$\sE \norm{ \uhomo_n - \sE \uhomo_n }^2 + \sE \norm{ \vhomo_n^{(0)} - \sE \vhomo_n^{(0)} }^2 \le \sE \norm{ \uhomo_n }^2 + \sE \norm{ \vhomo_n^{(0)} }^2 = \gO(1)$ and
$ \left| \sE \inner{ \uhomo_n - \sE \uhomo_n }{ \vhomo_n^{(j)} - \sE \vhomo_n^{(j)} } \right|
\le \sE \norm{ \uhomo_n - \sE \uhomo_n } \norm{ \vhomo_n^{(j)} - \sE \vhomo_n^{(j)} } 
\le \frac{1}{2} \sE \norm{ \uhomo_n - \sE \uhomo_n }^2 + \frac{1}{2} \sE \norm{ \vhomo_n^{(j)} - \sE \vhomo_n^{(j)} }^2 = \gO(1) $ for $0 \le j \le k-1$.
Then from the definition of $\vhomo_{ (n+1)- }^{(0)} $ and Equation (\ref{eq:pf_degenerate_u_n}), we have the following relation
\begin{align*}
    & \quad \ \left| \sE \inner{ \uhomo_{n+1} - \sE \uhomo_{n+1} }{ \vhomo_{ (n+1)- }^{(0)} - \sE \vhomo_{(n+1)-}^{(0)} } \right| \\
    & = \bigg| \sE \bigg\langle \left[ \rmI_d - \eta_n \gP_{\mA^\bot} \left( \mS - \frac{ k (1-\beta) }{ \eta_0 } \mathbbm{1}_{ \{ \alpha=1 \} } \rmI_d  \right) \gP_{\mA^\bot} + o(\eta_n) \right] (\uhomo_n - \sE \uhomo_n) \\
    & \qquad \  - 
    \eta_n (\vhomo_n^{(0)} - \sE \vhomo_n^{(0)} ) + \eta_n^{1 - k (1-\beta) } \xi_n^{(1)}
    + o(\eta_n) ( \vhomo_n^{(0)} - \sE \vhomo_n^{(0)} ) , \\
    & \qquad \ 
    (1 - \gamma \eta_n^\beta 
    ) ( \vhomo_n^{(0)} - \sE \vhomo_n^{(0)} ) 
    - 
    \eta_n^\beta ( \vhomo_n^{(1)} - \sE \vhomo_n^{(1)} )
    + (1 - \gamma \eta_n^\beta ) \eta_n^{ 1 - k (1-\beta) } \gS_\bot \xi_n^{(2)} \\
    & \qquad \ + o(\eta_n^\beta) ( \uhomo_n - \sE \uhomo_n + \vhomo_n^{(0)} - \sE \vhomo_n^{(0)} + \vhomo_n^{(1)} - \sE \vhomo_n^{(1)} )
    \bigg \rangle \bigg| \\
    & \le \left| \sE \inner{ \uhomo_n - \sE \uhomo_n }{ \vhomo_n^{(0)} - \sE \vhomo_n^{(0)} } \right|
    + \eta_n^\beta \left| \inner{ \uhomo_n - \sE \uhomo_n }{ \vhomo_n^{(1)} - \sE \vhomo_n^{(1)} } \right| + o(\eta_n^\beta).
\end{align*}
Substituting yields
\begin{align*}
    & \quad \ \left| \sE \inner{ \uhomo_{n+1} - \sE \uhomo_{n+1} }{ \vhomo_{n+1}^{(0)} - \sE \vhomo_{n+1}^{(0)} } \right| \\
    & \le (1 - \gamma \eta_n^\beta) \left| \sE \inner{ \uhomo_{n} - \sE \uhomo_{n} }{ \vhomo_{n}^{(0)} - \sE \vhomo_{n}^{(0)} } \right|
    + \eta_n^\beta \left| \inner{ \uhomo_n - \sE \uhomo_n }{ \vhomo_n^{(1)} - \sE \vhomo_n^{(1)} } \right| + o(\eta_n^\beta),
\end{align*}
which implies the limit of $\left| \sE \inner{ \uhomo_{n} - \sE \uhomo_{n} }{ \vhomo_{n}^{(0)} - \sE \vhomo_{n}^{(0)} } \right|$ depends on that of $\left| \sE \inner{ \uhomo_n - \sE \uhomo_n }{ \vhomo_n^{(1)} - \sE \vhomo_n^{(1)} } \right|$.
Following a similar procedure, we can also obtain
\begin{align*}
    & \quad \ \left| \sE \inner{ \uhomo_{n+1} - \sE \uhomo_{n+1} }{ \vhomo_{n+1}^{(k-i-1)} - \sE \vhomo_{n+1}^{(k-i-1)} } \right| \\
    & \le (1 {-} \gamma \eta_n^\beta) \left| \sE \inner{ \uhomo_{n} {-} \sE \uhomo_{n} }{ \vhomo_{n}^{(k-i-1)} {-} \sE \vhomo_{n}^{(k-i-1)} } \right|
    {+} \eta_n^\beta \left| \sE \inner{ \uhomo_n {-} \sE \uhomo_n }{ \vhomo_n^{(k-i)} {-} \sE \vhomo_n^{(k-i)} } \right| {+} o(\eta_n^\beta)
\end{align*}
for $1 \le i \le k$, and
\begin{align*}
    & \quad \ \left| \sE \inner{ \uhomo_{n+1} - \sE \uhomo_{n+1} }{ \vhomo_{n+1}^{(k-1)} - \sE \vhomo_{n+1}^{(k-1)} } \right| \\
    & \le (1 {-} \gamma \eta_n^\beta) \left| \sE \inner{ \uhomo_{n} {-} \sE \uhomo_{n} }{ \vhomo_{n}^{(k-1)} {-} \sE \vhomo_{n}^{(k-1)} } \right|
    {+} \eta_n^\beta \left| \sE \inner{ \uhomo_n {-} \sE \uhomo_n }{ \gS_\bot \gS^{k} (\vv_n {-} \sE \vv_n) } \right| {+} o(\eta_n^\beta).
\end{align*}
By Theorems~\ref{thm:cov_dgnrt_bias} and \ref{thm:converge} as well as Cauchy's inequality, we have
\[\left| \sE \inner{ \uhomo_n {-} \sE \uhomo_n }{ \gS_\bot \gS^{k} (\vv_n {-} \sE \vv_n) } \right| 
\preceq \sE \norm{ \uhomo_n - \sE \uhomo_n } \norm{ \vv_n - \sE \vv_n } 
\le \sqrt{ \sE \norm{\uhomo_n}^2  \sE \norm{ \vv_n }^2 }
= o(1).\] Then Lemma~\ref{lem:converge_r_t_3} implies 
$\left| \sE \inner{ \uhomo_{n} - \sE \uhomo_{n} }{ \vhomo_{n}^{(k-1)} - \sE \vhomo_{n}^{(k-1)} } \right| = o(1)$.
As a consequence, we have
$\left| \sE \inner{ \uhomo_{n} - \sE \uhomo_{n} }{ \vhomo_{n}^{(k-i-1)} - \sE \vhomo_{n}^{(k-i-1)} } \right| = o(1) $ for $1 \le i \le k-1$.
Then setting $i=k-1$ yields $ \left| \sE \inner{ \uhomo_n - \sE \uhomo_n }{ \vhomo_n^{(0)} - \sE \vhomo_n^{(0)} } \right| = o(1) $.

Now we can derive the recursive relation of $\sE \norm{ \uhomo_n - \sE \uhomo_n }^2$.
Recall that we have
$\sE \norm{ \uhomo_n - \sE \uhomo_n }^2 + \sE \norm{ \vhomo_n^{(0)} - \sE \vhomo_n^{(0)} }^2 
= \gO(1)$.
By (\ref{eq:pf_degenerate_u_n}), we have
\begin{align*}
    & \quad \ \sE \norm{ \uhomo_{n+1} - \sE \uhomo_{n+1} }^2 \\
    & = \sE \bigg\| \left[ \rmI_d - \eta_n \gP_{\mA^\bot} \left( \mS - \frac{ k (1-\beta) }{ \eta_0 } \mathbbm{1}_{ \{ \alpha=1 \} } \rmI_d  \right) \gP_{\mA^\bot} + o(\eta_n) \right] (\uhomo_n - \sE \uhomo_n ) \\
    & \qquad \
    - 
    \eta_n ( \vhomo_n^{(0)} - \sE \vhomo_n^{(0)} ) + \eta_n^{ 1 - k (1 - \beta) } \xi_n^{(1)}
    + o(\eta_n) ( \vhomo_n^{(0)} - \sE \vhomo_n^{(0)} ) \bigg\|^2 \\
    & \le \left( 1 - \Theta(\eta_n) \right) \sE \norm{ \uhomo_n - \sE \uhomo_n }^2 + \eta_n \left| \sE \inner{ \uhomo_n - \sE \uhomo_n }{ \vhomo_n^{(0)}  - \sE \vhomo_n^{(0)} } \right| + o(\eta_n) \\
    & = \left( 1 - \Theta(\eta_n) \right) \sE \norm{ \uhomo_n - \sE \uhomo_n }^2 + o(\eta_n).
\end{align*}
Combining this with Lemma~\ref{lem:converge_r_t_3}, we have $\sE \norm{ \uhomo_n - \sE \uhomo_n }^2 = o(1)$.
With $\vmu$ defined in (\ref{eq:pf_degenerate_mu}),
we obtain that $\sE \norm{ \uhomo_n - \vmu }^2 = \sE \norm{ \uhomo_n - \sE \uhomo_n }^2 + \norm{ \sE \uhomo_n - \vmu }^2 \rightarrow 0$ as $n \rightarrow \infty$.

\textbf{Case 3: $\beta = 1 - \frac{1}{2k}$.} The result can be obtained by combining the proof procedure of Theorem~\ref{thm:noncen_u} and the intermediate results in Case 2.
\end{proof}

\section{Proof of Section~\ref{sec:dlpsa}}
\subsection{Preliminary}
Foremost, we write down the iterative forms of the primal sequence $\{\vu_n\}$ and the auxiliary sequence $\{\vv_n\}$. For simplicity, we shorthand $\nabla f(\vx_n, \zeta_n^{(a)}) - \nabla f(\vx_n)$ 
and 
$\nabla f\left( \vx_n + \gamma^{-1}\eta_n^{1-\beta}\nabla f(\vx_n, \zeta_n^{(a)}) \right) - \nabla f\left( \vx_n + \gamma^{-1}\eta_n^{1-\beta}\nabla f(\vx_n, \zeta_n^{(a)}), \zeta_n^{(b)} \right)$ as $\xi_n^{(a)}$ and $\xi_n^{(b)}$ respectively, and we use $C$ to denote all universal constants that are independent of the parameters involved in the Assumptions.  Then, we have
\begin{equation}
\begin{aligned}
    &\vu_n - \vx^\star = \gP_{\mA^\bot}\left( \vx_{n+1} - \vx^\star\right)\\
    &= \left( \gP_{\mA^\bot}\vx_n - \vx^\star\right) - \eta_n \gP_{\mA^\bot}\nabla f\left( \vx_n + \gamma^{-1}\eta_n^{1-\beta}\nabla f(\vx_n) + \gamma^{-1}\eta_n^{1-\beta}\xi_n^{(a)} \right) + \eta_n\gP_{\mA^\bot}\xi_n^{(b)}\\
    & = (\vu_n - \vx^\star) - \eta_n \gP_{\mA^\bot}\left\{ \gQ_n^{(1)} + \gQ_n^{(2)} + \gQ_n^{(3)} + \gQ_n^{(4)} \right\} + \eta_n \gP_{\mA^\bot} \xi_n^{(b)},
\end{aligned}
\end{equation}
and
\begin{equation}
\begin{aligned}
    &\vv_{n+1} - \gamma^{-1}\eta_n^{1-\beta}\nabla f(\vx^\star)\\
    =& 
    \left\{
    \begin{aligned}
        &\left(\vv_n - \gamma^{-1}\eta_n^{1-\beta}\nabla f(\vx^\star)\right) - \eta_n\nabla f(\vx^\star) + \gamma^{-1}\gO\left(\frac{\eta_n^{1-\beta}}{n}\right)\nabla f(\vx^\star)\\
        & -\eta_n\gP_\mA\left\{\gQ_n^{(1)} + \gQ_n^{(2)} + \gQ_n^{(3)} + \gQ_n^{(4)}\right\}+ \eta_n\gP_\mA \xi_n^{(b)}; &&  1 - p_n \\
        & - \gamma^{-1}\eta_{n+1}^{1-\beta}\nabla f(\vx^\star); &&  p_n,
    \end{aligned}
    \right.
\end{aligned}
\end{equation}
where
\begin{equation}
\begin{aligned}
    &\gQ_n^{(1)}&&= \nabla f\left( \vu_n + 
    (\vv_n + \gamma^{-1}\eta_n^{1-\beta}\nabla f(\vx^\star))
    - \gamma^{-1}\eta_n^{1-\beta}(\nabla f(\vx^\star) - \nabla f(\vx_n)) + \gamma^{-1}\eta_n^{1-\beta}\xi_n^{(a)}\right)\\
    & &&- \nabla f\left( \vu_n + 
    (\vv_n + \gamma^{-1}\eta_n^{1-\beta}\nabla f(\vx^\star))
    - \gamma^{-1}\eta_n^{1-\beta}(\nabla f(\vx^\star) - \nabla f(\vx_n))\right)\\
    &\gQ_n^{(2)} &&= \nabla f\left( \vu_n + 
    (\vv_n + \gamma^{-1}\eta_n^{1-\beta}\nabla f(\vx^\star))
    - \gamma^{-1}\eta_n^{1-\beta}(\nabla f(\vx^\star) - \nabla f(\vx_n))\right)\\
    & &&- \nabla f\left( \vu_n + 
    (\vv_n + \gamma^{-1}\eta_n^{1-\beta}\nabla f(\vx^\star))\right)\\
    &\gQ_n^{(3)} && = \nabla f\left( \vu_n + 
    (\vv_n + \gamma^{-1}\eta_n^{1-\beta}\nabla f(\vx^\star))\right) - \nabla f\left( \vu_n\right)\\
    &\gQ_n^{(4)} &&= \nabla f(\vu_n) - \nabla f(\vx^\star).
\end{aligned}
\end{equation}
Before entering the next step of proof, we give these four terms, i.e., $\gQ_n^{(i)}$ with $i=1\cdots, 4$ more specific description.

For $\gQ_n^{(1)}$, by mean value theorems for definite integrals,
\begin{equation}
\begin{aligned}
    &\gQ_n^{(1)}= \\
    & \left(\int_0^1 \nabla^2 f\left(\vu_n 
    + (\vv_n + \frac{\eta_n^{1-\beta}}{\gamma}\nabla f(\vx^\star)) 
    - \frac{\eta_n^{1-\beta}}{\gamma}(\nabla f(\vx^\star) - \nabla f(\vx_n))
    + \lambda\frac{\eta_n^{1-\beta}}{\gamma}\xi_n^{(a)}
    \right)d\lambda \right)\\
    &\times\frac{\eta_n^{1-\beta}}{\gamma}\xi_n^{(a)}\\
    &= \nabla^2 f\left(\vu_n 
    + (\vv_n + \frac{\eta_n^{1-\beta}}{\gamma}\nabla f(\vx^\star)) 
    - \frac{\eta_n^{1-\beta}}{\gamma}(\nabla f(\vx^\star) - \nabla f(\vx_n))
    \right)\frac{\eta_n^{1-\beta}}{\gamma}\xi_n^{(a)}\\
    &+ \mathbf{Res}(\gQ_n^{(1)}),
\end{aligned}
\end{equation}
where $\mathbf{Res}(\gQ_n^{(1)}) = \left(\int_0^1 \nabla^2 f(\vh_n + \lambda\frac{\eta_n^{1-\beta}}{\gamma}\xi_n^{(a)}) - \nabla^2 f(\vh_n)\right)\frac{\eta_n^{1-\beta}}{\gamma}\xi_n^{(a)}$ with $\vh_n = \vu_n 
    + (\vv_n + \frac{\eta_n^{1-\beta}}{\gamma}\nabla f(\vx^\star)) 
    - \frac{\eta_n^{1-\beta}}{\gamma}(\nabla f(\vx^\star) - \nabla f(\vx_n))$.
Hence, by invoking Assumption~\ref{asp:hess_lip}, one can figure out that $\norm{\mathbf{Res}(\gQ_n^{(1)})} \le \Tilde{L}\norm{\frac{\eta_n^{1-\beta}}{\gamma}\xi_n^{(a)}}^2$.
So
\begin{equation}\label{eq:gq1_bd}
\begin{aligned}
    \gQ_n^{(1)} &= \nabla^2 f\left(\vu_n 
    + (\vv_n + \frac{\eta_n^{1-\beta}}{\gamma}\nabla f(\vx^\star)) 
    - \frac{\eta_n^{1-\beta}}{\gamma}(\nabla f(\vx^\star) - \nabla f(\vx_n))
    \right)\frac{\eta_n^{1-\beta}}{\gamma}\xi_n^{(a)}\\
    &+ \gO \left(\frac{\Tilde{L}\eta_n^{2-2\beta}}{\gamma^2}\norm{\xi_n^{(a)}}^2\right).
\end{aligned}
\end{equation}

Here we need to explain that, for a sequence of vector $\{\vartheta_n\}$ and a corresponding sequence of scalar $\{c_n\}$, the notation $\vartheta_n = \gO(c_n)$ means that there is a universal constant $C$ satisfies $\norm{\vartheta_n} \le C \cdot c_n$ for all $n \ge 1$. For $\gQ_n^{(2)}$, owing to the Lipschitzness of $\nabla f(\cdot)$, the following fact is easy to get,
\begin{equation}\label{eq:gq2_bd}
    \norm{\gQ_n^{(2)}} \le \frac{L\eta_n^{1-\beta}}{\gamma}\norm{\vx_n - \vx^\star}.
\end{equation}

For $\gQ_n^{(3)}$, by using mean value theorems for definite integrals and Assumption~\ref{asp:hess_lip} again, we have
\begin{equation}\label{eq:gq3_bd}
\begin{aligned}
    &\gQ_n^{(3)} = \nabla^2 f(\vu_n)\left(\vv_n + \frac{\eta_n^{1-\beta}}{\gamma} \nabla f(\vx^\star)\right)\\
    &+ \left(\int_0^1 \nabla^2 f\left(\vu_n + \lambda(\vv_n + \frac{\eta_n^{1-\beta}}{\gamma}\nabla f(\vx^\star))\right)d\lambda - \nabla^2 f(\vu_n)\right)\left(\vv_n + \frac{\eta_n^{1-\beta}}{\gamma} \nabla^2 f(\vx^\star)\right)\\
    &= \nabla^2 f(\vx^\star)\left(\vv_n + \frac{\eta_n^{1-\beta}}{\gamma} \nabla f(\vx^\star)\right) + (\nabla^2 f(\vu_n) - \nabla^2 f(\vx^\star))\left(\vv_n + \frac{\eta_n^{1-\beta}}{\gamma} \nabla f(\vx^\star)\right)\\
    &+ \gO\left( \Tilde{L}\norm{\vv_n + \frac{\eta_n^{1-\beta}}{\gamma}\nabla f(\vx^\star)}^2 \right)\\
    &= \nabla^2 f(\vx^\star)\left(\vv_n + \frac{\eta_n^{1-\beta}}{\gamma} \nabla f(\vx^\star)\right) + \gO\left( \Tilde{L}\norm{\vu_n - \vx^\star}\norm{\vv_n + \frac{\eta_n^{1-\beta}}{\gamma}\nabla f(\vx^\star)} \right)\\
    &+ \gO\left( \Tilde{L}\norm{\vv_n + \frac{\eta_n^{1-\beta}}{\gamma}\nabla f(\vx^\star)}^2 \right)
\end{aligned}
\end{equation}

The following lemma, whose proof is deferred to Appendix~\ref{sec:prf_4th_mmt_cov_rt}, will help us to derive the precise convergence rate of the new debiased algorithm.
\begin{lemma}\label{lem:4th_mmt_cov_rt}
Under Assumption~\ref{asp:smooth} - \ref{asp:noi_lip} and \ref{asp:lip_4th_mmt}, the following inequalities are true.
\begin{align*}
    &\EB\norm{\vu_n - \vx^\star}^4 \precsim 
    \frac{L^4 \eta_n^{4(1-\beta)}}{(\mu\gamma)^4}
    (1 + \norm{\nabla f(\vx^\star)}^4) + \frac{\Sigma^\prime\eta_n^2}{\mu^2}\\
    &\EB\norm{\vv_n}^4 \precsim \frac{\eta_n^{4(1-\beta)}}{\gamma^4}(1 + \norm{\nabla f(\vx^\star)}^4).
\end{align*}
\end{lemma}

\subsection{Proof of Theorem~\ref{thm:dlpsa_cov_rt}}
\begin{proof}

For the sequence $\EB \norm{\vu_n - \vx^\star}^2,~n\ge 1$, we have
\begin{equation}\label{eq:u_dblp_eq}
\begin{aligned}
&\EB \norm{\vu_{n+1} - \vx^\star}^2 = 
\EB 
\norm{\vu_n - \vx^\star}^2 - 
2\eta_n \EB \inner{\vu_n - \vx^\star}{\gP_{\mA^\bot}\ssum{i}{1}{4}\gQ_n^{(i)}}
\\
&+ \eta_n^2 \EB \left(
\norm{\ssum{i}{1}{4} \gQ_n^{(i)}}^2 + \norm{\gP_{\mA^\bot} \xi_n^{(b)}}^2
\right) + 2\eta_n \EB \inner{(\vu_n - \vx^\star) + \ssum{i}{1}{4}\gQ_n^{(i)}}{\gP_{\mA^\bot}\xi_n^{(b)}}\\
&= \EB \norm{\vu_n - \vx^\star}^2 +
+\eta_n^2 \EB \left(
\norm{\ssum{i}{1}{4}\gQ_n^{(i)}}^2 + \norm{\gP_{\mA^\bot}\xi_n^{(b)}}^2
\right)\\
&- 2\eta_n \EB \inner{\vu_n - \vx^\star}{\nabla f(\vu_n) - \nabla f(\vx^\star)} - 2\eta_n \ssum{i}{1}{3}\inner{\vu_n - \vx^\star}{\gQ_n^{(i)}}\\
&\le (1 - 2\mu\eta_n)\EB \norm{\vu_n - \vx^\star}^2 + 2L\eta_n \EB \norm{\vu_n - \vx^\star}\norm{\frac{\eta_n^{1-\beta}}{\gamma}
(\nabla f(\vx_n) - \nabla f(\vx^\star))}\\
&-2\eta_n \sum\limits_{i\in\{1,3\}} 
\EB\inner{\vu_n - \vx^\star}{\gQ_n^{(i)}} + 
\eta_n^2 \EB \left(
\norm{\ssum{i}{1}{4}\gQ_n^{(i)}}^2 + \norm{\xi_n^{(b)}}^2
\right)\\
&\le \left(1 - 2\mu\eta_n\right) \EB \norm{\vu_n - \vx^\star}^2 + \frac{2L^2 \eta_n^{2-\beta}}{\gamma}
\EB \norm{\vu_n - \vx^\star}\norm{\vx_n - \vx^\star}
\\
& -2\eta_n \sum\limits_{i\in\{1,3\}} 
\EB\inner{\vu_n - \vx^\star}{\gQ_n^{(i)}} + 
\eta_n^2 \EB \left(
\norm{\ssum{i}{1}{4}\gQ_n^{(i)}}^2 + \norm{\xi_n^{(b)}}^2
\right),
\end{aligned}
\end{equation}
where the second equality follows from $\EB \left[\xi_n^{(b)}| \sigma(\{\vu_k, \vv_k\}_{k\le n}, \zeta_n^{(a)})\right] = 0$. Now we aim to analysis every term in the last formula. Actually, for the second term, it holds that
\begin{align*}
    &\frac{2L^2\eta_n^{2-\beta}}{\gamma} \EB
    \norm{\vu_n - \vx^\star}\norm{\vx_n - \vx^\star} \le \frac{2L^2\eta_n^{2-\beta}}{\gamma}\left(
    \EB \norm{\vu_n - \vx^\star}^2 +
    \EB \norm{\vu_n - \vx^\star} \norm{\vv_n}
    \right)\\
    &\le \frac{\mu\eta_n}{16} \EB \norm{\vu_n - \vx^\star}^2 + \frac{64L^4\eta_n^{3-2\beta}}{\mu\gamma^2}\EB \norm{\vv_n}^2 + \frac{2L^2 \eta_n^{2-\beta}}{\gamma}\EB \norm{\vu_n - \vx^\star}^2\\
    &\le \frac{\mu\eta_n}{8}\EB \norm{\vu_n - \vx^\star}^2 + \frac{CL^4\eta_n^{3-2\beta}}{\gamma^2}\EB \norm{\vv_n}^2.
\end{align*}
For $-2\eta_n\EB \inner{\vu_n - \vx^\star}{\gQ_n^{(1)}}$, using (\ref{eq:gq1_bd}) yields
\begin{align*}
    &-2\eta_n \EB \inner{\vu_n - \vx^\star}{
    \gQ_n^{(1)}}\\
    &= -2\eta_n \EB \inner{\vu_n - \vx^\star}{
    \nabla^2 f(\vh_n) \frac{\eta_n^{1-\beta}}{\gamma} \xi_n^{(a)}} - 
    2\eta_n \EB \inner{\vu_n - \vx^\star}{
    \gO\left(
    \frac{\Tilde{L}\eta_n^{2-2\beta}}{\gamma^2}\norm{\xi_n^{(a)}}^2
    \right)}\\
    &\overset{(a)}{=} - 2\eta_n \EB \inner{\vu_n - \vx^\star}{
    \gO\left(
    \frac{\Tilde{L}\eta_n^{2-2\beta}}{\gamma^2}\norm{\xi_n^{(a)}}^2
    \right)}\\
    &\le \frac{2\Tilde{L}\eta_n^{3-2\beta}}{\gamma^2}\EB \norm{\vu_n - \vx^\star} \norm{\xi_n^{(a)}}^2
    \overset{(b)}{=} 
    \frac{2\Tilde{L}\eta_n^{3-2\beta}}{\gamma^2}
    \EB \norm{\vu_n - \vx^\star}\tr(\Sigma(\vx_n))\\
    &\le \frac{2\Tilde{L}\tr (\Sigma)\eta_n^{3-2\beta}}{\gamma^2}\EB 
    \norm{\vu_n - \vx^\star} +
    \frac{2\Tilde{L}L\eta_n^{3-2\beta}}{\gamma^2}
    \EB \norm{\vu_n - \vx^\star} \norm{\vx_n - \vx^\star}\\
    &\le \frac{\mu\eta_n}{16}\EB \norm{\vu_n - \vx^\star}^2 + \frac{C\Tilde{L}^2 \tr(\Sigma)^2}{\mu\gamma^4} \eta_n^{5-4\beta} +
    \frac{2\Tilde{L}L\eta_n^{3-2\beta}}{\gamma^2}(
    \EB \norm{\vu_n - \vx^\star}^2 +
    \EB \norm{\vu_n - \vx^\star} \norm{\vv_n}
    )
    \\
    &\le \left(
    \frac{\mu}{16} + \frac{2\Tilde{L}L\eta_n^{2-2\beta}}{\gamma^2}
    \right) \eta_n \EB \norm{\vu_n - \vx^\star}^2 +
    \frac{C\Tilde{L}^2 \tr(\Sigma)^2}{\mu \gamma^4}\eta_n^{5-4\beta} \\
    &\quad + \frac{\mu}{32}\eta_n \EB \norm{\vu_n - \vx^\star}^2 + \frac{C\Tilde{L}^2L^2}{\mu\gamma^4} \eta_n^{5-4\beta} \EB \norm{\vv_n}^2\\
    &\le \frac{\mu}{8}\eta_n \EB \norm{\vu_n - \vx^\star}^2 + \frac{C\Tilde{L}^2 \tr(\Sigma)^2}{\mu \gamma^4}\eta_n^{5-4\beta} +
    \frac{C\Tilde{L}^2L^2}{\mu\gamma^4} \eta_n^{5-4\beta} \EB \norm{\vv_n}^2,
\end{align*}
where $\vh_n = \vu_n + (\vv_n + \frac{\eta_n^{1-\beta}}{\gamma}\nabla f(\vx^\star)) - \frac{\eta_n^{1-\beta}}{\gamma}(\nabla f(\vx^\star) - \nabla f(\vx_n))$. In the derivation process, (a) is resulted from the martingale difference property of $\xi_n^{(a)}$ with respect to the natural filtration $\gF_n$, and (b) holds for Assumption~\ref{asp:noi_lip}. 
For $-2\eta_n \EB \inner{\vu_n - \vx^\star}{\gQ_n^{(3)}}$, invoking (\ref{eq:gq3_bd}) leads to
\begin{align*}
    &-2\eta_n \EB \inner{\vu_n - \vx^\star}{
    \gQ_n^{(3)}} = -2\eta_n \EB 
    \inner{\vu_n - \vx^\star}{ 
    \nabla^2 f(\vx^\star) \left(
    \vv_n + \frac{\eta_n^{1-\beta}}{\gamma}\nabla f(\vx^\star)
    \right)}\\
    &+ \eta_n\EB \inner{\vu_n - \vx^\star}{\gO
    \left(
    \Tilde{L}\norm{\vu_n - \vx^\star} 
    \norm{\vv_n + \frac{\eta_n^{1-\beta}}{\gamma}\nabla f(\vx^\star)}
    \right)}\\
    &+ \eta_n\EB \inner{\vu_n - \vx^\star}{\gO\left(
    \Tilde{L}\norm{\vv_n + \frac{\eta_n^{1-\beta}}{\gamma}\nabla f(\vx^\star)}^2
    \right)}\\
    &\le -2\eta_n \EB 
    \inner{\vu_n - \vx^\star}{ 
    \nabla^2 f(\vx^\star) \left(
    \vv_n + \frac{\eta_n^{1-\beta}}{\gamma}\nabla f(\vx^\star)
    \right)}\\
    &+ C\Tilde{L}\eta_n\EB \norm{\vu_n - \vx^\star}^2 \norm{\vv_n + \frac{\eta_n^{1-\beta}}{\gamma}\nabla f(\vx^\star)} + C\Tilde{L}\eta_n \EB \norm{\vu_n - \vx^\star} \norm{\vv_n + \frac{\eta_n^{1-\beta}}{\gamma}
    \nabla f(\vx^\star)}^2\\
    &\le -2\eta_n \EB 
    \inner{\vu_n - \vx^\star}{ 
    \nabla^2 f(\vx^\star) \left(
    \vv_n + \frac{\eta_n^{1-\beta}}{\gamma}\nabla f(\vx^\star)
    \right)} + \frac{C\Tilde{L}\eta_n^{2-\beta}}{\gamma}\norm{\nabla f(\vx^\star)}\EB \norm{\vu_n - \vx^\star}^2\\
    &+ C\Tilde{L}\eta_n \EB \norm{\vu_n - \vx^\star}^2 \norm{\vv_n} + C\Tilde{L}\eta_n \EB \norm{\vu_n - \vx^\star} \norm{\vv_n + \frac{\eta_n^{1-\beta}}{\gamma}\nabla f(\vx^\star)}^2\\
    &\le -2\eta_n \EB 
    \inner{\vu_n - \vx^\star}{ 
    \nabla^2 f(\vx^\star) \left(
    \vv_n + \frac{\eta_n^{1-\beta}}{\gamma}\nabla f(\vx^\star)
    \right)}\\
    &+ \left(\frac{\mu\eta_n}{16} + \frac{C\Tilde{L}\eta_n^{2-\beta}}{\gamma}\right)
    \EB \norm{\vu_n - \vx^\star}^2 + C\Tilde{L}\eta_n \left(
    \EB \norm{\vu_n - \vx^\star}^4 \EB \norm{\vv_n}^2
    \right)^{\frac{1}{2}}\\
    &+ \frac{C\Tilde{L}^2}{\mu} \left\{
    \EB \norm{\vv_n}^4 + \frac{\eta_n^{4-4\beta}}{\gamma^4} \norm{\nabla f(\vx^\star)}^4
    \right\}\\
    &\le -2\eta_n \EB 
    \inner{\vu_n - \vx^\star}{ 
    \nabla^2 f(\vx^\star) \left(
    \vv_n + \frac{\eta_n^{1-\beta}}{\gamma}\nabla f(\vx^\star)
    \right)} + \frac{\mu\eta_n}{8}\EB \norm{\vu_n - \vx^\star}^2\\
    &+ C\Tilde{L}\eta_n \left(
    \EB \norm{\vu_n - \vx^\star}^4 \EB \norm{\vv_n}^2
    \right)^{\frac{1}{2}} +
    \frac{C\Tilde{L}^2\eta_n}{\mu} \left\{
    \EB \norm{\vv_n}^4 + \frac{\eta_n^{4-4\beta}}{\gamma^4} \norm{\nabla f(\vx^\star)}^4
    \right\}.
\end{align*}
For the last term of (\ref{eq:u_dblp_eq}), by leveraging the Lipschitzness of $\nabla f(\cdot)$ and $\Sigma(\vx)$ defined in Assumption~\ref{asp:noi_lip}, we have
\begin{align*}
    &\EB \left(
    \ssum{i}{1}{4}\norm{\gQ_n^{(i)}}^2 + 
    \norm{\xi_n^{(b)}}^2
    \right) \precsim
    L^2 \left\{
    \frac{\eta_n^{2-2\beta}}{\gamma^2} \EB \norm{\xi_n^{(a)}}^2 + \frac{\eta_n^{2-2\beta}}{\gamma^2}\EB
    \norm{\nabla f(\vx_n) - \nabla f(\vx^\star)}^2
    \right.\\
    &\left.+ 
    \EB \norm{\vv_n + \frac{\eta_n^{1-\beta}}{\gamma}\nabla f(\vx^\star)}^2 + \EB \norm{\vu_n - \vx^\star}^2
    \right\} + \EB \tr \left(
    \Sigma\left(
    \vx_n + \frac{\eta_n^{1-\beta}}{\gamma}\nabla f(\vx_n) + \frac{\eta_n^{1-\beta}}{\gamma}\xi_n^{(a)}
    \right)
    \right)\\
    &\precsim L^2 \left( 1 + \frac{L^2\eta_n^{2-2\beta}}{\gamma^2}\right)
    \left( \EB \norm{\vu_n - \vx^\star}^2 
    + \EB \norm{\vv_n}^2\right) + \frac{L^2\eta_n^{2-2\beta}}{\gamma^2}\norm{
    \nabla f(\vx^\star)}^2\\
    &+ \tr(\Sigma) + L\left\{
    \left( 1 + \frac{L\eta_n^{1-\beta}}{\gamma}\right)\EB \norm{\vx_n - \vx^\star} + \frac{\eta_n^{1-\beta}}{\gamma}(\norm{\nabla f(\vx^\star)} + \EB \norm{\xi_n^{(a)}})
    \right\}\\
    &\precsim L^2 \EB \norm{\vu_n - \vx^\star}^2 + \tr(\Sigma).
\end{align*}
As a consequence of taking these back to (\ref{eq:u_dblp_eq}), we get
\begin{equation}\label{eq:u_dblp_ieq}
\begin{aligned}
&\EB \norm{\vu_{n+1} - \vx^\star}^2 \le 
(1 - 2\mu\eta_n)\EB \norm{\vu_n - \vx^\star}^2 +
\frac{\mu\eta_n}{8} \EB \norm{\vu_n - \vx^\star}^2\\
&+ \frac{CL^4 \eta_n^{3-2\beta}}{\gamma^2}\EB 
\norm{\vv_n}^2 + \frac{\mu\eta_n}{8} \EB
\norm{\vu_n - \vx^\star}^2 +
\frac{C\Tilde{L}^2 \tr(\Sigma)^2}{\mu\gamma^4} \eta_n^{5-4\beta}\\
&+ \frac{C\Tilde{L}^2L^2}{\mu\gamma^4}\eta_n^{5-4\beta} \EB \norm{\vv_n}^2 - 
2\eta_n \EB \inner{\vu_n - \vx^\star}{\nabla^2 f(\vx^\star)\left(
\vv_n + \frac{\eta_n^{1-\beta}}{\gamma}\nabla f(\vx^\star)
\right)}\\
&+ \frac{\mu\eta_n}{8} \EB \norm{\vu_n - \vx^\star}^2 + CL^2\eta_n^2 \EB 
\norm{\vu_n - \vx^\star}^2 + C\tr(\Sigma) \eta_n^2\\
&+ C\Tilde{L}\eta_n \left(
    \EB \norm{\vu_n - \vx^\star}^4 \EB \norm{\vv_n}^2
    \right)^{\frac{1}{2}} +
    \frac{C\Tilde{L}^2\eta_n}{\mu} \left\{
    \EB \norm{\vv_n}^4 + \frac{\eta_n^{4-4\beta}}{\gamma^4} \norm{\nabla f(\vx^\star)}^4
    \right\}\\
    &\overset{(a)}{\le} \left(1 - \frac{5\mu}{4}\eta_n\right)
    \EB \norm{\vu_n - \vx^\star}^2 - 
    2\eta_n \EB \inner{\vu_n - \vx^\star}{
    \nabla^2 f(\vx^\star) \left(
    \vv_n + \frac{\eta_n^{1-\beta}}{\gamma}
    \nabla f(\vx^\star)\right)}\\
    &+ \frac{CL^4 \eta_n^{1+4(1-\beta)}}{\gamma^6}
    \norm{\nabla f(\vx^\star)}^2 + 
    \frac{C\Tilde{L}^2\tr(\Sigma)^2}{\mu\gamma^4}
    \eta_n^{1+4(1-\beta)} +
    \frac{C\Tilde{L}^2\eta_n^{1+4(1-\beta)}}{\mu\gamma^4}\norm{\nabla f(\vx^\star)}^4\\
    &+ C\Tilde{L}\eta_n \left(
    \frac{L^4\eta_n^{4(1-\beta)}}{(\mu\gamma)^4}
    \norm{\nabla f(\vx^\star)}^4 + 
    \frac{\Sigma^\prime}{\mu^2}\eta_n^2
    \right)^{\frac{1}{2}}\frac{\eta_n^{1-\beta}}{
    \gamma} \norm{\nabla f(\vx^\star)}
    + C\tr(\Sigma)\eta_n^2\\
    &\le \left(1 - \frac{5\mu}{4}\eta_n\right)
    \EB \norm{\vu_n - \vx^\star}^2 - 
    2\eta_n \EB \inner{\vu_n - \vx^\star}{
    \nabla^2 f(\vx^\star) \left(
    \vv_n + \frac{\eta_n^{1-\beta}}{\gamma}
    \nabla f(\vx^\star)\right)}\\
    &+ \frac{C\Tilde{L}L^2 \eta_n^{1+3(1-\beta)}}{
    \mu^2 \gamma^3} + \frac{C\Tilde{L}\Sigma^\prime}{\mu^2\gamma}
    \eta_n^{3-\beta}\norm{\nabla f(\vx^\star)} + 
    C\tr(\Sigma)\eta_n^2\\
    &\le \left(1 - \frac{5\mu}{4}\eta_n\right)
    \EB \norm{\vu_n - \vx^\star}^2
    + \frac{C\Tilde{L}L^2 \eta_n^{1+3(1-\beta)}}{
    \mu^2 \gamma^3} + 
    C\tr(\Sigma)\eta_n^2\\
    &+ 2\eta_n \left|\EB \inner{\vu_n - \vx^\star}{
    \nabla^2 f(\vx^\star) \left(
    \vv_n + \frac{\eta_n^{1-\beta}}{\gamma}
    \nabla f(\vx^\star)\right)}\right|.
\end{aligned}
\end{equation}
Here (a) holds in view of Lemma~\ref{lem:4th_mmt_cov_rt}. And what remain to be handled is the inner production term $\left|\EB \inner{\vu_n - \vx^\star}{\nabla^2 f(\vx^\star) \left(\vv_n + \frac{\eta_n^{1-\beta}}{\gamma}
\nabla f(\vx^\star)\right)}
\right|$ which is shorthanded for $\gU_n$. Then we are going to find an iterative inequality for $\gU_n$ just like (\ref{eq:u_dblp_ieq}) for $\EB \norm{\vu_n - \vx^\star}^2$. And we denote $\nabla^2 f(\vx^\star)\nabla f(\vx^\star)$ as $\vv$ in the rest of the proof. We start by doing one-step iteration of $\vv_n$,
\begin{align*}
&\gU_{n+1} = \left|
\EB \inner{\vu_{n+1} - \vx^\star}{
\nabla^2 f(\vx^\star) \left(
\vv_{n+1} + \frac{\eta_{n+1}^{1-\beta}}{\gamma}
\nabla f(\vx^\star)
\right)}
\right|\\
&= \left|
(1 - \gamma\eta_n^\beta)\EB \inner{\vu_{n+1}- \vx^\star}{\nabla^2 f(\vx^\star) \vv_{n+\frac{1}{2}}} + \frac{\eta_n^{1-\beta}}{\gamma} \EB \inner{
\vu_{n+1} - \vx^\star}{\vv}
\right|\\
&= \left|
(1-\gamma\eta_n^\beta)\EB 
\inner{\vu_{n+1} - \vx^\star}{\nabla^2 f(\vx^\star)
\left(
\vv_n - \eta_n \nabla f(\vx^\star) -
\eta_n \ssum{i}{1}{4}\gQ_n^{(i)} + \eta_n \gP_{\mA} \xi_n^{(b)}
\right)}
\right.\\
&\quad \left.
+ \frac{\eta_{n+1}^{1-\beta}}{\gamma}\EB \inner{\vu_{n+1} - \vx^\star}{\vv}
\right|\\
&\overset{(a)}{=} \bigg|
(1 - \gamma\eta_n^\beta) \EB \inner{\vu_{n+1} - \vx^\star}{\nabla^2 f(\vx^\star) \vv_n} - 
\eta_n(1-\gamma\eta_n^{1-\beta}) \EB 
\inner{\vu_{n+1} - \vx^\star}{\vv}\\
&\left. \quad\quad- \eta_n \EB \inner{\vu_{n+1} - \vx^\star}{
\ssum{i}{1}{4}\gQ_n^{(i)}} +
\frac{\eta_{n+1}^{1-\beta}}{\gamma}\EB 
\inner{\vu_{n+1} - \vx^\star}{\vv}\right|\\
&\le (1 - \gamma \eta_n^{\beta}) \left|
\inner{\vu_{n+1} - \vx^\star}{\nabla^2 
f(\vx^\star) \left(
\vv_n + \frac{\eta_n^{1-\beta}}{\gamma}\nabla f(\vx^\star)
\right)}
\right| \\
&\quad+ \eta_n(1 - \gamma \eta_n^\beta)
\left|
\EB \inner{\vu_{n+1} - \vx^\star}{
\nabla^2 f(\vx^\star)\ssum{i}{1}{4}\gQ_n^{(i)}}
\right|\\
&\quad + \left|
\left(
\frac{\eta_{n+1}^{1-\beta}}{\gamma} -
\eta_n (1 - \gamma \eta_n^{1-\beta}) -
\frac{\eta_n^{1-\beta}}{\gamma}(1 - \gamma \eta_n^\beta)
\right) \EB \inner{\vu_{n+1} - \vx^\star}{\vv}
\right|\\
&\le (1 - \gamma \eta_n^{\beta}) \left|
\EB \inner{\vu_{n+1} - \vx^\star}{\nabla^2 
f(\vx^\star) \left(
\vv_n + \frac{\eta_n^{1-\beta}}{\gamma}\nabla f(\vx^\star)
\right)}
\right| \\
&\quad + \eta_n(1 - \gamma \eta_n^\beta)
\left|
\EB \inner{\vu_{n+1} - \vx^\star}{
\nabla^2 f(\vx^\star)\ssum{i}{1}{4}\gQ_n^{(i)}}
\right|\\
&\quad + \eta_n^{1-\beta}(\eta_n + \frac{1}{\gamma n})\left| \EB \inner{\vu_{n+1} - \vx^\star}{\vv}
\right|,
\end{align*}
where (a) is guaranteed by the zero mean nature of $\xi_n^{(b)}$ which is conditional on $\sigma(\gF_n , \zeta_n^{(a)})$. We take the update rule of $\vu_{n+1}$ into the equation sequentially and obtain
\begin{align*}
    &\gU_{n+1} \le 
     (1 - \gamma \eta_n^\beta)\gU_n + 
\eta_n(1-\gamma\eta_n^\beta) \left| 
\EB \inner{\ssum{i}{1}{4}\gQ_n^{(i)}}{
\nabla^2 f(\vx^\star) \left(
\vv_n + \frac{\eta_n^{1-\beta}}{\gamma}
\nabla f(\vx^\star)
\right)}
\right|\\
&\quad + \eta_n \left|
\EB \inner{\vu_n - \vx^\star}{\nabla^2 f(\vx^\star)
\ssum{i}{1}{4}\gQ_n^{(i)}}
\right| + L\eta_n \EB \norm{\ssum{i}{1}{4}\gQ_n^{(i)}}^2\\
&\quad + (1 + \frac{\mathbbm{1}_{\alpha = 1}}{\gamma}) \eta_n^{2-\beta} \left| \EB \inner{
\vu_n - \vx^\star}{\vv}\right| + 
(1 + \frac{\mathbbm{1}_{\alpha}}{\gamma}) \eta_n^{3-\beta} \left| \EB \inner{\ssum{i}{1}{4}\gQ_n^{(i)}}{\vv} \right|\\
&\le (1 - \gamma \eta_n^\beta) \gU_n + L\eta_n 
\left\{
\EB \norm{\ssum{i}{1}{4}\gQ_n^{(i)}}^2 +
 \EB \norm{\vv_n}^2 + \frac{\eta_n^{2-2\beta}}{\gamma^2} \norm{\nabla f(\vx^\star)}^2
\right\}\\
&\quad + L\eta_n \EB\norm{\vu_n - \vx^\star}^2 +
2L\eta_n \EB \norm{\ssum{i}{1}{4}\gQ_n^{(i)}}^2 +
L\eta_n \EB \norm{\vu_n - \vx^\star}^2\\
&\quad + 
CL(1+ \frac{\mathbbm{1}_{\alpha = 1}}{\gamma^2}) \eta_n^{3-2\beta} \norm{\nabla f(\vx^\star)}^2
+ CL\norm{\nabla f(\vx^\star)}(1+\frac{\mathbbm{1}_{\alpha = 1}}{\gamma}) \eta_n^{3-\beta} \EB \norm{\ssum{i}{1}{4}\gQ_n^{(i)}}\\
&\le (1 - \gamma \eta_n^\beta) \gU_n + 3L\eta_n \EB \norm{\vu_n - \vx^\star}^2 + 4L\eta_n \EB \norm{\ssum{i}{1}{4}\gQ_n^{(i)}}^2\\
&\quad + CL\left(
\left(1 + \frac{\mathbbm{1}_{\alpha = 1}}{\gamma^2}\right) \eta_n^{5-2\beta} + \frac{1}{\gamma^2}\eta_n^{3-2\beta}
\right)\norm{\nabla f(\vx^\star)}^2.
\end{align*}
Recall the bound of $\EB \norm{\ssum{i}{1}{4}\gQ_n^{(i)}}^2$ that we have analyzed before,
\begin{align*}
    &\EB \norm{\ssum{i}{1}{4}\gQ_n^{(i)}}^2\le 
    \ssum{i}{1}{4} \EB \norm{\gQ_n^{(i)}}^2\\
    &\precsim L^2\left(
    1 + \frac{L^2 \eta_n^{2-2\beta}}{\gamma^2}
    \right)(\EB \norm{\vu_n - \vx^\star}^2 + 
    \EB \norm{\vv_n}^2) + 
    \frac{L^2 \eta_n^{2-2\beta}}{\gamma^2}(\tr(\Sigma) + \norm{\nabla f(\vx^\star)}^2)\\
    &\precsim L^2 \EB \norm{\vu_n - \vx^\star}^2
    + \frac{L^2 \eta_n^{2-2\beta}}{\gamma}
    (\tr(\Sigma) + \norm{\nabla f(\vx^\star)}^2).
\end{align*}
Bring this back, we finally get
\begin{equation}\label{eq:v_dblp_ieq}
\begin{aligned}
&\gU_{n+1} \le (1 - \gamma \eta_n^\beta) \gU_n +
3L\eta_n \EB \norm{\vu_n - \vx^\star}^2 +
\frac{CL\eta_n^{3-2\beta}}{\gamma^2}\eta_n^{3-2\beta} \norm{\nabla f(\vx^\star)}^2\\
&\quad + CL^3 \eta_n \EB \norm{\vu_n - \vx^\star}^2 + CL^3\eta_n\times  \frac{\eta_n^{2-2\beta}}{\gamma^2}
(\tr(\Sigma) + \norm{\nabla f(\vx^\star)}^2)\\
&\le (1 - \gamma \eta_n^\beta) \gU_n + CL^3 \eta_n \EB \norm{\vu_n - \vx^\star}^2 +
\frac{CL^3\eta_n^{3-2\beta}}{\gamma^2}
(\tr(\Sigma) + \norm{\nabla f(\vx^\star)}^2).
\end{aligned}
\end{equation}
Consider the combination of (\ref{eq:u_dblp_ieq}) and $\ref{eq:v_dblp_ieq}$. It can be derived that
\begin{align*}
    &\EB \norm{\vu_{n+1} - \vx^\star}^2 + \frac{\mu}{CL^3}\gU_{n+1}\\
    &\le 
    \left(1 - \frac{5\mu}{4}\eta_n\right)\EB 
    \norm{\vu_n - \vx^\star}^2 + 2\eta_n\gU_n
    + \frac{C\Tilde{L}L^2\eta_n^{1+3(1-\beta)}}{
    \mu^2\gamma^3} + C\tr(\Sigma)\eta_n^2\\
    &\quad +
    (1 - \gamma \eta_n^\beta)\frac{\mu}{CL^3}\gU_n
    + \frac{\mu\eta_n}{4}\EB \norm{\vu_n - \vx^\star}^2 + 
    \frac{\mu\eta_n^{3-2\beta}}{\gamma^2}
    (\tr(\Sigma) + \norm{\nabla f(\vx^\star)}^2)\\
    &\le (1 - \mu\eta_n) \EB \norm{\vu_n - \vx^\star}^2 + (1 - \frac{\gamma}{2}\eta_n^\beta)
    \frac{\mu}{CL^3}\gU_n + C\tr(\Sigma)\eta_n^2\\
    &\quad + \frac{\mu\eta_n^{3-2\beta}}{\gamma^2}
    (\tr(\Sigma) + \norm{\nabla f(\vx^\star)}^2)\\
    &\le (1 - \mu\eta_n) \left(\EB \norm{\vu_n - \vx^\star}^2 +
    \frac{\mu}{CL^3}\gU_n \right) + C\tr(\Sigma)\eta_n^2\\
    &\quad + \frac{\mu\eta_n^{3-2\beta}}{\gamma^2}
    (\tr(\Sigma) + \norm{\nabla f(\vx^\star)}^2)
\end{align*}
In view of Lemma~\ref{lem:converge_r_t}, it holds that
\[
\EB \norm{\vu_{n} - \vx^\star}^2 + \frac{\mu}{CL^3}\gU_{n} \precsim \frac{\eta_n^{2(1-\beta)}}{\gamma^2}(\tr(\Sigma) + \norm{\nabla f(\vx^\star)}^2) + \frac{\tr(\Sigma)}{\mu}\eta_n.
\]
By substituting $\EB \norm{\vu_n - \vx^\star}^2$ in (\ref{eq:v_dblp_ieq}) for the bound, it follows that
\begin{align*}
    \gU_{n+1} \le (1 - \gamma\eta_n^\beta)\gU_n
    + \frac{CL^3\eta_n^{3-2\beta}}{\gamma}(\tr(\Sigma) + \norm{\nabla f(\vx^\star)^2}) + \frac{CL^3\tr(\Sigma)}{\mu}\eta_n^2,
\end{align*}
which implies
\[
\gU_n \precsim \frac{L^3\eta_n^{3(1-\beta)}}{\gamma^2}
(\tr(\Sigma) + \norm{\nabla f(\vx^\star)}^2)+
\frac{L^3\tr(\Sigma)}{\mu\gamma}\eta_n^{2-\beta}.
\]
Replacing this result into (\ref{eq:u_dblp_ieq}) yields
\begin{align*}
    &\EB \norm{\vu_{n+1} - \vx^\star}^2 \le 
    (1 - \mu\eta_n) \EB \norm{\vu_n - \vx^\star}^2
    + \frac{C\Tilde{L}L^2\eta_n^{1 + 3(1-\beta)}}{\mu^2\gamma^3} + C\tr(\Sigma)\eta_n^2\\
    &+ \frac{L^3\eta_n^{1 + 3(1-\beta)}}{\gamma^2}
    (\tr(\Sigma) + \norm{\nabla f(\vx^\star)}^2)
    + \frac{L^3\tr(\Sigma)}{\mu\gamma}\eta_n^{3-\beta}.
\end{align*}
Invoking Lemma~\ref{lem:converge_r_t} leads to
\[
\EB\norm{\vu_n - \vx^\star}^2 \precsim
\left(\frac{\Tilde{L}}{\mu^2 \gamma} + L(\tr(\Sigma) + \norm{\nabla f(\vx^\star)}^2)\right)\frac{L^2\eta_n^{3(1-\beta)}}{\mu\gamma^2} + \frac{\tr(\Sigma)}{\mu}\eta_n.
\]
Up to now, the proof is complete.
\end{proof}

\subsection{Proof of Lemma~\ref{lem:4th_mmt_cov_rt}}\label{sec:prf_4th_mmt_cov_rt}

\begin{proof}[Proof of Lemma~\ref{lem:4th_mmt_cov_rt}]
By direct expanding, we have an iterative equation
\begin{equation}\label{eq:u_4th_mmt_eq}
\begin{aligned}
&\EB \norm{\vu_{n+1} - \vx^\star}^4\\
=& \EB\left\{ \norm{\vu_n - \vx^\star}^2 - 2\eta_n
\left\langle \vu_n - \vx^\star, \gP_{\mA^\bot}\ssum{i}{1}{4}\gQ_n^{(i)}\right\rangle \right.\\
&\left.+ 2\eta_n \langle \vu_n - \vx^\star, \gP_{\mA^\bot}\xi_n^{(b)} \rangle + \eta_n^2 \norm{\gP_{\mA^\bot}\ssum{i}{1}{4}\gQ_n^{(i)} + \gP_{\mA^\bot}\xi_n^{(b)}}^2
\right\}^2\\
\le& \EB \norm{\vu_n - \vx^\star}^4 - 4\eta_n \EB\norm{\vu_n - \vx^\star}^2 \left\langle \vu_n - \vx^\star, \gP_{\mA^\bot}\ssum{i}{1}{4}\gQ_n^{(i)} \right\rangle\\
&+ 4\eta_n \EB \norm{\vu_n - \vx^\star}^2 \langle \vu_n - \vx^\star, \gP_{\mA^\bot} \xi_n^{(b)} \rangle
 + 6\eta_n^2 \EB \norm{\vu_n - \vx^\star}^2 \norm{\gP_{\mA^\bot}\left( \ssum{i}{1}{4}\gQ_n^{(i)} + \xi_n^{(b)} \right)}^2\\
 &+ 8\eta_n \EB\inner{\vu_n - \vx^\star}{\gP_{\mA^\bot}\ssum{i}{1}{4}\gQ_n^{(i)}}  \inner{\vu_n - \vx^\star}{\gP_{\mA^\bot} \xi_n^{(b)}}\\
 & - 2\eta_n^3 \EB \inner{\vu_n - \vx^\star}{\gP_{\mA^\bot}\ssum{i}{1}{4}\gQ_n^{(i)} + \gP_{\mA^{\bot}}\xi_n^{(b)}}\norm{\gP_{\mA^\bot}\left(  \ssum{i}{1}{4}\gQ_n^{(i)} + \xi_n^{(b)}\right)}^2\\
 &+ \eta_n^4 \EB \norm{\gP_{\mA^{\bot}}\left( \ssum{i}{1}{4}\gQ_n^{(i)} + \xi_n^{(4)}\right)}^4.
\end{aligned}
\end{equation}
Since $\EB [\xi_n^{(a)}| \sigma\left( \{\vu_j, \vv_j\}_{j = 1}^n \right)] = 0$ and $\EB [\xi_n^{(b)}| \sigma\left( \{\vu_j, \vv_j\}_{j = 1}^n, \xi_n^{(a)} \right)] = 0$, we have
\begin{align*}
    &\EB \norm{\vu_n - \vx^\star}^2 \inner{\vu_n - \vx^\star}{\gP_{\mA^\bot}\xi_n^{b}} = 0,\\
    &\EB\inner{\vu_n - \vx^\star}{\gP_{\mA^\bot}\ssum{i}{1}{4}\gQ_n^{(i)}}\inner{\vu_n - \vx^\star}{\gP_{\mA^\bot}\xi_n^{(b)}}=0.
\end{align*}
Then we combine the description of $\gQ_n^{(i)}$, i.e., (\ref{eq:gq1_bd}, \ref{eq:gq2_bd}) and $\norm{\gQ_n^{(3)}} \le L \norm{\vv_n - \frac{\eta_n^{1-\beta}}{\gamma}\nabla f(\vx^\star)}$, with Young's inequality and get
\begin{align*}
    &- \EB \norm{\vu_n - \vx^\star}^2 \inner{\vu_n - \vx^\star}{\gP_{\mA^\bot}\ssum{i}{1}{4}\gQ_n^{(i)}}
    = -\ssum{i}{1}{4} \EB \norm{\vu_n - \vx^\star}^2 \inner{\vu_n - \vx^\star}{ \gQ_n^{(i)}}\\
    \le& \left| \EB \norm{\vu_n - \vx^\star}^2\inner{\vu_n - \vx^\star}{\nabla^2 f\left( \vh_n\right)\frac{\eta_n^{1-\beta}}{\gamma}\xi_n^{(a)} + \gO\left( \frac{L\eta_n^{2-2\beta}}{\gamma^2}\norm{\xi_n^{(a)}}^2 \right)} \right|\\
    +& \frac{L\eta_n^{1-\beta}}{\gamma}\EB \norm{\vu_n - \vx^\star}^3\norm{\vx_n - \vx^\star} 
    + L\EB \norm{\vu_n - \vx^\star}^3\norm{\vv_n + \frac{\eta_n^{1-\beta}}{\gamma}\nabla f(\vx^\star)}\\
    -& \EB \norm{\vu_n - \vx^\star}^2 \inner{\vu_n - \vx^\star}{\nabla f(\vu_n) - \nabla f(\vx^\star)}\\
    \overset{(a)}{\le}& \frac{L\eta_n^{2-2\beta}}{\gamma}\EB \norm{\vu_n - \vx^\star}^3\norm{\xi_n^{(a)}}^2 + \frac{L\eta_n^{1-\beta}}{\gamma}\EB \norm{\vu_n - \vx^\star}^3\norm{\vx_n - \vx^\star}\\
    +& L\EB \norm{\vu_n - \vx^\star}^3 \norm{\vv_n + \frac{\eta_n^{1-\beta}}{\gamma}\nabla f(\vx^\star)} - \EB \norm{\vu_n - \vx^\star}^2 \inner{\vu_n - \vx^\star}{\nabla f(\vu_n) - \nabla f(\vx^\star)}\\
    \overset{(b)}{=}& \frac{L\eta_n^{2-2\beta}}{\gamma}\EB \left\{\norm{\vu_n - \vx^\star}^3 \tr \left(\Sigma(\vx_n)\right)\right\} + \frac{L\eta_n^{1-\beta}}{\gamma}\EB \norm{\vu_n - \vx^\star}^3 \norm{(\vu_n - \vx^\star) + \vv_n}\\
    +& L\EB\norm{\vu_n - \vx^\star}^3\norm{\vv_n + \frac{\eta_n^{1-\beta}}{\gamma}\nabla f(\vx^\star)} - \EB\norm{\vu_n - \vx^\star}^2 
    \inner{\vu_n - \vx^\star}{\nabla f(\vu_n) - \nabla f(\vx^\star)}\\
    \le& \frac{L^2\eta_n^{2-2\beta}}{\gamma}\EB \norm{\vu_n - \vx^\star}^4 + \left( \frac{\mu}{8}\EB \norm{\vu_n - \vx^\star}^4 + \frac{54\tr(\Sigma)^4L^4 \eta_n^{8-8\beta}}{\mu^3 \gamma^4} \right) + \frac{L\eta_n^{1-\beta}}{\gamma}\EB \norm{\vu_n - \vx^\star}^4\\
    +& \left( \frac{\mu}{8}\EB\norm{\vu_n - \vx^\star}^4
    + \frac{54L^4\eta_n^{4-4\beta}}{\mu^3 \gamma^4}\EB \norm{\vv_n}^4 \right)\\
    +& \left( \frac{\mu}{8}\EB\norm{\vu_n - \vx^\star}^4 + \frac{54L^4}{\mu^3} \EB \norm{\vv_n + \frac{\eta_n^{1-\beta}}{\gamma}\nabla f(\vx^\star)}^4 \right) - \mu \EB \norm{\vu_n - \vx^\star}^4\\
    \le& \left(-\mu + \frac{3}{8}\mu + \frac{L\eta_n^{1-\beta}}{\gamma} + \frac{L^2 \eta_n^{2-2\beta}}{\gamma}\right)\EB\norm{\vu_n - \vx^\star}^4 + \left( \frac{432L^4}{\mu^3} + \frac{54L^4\eta_n^{4-4\beta}}{\mu^3\gamma^4} \right)\EB\norm{\vv_n}^4\\
    +& \left(\frac{432L^4\eta_n^{4-4\beta}}{\mu^3\gamma^4}\norm{\nabla f(\vx^\star)}^4 + \frac{54\tr(\Sigma)^4L^4\eta_n^{8-8\beta}}{\mu^3\gamma^4}\right)\\
    \le& -\frac{\mu}{2}\EB\norm{\vu_n - \vx^\star}^4 + \frac{864L^4}{\mu^3}\EB\norm{\vv_n}^4
    + \frac{864L^4\eta_n^{4-4\beta}}{\mu^3 \gamma^4}\norm{\nabla f(\vx^\star)}^4,
\end{align*}
where (a) holds for $\EB \norm{\vu_n - \vx^\star}^2 \inner{\vu_n - \vx^\star}{\nabla^2 f(\vh_n) \xi_n^{(a)}} = 0$, and (b) holds by Assumption~\ref{asp:noi_lip}. 
\begin{align*}
    &\EB\norm{\gP_{\mA^\bot}\left( \ssum{i}{1}{4}\gQ_n^{(i)} + \xi_n^{(b)} \right)}^4 \le 125\left\{ \ssum{i}{1}{4}\EB\norm{\gQ_n^{(i)}}^4 + \EB\norm{\xi_n^{(b)}}^4 \right\}\\
    &\le 125\left\{ \frac{L^4\eta_n^{4-4\beta}}{\gamma^4}\EB\left(\norm{\xi_n^{(a)}}^4 + L^4\norm{\vx_n - \vx^\star}^4\right) + L^4\EB\norm{\vv_n + \frac{\eta_n^{1-\beta}}{\gamma}\nabla f(\vx^\star)}^4\right\}\\
    &~~+ 125 \left\{ L^4\EB\norm{\vu_n - \vx^\star}^4 + \EB \norm{\xi_n^{(b)}}^4\right\}\\
    &\overset{(a)}{\le} 125 \left\{ \frac{L^8\eta_n^{4-4\beta}}{\gamma^4} \EB \norm{\vx_n - \vx^\star}^4 + 
    L^4\EB \norm{\vv_n + \frac{\eta_n^{1- \beta}}{\gamma}\nabla f(\vx^\star)}^4 + L^4\EB\norm{\vu_n - \vx^\star}^4 \right\}\\
    &~~+ 125\Sigma^\prime \frac{L^4\eta_n^{4-4\beta}}{\gamma^4}(1+\EB\norm{\vx_n - \vx^\star}^4)\\
    &~~+ 125\Sigma^\prime \left(1+\EB\norm{\vx_n - \vx^\star + \frac{\eta_n^{1-\beta}}{\gamma}\nabla f(\vx_n) +
    \frac{\eta_n^{1-\beta}}{\gamma}\xi_n^{(a)}
    }^4 \right)\\
    &\overset{(b)}{\le} 125\left(L^4 + 64 \Sigma^\prime + \frac{8L^8\eta_n^{4-4\beta}}{\gamma^4}(1+ \Sigma^\prime) \right) \EB \norm{\vu_n - \vx^\star}^4\\
    &~~+ 1000\left(L^4 + 8\Sigma^\prime + \frac{L^8\eta_n^{4-4\beta}}{\gamma^4}(1+\Sigma^\prime)\right)\EB\norm{\vv_n}^4\\
    &~~+ 125\left\{ \frac{L^4\eta_n^{4-4\beta}}{\gamma^4}
    \left(8\norm{\nabla f(\vx^\star)}^4 + 1 \right) + \Sigma^\prime + 64\frac{\Sigma^\prime \eta_n^{4-4\beta}}{\gamma^4}\left(\EB\norm{\nabla f(\vx_n)}^4 + \EB\norm{\xi_n^{(a)}}^4\right)\right\}\\
    &\precsim (L^4 + \Sigma^\prime)\EB \norm{\vu_n - \vx^\star}^4
    + (L^4 + \Sigma^\prime)\EB \norm{\vv_n}^4
    + \Sigma^\prime,
\end{align*}
where (a) holds by Assumption~\ref{asp:lip_4th_mmt} and (b) holds by using Jensen's inequality and rearranging. Finally, we use Young's inequality and Jensen's inequality to simplify rest terms in (\ref{eq:u_4th_mmt_eq}).
\begin{align*}
    &6\eta_n^2 \EB \norm{\vu_n - \vx^\star}^2 \norm{\gP_{\mA^\bot}\left( \ssum{i}{1}{4}\gQ_n^{(i)} + \xi_n^{(b)}\right)}^2\\
    \le& \frac{\mu\eta_n}{4}\EB\norm{\vu_n - \vx^\star}^4
    + \frac{144}{\mu}\eta_n^3 \EB \norm{\gP_{\mA^\bot}\left(\ssum{i}{1}{4}\gQ_n^{(i)} + \xi_n^{(b)}\right)}^4;\\
    & -2 \eta_n^3 \EB \inner{\vu_n - \vx^\star}{ \gP_{\mA^\bot}\left(\ssum{i}{1}{4}\gQ_n^{(i)} + 
    \xi_n^{(b)}
    \right)}
    \norm{\gP_{\mA^\bot}\left( \ssum{i}{1}{4}\gQ_n^{(i)} + \xi_n^{(b)} \right)}^2\\
    \le& 2\eta_n^3 \EB \norm{\vu_n - \vx^\star}
    \norm{\gP_{\mA^\bot}\left(\ssum{i}{1}{4}\gQ_n^{(i)} + 
    \xi_n^{(b)}
    \right)}^3\\
    \le& \frac{\mu\eta_n}{2}\EB \norm{\vu_n - \vx^\star}^4 +  \frac{3}{2\mu^{1/3}}\eta_n^{\frac{11}{3}}\EB \norm{\gP_{\mA^\bot}\left(\ssum{i}{1}{4}\gQ_n^{(i)} + 
    \xi_n^{(b)}
    \right)}^4
\end{align*}
Taking all the results we got above into (\ref{eq:u_4th_mmt_eq}) can generate an iterative inequality,
\begin{equation}\label{eq:u_4th_mmt_ieq}
\begin{aligned}
&\EB\norm{\vu_{n+1} - \vx^\star}^4 
=\EB \norm{\vu_n - \vx^\star}^4 - 4 \eta_n \EB \norm{\vu_n - \vx^\star}^2\inner{\vu_n - \vx^\star}{ \gP_{\mA^\bot}\ssum{i}{1}{4}\gQ_n^{(i)}}\\
&- 2\eta_n^3 \EB \inner{\vu_n - \vx^\star}{\gP_{\mA^\bot}\left( \ssum{i}{1}{4}\gQ_n^{(i)} + \xi_n^{(b)} \right)}
\norm{\gP_{\mA^\bot}\left( \ssum{i}{1}{4}\gQ_n^{(i)} + \xi_n^{(b)} \right)}^2\\
&+ 6\eta_n^2 \EB \norm{\vu_n - \vx^\star}^2\norm{\gP_{\mA^\bot}\left( \ssum{i}{1}{4}\gQ_n^{(i)} + \xi_n^{(b)} \right)}^2
+ \eta_n^4 \EB \norm{\gP_{\mA^\bot}\left( \ssum{i}{1}{4}\gQ_n^{(i)} + \xi_n^{(b)} \right)}^4\\
&\le \EB \norm{\vu_n - \vx^\star}^4 - 4\eta_n\left(\frac{\mu}{2}\EB \norm{\vu_n - \vx^\star}^4 + \frac{864L^4}{\mu^3}\EB \norm{\vv_n}^4 + \frac{864L^4 \eta_n^{4-4\beta}}{\mu^3 \gamma^4}\norm{\nabla f(\vx^\star)}^4\right)\\
&+ \frac{\eta_n}{4}\EB \norm{\vu_n - \vx^\star}^4 + \frac{144}{\mu}\eta_n^3 \EB \norm{\gP_{\mA^\bot}\left( \ssum{i}{1}{4}\gQ_n^{(i)} + \xi_n^{(b)} \right)}^4\\
&+ \frac{\eta_n}{2} \EB \norm{\vu_n - \vx^\star}^4 + \left(\frac{3}{2\mu^{1/3}}\eta_n^{\frac{11}{3}} + \eta_n^4\right)\EB \norm{\gP_{\mA^\bot}\left( \ssum{i}{1}{4}\gQ_n^{(i)} + \xi_n^{(b)} \right)}^4\\
&\le (1 - \mu\eta_n) \EB \norm{\vu_n - \vx^\star}^4 + \frac{CL^4 \eta_n}{\mu^3}\EB \norm{\vv_n}^4 + \frac{CL^4 \eta_n^{5-4\beta}}{\mu^3 \gamma^4}\norm{\nabla f(\vx^\star)}^4\\
&+ C(\mu^{-1}\eta_n^3 + \mu^{-1/3}\eta_n^{11/3} + \eta_n^4)\EB \norm{\gP_{\mA^\bot}\left( \ssum{i}{1}{4}\gQ_n^{(i)} + \xi_n^{(b)} \right)}^4\\
&\le (1 - \frac{5}{4}\mu\eta_n) \EB \norm{\vu_n - \vx^\star}^4 + \frac{CL^4 \eta_n}{\mu^3}\EB \norm{\vv_n}^4 + \frac{CL^4 \eta_n^{5-4\beta}}{\mu^3 \gamma^4}\norm{\nabla f(\vx^\star)}^4\\
&+ C\frac{\eta_n^3}{\mu}\left\{ (L^4 + \Sigma^\prime)\EB \norm{\vu_n - \vx^\star}^4
+ (L^4 + \Sigma^\prime)\EB \norm{\vv_n}^4 + \Sigma^\prime
\right\}\\
&\le (1 - \mu \eta_n)\EB \norm{\vu_n - \vx^\star}^4 + \frac{CL^4\eta_n}{\mu^3}\EB \norm{\vv_n}^4 + \frac{CL^4 \eta_n^{5-4\beta}}{\mu^3\gamma^4}\norm{\nabla f(\vx^\star)}^4 + \frac{C\eta_n^3}{\mu}\Sigma^\prime.
\end{aligned}
\end{equation}
Here we use $C$ to refer to all universal constant.
After this analysis on $\EB \norm{\vu_n - \vx^\star}$, we start to deal with the similar property of $\EB\norm{\vv_n}^4$. Actually, we have
\begin{equation}\label{eq:v_4th_mmt_eq}
\begin{aligned}
&\EB \norm{\vv_n}^4
= (1-\gamma\eta_n^\beta) \EB \norm{\vv_n - \eta_n\nabla f(\vx^\star) - \eta_n \gP_{\mA}\ssum{i}{1}{4}\gQ_n^{(i)} + \eta_n\gP_{\mA}\xi_n^{(b)}}^4\\
&= (1-\gamma \eta_n^\beta) \EB \left\{
\norm{\vv_n}^2 - 2\eta_n\inner{\vv_n}{\nabla f(\vx^\star)} - 2\eta_n \inner{\vv_n}{\Omega_n} + \eta_n^2 \norm{\nabla f(\vx^\star) + \Omega_n}^2\right\}^2\\
&= (1 - \gamma \eta_n^\beta) \EB
\left\{
\norm{\vv_n}^4 - 4\eta_n \norm{\vv_n}^2 \inner{\vv_n}{\nabla f(\vx^\star)}
- 4\eta_n \norm{\vv_n}^2 \inner{\vv_n}{\Omega_n}\right.\\
&\quad + 2\eta_n^2 \norm{\vv_n}^2 \norm{\nabla f(\vx^\star) + \Omega_n}^2
+ 4\eta_n^2 \inner{\vv_n}{\nabla f(\vx^\star)}^2
+ 8\eta_n^2 \inner{\vv_n}{\nabla f(\vx^\star)} \inner{\vv_n}{\Omega}\\
&\quad - 4\eta_n^3 \inner{\vv_n}{\nabla f(\vx^\star) + \Omega_n} \norm{\nabla f(\vx^\star) + \Omega_n}^2
+ 4\eta_n^2 \inner{\vv_n}{\Omega_n}^2\\
&\left.\quad + \eta_n^4 \norm{\nabla f(\vx^\star) + \Omega_n}^4\right\},
\end{aligned}
\end{equation}
where $\Omega_n$ is used to representing $\gP_{\mA}\left(\ssum{i}{1}{4}\gQ_n^{(i)} - \xi_n^{(b)}\right)$. As what we did above, we start to control every term in (\ref{eq:v_4th_mmt_eq}) respectively.
\begin{align*}
    &\EB -4\eta_n \norm{\vv_n}^2\inner{\vv_n}{\nabla f(\vx^\star)} \le \EB 4\eta_n\norm{\vv_n}^3\norm{\nabla f(\vx^\star)}\\
    \le& \frac{\gamma}{16}\eta_n^\beta \EB \norm{\vv_n}^4 + \frac{C}{\gamma^3}\eta_n^{4-3\beta}\norm{\nabla f(\vx^\star)}^4;\\\\
    &\EB -4\eta_n\norm{\vv_n}^2\inner{\vv_n}{\Omega_n}
    = \EB -4\eta_n \norm{\vv_n}\inner{\vv_n}{\ssum{i}{1}{4}\gQ_n^{(i)}}\\
    \le& \EB 4\eta_n \norm{\vv_n}^3 \norm{\Omega_n}
    \le \frac{\gamma}{16}\eta_n^\beta \norm{\vv_n}^4 + \frac{C}{\gamma^3}\eta_n^{4-3\beta}\ssum{i}{1}{4}\norm{\gQ_n^{(i)}}^4;\\\\
    & \EB 2\eta_n^2 \norm{\vv_n}^2 \norm{\nabla f(\vx^\star) + \Omega_n}^4\\
    \le& \frac{\gamma\eta_n^\beta}{16}\EB \norm{\vv_n}^4 + \frac{C\eta_n^{4-\beta}}{\gamma}\left\{ 
    \norm{\nabla f(\vx^\star)}^4 + \ssum{i}{1}{4}\EB\norm{\gQ_n^{(i)}}^4 + \EB\norm{\xi_n^{(b)}}^4
    \right\};\\\\
    & \EB 4\eta_n^2 \inner{\vv_n}{\nabla f(\vx^\star)}^2 \le \EB 4\eta_n^2 \norm{\vv_n}^2 \norm{\nabla f(\vx^\star)}^2\\
    \le& \frac{\gamma\eta_n^\beta}{16}\EB \norm{\vv_n}^4 + \frac{C\eta_n^{4-\beta}}{\gamma}\norm{\nabla f(\vx^\star)}^4;\\\\
    &\EB 8\eta_n^2 \inner{\vv_n}{\nabla f(\vx^\star)} \inner{\vv_n}{\Omega_n}
    = \EB 8\eta_n^2 \inner{\vv_n}{\nabla f(\vx^\star)} \inner{\vv_n}{\ssum{i}{1}{4}\gQ_n^{(i)}}\\
    \le& \ssum{i}{1}{4}\EB 8\eta_n^2 \norm{\vv_n}^2 \norm{\nabla f(\vx^\star)} \norm{\gQ_n^{(i)}}\\
    \le& \ssum{i}{1}{4}\EB \left\{
    \frac{\gamma\eta_n^{\beta}}{64}\norm{\vv_n}^4
    + \frac{C\eta_n^{4-\beta}}{\gamma}\norm{\nabla f(\vx^\star)}^2 \norm{\gQ_n^{(i)}}^2\right\}\\
    \le& \frac{\gamma\eta_n^\beta}{16}\EB \norm{\vv_n}^4 + \frac{C\eta_n^{4-\beta}}{\gamma}\left\{\norm{\nabla f(\vx^\star)}^4 + \ssum{i}{1}{4}\EB\norm{\gQ_n^{(i)}}^4\right\};\\\\
    &\EB -4\eta_n^3 \inner{\vv_n}{\nabla f(\vx^\star) + \Omega_n}\norm{\nabla f(\vx^\star) + \Omega_n}^2
    \le \EB 4\eta_n^3 \norm{\vv_n} \norm{\nabla f(\vx^\star) + \Omega_n}^3\\
    \le& \frac{3}{4}\EB \left( \left(\frac{\gamma}{12}\right)^{\frac{1}{4}}
    \eta_n^{\frac{\beta}{4}} \norm{\vv_n}
    \right)^4 + 
    C\left(
    \gamma^{-\frac{1}{4}} \eta_n^{3 - \frac{\beta}{4}} \norm{\nabla f(\vx^\star) + \Omega_n}^3
    \right)^{\frac{4}{3}}\\
    =& \frac{\gamma\eta_n^\beta}{16}\EB \norm{\vv_n}^4 + \frac{C\eta_n^{4-\beta/3}}{\gamma^{1/3}}\left\{
    \norm{\nabla f(\vx^\star)}^4 + \ssum{i}{1}{4} \EB\norm{\gQ_n^{(i)}}^4 + \EB\norm{\xi_n^{(b)}}^4
    \right\}.
\end{align*}
Integrating these results into (\ref{eq:v_4th_mmt_eq}) yields
\begin{align*}
\EB \norm{\vv_{n+1}}^4 &\le \left( 1 - \gamma \eta_n^\beta + \frac{7}{16}\gamma\eta_n^\beta\right) \EB \norm{\vv_n}^4 + C\left\{ \frac{\eta_n^{4-\beta}}{\gamma} + \frac{\eta_n^{4-\beta/3}}{\gamma^{1/3}} \right\}\EB \norm{\xi_n^{(b)}}^4\\
&+ C\left\{
\frac{\eta_n^{4-3\beta}}{\gamma^3}
+ \frac{\eta_n^{4-\beta}}{\gamma}
+ \frac{\eta_n^{4-\beta/3}}{\gamma^{1/3}}
\right\} \left\{\norm{\nabla f(\vx^\star)}^4 
+ \ssum{i}{1}{4}\EB \norm{\gQ_n^{(i)}}^4
\right\}\\
\overset{(a)}{\le}& \left( 1- \frac{9\gamma}{16}\eta_n^\beta\right)\EB \norm{\vv_n}^4 + \frac{C\Sigma^\prime\eta_n^{4-\beta}}{\gamma}\left(1+ \EB\norm{\vx_n - \vx^\star + \frac{\eta_n^{1-\beta}}{\gamma}(\nabla f(\vx_n) + \xi_n^{(a)})}^4\right)\\
&+ \frac{C\eta_n^{4-3\beta}}{\gamma^3} \left\{\norm{\nabla f(\vx^\star)}^4 
+ \ssum{i}{1}{4}\EB \norm{\gQ_n^{(i)}}^4
\right\},
\end{align*}
where (a) follows from Assumption~\ref{asp:lip_4th_mmt}. For the second term, we have
\begin{align*}
    &\frac{\Sigma^\prime\eta_n^{4-\beta}}{\gamma}\left(1+ \EB\norm{\vx_n - \vx^\star + \frac{\eta_n^{1-\beta}}{\gamma}(\nabla f(\vx_n) + \xi_n^{(a)})}^4\right)\\
    \le& 
    \frac{\Sigma^\prime \eta_n^{4-\beta}}{\gamma} \left\{
    1 + \norm{\nabla f(\vx^\star)}^4
    + \frac{\eta_n^{4-4\beta}}{\gamma^4}\EB \norm{\xi_n^{(a)}}^4\right.\\
    &\left.
    + \left(1+ \frac{L^4\eta_n^{4-4\beta}}{\gamma^4}\right)(\EB \norm{\vu_n - \vx^\star}^4 + \EB \norm{\vv_n}^4)
    \right\}.
\end{align*}
For the last term, using Assumption~\ref{asp:smooth} and the definition of $\gQ_n^{(i)}$ leads to
\begin{align*}
    &\frac{\eta_n^{4-3\beta}}{\gamma^3}\left\{
    \norm{\nabla f(\vx^\star)}^4 +
    \ssum{i}{1}{4} \EB\norm{\gQ_n^{(i)}}^4
    \right\}\\
    \le& \frac{\eta_n^{4-3\beta}}{\gamma^3}\norm{\nabla f(\vx^\star)}^4
    + \frac{L^4\eta_n^{4-3\beta}}{\gamma^3}
    \left\{
    \frac{\eta_n^{4-4\beta}}{\gamma^4}\EB\norm{\xi_n^{(a)}}^4 +
    \frac{L^4\eta_n^{4-4\beta}}{\gamma^4}
    (\EB \norm{\vu_n - \vx^\star}^4 + \EB \norm{\vv_n}^4)
    \right.\\
    &\left.
    + \left(
    \EB \norm{\vv_n}^4 + \frac{\eta_n^{4-4\beta}}{\gamma^4}\norm{\nabla f(\vx^\star)}^4
    \right) + \norm{\vu_n - \vx^\star}^4
    \right\}.
\end{align*}
Combining these facts and rearranging implies
\begin{equation}\label{eq:v_4th_mmt_ieq}
\begin{aligned}
&\EB \norm{\vv_{n+1}}^4 \le
\left[
1 - \frac{9\gamma}{16}\eta_n^\beta + 
\left(
\frac{C\Sigma^\prime\eta_n^{4-\beta}}{\gamma} + 
\frac{CL^4 \eta_n^{4-3\beta}}{\gamma^3}
\right) \left(
1 + \frac{L^4\eta_n^{4-4\beta}}{\gamma^4}
\right)
\right] \EB \norm{\vv_n}^4\\
&+ \left[ 
\frac{C\Sigma^\prime\eta_n^{4-\beta}}{\gamma}
\left(1+ \frac{L^4\eta_n^{4-4\beta}}{\gamma^4}\right) + 
\frac{CL^4 \eta_n^{4-3\beta}}{\gamma^3}\left(
1 + \frac{L^4\eta_n^{4-4\beta}}{\gamma^4}
\right)
\right] \EB \norm{\vu_n - \vx^\star}^4\\
&+ \left[ 
\frac{C\Sigma^\prime\eta_n^{4-\beta}}{\gamma} +
\frac{C\eta_n^{4-3\beta}}{\gamma^3} +
\frac{CL^4 \eta_n^{4-3\beta}}{\gamma^3} \times
\frac{\eta_n^{4-4\beta}}{\gamma^4}
\right] ( 1 + \norm{\nabla f(\vx^\star)}^4)\\
&+ \left(
\frac{C\Sigma^\prime \eta_n^{4-\beta}}{\gamma} + 
\frac{CL^4 \eta_n^{4-3\beta}}{\gamma^3}
\right)\frac{\eta_n^{4-4\beta}}{\gamma^4} \EB \norm{\xi_n^{(a)}}^4\\
&\overset{(a)}{\le} \left(
1 - \frac{9\gamma}{16}\eta_n^{\beta} + \frac{\gamma}{32}\eta_n^\beta
\right)\EB \norm{ \vv_n }^4 +
\frac{CL^4 \eta_n^{4-3\beta}}{\gamma^3} \EB \norm{\vu_n - \vx^\star}^4\\
&+ \frac{C\eta_n^{4-3\beta}}{\gamma^3} (1+\norm{\nabla f(\vx^\star)}^4) + 
\frac{CL^4 \eta_n^{1+7(1-\beta)}}{\gamma^7} \Sigma^\prime (1+ \EB \norm{\vx_n - \vx^\star}^4)\\
&\le \left(
1 - \frac{17\gamma}{32}\eta_n^\beta +
\frac{C\Sigma^\prime L^4 \eta_n^{1+7(1-\beta)}}{\gamma^7}
\right) \EB \norm{\vv_n}^4 +
\frac{C\eta_n^{4-3\beta}}{\gamma^3} ( 1 + \norm{\nabla f(\vx^\star)}^4)
\\
&+ 
\frac{CL^4 \eta_n^{4-3\beta}}{\gamma^3}\left(
1 + \frac{\Sigma^\prime \eta_n^{4-4\beta}}{\gamma^4}
\right) \EB \norm{ \vu_n - \vx^\star }^4 +
\frac{C\Sigma^\prime L^4 \eta_n^{1+7(1-\beta)}}{\gamma^7}\\
&\le \left( 1 - \frac{\gamma}{2}\eta_n^{\beta}\right) \EB \norm{\vv_n}^4 + \frac{CL^4\eta_n^{4-3\beta}}{\gamma^3}\EB \norm{ \vu_n - \vx^\star}^4 + \frac{C\eta_n^{4-3\beta}}{\gamma^3} (1 + \norm{\nabla f(\vx^\star)}^4).
\end{aligned}
\end{equation}
Here (a) is also guaranteed by Assumption~\ref{asp:lip_4th_mmt}.
Summing (\ref{eq:u_4th_mmt_ieq}) and (\ref{eq:v_4th_mmt_ieq}) produces
\begin{align*}
    &(\EB \norm{\vu_{n+1} - \vx^\star}^4 + \EB \norm{\vv_{n+1}}^4)\le
    \left(
    1 - \frac{\gamma}{2} \eta_n^\beta + \frac{CL^4\eta_n}{\mu^3}
    \right) \EB \norm{ \vv_n}^4\\
    &+ \left(
    1 - \mu\eta_n + \frac{CL^4 \eta_n^{4-4\beta}}{\mu^3}
    \right) \EB \norm{\vu_n - \vx^\star}^4 +
    \frac{C\eta_n^{4-3\beta}}{\gamma^3}(1+\norm{\nabla f(\vx^\star)}^4) + \frac{C\Sigma^\prime \eta_n^3}{\mu}\\
    &\le \left(1 - \frac{\mu}{2}\eta_n\right) 
    \EB \norm{\vu_n - \vx^\star}^4 + 
    \left( 1 - \frac{\gamma}{3}\eta_n^\beta\right) \EB \norm{\vv_n}^4\\
    &+ \frac{C\eta_n^{4-3\beta}}{\gamma^3}(1+ \norm{\nabla f(\vx^\star)}^4) +
    \frac{C\Sigma^\prime\eta_n^3}{\mu}\\
    &\le \left(
    1 - \frac{\mu}{2}\eta_n
    \right) (\EB \norm{\vu_n - \vx^\star}^4 + \EB \norm{\vv_n}^4) + \frac{C\eta_n^{4-3\beta}}{\gamma^3}(1+ \norm{\nabla f(\vx^\star)}^4) +
    \frac{C\Sigma^\prime\eta_n^3}{\mu}
\end{align*}
By leveraging Lemma~\ref{lem:converge_r_t}, this implies
\[
\EB \norm{\vu_n - \vx^\star}^4 +
\EB \norm{\vv_n}^4 \precsim
\frac{\eta_n^{3(1-\beta)}}{\mu\gamma^3}(1+\norm{\nabla f(\vx^\star)}^4) + \frac{\Sigma^\prime \eta_n^2}{\mu^2}.
\]
Especially, the left two terms of this both satisfy the bound. What can be deduced is
\begin{align*}
\EB \norm{\vv_{n+1}}^4 &\le 
(1 - \frac{\gamma}{2}\eta_n^{\beta})\EB \norm{\vv_n}^4 + \frac{C\eta_n^{4-3\beta}}{\gamma^3}(1 + \norm{\nabla f(\vx^\star)}^4) 
\\
&+ \frac{CL^4 \eta_n^{4-3\beta}}{\mu\gamma^3}
\left\{
\frac{\eta_n^{3(1-\beta)}}{\gamma^3}(1+\norm{\nabla f(\vx^\star)}^4) + \frac{\Sigma^\prime \eta_n^2}{\mu^2}\right\}\\
&\le ( 1- \frac{\gamma}{2}\eta_n^\beta) \EB \norm{\vv_n}^4 + \frac{C\eta_n^{4-3\beta}}{\gamma^3}
\left\{
\left( 1 + \frac{L^4\eta_n^{3(1-\beta)}}{\mu\gamma^3} \right)(1+\norm{\nabla f(\vx^\star)}^4) + \frac{\Sigma^\prime\eta_n^2}{\mu^2}
\right\}\\
&\le (1 - \frac{\gamma}{2}\eta_n^\beta) \EB \norm{\vv_n}^4 + \frac{C\eta_n^{4-3\beta}}{\gamma^3}
(1 + \norm{\nabla f(\vx^\star)}^4).
\end{align*}
By Lemma~\ref{lem:converge_r_t}, it follows that
\[
\EB \norm{\vv_n}^4 \precsim \frac{\eta_n^{4(1-\beta)}}{\gamma^4}(1 + \norm{\nabla f(\vx^\star)}^4).
\]
Taking this result back to (\ref{eq:u_4th_mmt_ieq}) yields
\begin{align*}
    &\EB \norm{ \vu_n - \vx^\star }^4
    \le (1 - \mu\eta_n) \EB \norm{\vu_n - \vx^\star}^4 + \frac{CL^4\eta_n^{5-4\beta_n}}{\mu^3\gamma^4}\norm{\nabla f(\vx^\star)}^4 + \frac{C\Sigma^\prime\eta_n^3}{\mu}\\
    &+ \frac{CL^4 \eta_n}{\mu^3} \times \frac{\eta_n^{4(1-\beta)}}{\gamma^4}(1 + \norm{\nabla f(\vx^\star)}^4)\\
    &\le (1 - \mu\eta_n) \EB \norm{\vu_n - \vx^\star}^4 + \frac{CL^4\eta_n^{5-4\beta_n}}{\mu^3\gamma^4}(1 + \norm{\nabla f(\vx^\star)}^4) + \frac{C\Sigma^\prime\eta_n^3}{\mu}.
\end{align*}
Finally, owing to Lemma~\ref{lem:converge_r_t}, it holds that
\[
\EB \norm{\vu_n - \vx^\star}^4 \precsim
\frac{L^4 \eta_n^{4(1-\beta)}}{(\mu\gamma)^4}
(1 + \norm{\nabla f(\vx^\star)}^4) + \frac{\Sigma^\prime \eta_n^2}{\mu^2}.
\]
Up to now, the proof is concluded.
\end{proof}

\section{Experimental Details}\label{sec:append_expe}
This section, similar to Section~\ref{sec:exp}, is also devoted to experiments.
We first present the choice of parameters in Section~\ref{sec:exp} and provide some complementary results on problem~\eqref{eq:proj_opt_quadratic}.
Then we construct an example of distributed learning and show the experimental results of this problem.
Finally, we focus on synthetic classification problems and plot the graphs of MSEs to show the convergence rate.

\paragraph{Parameters}
For $\alpha = 1$, we take the initial step size $\eta_0$ as $1$; for $\alpha < 1$, we take $\eta_0$ as $0.2$.
For $\beta = 0$, we take the initial value of the projection probability $\gamma$ as $0.1 $; for $\beta > 1$, we take $\gamma$ as $0.5$.
In Appendix~\ref{sec:append_expe_dis}, the choices of parameters are the same without additional notes.

Note that for the problem~\eqref{eq:proj_opt_quadratic},
the objective is $\mu$-strongly convex with $\mu=1$.
Then our theorem results require $\eta_0 > 2 / \mu = 2 $ for $\alpha = 1$.
It seems that taking $\eta_0 = 1$ violates this condition.
In fact, the condition $\eta_0 > 2 / \mu $ is just for the convenience of proof and the factor $2$ is not essential.
It can be replaced by any real number greater than $1/2$.

\begin{figure}[t!]
    \centering
    \includegraphics[width=\textwidth]{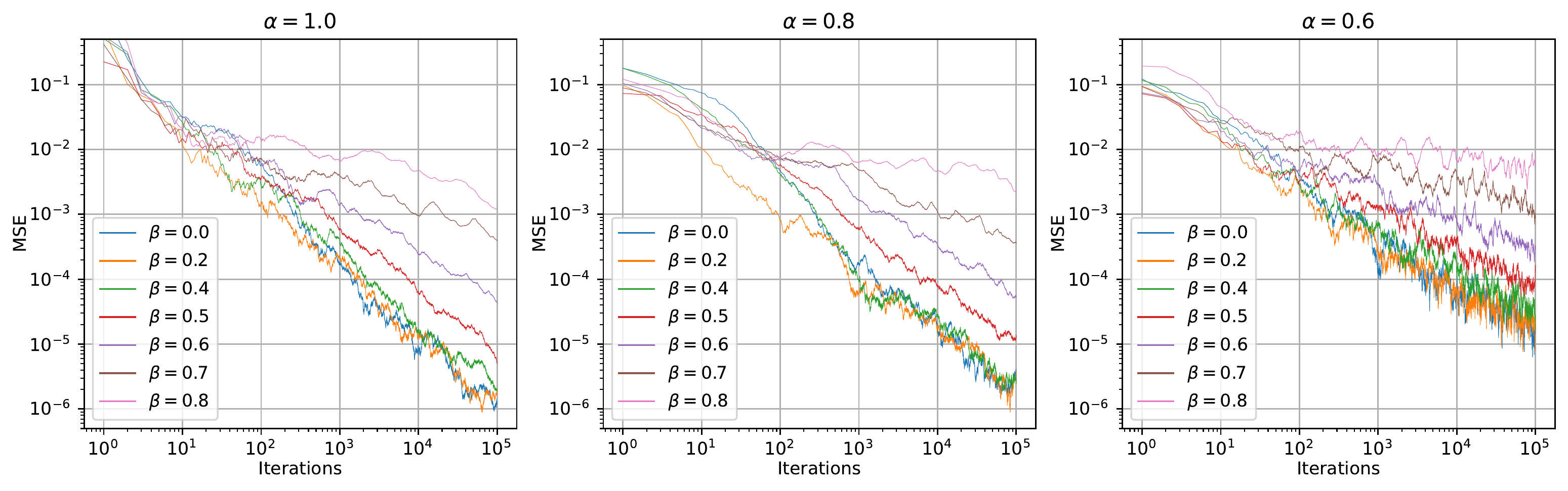}
    \caption{The log-log scale graphs of averaged MSEs 
    on problem~\eqref{eq:dis_opt_quadratic} 
    over $10$ repetition vs iterations.}    \label{fig:fl_converge_rate}
\end{figure}

\subsection{Complementary Results for Section~\ref{sec:exp}}\label{sec:append_expe_complement}
In this subsection, we present some complementary results on problem~\eqref{eq:proj_opt_quadratic}.

\begin{figure}[t!]
    \centering
    \includegraphics[width=\textwidth]{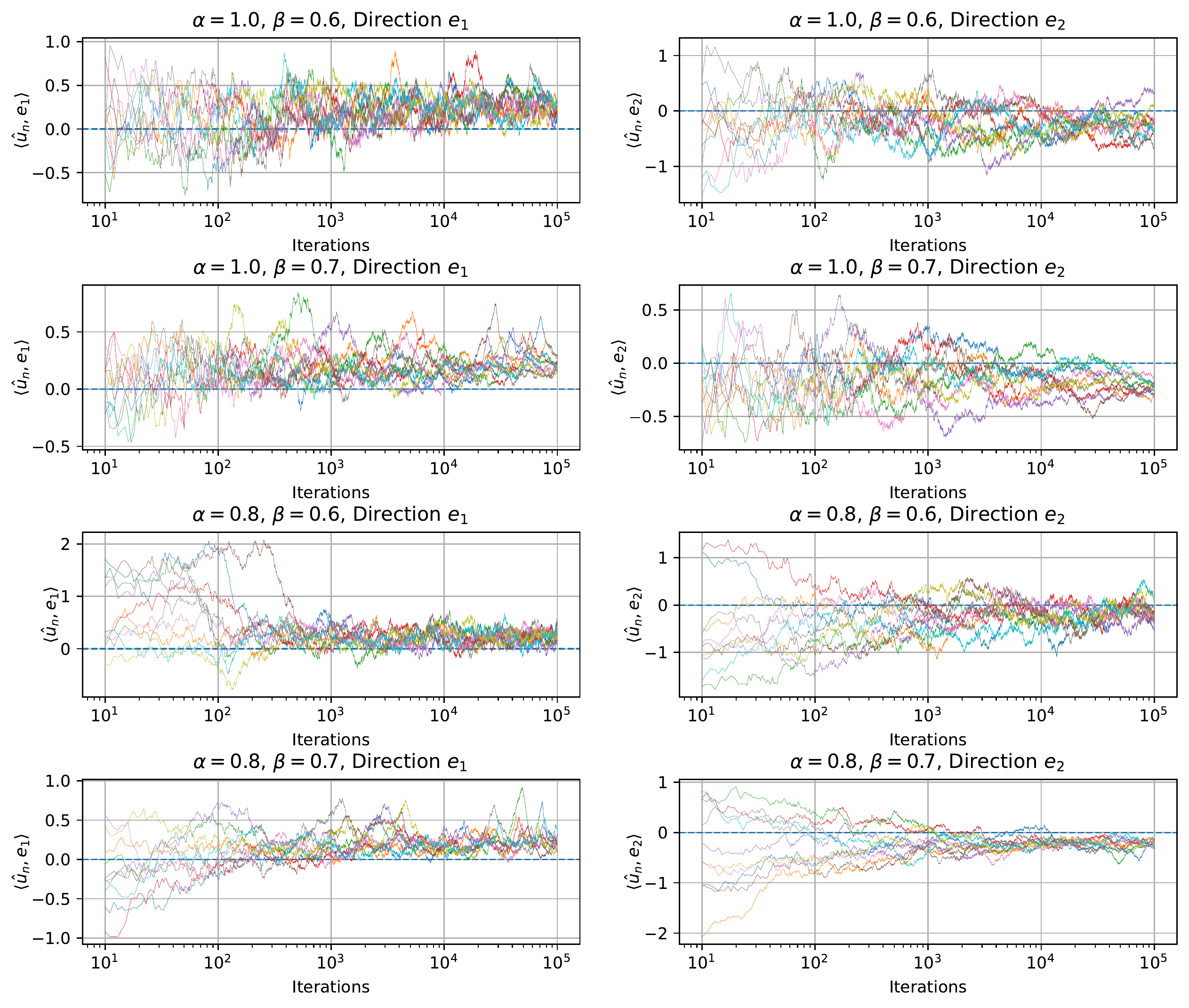}
    \caption{Trajectories of $\hat{\vu}_n$ along two orthogonal directions
    on problem~\eqref{eq:proj_opt_quadratic} 
    over $10$ repetition vs iterations.}    \label{fig:lc_traj_u}
\end{figure}

\paragraph{Trajectories of $\hat{\vu}_n$}
Recall that in Section~\ref{sec:exp}, we set $d_1 = 5$ and $d_2 = 2$. Then one can show the $\hat{\vu}_n$ sequence lies in a two-dimensional subspace of $\sR^5$.
For $\alpha \in \{ 1, 0.8\}$ and $\beta \in \{ 0.6, 0.7\}$, we plot the trajectories of $\hat{\vu}_n$ on problem~\eqref{eq:proj_opt_quadratic} along two orthogonal directions $\ve_1$ and $\ve_2$ in Fig.~\ref{fig:lc_traj_u}.
Different colors indicate different repetitions.
Fig.~\ref{fig:lc_traj_u} can be analyzed together with Fig.~\ref{fig:lc_bias_mse}.
From Fig.~\ref{fig:lc_traj_u}, we find $\hat{\vu}_n$ will not converge to $\vzero$.
In fact, it seems that the sequence can never converge to any constant vector because it keeps fluctuating, This is consistent with the results in Fig.~\ref{fig:lc_bias_mse}, where the MSEs between $\hat{\vu}_n$ and the bias vector $\vmu$ have a non-negligible magnitude.
Theorem~\ref{thm:asy_bias_u} shows $\hat{\vu}_n$ will converge eventually.
We believe this will happen if we let the number of iterations sufficiently large, which, however, is beyond our computing power.

\begin{figure}[t!]
    \centering
    \includegraphics[width=\textwidth]{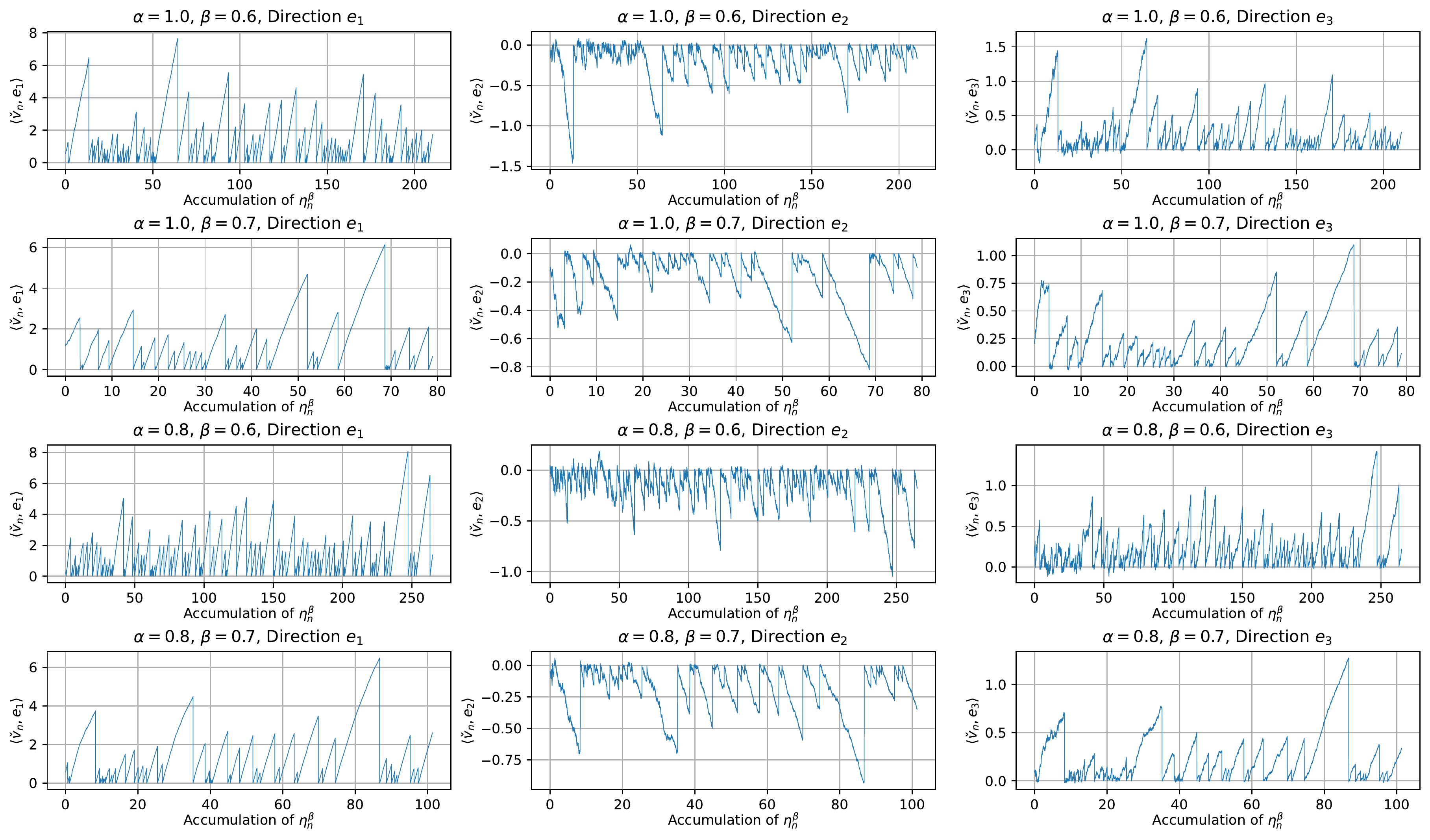}
    \caption{Trajectories of $\hat{\vv}_n$ along three orthogonal directions
    on problem~\eqref{eq:proj_opt_quadratic} 
    vs the accumulation of $\eta_n^\beta$.}    \label{fig:lc_traj_v}
\end{figure}

\paragraph{Trajectories of $\hat{\vv}_n$}
Since in Section~\ref{sec:exp}, we set $d_1 = 5$ and $d_2 = 2$, one can show the $\hat{\vv}_n$ sequence lies in a three-dimensional subspace of $\sR^5$.
For $\alpha \in \{1, 0.8\}$ and $\beta \in \{ 0.6, 0.7 \}$, we plot the trajectories of $\check{\vv}_n$ along three orthogonal directions $\ve_1$, $\ve_2$ and $\ve_3$ vs accumulation of $\eta_n$ in Fig.~\ref{fig:lc_traj_v}. 
The value of the horizontal coordinate is $\sum_{i=s}^n \eta_i^\beta$, where $s$ is the start point of the trajectory that aims to eliminate the irregular behavior in the early stage of the optimization process.
We set $s = 1000$.
Observe that the trajectories in Fig.~\ref{fig:lc_traj_v} come in a jagged manner and the peak value does not vanish or explode. This is because we have chosen a suitable rescaled version $\check{\vv}_n$ of $\vv_n$ and such a behavior can be captured by Theorem~\ref{thm:jump_approx}.

\begin{figure}[t!]
    \centering
    \includegraphics[width=0.8\textwidth]{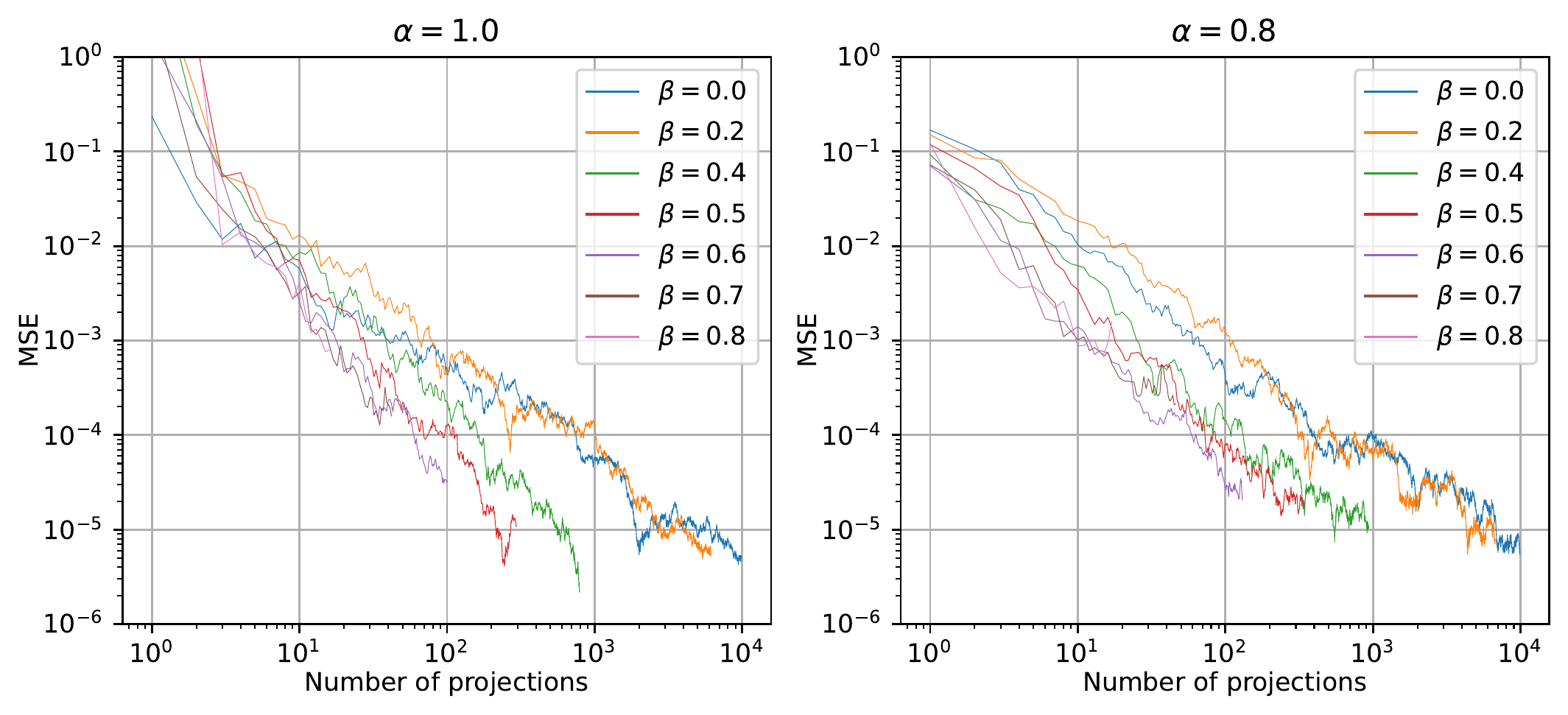}
    \caption{The log-log scale graphs of averaged MSE on problem~\eqref{eq:proj_opt_quadratic} over $10$ repetitions vs the number of projections.}
    \label{fig:lc_mse_over_proj}
\end{figure}

\paragraph{Convergence Rates in terms of the Number of Projections}
Recall that in Section~\ref{sec:converge}, 
and Section~\ref{sec:exp}, 
we establish the convergence rates in terms of the number of iterations and provide the log-log scale graphs of averaged MSEs vs the number of iterations. To better capture the influence of projections, we consider the convergence rates in terms of the number of projections and plot the log-log scale graphs of averaged MSEs vs the number of projections.

In our method, at the $n$-th iteration, the projection probability is $p_n = \min\{ \gamma \eta_n^\beta, 1 \} = \Theta (n^{-\alpha \beta})$.
As a result, after $n$ steps of iterations, the number of projections $m$ should be of the order $\Theta (n ^{1 - \alpha \beta} )$.
Suppose that after $m$ steps of projections, we obtain the variable $\vx_m$ and $\vu_m = \gP_{\mA^\bot} (\vx_m)$.
By Theorem \ref{thm:converge}, we have 
$\sE \norm{ \vu_m - \vx^\star }^2 = \gO \left( m^ { - \frac{ \alpha \min\{ 1, 2 - 2 \beta \} }{ 1 - \alpha \beta } } \right)$.
For $0 \le \beta < 0.5$, the rate is of the order $\gO \left( m^{- \frac{\alpha}{1 - \alpha \beta}} \right)$ and a larger $\beta$ leads to a faster rate.
For $0.5 \le \beta < 1$, the rate is of the order $ \gO \left( m^{ - \frac{2 \alpha (1 - \beta)}{1 - \alpha \beta} } \right) = \gO \left( m^{ - 2 \left( 1 - \frac{1 - \alpha}{1 - \alpha \beta} \right)} \right).$
If $\alpha = 1$, the choice of $\beta$ does not influence the convergence rate; if $\alpha < 1$, a larger $\beta$ leads to a slower rate.

To conclude, if we only focus on the complexity of projection steps and ignore the cost of gradient computation, $\beta=0.5$ is always the best choice. 

Then for $\alpha \in \{1.0, 0.8 \}$ and $\beta \in \{0, 0.2, 0.4, 0.5, 0.6, 0.7, 0.8 \}$, we plot the log-log scale graphs of averaged MSEs vs the number of projections over $10$ repetitions in Fig.~\ref{fig:lc_mse_over_proj}.
For each repetition, we run $100000$ steps of LPSA.
We find that for $\alpha = 1$, the line of convergence rates becomes steeper as $\beta$ increases, until $\beta \ge 0.5$, when the lines of convergence rates are almost parallel.
For $\alpha = 0.8$, when $\beta$ is closer to $0.5$, the lines of convergence rates are steeper. 
These are consistent with our analysis above.

\subsection{Results on a Distributed Example}\label{sec:append_expe_dis}

\paragraph{A Distributed Example}
Recall that in Example~\ref{eg:distributed}, we have shown that typical distributed optimization problems can be formulated as a linear equality-constrained optimization problem \eqref{eq:proj_opt}.
Here we focus on the quadratic objective of the following form
\begin{align}
\label{eq:dis_opt_quadratic}
    \min_{\vx^{(1)}, \vx^{(2)}, \dots, \vx^{(N)}} 
    \frac{1}{N} \sum_{i=1}^N \left[ \frac{1}{2} (\vx^{(i)})^\top \mS_i \vx^{(i)} - \vb_i^\top \vx^{(i)} + (\zeta^i)^\top \vx^{(i)} \right]
    ~~\text{ s.t. } \vx^{(1)}=\cdots=\vx^{(N)},
\end{align}
where $\mS_i \in \sR^{d \times d} $ is the local Hessian matrix at client $i$, $\vb_i \in \sR^d $ is a local vector at client $i$, and $\zeta_i \in \sR^d$ is the local randomness at client $i$. 
And $\mS_i$, $\vb_i$ and $\zeta_i$
are generated similarly as we generate $\mS, \vb$ and $\zeta$ in the \textit{Experimental Setup} part of Section~\ref{sec:exp}, except that we replace $d_1$ there with $d$ here.
Besides, the tuples ${(\mS_i, \vb_i, \zeta_i)}_{i=1}^N$ are mutually independent.
In our experiments, we set $d = 2$ and $N = 5$.


\paragraph{Convergence rates}
We plot the log-log scale graphs of averaged MSEs on problem~\eqref{eq:dis_opt_quadratic} over $10$ repetitions vs iterations in Fig.~\ref{fig:fl_converge_rate}.
The value of $\alpha $ is set as $\{1, 0.8, 0.6 \}$ and the value of $\beta$ is from $\{0, 0.2, 0.4, 0.5, 0.6, 0.7, 0.8 \}$.
For each repetition, we run $100000$ steps of LPSA.
The results in Fig.~\ref{fig:fl_converge_rate} satisfy the convergence rate in Theorem~\ref{thm:converge} and are consistent with those in Fig.~\ref{fig:lc_converge_rate}.
For $\beta \in [0, 1/2)$, the value of $\beta$ does not affect the slope, while for $\beta \in  (1/2, 1)$, larger $\beta$ and smaller $\alpha$ both lead to smoother lines. $\beta = 1/2$ is the point where the phase transition occurs.

\begin{figure}[t!]
    \centering
    \includegraphics[width=\linewidth]{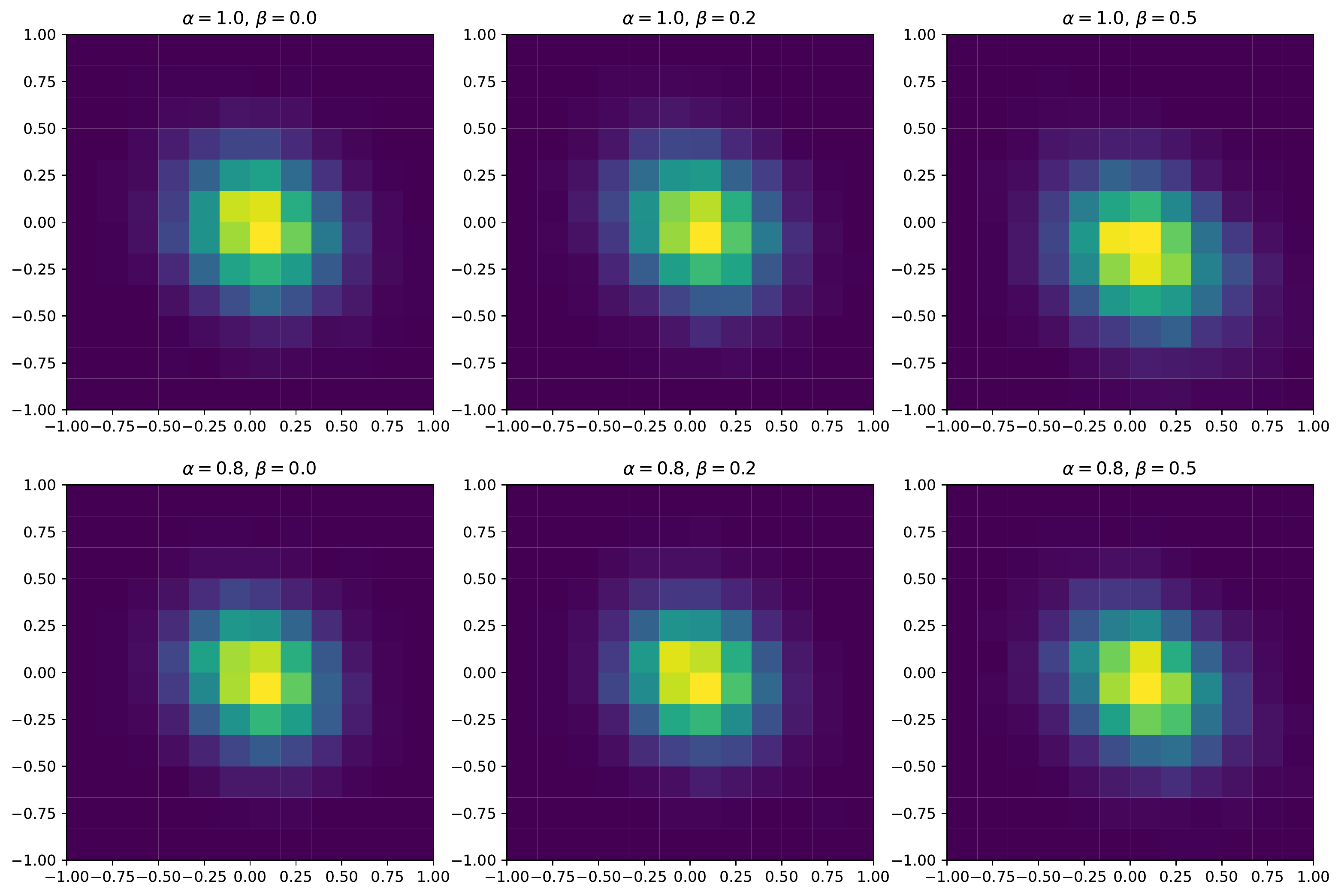}
        \caption{The heatmaps of (centralized) $\check{\vu}_{n}$ ($n=10001$) along two orthogonal directions over $10000$ repetitions
        on problem~\eqref{eq:dis_opt_quadratic}. 
        }
        \label{fig:fl_heatmap}
\end{figure}

\paragraph{Frequent Projection}
For $\alpha=1$ and $\beta \in \{ 0, 0.2, 0.5 \}$, we run $10000$ steps of LPSA on problem~\eqref{eq:dis_opt_quadratic} over $10000$ repetitions 
and pick up the last 
iterates ${\vu}_{10001} $. 
For these iterates, we compute the rescaled vectors $\check{\vu}_{10001}$ (as defined in Section~\ref{sec:freq_proj}).
One can check that the $\check{\vu}_n$ sequence lies in a two-dimensional subspace of $\sR^5$.
In fact, it suffices to pick up the first two dimensions of $\check{\vu}_n$ because all the odd dimensions share the same value, and so do all the even dimensions.
Then for $\beta \in \{0, 0.2 \}$, we plot the heatmaps of these $\check{\vu}_{10001}$s across two orthogonal directions of the subspace in the left two columns of Fig.~\ref{fig:fl_heatmap}.
For $\beta = 0.5$, we first centralize $\check{\vu}_{10001}$, that is, subtracting $\check{\vu}_{10001}$ by the bias vector $\vmu$ in Theorem~\ref{thm:noncen_u},
and then plot the heatmaps as we did before
in the last column of Fig.~\ref{fig:fl_heatmap}.
We find the heatmaps in Fig.~\ref{fig:fl_heatmap} share a similar form as those in Fig.~\ref{fig:lc_heatmap},
and agree with Theorems~\ref{thm:diff_approx} and \ref{thm:noncen_u}, where the limiting distribution of the (centralized) $\check{\vu}_n$ is Gaussian. 

\begin{figure}[t!]
    \centering
    \begin{minipage}{0.49\textwidth}
        \centering
        \includegraphics[width=\linewidth]{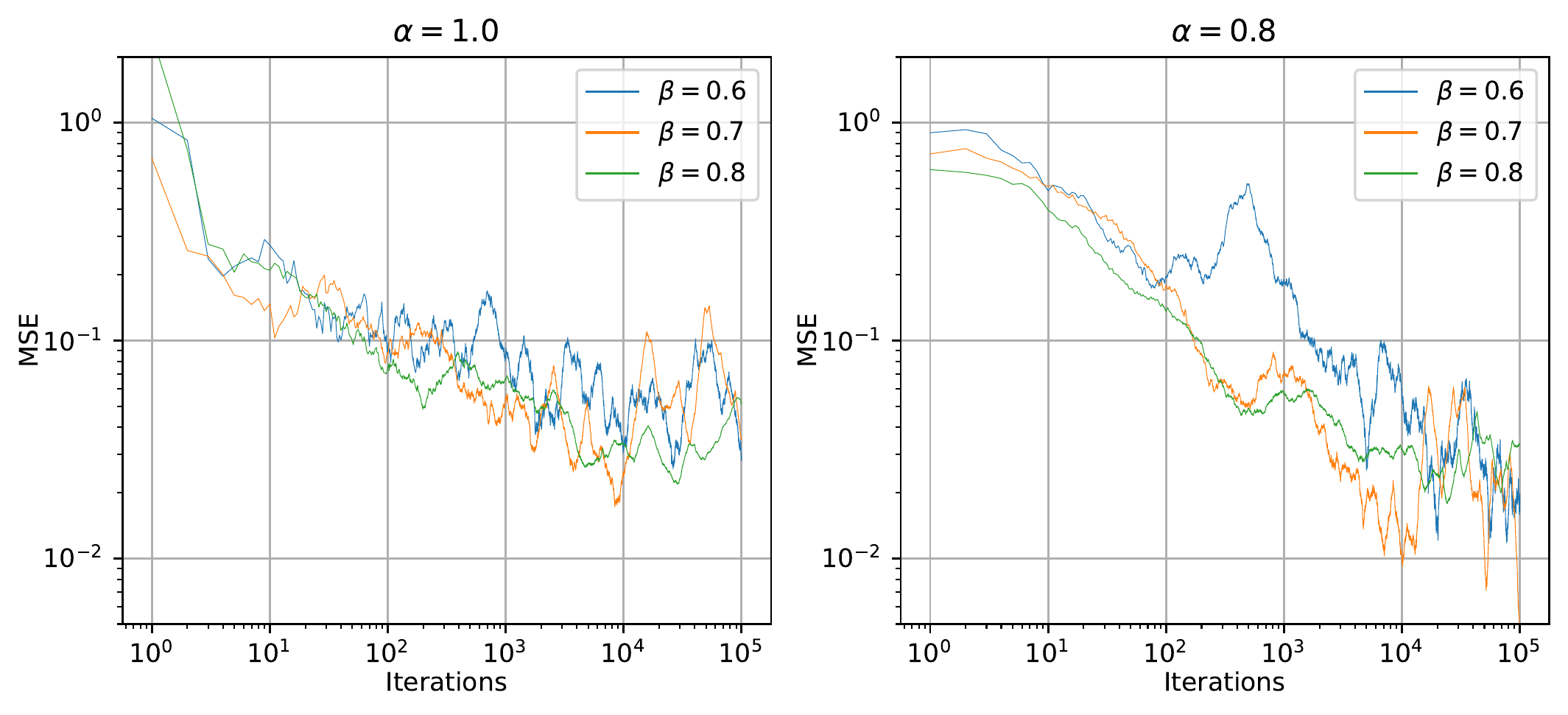}
        \caption{ The log-log scale graphs of averaged MSEs between $\check{\vu}_n$ and the bias vector $\vmu$ (defined in Theorem~\ref{thm:jump_approx})on problem~\eqref{eq:dis_opt_quadratic} over $10$ repetitions vs iterations .}
        \label{fig:fl_bias_mse}
    \end{minipage}
    \hfill
    \begin{minipage}{0.49\textwidth}
        \centering
        \includegraphics[width=\linewidth]        {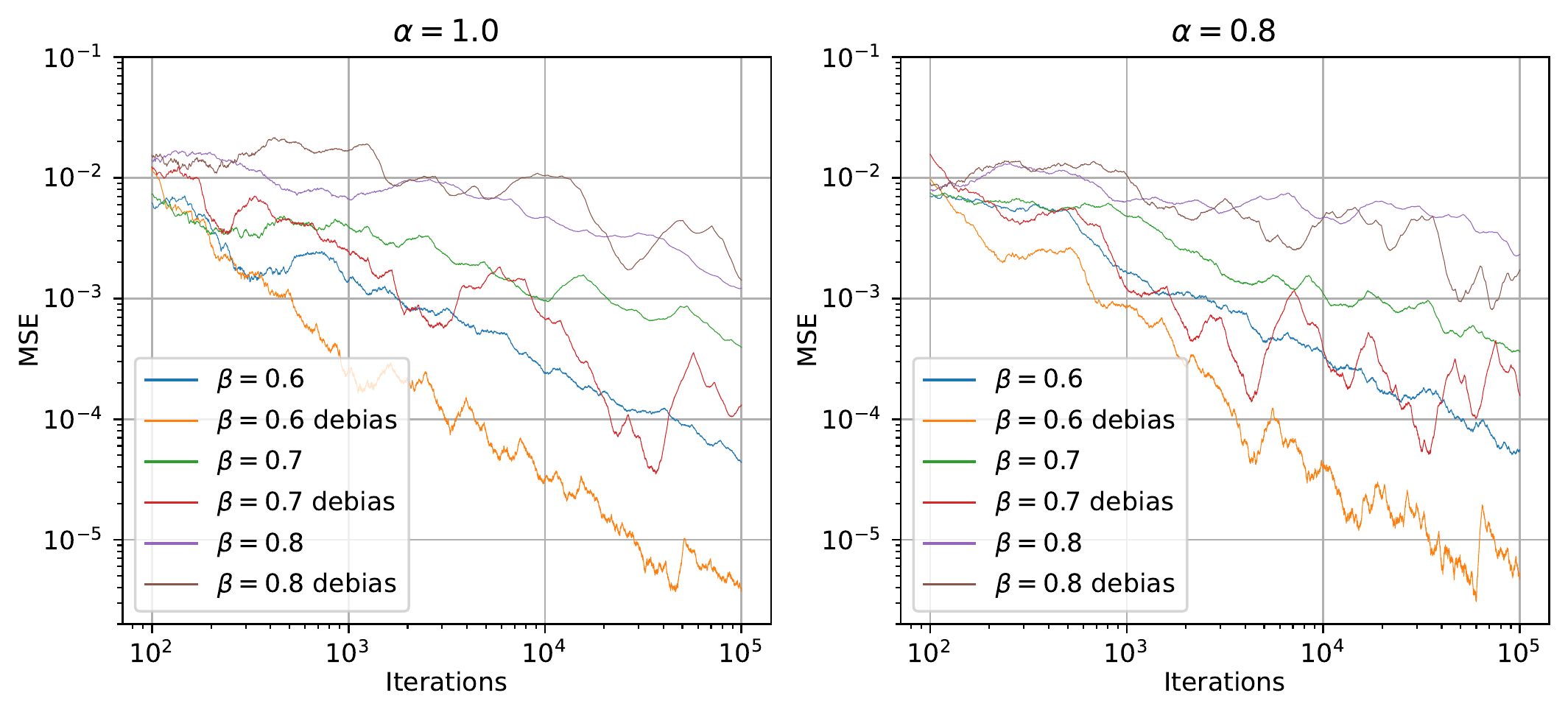}
        \caption{
        The log-log scale graphs of averaged MSEs of LPSA and DLPSA on problem~\eqref{eq:dis_opt_quadratic} over $10$ repetitions vs iterations.
        The `debias' in the legend indicates the line for DLPSA.
        }
        \label{fig:fl_debias}
    \end{minipage}
\end{figure}

\begin{figure}[t!]
    \centering
    \includegraphics[width=\textwidth]{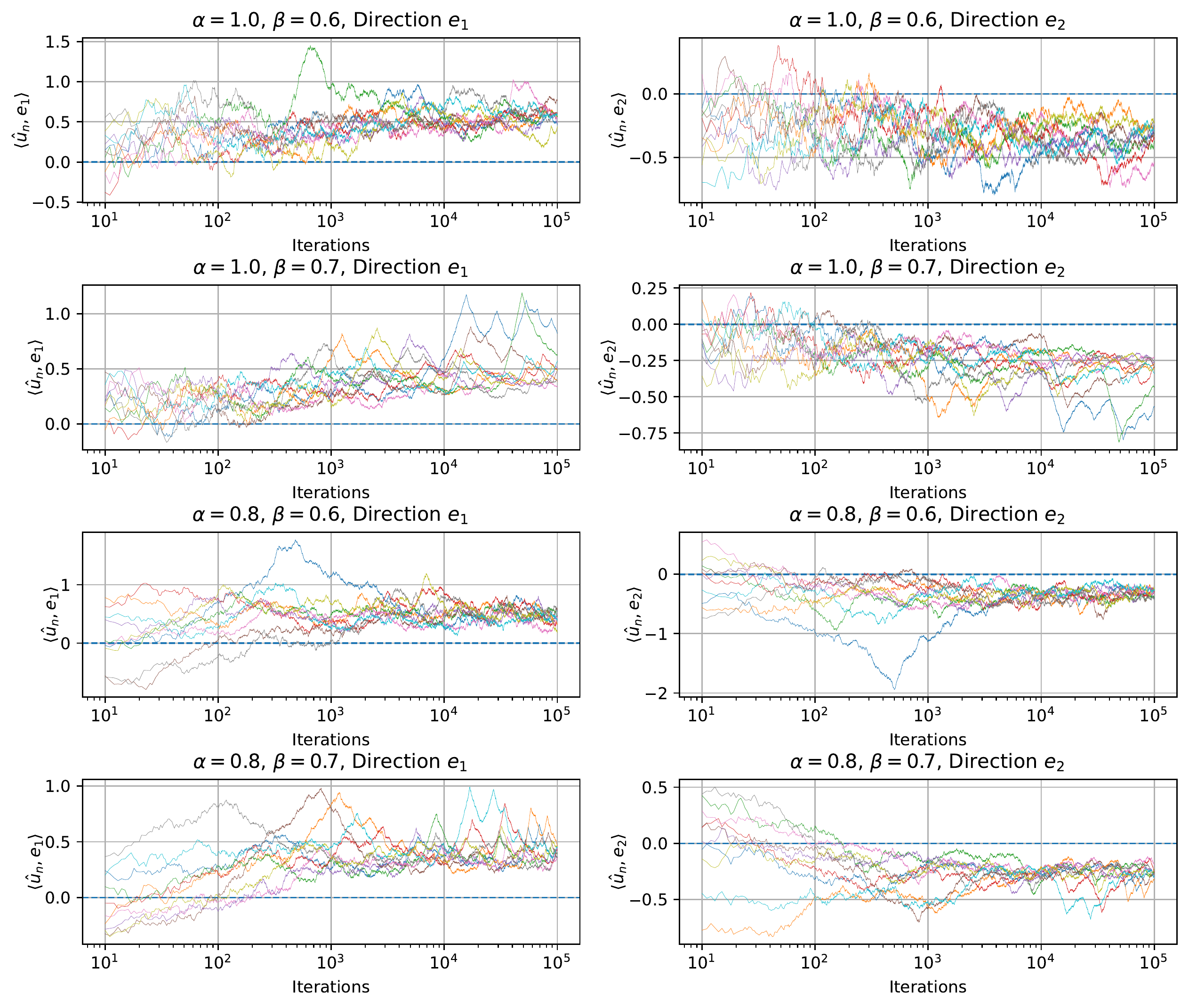}
    \caption{Trajectories of $\hat{\vu}_n$ along two orthogonal directions
    on problem~\eqref{eq:dis_opt_quadratic} 
    over $10$ repetition vs iterations.}    \label{fig:fl_traj_u}
\end{figure}

\paragraph{Occasional Projection}
For $\alpha \in \{ 1.0, 0.8 \}$ and $\beta \in \{ 0.6, 0.7, 0.8 \}$,
we run $100000$ steps of LPSA on problem~\eqref{eq:dis_opt_quadratic} over $10$ repetition, 
and plot the average MSEs between $\hat{\vu}_n$ (defined in Section~\ref{sec:occa_proj}) and the bias vector $\vmu$ (defined in Theorem~\ref{thm:asy_bias_u}) in Fig.~\ref{fig:fl_bias_mse}.
Similar to Fig.~\ref{fig:lc_bias_mse},
the MSEs have a non-negligible magnitude,
but gradually decrease as the number of iterations increases.


\paragraph{Trajectories of $\hat{\vu}_n$}
Since we set $d=2$ and $N=5$, one can check the $\hat{\vu}_n$ sequence lies in a two-dimensional subspace of $\sR^10$.
For $\alpha \in \{ 1, 0.8\}$ and $\beta \in \{ 0.6, 0.7\}$, we plot the trajectories of $\hat{\vu}_n$ on problem~\eqref{eq:dis_opt_quadratic} along two orthogonal directions $\ve_1$ and $\ve_2$ in Fig.~\ref{fig:fl_traj_u}.
Different colors indicate different repetitions.
Together with Fig.~\ref{fig:fl_bias_mse},
Fig.~\ref{fig:fl_traj_u} shows $\hat{\vu}_n$ keeps fluctuating and will not converge to $\vzero$.
Theorem~\ref{thm:asy_bias_u} shows $\hat{\vu}_n$ will converge eventually.
We believe this will happen if we let the number of iterations sufficiently large, which, however, is beyond our computing power.

\paragraph{Debiased algorithm}
Recall that in Section~\ref{sec:dlpsa}, we propose a debiased algorithm DLPSA.
Theorem~\ref{thm:dlpsa_cov_rt} shows that compared to LPSA, DLPSA enjoys superior convergence rates $\gO (n^{ - \alpha \min \{1, 3(1 - \beta) \}})$ for $\beta > 0.5$. 
Thus for $\alpha \in \{ 1.0, 0.8 \}$ and $\beta \in \{ 0.6, 0.7, 0.8 \}$, we run $100000$ steps of DLPSA on problem~\eqref{eq:dis_opt_quadratic} and compare the averaged MSEs over $10$ repetitions between LPSA and DLPSA in Fig.~\ref{fig:fl_debias}.
Since
the approximation of $\nabla f(\vx^\star)$ by $\nabla f(\vx_n, \zeta_n')$ is only valid for large $n$,
we replace the first $100$ steps of DLPSA with LPSA as a warm-up stage.
In Fig.~\ref{fig:fl_debias}, we only plot the MSEs after this stage.
The results in Fig.~\ref{fig:fl_debias} exhibit similar behavior as those in Fig.~\ref{fig:lc_debias}, which implies DLPSA indeed converges faster than LPSA for $\beta > 0.5$.

\begin{figure}[t!]
    \centering
    \includegraphics[width=\textwidth]{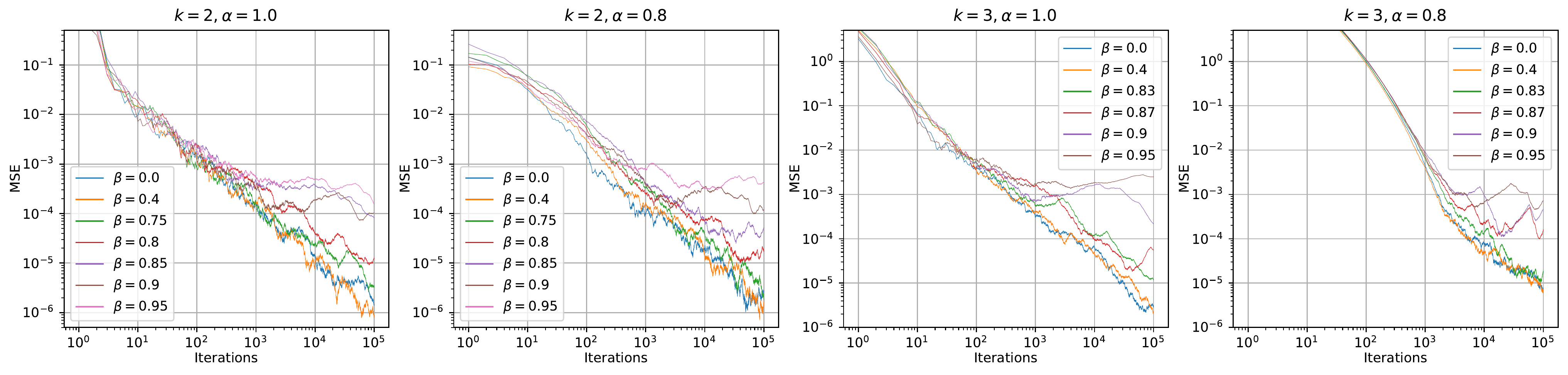}
    \caption{The log-log scale graphs of averaged MSEs on degenerate cases on problem~\eqref{eq:dis_opt_quadratic} over $10$ repetition vs iterations.}    \label{fig:fl_degenerate}
\end{figure}

\paragraph{Degenerate cases}
Generally, the problem~\eqref{eq:dis_opt_quadratic} constructed before  does not satisfy the degenerate condition $\gP_{\mA^\bot} \nabla^2 f(\vx^\star) \nabla f(\vx^\star) = \vzero$.
Similar to Lemma~\ref{lem:well_def_dgnt}, 
we can also find a degenerated example for the distributed setting.
To illustrate the performance of LPSA on degenerate cases, we construct two problems corresponding to $k = 2$ and $3$ respectively, with $k$ the degenerate order defined in Assumption~\ref{asp:degenerate}.
For $k=2$, we still set $d= 2$ and $N = 5$, and choose $\alpha \in \{ 1.0, 0.8 \}$ and $\beta \in \{ 0, 0.4, 0.75, 0.8, 0.85, 0.9, 0.95 \}$.
For $k=3$, we set $d = 4$ and $N = 5$,
and choose $\alpha \in \{ 1.0, 0.8 \}$ and $\beta \in \{ 0, 0.4, 0.83, 0.7, 0.9, 0.95 \}$.
For each pair of $(\alpha, \beta)$, we 
run $100000$ steps of LPSA and plot the log-log scale graphs of averaged MSEs over $10$ repetitions in Fig.~\ref{fig:fl_degenerate}\footnote{For $\alpha < 1$, we set the initial step size as $\eta_0 = 0.3$ instead of $0.2$ before}.
Theorem~\ref{thm:degenerate_bias_u} implies the convergence rate is of the order $ \gO( n^{-\alpha \min\{ 1, 2k(1-\beta) \}} )$.
Fig.~\ref{fig:lc_degenerate} captures the property that the phase transition occurs when $\beta$ crosses $1 - 1 / (2k)$.

\begin{figure}[t!]
    \centering
    \includegraphics[width=0.8\textwidth]{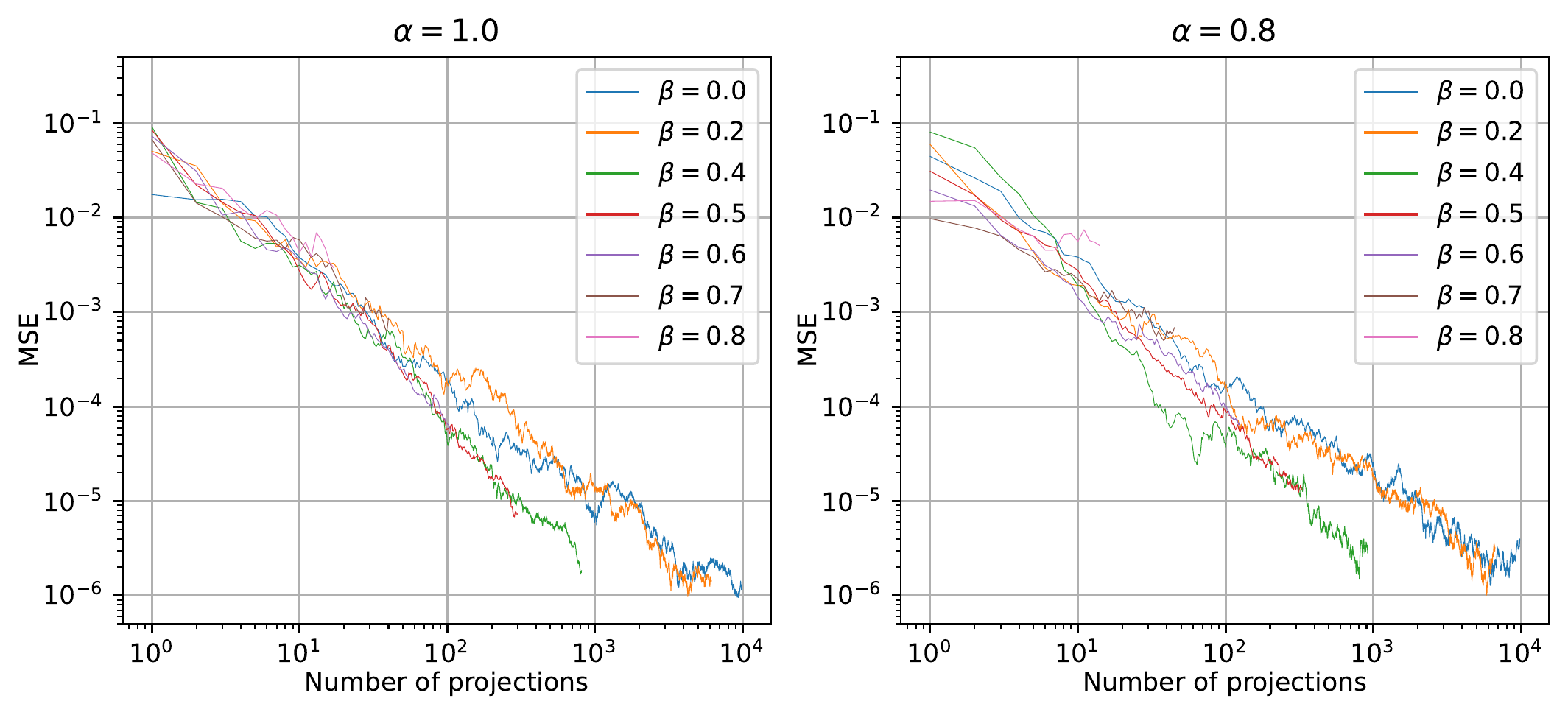}
    \caption{The log-log scale graphs of averaged MSE on problem~\eqref{eq:dis_opt_quadratic} over $10$ repetitions vs the number of projections.}
    \label{fig:fl_mse_over_proj}
\end{figure}

\paragraph{Convergence Rates in terms of the Number of Projections}

For $\alpha \in \{1.0, 0.8 \}$ and $\beta \in \{0, 0.2, 0.4, 0.5, 0.6, 0.7, 0.8 \}$, we plot the log-log scale graphs of averaged MSEs vs the number of projections over $10$ repetitions in Fig.~\ref{fig:fl_mse_over_proj}.
For each repetition, we run $100000$ steps of LPSA.
The results in Fig.~\ref{fig:fl_mse_over_proj} are similar to those in Fig.~\ref{fig:lc_mse_over_proj} and thus consistent with our analysis in Appendix~\ref{sec:append_expe_complement}.

\subsection{Convergence Rates on Synthetic Classification Problems}
In this subsection, we focus on classification problems with cross-entropy loss, and plot graphs of MSEs to show the convergence rates.
$\ell_2^2$ regularization is imposed to ensure the strong convexity of the objective function.

The synthetic datasets are generated by following \cite{li2018federated}.
There are $K$ clients and
the sample $(\vx_k, z_k)$ on the $k$-th client is modeled as $\vx_k \sim \gN (\nu_k, \Lambda) $ and $z_k = \mathrm{argmax}( \mathrm{softmax} (\mW_k \vx_k + \vb_k) )$  where $\Lambda \in \sR^{d \times d}$ is diagonal with the  entry $(j,j)$ equal to $j^{-1.2}$, $\mW_k \in \sR^{C \times d}$ and $\vb_k \in \sR^C$. 
Specially, we consider the following three datasets.

\paragraph{ \texttt{IID} }
All the clients share the same $\mW_k \in \sR^{C \times d}$ and $\vb_k \in \sR^C$ and their entries are modeled as $\gN (0, 1)$. 
$\nu_k \sim \gN(\vzero, \rmI_d)$.
We set $K=100$, $d=60$ and $C=10$.
For this dataset, there is no heterogeneity between the optimal local parameters.
The heterogeneity is all from the diversity of the distributions of $\vx_k$.
For each client, the sample size is around $100$.

The parameter of $\ell_2^2$ regularization is $0.005$.
For $\alpha=1$, we set $\eta_0 = 200$;
for $\alpha = 0.8$, we set $\eta_0 = 40$;
for $\alpha = 0.6$, we set $\eta_0 = 20$.

\paragraph{ \texttt{Synthetic} ($a,b$) }
Each entry of $\mW_k$ and $\vb_k$ is modeled as $\gN (\mu_k, 1)$ with $\mu_k \sim \gN (0, a)$ and $\nu_k \sim \gN(\zeta_k, \rmI_d) $ with $\zeta_k \sim \gN(\vzero, b \rmI_d)$. We set $K=20$, $d=10$ and $C=5$.
$a$ controls how many local models differ from each other and $b$ controls how much the local data for each client differs from that of other clients. They are the two sources of heterogeneity.
For each client, the sample size is around $50$. In this paper, we let $a = b = 1$.

The parameter of $\ell_2^2$ regularization is $0.5$.
For $\alpha=1$, we set $\eta_0 = 1$;
for $\alpha = 0.8$, we set $\eta_0 = 0.3$;
for $\alpha = 0.6$, we set $\eta_0 = 0.1$.

\paragraph{ \texttt{Lincons} }
The data are generated by the same way in \texttt{IID}. Since in \texttt{IID}, all the clients share the same $\mW_k$ and $\vb_k$, we can combine all the samples and obtain the dataset \texttt{Lincons}.
Then we generate the matrix $\mA \in \sR^{610 \times 400}$ whose entries are independent and modeled as $\gN (0, 1)$.

The parameter of $\ell_2^2$ regularization is $0.05$.
For $\alpha=1$, we set $\eta_0 = 8$;
for $\alpha = 0.8$, we set $\eta_0 = 2$;
for $\alpha = 0.6$, we set $\eta_0 = 0.8$.

For all the three datasets, the loss function is defined as the sum of cross entropy loss and $\ell_2^2$ regularization.
The mini-batch size is $4$. As for the probability $p_n$, we reparameterize it as $p_0 n^{-\alpha \beta}$ with $p_0 < 1$. The value of $\alpha$ is from $\{1, 0.8, 0.6 \}$ and the value of $\beta$ is from $\{0, 0.2, 0.4, 0.6, 0.8 \}$. For $\beta = 0$, we set $p_0 = 0.2$; for $\beta > 0$, we set $p_0 = 0.5$. And we run gradient descent $1000$ steps to obtain the value of $\vx^\star$.

\begin{figure}[t!]
    \centering
    \includegraphics[width=\textwidth]{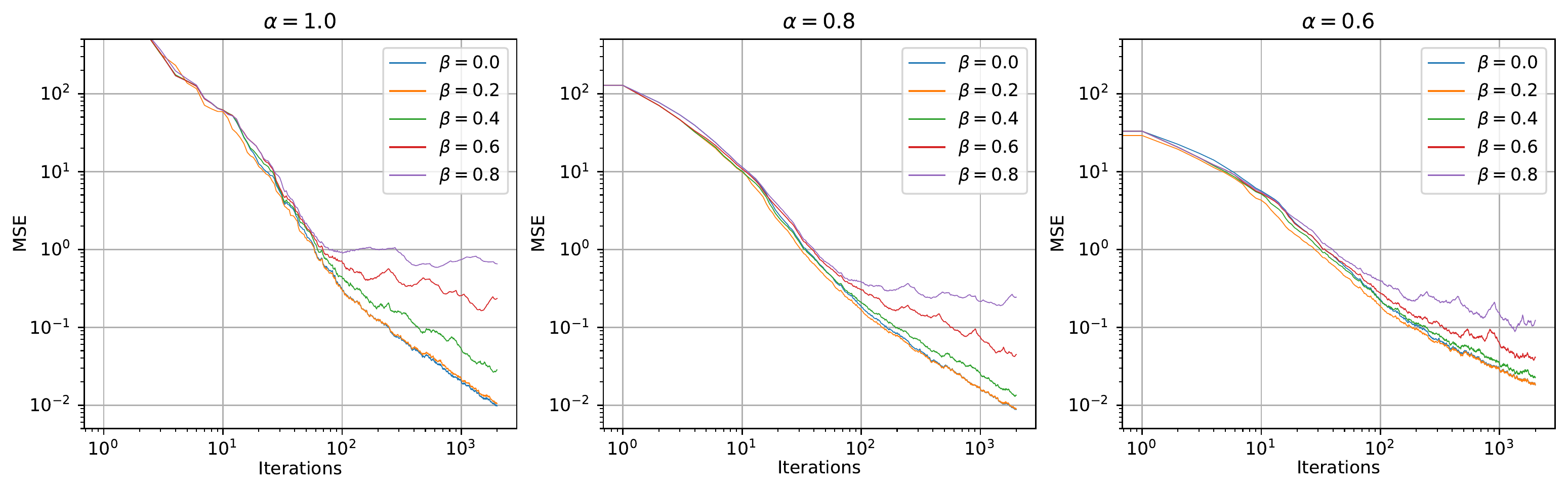}
    \caption{The log-log scale graphs of averaged MSE on \texttt{IID} over $5$ repetition vs iterations.}
    \label{fig:main_converge_rate}
\end{figure}

\begin{figure}[t!]
    \centering
    \includegraphics[width=\textwidth]{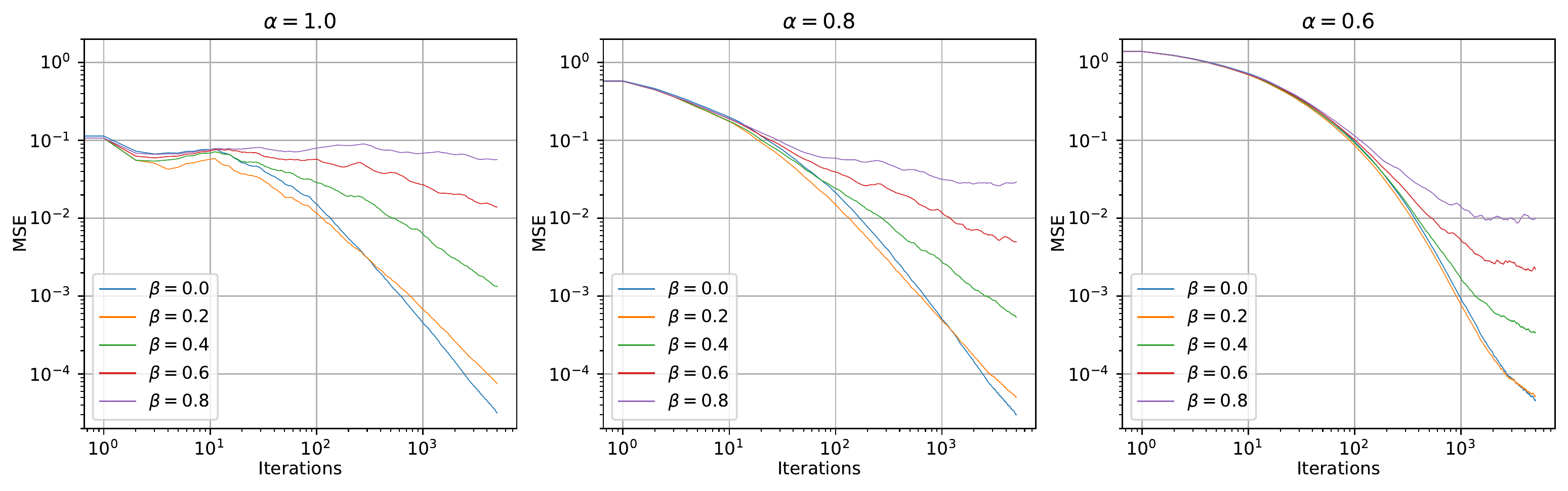}
    \caption{The log-log scale graphs of averaged MSE on \texttt{Synthetic} ($1,1$) over $10$ repetitions vs iterations.}
    \label{fig:converge_syn}
\end{figure}

\begin{figure}[t!]
    \centering
    \includegraphics[width=\textwidth]{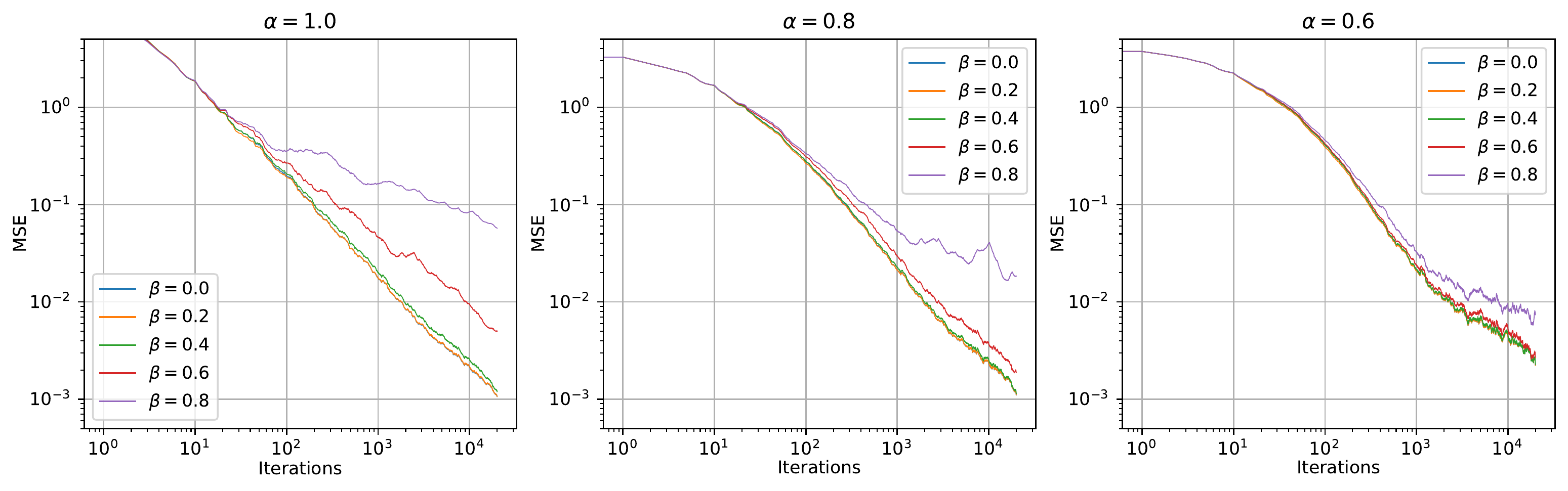}
    \caption{The log-log scale graphs of averaged MSE on \texttt{Lincons} over $10$ repetitions vs iterations.}
    \label{fig:converge_cons}
\end{figure}

We plot the log-log scale graphs of averaged MSEs over $5$ repetitions on \texttt{IID} vs iterations in Figs. \ref{fig:main_converge_rate}, \ref{fig:converge_syn} and \ref{fig:converge_cons}.

\paragraph{Convergence Rates}

We plot the log-log scale graphs of averaged MSEs over $5$ repetitions on \texttt{IID} vs iterations in Fig. \ref{fig:main_converge_rate}
and the log-log scale graphs of averaged MSEs over $10$ repetitions on \texttt{IID} and \texttt{Lincons} vs iterations in Figs. \ref{fig:converge_syn} and \ref{fig:converge_cons}. 
When $\beta < 1/2$, the value of $\beta$ hardly affects the convergence rate; when $\beta > 1/2$, both larger $\beta$ and smaller $\alpha$ lead to a slower convergence rate. This is consistent with the result of Theorem \ref{thm:converge}.

\end{document}